\documentclass[a4paper,12pt,amsart,frenchb]{article}

\usepackage{amsmath,amsbsy,amsfonts,amssymb}
\usepackage[french]{babel}
\newcommand{\ass}[2]{\vskip0.3cm\noindent
{\bf {#1}}. { \sl {#2}}\vskip0.3cm\noindent
}
  
\begin{document}
 
   \title{  Stabilisation de la formule des traces tordue II:   int\'egrales orbitales et endoscopie sur un corps local non-archim\'edien; d\'efinitions et \'enonc\'es des r\'esultats}
\author{J.-L. Waldspurger}
\date{23 janvier 2014}
\maketitle

 {\bf Introduction}
 
 Ceci est le deuxi\`eme d'une s\'erie d'articles, en collaboration avec C. Moeglin, visant \`a \'etablir la stabilisation de la formule des traces tordue. On y donne les d\'efinitions des termes locaux intervenant dans la partie g\'eom\'etrique de cette formule, sous la restriction que le corps local de base $F$ est suppos\'e non-archim\'edien. On \'enonce les principaux r\'esultats concernant ces objets. Les deux plus importants, \`a savoir les th\'eor\`emes 1.10 et 1.16, ainsi que quelques autres, ne seront pas d\'emontr\'es ici mais seulement \'enonc\'es comme assertions \`a prouver. Le th\'eor\`eme 1.10 sera d\'eduit dans l'article suivant des r\'esultats d'Arthur. La preuve du th\'eor\`eme 1.16 n\'ecessite un argument global, elle ne sera donn\'ee que beaucoup plus tard. 
 
 On utilise les notations introduites dans le premier article [I]. Consid\'erons un triplet $(G,\tilde{G},{\bf a})$ comme dans celui-ci. Le terme $G$ est un groupe r\'eductif connexe sur $F$, $\tilde{G}$ est un espace tordu sous $G$ et ${\bf a}$ est un \'el\'ement de $H^1(\Gamma_{F};Z(\hat{G}))$, qui d\'etermine un caract\`ere $\omega$ de $G(F)$. Dans la premi\`ere section, on commence par d\'efinir les objets de base, \`a savoir, pour un espace de Levi $\tilde{M}$ de $\tilde{G}$, les int\'egrales orbitales pond\'er\'ees $J_{\tilde{M}}^{\tilde{G}}(\gamma,\omega,f)$ et leurs avatars $\omega$-\'equivariants $I_{\tilde{M}}^{\tilde{G}}(\gamma,\omega,f)$. Pour cela, nous suivons bien s\^ur Arthur mais nous modifions un peu ses d\'efinitions. Expliquons cela en consid\'erons le cas $\tilde{G}=G=SO(7)$, ${\bf a}=1$, $\tilde{M}=M=GL(2)\times SO(3)$. Nous ne changeons rien aux d\'efinitions d'Arthur pour un \'el\'ement $\gamma\in M(F)$ qui est $G$-\'equisingulier, c'est-\`a-dire tel que $G_{\gamma}=M_{\gamma}$ (on note par exemple $G_{\gamma}$ la composante neutre du centralisateur de $\gamma$ dans $G$). Le changement concerne les \'el\'ements non-\'equisinguliers, par exemple l'\'el\'ement $\gamma=1$. Soit $(B,T)$ une paire de Borel de $G$ d\'efinie sur $F$, telle que $M$ soit standard pour cette paire. Notons $\mathfrak{g}$ et $\mathfrak{t}$ les alg\`ebres de Lie de $G$ et $T$ et notons $\Sigma(T)$ l'ensemble des racines de $T$ dans $\mathfrak{g}$.  On identifie $\mathfrak{t}(F)$ \`a $F^3$ de sorte que $\Sigma(T)$ s'identifie \`a l'ensemble $\{\pm\alpha_{i,\pm j}; 1\leq i<j\leq 3\}\cup\{\pm\alpha_{i};1\leq i\leq 3\}$ 
 de formes lin\'eaires sur $\mathfrak{t}(F)$, o\`u

 $\bullet$ $\alpha_{i,j}(x_{1},x_{2},x_{3})=x_{i}-x_{j}$;
 
 $\bullet$ $\alpha_{i,-j}(x_{1},x_{2},x_{3})=x_{i}+x_{j}$;
 
 $\bullet$ $\alpha_{i}(x_{1},x_{2},x_{3})=x_{i}$.
 
 Les racines dans $M$ sont $\pm\alpha_{1,2}$ et $\pm\alpha_{3}$. Introduisons l'alg\`ebre de Lie $\mathfrak{a}_{M}$ du centre de $M$. Alors $\mathfrak{a}_{M}(F)=\{H(x);x\in F\}$, o\`u $H(x)=(x,x,0)$.  Les racines ci-dessus se restreignent \`a ce sous-espace en $0$, $\pm\beta$ ou $\pm2\beta$, o\`u $\beta(H(x))=x$. On pose
 $$l_{\beta}(x)=\vert e^{\beta(H(x))}-e^{-\beta(H(x))}\vert _{F},\,\,l_{2\beta}(x)=\vert e^{2\beta(H(x))}-e^{-2\beta(H(x))}\vert _{F}.$$
  Pour $x\not=0$ mais assez voisin de $0$, l'\'el\'ement $exp(H(x))\in M(F)$ est $G$-\'equisingulier. Pour $f\in C_{c}^{\infty}(G(F))$, l'int\'egrale orbitale pond\'er\'ee $J_{M}^G(exp(H(x)),f)$ est bien d\'efinie, ainsi que $J_{G}^G(exp(H(x)),f)=I^G(exp(H(x)),f)$. Arthur montre qu'il existe un r\'eel $k_{M}^G$, n\'ecessairement unique, de sorte que, pour tout $f$, l'expression
 $$J_{M}^G(exp(H(x)),f)+k_{M}^Gl_{\beta}(x)I^G(exp(H(x)),f)$$
 ait une limite quand $x$ tend vers $0$. Il d\'efinit $J_{M}^G(1,f)$ comme \'etant cette limite. Pour des raisons de compatibilit\'e \`a l'induction, il nous semble pr\'ef\'erable de ne pas privil\'egier la racine indivisible $\beta$ mais de r\'epartir plut\^ot le coefficient $k_{M}^G$ sur les deux racines $\beta$ et $2\beta$. Pour cela, consid\'erons l'ensemble des \'el\'ements de $\Sigma(T)$ qui se restreignent en un multiple entier de $2\beta$. Ce sont $\pm\alpha_{1,2}$, $\pm\alpha_{1,-2}$ et $\pm\alpha_{3}$. C'est le syst\`eme de racines d'un sous-groupe $G_{2\beta}=SO(4)\times SO(3)$ de $G$ qui contient $M$. On a de m\^eme un nombre r\'eel $k_{M}^{G_{2\beta}}$. Consid\'erons l'expression
$$J_{M}^G(exp(H(x)),f)+\left((k_{M}^G-k_{M}^{G_{2\beta}})l _{\beta}(x)+k_{M}^{G_{2\beta}}l_{2\beta}(x)\right)I^G(exp(H(x)),f).$$
Elle a encore une limite quand $x$ tend vers $0$, qui est \'egale \`a la pr\'ec\'edente, plus $k_{M}^{G_{2\beta}}log(\vert 2\vert _{F})I^G(1,f)$. C'est cette limite que nous notons $J_{M}^G(1,f)$.

On s'aper\c{c}oit que le proc\'ed\'e admet diverses variantes. Par exemple, fixons un rationnel $b>0$. On peut remplacer dans la formule ci-dessus $ l_{\beta}(x)$ et $l_{2\beta}(x)$ par $ l_{\beta}(bx)$ et $l_{2\beta}(bx)$. Il y a encore une limite, \'egale \`a la pr\'ec\'edente plus $k_{M}^Glog(\vert b\vert _{F})I^G(1,f)$. On la note $J_{M}^G(1,B,f)$, $B$ d\'esignant la fonction constante sur $\Sigma(T)$ de valeur $b$.  Plus subtilement, consid\'erons la fonction $B$ sur $\Sigma(T)$ d\'efinie par $B(\pm \alpha_{i,j})=1$, $B(\pm \alpha_{i,-j})=1$, $B(\pm \alpha_{i})=1/2$. Cette fonction est proportionnelle au carr\'e de la longueur usuelle. Pour $\alpha\in \Sigma(T)$, $\frac{\alpha}{B(\alpha)}$ est encore une forme lin\'eaire sur $\mathfrak{t}(F)$. Les restrictions de ces formes \`a $\mathfrak{a}_{M}(F)$ sont encore $0$, $\pm \beta$, $\pm 2\beta$. Consid\'erons l'ensemble des $\alpha$ telles que la restriction de $\frac{\alpha}{B(\alpha)}$ soit un multiple entier de $2\beta$. Il est form\'e de $\pm \alpha_{1,2}$, $\pm\alpha_{1,-2}$, $\alpha_{1}$, $\alpha_{2}$, $\alpha_{3}$. C'est le syst\`eme de racines d'un groupe $G_{2\beta,B}=SO(5)\times SO(3)$. Ce n'est plus un sous-groupe de $G$, mais il contient encore $M$. L'expression
$$J_{M}^G(exp(H(x)),f)+\left((k_{M}^G-k_{M}^{G_{2\beta,B}})l _{\beta}(x)+k_{M}^{G_{2\beta,B}}l_{2\beta}(x)\right)I^G(exp(H(x)),f)$$
a encore une limite quand $x$ tend vers $0$. On la note $J_{M}^G(1,B,f)$.

Ainsi, pour certaines fonctions $B$ sur $\Sigma(T)$, on peut d\'efinir $J_{M}^G(1,B,f)$, ainsi que son avatar invariant $I_{M}^G(1,B,f)$. La consid\'eration de ces diverses d\'efinitions est utile pour notre propos. Expliquons pourquoi en revenant au cas g\'en\'eral. Consid\'erons une donn\'ee endoscopique ${\bf G}'=(G',{\cal G}',\tilde{s})$ de $(G,\tilde{G},{\bf a})$. Soient $\eta$ un \'el\'ement semi-simple de $\tilde{G}(F)$, $\epsilon$ un \'el\'ement semi-simple de $\tilde{G}'(F)$, supposons que ces deux \'el\'ements se correspondent par la correspondance endoscopique usuelle. Fixons des formes quasi-d\'eploy\'ees $G_{\eta}^*$ et $G_{\epsilon}^{_{'}*}$ de $G_{\eta}$ et $G'_{\epsilon}$ et des paires de Borel d\'efinies sur $F$ dans ces deux groupes, dont on note les tores $T_{\eta}$ et $T'_{\epsilon}$. On note $\Sigma(T_{\eta})$ et $\Sigma(T'_{\epsilon})$ les ensembles de racines de $T_{\eta}$ dans $G_{\eta}^*$ et de $T'_{\epsilon}$ dans $G^{_{'}*}_{\epsilon}$.  Il y a un isomorphisme naturel $\mathfrak{t}_{\eta}\simeq \mathfrak{t}'_{\epsilon}$. Mais il n'identifie pas $\Sigma(T'_{\epsilon})$ \`a un sous-ensemble de $\Sigma(T_{\eta})$. Par contre, il existe une fonction $B^{\tilde{G}}_{\epsilon}:\Sigma(T'_{\epsilon})\to {\mathbb Q}_{>0}$ telle que l'ensemble $\{\frac{\alpha}{B^{\tilde{G}}_{\epsilon}(\alpha)}; \alpha\in \Sigma(T'_{\epsilon})\}$ s'identifie \`a un sous-ensemble de $\Sigma(T_{\eta})$.  Soient $\tilde{M}$ et $\tilde{M}'$ des espaces de Levi de $\tilde{G}$ et $\tilde{G}'$ qui se correspondent, supposons $\eta\in \tilde{M}(F)$ et $\epsilon\in \tilde{M}'(F)$. Soit enfin $\gamma\in \tilde{M}(F)$ de partie semi-simple $\eta$.   "Stabiliser"  la distribution $f\mapsto I_{\tilde{M}}^{\tilde{G}}(\gamma,\omega,f)$ revient \`a \'etablir une relation entre celle-ci et  d'autres distributions vivant sur des espaces endoscopiques. Parmi ces derni\`eres, il y a en premi\`ere approximation les distributions $f'\mapsto I_{\tilde{M}'}^{\tilde{G}'}(\delta,f')$, o\`u $\delta$ est un \'el\'ement de $\tilde{M}'(F)$ de partie semi-simple $\epsilon$.  Il s'av\`ere qu'il  est plus pertinent d'utiliser la distribution $f'\mapsto I_{\tilde{M}'}^{\tilde{G}'}(\delta,B_{\epsilon}^{\tilde{G}},f')$.

 Dans la suite de la premi\`ere section, on d\'efinit les avatars stables et endoscopiques des int\'egrales orbitales pond\'er\'ees $\omega$-\'equivariantes. Le th\'eor\`eme 1.10, qui ne concerne que le cas o\`u $(G,\tilde{G},{\bf a})$ est quasi-d\'eploy\'e et \`a torsion int\'erieure, affirme que les avatars stables sont bel et bien stables. Le th\'eor\`eme 1.16 affirme l'\'egalit\'e des int\'egrales orbitales pond\'er\'ees $\omega$-\'equivariantes avec leurs avatars endoscopiques, c'est-\`a-dire une \'egalit\'e $I^{\tilde{G}}_{\tilde{M}}(\boldsymbol{\gamma},{\bf f})=I_{\tilde{M}}^{\tilde{G},{\cal E}}(\boldsymbol{\gamma},{\bf f})$ avec des notations proches de celles d'Arthur. Encore une fois, ces th\'eor\`emes ne sont ici qu'\'enonc\'es comme des assertions \`a prouver.
 
 Dans la deuxi\`eme section, on d\'eveloppe pour nos int\'egrales la th\'eorie des germes de Shalika. On pourrait esp\'erer que ceux-ci permettent de ramener les th\'eor\`emes 1.10 et 1.16 aux m\^emes th\'eor\`emes restreints aux distributions \`a support fortement r\'egulier dans $\tilde{G}(F)$ (c'est-\`a-dire de prouver que, si ces th\'eor\`emes sont v\'erifi\'es pour de telles distributions, ils sont vrais pour toute distribution). Cet espoir est vain, pour autant que je le sache, car on n'a pas assez de renseignements sur les germes. Ceux-ci permettent toutefois de prouver que les th\'eor\`emes, restreints aux distributions \`a support fortement r\'egulier dans $\tilde{G}(F)$, entra\^{\i}nent les m\^emes th\'eor\`emes pour les distributions \`a support seulement $\tilde{G}$-\'equisingulier. 
 
 Dans la troisi\`eme section, on \'etudie plus finement la d\'efinition des int\'egrales orbitales pond\'er\'ees $\omega$-\'equivariante. Par d\'efinition, une int\'egrale $I_{\tilde{M}}^{\tilde{G}}(\boldsymbol{\gamma},{\bf f})$  est limite de combinaisons lin\'eaires d'int\'egrales $I_{\tilde{L}}^{\tilde{G}}(a\boldsymbol{\gamma},{\bf f})$, o\`u $a\in A_{\tilde{M}}(F)$ est en position g\'en\'erale et tend vers $1$ et $\tilde{L}$ est un espace de Levi contenant $\tilde{M}$. On change l\'eg\`erement de point de vue et  on \'etudie plut\^ot le germe en $1$ de la fonction $a\mapsto I_{\tilde{M}}^{\tilde{G}}(a\boldsymbol{\gamma},{\bf f})$. On obtient un d\'eveloppement de cette fonction en termes de fonctions assez \'el\'ementaires de $a$. C'est ce d\'eveloppement qui, dans l'article suivant, nous permettra de ramener les th\'eor\`emes 1.10 et 1.16 aux m\^emes th\'eor\`emes restreints aux distributions \`a support fortement r\'egulier dans $\tilde{G}(F)$.

 Dans la derni\`ere section, on traite le cas non ramifi\'e, o\`u on \'etudie seulement les int\'egrales orbitales pond\'er\'ees non $\omega$-\'equivariantes de la fonction caract\'eristique d'un espace hypersp\'ecial. Le principal r\'esultat est que le lemme fondamental pond\'er\'e, qui est connu gr\^ace \`a Ngo Bao Chau pour les distributions \`a support fortement r\'egulier dans $\tilde{G}(F)$, est v\'erifi\'e pour toute distribution. Cela utilise les r\'esultats des sections 2 et 3.
 
 Pour conclure cette introduction, il faut dire que cet article doit tout aux travaux ant\'erieurs d'Arthur sur ce sujet et que, si on ne le cite pas \`a chaque ligne, c'est seulement pour ne pas lasser le lecteur.

\section{Int\'egrales orbitales pond\'er\'ees}

\bigskip

\subsection{Les hypoth\`eses}
Dans tout l'article, le corps de base $F$ est local, de caract\'eristique nulle et non-archim\'edien. On note $p$ la caract\'eristique r\'esiduelle de $F$.  
 On consid\`ere des triplets $(G,\tilde{G},{\bf a})$ comme dans [I]. Le terme $G$ est un  groupe r\'eductif connexe d\'efini sur $F$, $\tilde{G}$ est un espace tordu sur $G$, ${\bf a}$ est un \'el\'ement de $H^1(W _{F},Z(\hat{G}))$ qui d\'etermine  un caract\`ere $\omega$ de $G(F)$.  On suppose
 
$ \bullet$ $\tilde{G}(F)\not=\emptyset$;
 
$ \bullet$ l'automorphisme $\theta$ de $Z(G)$ est d'ordre fini;

$\bullet$ le caract\`ere $\omega$ est unitaire.

On aura \`a prouver des assertions concernant un tel triplet. On raisonne par r\'ecurrence sur l'entier $dim(G_{SC})$. 

Pour d\'emontrer une assertion concernant un triplet $(G,\tilde{G},{\bf a})$   quasi-d\'eploy\'e et \`a torsion int\'erieure, on suppose connues toutes les assertions concernant des triplets  $(G',\tilde{G}',{\bf a}')$ quasi-d\'eploy\'es et \`a torsion int\'erieure tels que $dim(G'_{SC})<dim(G_{SC})$. 

Pour d\'emontrer une assertion concernant un triplet $(G,\tilde{G},{\bf a})$  qui n'est pas quasi-d\'eploy\'e et \`a torsion int\'erieure, on suppose connues toutes les assertions concernant des triplets  $(G',\tilde{G}',{\bf a}')$ quasi-d\'eploy\'es et \`a torsion int\'erieure tels que $dim(G'_{SC})\leq dim(G_{SC})$. On suppose connues toutes les assertions concernant des triplets $(G',\tilde{G}',{\bf a}')$ quelconques tels que $dim(G'_{SC})<dim(G_{SC})$. 

Beaucoup d'assertions concernant un triplet $(G,\tilde{G},{\bf a})$ sont relatives \`a un espace de Levi $\tilde{M}$ de $\tilde{G}$. On supposera connues toutes les assertions concernant ce m\^eme triplet $(G,\tilde{G},{\bf a})$, relatives \`a un espace de Levi $\tilde{L}\in {\cal L}(\tilde{M})$ tel que $\tilde{L}\not=\tilde{M}$.

\bigskip

 \subsection{D\'efinition des int\'egrales pond\'er\'ees   d'apr\`es Arthur}

Soit $(G,\tilde{G},{\bf a})$ un triplet comme en 1.1. 
  Soit $\tilde{M}$ un espace de Levi de $\tilde{G}$. On doit fixer une mesure sur ${\cal A}_{\tilde{M}}^{\tilde{G}}$. Pour ce faire, introduisons $\underline{la}$ paire de Borel \'epingl\'ee ${\cal E}^*=(B^*,T^*,(E^*_{\alpha})_{\alpha\in \Delta})$ de $G$, cf. [I] 1.2. Elle est munie d'une action $\sigma\mapsto \sigma_{G^*}$ du groupe de Galois $\Gamma_{F}$ et d'un automorphisme $\theta^*$. Le groupe de Weyl $W$ agit sur $T^*$. On fixe une  forme quadratique d\'efinie positive sur $X_{*}(T^*)\otimes {\mathbb R}$, invariante par les actions de $W$ et  de $\Gamma_{F}$ et par  $\theta^*$. En fixant $\tilde{P}\in {\cal P}(\tilde{M})$ et en identifiant ${\cal E}^*$ \`a une paire de Borel \'epingl\'ee contenue dans $\tilde{ P}$ et dont le tore est contenu dans $\tilde{M}$, ${\cal A}_{\tilde{M}}^{\tilde{G}}$ s'identifie \`a un sous-espace de  $X_{*}(T^*)\otimes {\mathbb R}$. Par restriction, on obtient une forme quadratique d\'efinie positive sur ${\cal A}_{\tilde{M}}^{\tilde{G}}$. Elle ne d\'epend pas des choix. De cette forme se d\'eduit la mesure cherch\'ee. On fixe un sous-groupe compact maximal sp\'ecial $K$ de $G(F)$ en bonne position relativement \`a $M$. On fixe aussi une mesure de Haar sur $G(F)$.  
  
  Introduisons  la notion d'\'el\'ement $\tilde{G}$-\'equisingulier de $\tilde{M}$. Soit $\gamma\in \tilde{M}$, notons $\eta$ sa partie semi-simple. On a
  
  (1) les \'egalit\'es $M_{\gamma}=G_{\gamma}$ et $M_{\eta}=G_{\eta}$ sont \'equivalentes.
  
  Preuve. Ces \'egalit\'es sont  \'equivalentes aux inclusions $G_{\gamma}\subset M$, resp. $G_{\eta}\subset M$. Ecrivons $\gamma=u\eta$, o\`u $u$ est un \'el\'ement unipotent de $M_{\eta}$. On a les \'egalit\'es $G_{\gamma}=(G_{\eta})_{u}$ et $M_{\gamma}=(M_{\eta})_{u}$. Si $G_{\eta}\subset M$, on a $G_{\gamma}\subset G_{\eta}\subset M$, donc  $G_{\gamma}=M_{\gamma}$. Inversement, supposons $M_{\gamma}=G_{\gamma}$. Posons $H=G_{\eta}$ et $L=M_{\eta}$. Alors $L$ est un Levi de $H$ et $u$ est un \'el\'ement unipotent de $L$ tel que $H_{u}\subset L$. On veut en d\'eduire que $L=H$. Mais soit $Q\in {\cal P}(L)$. Si $L\not=H$, le radical unipotent $U_{Q}$ est non trivial. L'automorphisme $ad_{u}$ agit de fa\c{c}on unipotente sur ce radical, ce qui implique que son ensemble de points fixes dans $U_{Q}$ est non trivial. Cet ensemble est inclus dans $H_{u}$, ce qui contredit l'inclusion $H_{u}\subset L$. $\square$
  
  On appelle \'el\'ement $\tilde{G}$-\'equisingulier  de $\tilde{M}$ un \'el\'ement $\gamma$ v\'erifiant les \'egalit\'es (1). 
   
  Soit $\gamma\in \tilde{M}(F)$. Fixons une mesure de Haar sur le groupe $M_{\gamma}(F)$. Arthur d\'efinit  dans [A1] une distribution $f\mapsto J_{\tilde{M}}^{\tilde{G}}(\gamma,\omega,f)$ sur $C_{c}^{\infty}(\tilde{G}(F))$. On va rappeler  sa d\'efinition.  Nous la modifierons dans le paragraphe suivant, c'est pourquoi nous affecterons des exposants $Art$ \`a certains objets d\'efinis par Arthur.
   
   Si $\omega$ n'est pas trivial sur $M_{\gamma}(F)$, on pose $J_{\tilde{M}}^{\tilde{G}}(\gamma,\omega,f)=0$ pour tout $f$. On suppose d\'esormais $\omega$ trivial sur $M_{\gamma}(F)$.

 Premier cas: on suppose que $\gamma$ est $\tilde{G}$-\'equisingulier.   Arthur d\'efinit pour tout $g\in G(F)$ une $(\tilde{G},\tilde{M})$-famille $(v_{\tilde{P}}(g;\lambda))_{\tilde{P}\in {\cal P}(\tilde{M})}$ ($\lambda$ est une variable dans $i{\cal A}_{\tilde{M}}^*$). Comme de toute $(\tilde{G},\tilde{M})$-famille, il s'en d\'eduit une fonction $v_{\tilde{M}}^{\tilde{G}}(g;\lambda)$. On pose $v_{\tilde{M}}^{\tilde{G}}(g)=v_{\tilde{M}}^{\tilde{G}}(g;0)$.  La fonction $g\mapsto v^{\tilde{G}}_{\tilde{M}}(g)$ est la fonction "poids". Pour $f\in C_{c}^{\infty}(\tilde{G}(F))$, on pose
 $$J_{\tilde{M}}^{\tilde{G}}(\gamma,\omega,f)=D^{\tilde{G}}(\gamma)^{1/2}\int_{M_{\gamma}(F)\backslash G(F)}\omega(g)f(g^{-1}\gamma g)v_{\tilde{M}}^{\tilde{G}}(g)\,dg.$$
 
 Cas g\'en\'eral. On \'ecrit $\gamma=u\eta$, o\`u $\eta$ est la partie semi-simple de $\gamma$ et $u$ est un unipotent dans $M_{\eta}(F)$. Notons $\Sigma(A_{\tilde{M}})$ l'ensemble des racines de $A_{\tilde{M}}$ dans $G$ (toutes les racines, pas seulement les indivisibles). Notons  $\Sigma_{ind}(A_{M_{\eta}})$ l'ensemble des racines indivisibles de $A_{M_{\eta}}$ dans $G_{\eta}$. La restriction d\'efinit une application naturelle $\beta\mapsto \beta_{\tilde{M}}$ de $\Sigma_{ind}(A_{M_{\eta}})$ dans $\Sigma(A_{\tilde{M}})\cup\{0\}$. Pour tout $\beta\in \Sigma_{ind}(A_{M_{\eta}})$, Arthur d\'efinit un r\'eel $\rho^{Art}(\beta,u)$ et une "coracine" $\check{\beta} \in {\cal A}_{M_{\eta}}$.   Pour  $\alpha\in \Sigma(A_{\tilde{M}})$, pour $a\in A_{\tilde{M}}(F)$ en position g\'en\'erale et pour $\lambda\in i{\cal A}_{\tilde{M},{\mathbb C}}$, posons
 $$r_{\alpha}^{Art}(\gamma,a;\lambda)=\prod_{\beta\in \Sigma_{ind}(A_{M_{\eta}}); \beta_{\tilde{M}}=\alpha}\vert \alpha(a)-\alpha(a)^{-1}\vert _{F}^{<\lambda,\rho^{Art}(\beta,u)\check{\beta}_{\tilde{M}}>}.$$
  On d\'efinit ensuite une $(\tilde{G},\tilde{M})$-famille $(r_{\tilde{P}}^{Art}(\gamma,a;\lambda))_{\tilde{P}\in {\cal P}(\tilde{M})}$ par
 $$r_{\tilde{P}}^{Art}(\gamma,a;\lambda)=\prod_{\alpha>_{P}0}r_{\alpha}^{Art}(\gamma,a;\lambda/2)$$
 pour $\lambda\in i{\cal A}^*_{\tilde{M}}$, o\`u $\alpha$ parcourt les \'el\'ements de $\Sigma(A_{\tilde{M}})$ qui sont "positifs" pour $P$.  On d\'eduit de cette $(\tilde{G},\tilde{M})$-famille une fonction $r_{\tilde{M}}^{\tilde{G},Art}(\gamma,a;\lambda)$ et on pose $r_{\tilde{M}}^{\tilde{G},Art}(\gamma,a)=r_{\tilde{M}}^{\tilde{G},Art}(\gamma,a;0)$.
 
   Pour $f\in C_{c}^{\infty}(\tilde{G}(F))$, consid\'erons la fonction
 $$(2) \qquad a\mapsto \sum_{\tilde{L}\in {\cal L}(\tilde{M})}r_{\tilde{M}}^{\tilde{L},Art}(\gamma,a)J_{\tilde{L}}^{\tilde{G}}(a\gamma,\omega,f).$$
 Pour $a$ en position g\'en\'erale, elle est bien d\'efinie: on a $G_{a\gamma}=M_{a\gamma}=M_{\gamma}$ et les int\'egrales orbitales pond\'er\'ees $J_{\tilde{L}}^{\tilde{G}}(a\gamma,\omega,f)$ sont d\'efinies d'apr\`es le premier cas ci-dessus. Arthur montre que la fonction (2) a une limite quand $a$ tend vers $1$ (Arthur traite le cas $\omega=1$ mais sa preuve s'\'etend sans changement au cas g\'en\'eral).  Notons $J_{\tilde{M}}^{\tilde{G},Art}(\gamma,\omega,f)$ la limite de la fonction (2). C'est l'int\'egrale orbitale pond\'er\'ee telle que d\'efinie par Arthur. Notons que, dans le cas o\`u $M_{\gamma}=G_{\gamma}$, on retrouve celle donn\'ee plus haut.
 
 {\bf Remarque.}  On v\'erifie que $J_{\tilde{M}}^{\tilde{G},Art}(\gamma,\omega,f)=0$ si $\omega$ n'est pas trivial sur $Z_{M}(\gamma;F)$ tout entier.  
 
 \bigskip
 
 On aura besoin d'un r\'esultat un peu plus pr\'ecis. On d\'efinit une distance $d$ au voisinage de $1$ dans $A_{\tilde{M}}(F)$ de la fa\c{c}on suivante. On fixe une norme $ \vert .\vert $ sur l'alg\`ebre de Lie $\mathfrak{a}_{\tilde{M}}(F)$. On fixe des voisinages $U$ de $1$ dans $A_{\tilde{M}}(F)$ et $\mathfrak{u}$ de $0$ dans $\mathfrak{a}_{\tilde{M}}(F)$ tels que l'exponentielle soit bijective de $\mathfrak{u}$ dans $U$. Pour $a\in U$, on \'ecrit $a=exp(H)$, avec $h\in \mathfrak{u}$, et on pose $d(a)=\vert H\vert $. On a alors
 
 (3) il existe $r>0$ tel que, pour tout $\gamma\in \tilde{M}(F)$, tout $f\in C_{c}^{\infty}(\tilde{G}(F))$, il existe $C>0$ de sorte que
 $$\vert J_{\tilde{M}}^{\tilde{G},Art}(\gamma,\omega,f)-\sum_{\tilde{L}\in {\cal L}(\tilde{M})}r_{\tilde{M}}^{\tilde{L},Art}(\gamma,a)J_{\tilde{L}}^{\tilde{G}}(a\gamma,\omega,f)\vert \leq Cd(a)^{r}$$
 pour tout $a\in A_{\tilde{M}}(F)$ en position g\'en\'erale et assez proche de $1$.
 
 Preuve. Un examen attentif de la preuve d'Arthur montre qu'il suffit d'am\'eliorer son lemme  6.1 de [A1]. Reprenons les notations de ce lemme dans la situation simplifi\'ee qui nous concerne: l'ensemble de places $S$ est r\'eduit \`a un \'el\'ement, on a $F_{S}=F$ et les $v$ disparaissent. On consid\`ere une famille d'\'el\'ements de l'espace ${\cal P}^+(\Omega)$ d\'ependant d'un param\`etre $a$ parcourant un voisinage de $1$ dans $A_{\tilde{M}}(F)$. On note $p[a]=\oplus_{\omega\in \Omega}p[a]_{\omega}$ l'\'el\'ement de cette famille param\'etr\'e par $a$. Le terme $p[a]_{\omega}$ est un polyn\^ome sur $F^d$ \`a valeurs dans un espace $V_{\omega}$ de dimension finie sur $F$, muni d'une norme $\vert\vert.\vert\vert$. On suppose que l'application $(a,x)\to p[a]_{\omega}(x)$  est la restriction (au voisinage de $a=1$) d'un polyn\^ome d\'efini sur $A_{\tilde{M}}(F)\times F^d$.  Pour $x\in F^d$, on pose
 $$\lambda_{p[a]}(x)=\prod_{\omega\in \Omega}\vert log(\vert\vert p[a]_{\omega}(x)\vert \vert)\vert .$$
  Arthur montre que, pour tout $\phi\in C_{c}^{\infty}(F^d )$, l'application
 $$a\mapsto\lambda_{p[a]}(\phi)= \int_{{\cal O}}\phi(x)\lambda_{p[a]}(x)\,dx$$
 est continue. Pour obtenir (3), on doit montrer qu'il existe $r>0$   et $C>0$ de sorte que
 $$\vert \lambda_{p[a]}(\phi)-\lambda_{p[1]}(\phi)\vert \leq Cd(a)^r.$$
 On veut de plus que $r$ ne d\'epende pas de $\phi$ et, si on fixe un entier $D$ et que l'on impose que tous les $p[a]_{\omega}$ sont de degr\'e au plus $D$, que $r$ ne d\'epende pas non plus de la famille de polyn\^omes. Il suffit pour cela de reprendre la fin de la preuve du lemme 6.1.  On pose $p^0=p[1]$. En choisissant un param\`etre auxiliaire $\epsilon$, Arthur montre que $\vert \lambda_{p[a]}(\phi)-\lambda_{p^0}(\phi)\vert $ est major\'e par la somme de trois expressions (7.2), (7.3) et (7.4). La relation (7.1) nous dit que le terme (7.3) est major\'e par $C_{1}\epsilon^{r_{1}}$, o\`u $C_{1}$ et $r_{1}$ v\'erifient les conditions requises. Le terme (7.2) v\'erifie une majoration analogue pourvu que l'on ait l'inclusion $\Gamma(p^0,\epsilon)\subset \Gamma(p[a],2\epsilon)$. Rappelons que  $\Omega$ est un ensemble fini, que $\Gamma$ est un sous-ensemble compact de $F^d$ et que $\Gamma(p[a],\epsilon)$ est la r\'eunion sur les $\omega\in \Omega$ des ensembles des $x\in \Gamma$ tels que $\vert \vert p[a]_{\omega}(x)\vert \vert <\epsilon$. Puisque les $p[a]_{\omega}$ sont polynomiaux en $a$, il existe $C_{2}$ tel que 
$$ \vert \vert \vert p[a]_{\omega}(x)\vert \vert -\vert \vert p^0_{\omega}(x)\vert \vert \vert <C_{2}d(a)$$
pour tout $a$ voisin de $1$, tout $x\in \Gamma$ et tout $\omega\in \Omega$. L'inclusion $\Gamma(p^0,\epsilon)\subset \Gamma(p[a],2\epsilon)$ est v\'erifi\'ee pourvu que $C_{2}d(a)<\epsilon$. Imposons plut\^ot $2C_{2}d(a)<\epsilon$. Le m\^eme calcul montre que l'on a l'inclusion en sens inverse $\Gamma(p[a],\epsilon/2)\subset \Gamma(p^0,\epsilon)$. Le terme (7.4) est de la forme
$$C_{3}\int_{\Gamma-\Gamma(p^0,\epsilon)}\vert \lambda_{p[a]}(x)-\lambda_{p^0}(x)\vert \,dx.$$
Sur le domaine d'int\'egration, on a $\vert \vert p^0_{\omega}(x)\vert\vert >\epsilon$ pour tout $\omega$ et, d'apr\`es l'inclusion ci-dessus, on a aussi $\vert\vert p[a]_{\omega}(x)\vert\vert >\epsilon/2$.  Ecrivons $\Omega=\{\omega_{1},...,\omega_{\ell}\}$. On peut \'ecrire
$$\lambda_{p[a]}(x)-\lambda_{p^0}(x)=\sum_{k=1,...,\ell}\left(\prod_{i=1,...,k-1}\vert log(\vert \vert p[a]_{\omega_{i}}(x)\vert \vert )\vert \right)$$
$$\left(\vert log(\vert \vert p[a]_{\omega_{i}}(x)\vert \vert )\vert -\vert log(\vert \vert p[a]_{\omega_{i}}(x)\vert \vert )\vert\right)\left(\prod_{j=k+1,...,\ell} \vert log(\vert \vert p^0_{\omega_{i}}(x)\vert \vert )\vert \right).$$
On  d\'eduit des in\'egalit\'es pr\'ec\'edentes que $\vert \lambda_{p[a]}(x)-\lambda_{p^0}(x)\vert$ est essentiellement born\'e par la somme sur les $\omega$ de
$$\vert log(\epsilon/2)\vert ^{\vert \Omega\vert -1}\vert log(\vert \vert p[a]_{\omega}(x)\vert\vert)-log(\vert\vert p^0_{\omega}(x)\vert\vert)\vert .$$
Le dernier terme est \'egal \`a la valeur absolue de 
$$log( \frac{\vert\vert p[a]_{\omega}(x)\vert\vert}{\vert\vert p^0_{\omega}(x)\vert\vert} ).$$
On peut  \'ecrire $p[a]_{\omega}(x)=p^0_{\omega}(x)+q(a,x)$, o\`u $q(a,x)$ est un polyn\^ome en $a$ et $x$ qui est nul en $a=1$. On a une majoration $\vert \vert q(a,x)\vert\vert \leq C_{4}d(a)$ pour tout $x\in \Gamma$. Puisque $\vert\vert p^0_{\omega}(x)\vert\vert >\epsilon$, on obtient 
$$\vert \frac{\vert\vert p[a]_{\omega}(x)\vert\vert}{\vert \vert p^0_{\omega}(x)\vert\vert}-1\vert \leq C_{4}d(a)\epsilon^{-1}.$$
Renfor\c{c}ons la minoration impos\'ee \`a $\epsilon$ en supposant $d(a)^{1/2}\leq\epsilon$ (c'est plus fort que $2C_{2}d(a)<\epsilon$ pour $a$ proche de $1$). Alors 
$$\vert \frac{\vert\vert p[a]_{\omega}(x)\vert \vert}{\vert\vert p^0_{\omega}(x)\vert\vert}-1\vert \leq C_{4}d(a)^{1/2}$$
d'o\`u 
$$\vert log( \frac{\vert\vert p[a]_{\omega}(x)\vert\vert}{\vert\vert p^0_{\omega}(x)\vert\vert})\vert \leq C_{5}d(a)^{1/2}$$
pour une constante $C_{5}$ convenable. Alors le terme (7.4) est essentiellement major\'e par
$$C_{3}C_{5}\vert \Omega\vert \vert log(\epsilon/2)\vert ^{\vert \Omega\vert -1}d(a)^{1/2}.$$
On fixe maintenant $\epsilon=d(a)^{1/2}$. Le terme ci-dessus est major\'e par
$$C_{6}d(a)^{r_{2}}$$
pour tout r\'eel $r_{2}<1/2$ et pour une constante $C_{6}$ convenable. Les majorations des termes (7.2) et (7.3) deviennent de la forme $C_{1}d(a)^{r_{1}/2}$. En prenant pour $r$ l'inf de $r_{2}$ et $r_{1}/2$, on a obtenu la majoration cherch\'ee. $\square$

 \bigskip
 
 \subsection{Propri\'et\'es des termes $\rho^{Art}(\beta,u)\check{\beta}$}
 On suppose dans ce paragraphe et le suivant qu'il n'y a pas de torsion, c'est-\`a-dire $\tilde{G}=G$. On consid\`ere un Levi $M$ de $G$ et un \'el\'ement unipotent $u\in M(F)$. Comme on l'a rappel\'e, Arthur d\'efinit pour toute racine $\beta\in \Sigma_{ind}(A_{M})$ un r\'eel $\rho^{Art}(\beta,u)$ et une coracine $\check{\beta}$. Fixons une  paire de Borel $(B,T)$ de $G$ telle que $M$ soit standard pour cette paire. Fixons une extension finie $F'$ de $F$ telle que $(B,T)$ soit d\'efinie sur $F'$ et que $G$ soit d\'eploy\'e sur $F'$. Pla\c{c}ons-nous sur le corps de base $F'$.  Le tore $Z(M)^0$ est alors l'analogue de $A_{M}$. Pour tout  $\beta'\in \Sigma_{ind}(Z(M)^0)$, on d\'efinit le r\'eel $\rho^{Art}(\beta',u)$ et une coracine $\check{\beta}'\in X_{*}(Z(M)^0)\otimes_{{\mathbb Z}}{\mathbb R}$.  De l'inclusion $A_{M}\subset Z(M)^0$ se d\'eduisent des applications de restriction
 $$\begin{array}{ccc}\Sigma_{ind}(Z(M)^0)&\to& \Sigma(A_{M})\\ \beta'&\mapsto& \beta'_{A_{M}}\\ \end{array},\,\, \begin{array}{ccc}X_{*}(Z(M)^0)\otimes_{{\mathbb Z}}{\mathbb R}&\to& {\cal A}_{M}\\ H&\mapsto&H_{{\cal A}_{M}}\\ \end{array}.$$
 Remarquons qu'une racine indivisible de $Z(M)^0$ ne se restreint pas forc\'ement en une racine indivisible.
 
 Pour $\beta\in \Sigma_{ind}(A_{M})$, on a l'\'egalit\'e
 
 (1) $\rho^{Art}(\beta,u)\check{\beta}=\sum_{n\geq1}\sum_{\beta';\beta'_{A_{M}} =n\beta}\rho^{Art}(\beta',u)\check{\beta}'_{{\cal A}_{M}}$.

Preuve. Soit $P\in {\cal P}(M)$ et soit $\omega$ un poids de $A_{M}$ qui est dominant pour $P$ (c'est la notation d'Arthur; il ne s'agit pas de notre caract\`ere $\omega$ que nous oublions pour un temps). Notons ${\cal U}$ l'orbite g\'eom\'etrique de $u$ dans $M$. Arthur d\'efinit une fonction $W_{\omega}(a,\pi)$ sur $A_{M}\times {\cal U}U_{P}$, \`a valeurs dans un espace de dimension finie sur $\bar{F}$(cf. [A1] 3.8; on consid\`ere ici le cas $P=\bar{P}_{1}$ avec les notations de cette r\'ef\'erence).  Le groupe $M$ agit sur cet espace et la fonction est \'equivariante pour l'action de $M$ par conjugaison sur ${\cal U}U_{P}$ et cette action sur l'espace d'arriv\'ee. Arthur montre que cette fonction est  polynomiale et n'est pas identiquement nulle en $a=1$ ([A1] corollaire 4.3). Elle est donc non nulle sur $\{1\}\times {\cal O}$, o\`u ${\cal O}$ est un ouvert de Zariski de ${\cal U}U_{P}$, qui est dense et invariant par conjugaison par $M$. En se pla\c{c}ant sur $F'$, on a de m\^eme une fonction $W'_{\omega}(a',\pi)$ sur $Z(M)^0\times {\cal U}U_{P}$, qui est  polynomiale et est non nulle sur $\{1\}\times {\cal O}'$, o\`u ${\cal O}'$ est un ouvert de Zariski de ${\cal U}U_{P}$, qui est dense et invariant par conjugaison par $M$.  Sa restriction \`a $A_{M}\times {\cal U}U_{P}$ v\'erifie donc la m\^eme propri\'et\'e. Or il r\'esulte de la d\'efinition (3.8) de [A1] que, pour $(a,\pi)\in A_{M}\times {\cal U}U_{P}$, on a l'\'egalit\'e
$$W_{\omega}(a,\pi)=W'_{\omega}(a,\pi)\frac{Q(a)}{Q'(a)},$$
o\`u
$$Q(a)=\prod_{\beta\in \Sigma_{ind}(A_{M}),\beta>_{P}0}(\beta(a)-\beta(a)^{-1})^{\rho^{Art}(\beta,u)<\omega,\check{\beta}>},$$
$$Q'(a)=\prod_{\beta'\in \Sigma_{ind}(Z(M)^0),\beta'>_{P}0}(\beta'(a)-\beta'(a)^{-1})^{\rho^{Art}(\beta',u)<\omega,\check{\beta}'>}.$$
Les propri\'et\'es des deux fonctions $W_{\omega}$ et $W'_{\omega}$ entra\^{\i}nent que  la fraction rationnelle $\frac{Q(a)}{Q'(a)}$ n'a ni z\'ero, ni p\^ole en $a=1$. 
 Remarquons que l'on peut r\'ecrire
$$Q'(a)=\prod_{\beta\in \Sigma_{ind}(A_{M}),\beta>_{P}0}\prod_{n\geq1}\prod_{\beta';\beta'_{A_{M}}= n\beta}(\beta'(a)-\beta'(a)^{-1})^{\rho^{Art}(\beta',u)<\omega,\check{\beta}'>}.$$
Si $\beta'_{A_{M}}= n\beta$, la fonction
$$\frac{\beta'(a)-\beta'(a)^{-1}}{\beta(a)-\beta(a)^{-1}}$$
n'a ni z\'ero, ni p\^ole en $a=1$. Posons
$$Q''(a)=\prod_{\beta\in \Sigma_{ind}(A_{M}),\beta>_{P}0}(\beta(a)-\beta(a)^{-1})^{<\omega,X(\beta)>},$$
o\`u $X(\beta)$ est le membre de droite de (1).  Alors $\frac{Q'(a)}{Q''(a)}$ n'a ni z\'ero, ni p\^ole en $a=1$. Donc $\frac{Q(a)}{Q''(a)}$ a la m\^eme propri\'et\'e. Cela \'equivaut \`a $\rho^{Art}(\beta,u)<\omega,\check{\beta}>=<\omega,X(\beta)>$ pour tout $\beta$. Cela \'etant vrai pour tout   poids dominant $\omega$, cela entra\^{\i}ne l'\'egalit\'e (1) cherch\'ee. $\square$

Soit $L\in {\cal L}(M)$. On sait d\'efinir la classe de conjugaison (g\'eom\'etrique) induite de $M$ \`a $L$ de la classe de conjugaison de $u$. Soit $R=MU_{R}\in {\cal P}^L(M)$. Alors l'intersection de cette classe de conjugaison et de $uU_{R}$ est Zariski-dense dans cet ensemble. Soit $u'$ dans cette classe.  On peut d\'efinir des termes $\rho^{Art}(\beta',u')$ et $\check{\beta}'\in {\cal A}_{L}$ pour $\beta'\in \Sigma_{ind}(A_{L})$. On a des applications de restriction
$$\begin{array}{ccc}\Sigma_{ind}(A_{M})&\to& \Sigma(A_{L})\cup\{0\}\\ \beta&\mapsto& \beta_{L}\ \end{array},\,\, \begin{array}{ccc} {\cal A}_{M}&\to& {\cal A}_{L}\\ H&\mapsto&H_{L}\\ \end{array}.$$
Remarquons que la restriction d'une racine indivisible de $A_{M}$  peut \^etre nulle ou divisible. Pour $\beta'\in \Sigma_{ind}(A_{L})$, on a l'\'egalit\'e

(2) $\rho^{Art}(\beta',u')\check{\beta}'=\sum_{n\geq1}\sum_{\beta\in \Sigma_{ind}(A_{M});\beta_{L} =n\beta'}\rho^{Art}(\beta,u)\check{\beta}_{L}$.

Preuve. On fixe un sous-groupe parabolique $P'\in {\cal P}(L)$ et un poids $\omega$ de $A_{L}$ qui est dominant pour $P'$. On  fixe un  sous-groupe parabolique $P\in {\cal P}(M)$ contenu dans $P'$. 
 On d\'efinit comme ci-dessus des fonctions $W_{\omega}$ sur $A_{M}\times {\cal U}U_{P}$ et $W'_{\omega}$ sur $A_{L}\times {\cal U}'U_{P'}$, o\`u ${\cal U}'$ est l'orbite g\'eom\'etrique de $u'$.   Comme plus haut, elles sont polynomiales et  non nulles  sur $\{1\}\times {\cal O}$, respectivement sur $\{1\}\times {\cal O}'$, o\`u ${\cal O}$ est  un ouvert de Zariski de ${\cal U}U_{P}$ qui est dense et invariant par conjugaison par $M$ et ${\cal O}'$ est un ouvert de Zariski de ${\cal U}'U_{P'}$, qui est dense et invariant par conjugaison par $L$. Il r\'esulte de la d\'efinition de $u'$ que ${\cal O}\cap{\cal O}'\not=\emptyset$.
Donc les deux fonctions  sont toutes  deux non nulles sur $\{1\}\times ({\cal O}\cap {\cal O}')$.  Pour $a\in A_{M}$, d\'efinissons $Q(a)$ comme ci-dessus et, pour $a'\in A_{L}$, posons
$$Q'(a')=\prod_{\beta'\in \Sigma_{ind}(A_{L}),\beta'>_{P'}0}(\beta'(a')-\beta'(a')^{-1})^{\rho^{Art}(\beta',u')<\omega,\check{\beta}'>}.$$
Le m\^eme argument que plus haut montre que la fraction rationnelle $\frac{Q(a')}{Q'(a')}$ sur $A_{L}$ n'a ni z\'ero, ni p\^ole en $a'=1$. Remarquons que l'on peut supprimer de la d\'efinition de $Q$ les $\beta$ dont la restriction \`a $A_{L}$ est nulle: pour ceux-l\`a, on a $<\omega,\check{\beta}>=0$. Comme plus haut, si $\beta$ se restreint en $n\beta'$, la fonction
$$\frac{\beta(a')-\beta(a')^{-1}}{\beta'(a')-\beta'(a')^{-1}}$$
n'a ni z\'ero, ni p\^ole en $a'=1$.  En notant $Y(\beta')$ la diff\'erence entre le membre de gauche de (2) et celui de droite, on obtient alors que la fonction $\frac{Q(a')}{Q'(a')}$ a la m\^eme singularit\'e en $a=1$ que la fonction
$$\prod_{\beta'\in \Sigma_{ind}(A_{L}),\beta'>_{P}0}(\beta'(a')-\beta'(a')^{-1})^{<\omega,Y(\beta')>}.$$
Donc ce produit n'a lui-m\^eme ni z\'ero, ni p\^ole en $a'=1$. Cela entra\^{\i}ne $Y(\beta')=0$ pour tout $\beta'$. $\square$

 \bigskip
 
 \subsection{D\'efinition d'un nouveau terme $\rho(\beta,u)$}
On conserve la situation du d\'ebut du paragraphe pr\'ec\'edent et on fixe une paire de Borel $(B,T)$ et une extension $F'$ comme alors. On va d\'efinir un \'el\'ement $\rho(\beta,u)\in  X_{*}(Z(M)^0)\otimes_{{\mathbb Z}}{\mathbb R}$, ou plus pr\'ecis\'ement $\rho^G(\beta,u)$,  pour toute racine $\beta\in \Sigma(Z(M)^0)$ et non plus seulement pour les racines indivisibles. La d\'efinition se fait par r\'ecurrence sur la dimension de $G_{SC}$. Soit $\beta\in \Sigma(Z(M)^0)$. On introduit le sous-groupe $G_{\beta}$ de $G$  engendr\'e par $M$ et les sous-groupes radiciels associ\'es aux racines  $n\beta$ pour $n\in {\mathbb Z}$. Si $dim(G_{\beta,SC})<dim(G_{SC})$, le terme $\rho^{G_{\beta}}(\beta,u)$ relatif \`a $G_{\beta}$ est d\'ej\`a d\'efini et on pose $\rho^G(\beta,u)=\rho^{G_{\beta}}(\beta,u)$. Supposons $dim(G_{\beta,SC})=dim(G_{SC})$. Dans ce cas, $M$ est un Levi maximal de $G$ et $\beta$ est une racine indivisible (une telle racine est unique au signe pr\`es). On pose
$$\rho^G(\beta,u)=\rho^{Art}(\beta,u)\check{\beta}-\sum_{n>1}\rho^G(n\beta,u),$$
avec la convention $\rho^G(n\beta,u)=0$ si $n\beta$ n'est pas une racine. 

On redescend \`a la situation d\'efinie sur $F$ de la fa\c{c}on suivante. On a encore des applications de restriction
 $$\begin{array}{ccc}\Sigma(Z(M)^0)&\to& \Sigma(A_{M})\\ \beta'&\mapsto& \beta'_{A_{M}}\\ \end{array},\,\, \begin{array}{ccc}X_{*}(Z(M)^0)\otimes_{{\mathbb Z}}{\mathbb R}&\to& {\cal A}_{M}\\ H&\mapsto&H_{{\cal A}_{M}}\\ \end{array}.$$
 Pour $\beta\in \Sigma(A_{M})$, on pose
 $$(1) \qquad \rho^G(\beta,u)=\sum_{\beta';\beta'_{A_{M}}= \beta}\rho^G(\beta',u)_{{\cal A}_{M}}.$$ 
 
On a, avec la m\^eme convention que ci-dessus,

(2) pour tout $\beta\in \Sigma_{ind}(A_{M})$,
$$\rho^{Art}(\beta,u)\check{\beta}=\sum_{n\geq1}\rho^G(n\beta,u).$$

Preuve. D'apr\`es 1.3(1), le membre de gauche est
$$\sum_{n\geq1}\sum_{\beta'\in \Sigma_{ind}(Z(M)^0), \beta'_{A_{M}}= n\beta}\rho^{Art}(\beta',u)\check{\beta}'_{{\cal A}_{M}}.$$
D'apr\`es (1) ci-dessus, le membre de droite est
$$\sum_{n\geq1}\sum_{\beta'\in \Sigma(Z(M)^0),\beta'_{A_{M}}= n\beta}\rho^{G}(\beta',u)_{{\cal A}_{M}}.$$
Les racines $\beta'$  qui interviennent dans la deuxi\`eme expression sont exactement les multiples  positifs de racines indivisibles intervenant dans la premi\`ere. Cela nous ram\`ene \`a prouver l'analogue suivant de l'assertion (2): pour $\beta'\in \Sigma_{ind}(Z(M)^0)$, on a l'\'egalit\'e
$$\rho^{Art}(\beta',u)\check{\beta}'=\sum_{n\geq1}\rho^G(n\beta',u).$$
Introduisons le groupe $G_{\beta'}$ comme plus haut. D'apr\`es les d\'efinitions, le membre de droite est l'analogue de $\rho^{Art}(\beta',u)\check{\beta}'$ quand on remplace le groupe ambiant $G$ par $G_{\beta'}$. Mais, $\beta'$ \'etant indivisible, le groupe $G_{\beta'}$ est un Levi (c'est un Levi minimal parmi ceux qui contiennent $M$). Il r\'esulte de la d\'efinition d'Arthur ([A1] paragraphe 3) que le terme  $\rho^{Art}(\beta',u)\check{\beta}'$ est le m\^eme, que le groupe ambiant soit $G$ ou $G_{\beta'}$. Cela prouve (2). $\square$

Pour $\beta\in \Sigma(A_{M})$, notons $G_{\beta}$ le sous-groupe de $G$ engendr\'e par $M$ et les sous-groupes radiciels associ\'es aux racines  $n\beta$ pour $n\in {\mathbb Z}$. On a

(3) $\rho^G(\beta,u)=\rho^{G_{\beta}}(\beta,u)$.

Cela r\'esulte comme (2) d'un d\'evissage facile. 

Soit $L\in {\cal L}(M)$. Comme en 1.3(2), soit $u'$ un \'el\'ement de l'orbite induite de $M$ \`a $L$ par l'orbite de $u$. Pour tout $\beta'\in \Sigma(A_{L})$, on a l'\'egalit\'e
$$(4) \qquad \rho^G(\beta',u')=\sum_{\beta\in \Sigma(A_{M}),\beta_{L}= \beta'}\rho^G(\beta,u)_{L}.$$

Preuve. On note $G_{\beta'}$ le sous-groupe de $G$ engendr\'e par $L$ et les sous-groupes radiciels associ\'es aux racines  $n\beta'$ pour $n\in {\mathbb Z}$. Pour $\beta\in \Sigma(A_{M})$, on d\'efinit $G_{\beta}$ comme ci-dessus. Les racines $\beta\in \Sigma(A_{M})$ qui se restreignent en $\beta'$ sont exactement les racines $\beta\in \Sigma^{G_{\beta'}}(A_{M})$ qui se restreignent en $\beta'$. Pour celles-ci, le groupe $G_{\beta}$ est contenu dans $G_{\beta'}$. La relation (3) appliqu\'ee pour $G$ et pour $G_{\beta'}$ entra\^{\i}ne que $\rho^G(\beta,u)=\rho^{G_{\beta'}}(\beta,u)$. Il en r\'esulte que le membre de droite de (4) ne change pas quand on remplace le groupe ambiant $G$ par $G_{\beta'}$. D'apr\`es (3), il en est de m\^eme du membre de gauche. Si $dim(G_{\beta',SC})<dim(G_{SC})$, on conclut en raisonnant par r\'ecurrence sur cette dimension. Reste le cas o\`u $G_{\beta'}=G$. Alors $L$ est un Levi propre maximal et $\beta'$ est une racine indivisible. Dans ce cas, on a d'apr\`es (2)
$$(5)\qquad \rho^G(\beta',u')=\rho^{Art}(\beta',u')\check{\beta}'-\sum_{n\geq2}\rho^G(n\beta',u').$$
On applique 1.3(2):
$$\rho^{Art}(\beta',u')\check{\beta}'=\sum_{m\geq1}\sum_{\beta\in \Sigma_{ind}(A_{M}),\beta_{L}= m\beta'}\rho^{Art}(\beta,u)\check{\beta}_{L}.$$
En utilisant (2), c'est aussi
$$\sum_{m\geq1}\sum_{\beta\in \Sigma_{ind}(A_{M}),\beta_{L}= m\beta'}\sum_{k\geq1}\rho^G(k\beta,u)_{L}.$$
La triple somme se simplifie: les racines $k\beta$ intervenant ici sont exactement les \'el\'ements de $\Sigma(A_{M})$ qui se restreignent en un multiple positif de $\beta'$. On obtient
$$\sum_{n\geq1}\sum_{\beta\in \Sigma(A_{M}),\beta_{L}= n\beta'}\rho^G(\beta,u)_{L}.$$
Pour $n\geq2$, on applique la relation (4) d\'ej\`a d\'emontr\'ee pour $n\beta'$. La somme ci-dessus devient
$$(\sum_{\beta\in \Sigma(A_{M}),\beta_{L}= \beta'}\rho^G(\beta,u)_{L})+\sum_{n\geq2}\rho^G(n\beta',u).$$
En glissant cette expression de $\rho^{Art}(\beta',u')\check{\beta}'$ dans l'expression (5), on obtient (4). $\square$

\bigskip
\subsection{Modification de la d\'efinition des int\'egrales orbitales pond\'er\'ees}

On revient \`a la situation g\'en\'erale de 1.2 dont on reprend les notations. Si $\omega$ n'est pas  trivial sur $M_{\gamma}(F)$, on pose encore $J_{\tilde{M}}^{\tilde{G}}(\gamma,\omega,f)=0$ pour tout $f$.

On suppose maintenant que $\omega$ est trivial sur $M_{\gamma}(F)$. On ne change rien dans le cas o\`u  $\gamma$ est $\tilde{G}$-\'equisingulier. Dans le cas g\'en\'eral, on d\'efinit pour tout $\alpha\in \Sigma(A_{\tilde{M}})$ un \'el\'ement $\rho(\alpha,\gamma)\in {\cal A}_{\tilde{M}}$, ou plus pr\'ecis\'ement $\rho^{\tilde{G}}(\alpha,\gamma)$, par la formule
$$\rho^{\tilde{G}}(\alpha,\gamma)=\sum_{\beta\in \Sigma(A_{M_{\eta}}),\beta_{\tilde{M}}=\alpha}\rho^{G_{\eta}}(\beta,u)_{\tilde{M}}.$$

  Pour  $\alpha\in \Sigma(A_{\tilde{M}})$, pour $a\in A_{\tilde{M}}(F)$ en position g\'en\'erale et pour $\lambda\in i{\cal A}_{\tilde{M},{\mathbb C}}$, posons
 $$r_{\alpha}(\gamma,a;\lambda)= \vert \alpha(a)-\alpha(a)^{-1}\vert _{F}^{<\lambda,\rho(\alpha,\gamma)>}.$$
  On d\'efinit ensuite une $(\tilde{G},\tilde{M})$-famille $(r_{\tilde{P}}(\gamma,a;\lambda))_{\tilde{P}\in {\cal P}(\tilde{M})}$ par
 $$r_{\tilde{P}}(\gamma,a;\lambda)=\prod_{\alpha>_{P}0}r_{\alpha}(\gamma,a;\lambda/2)$$
 pour $\lambda\in i{\cal A}^*_{\tilde{M}}$. Comme pr\'ec\'edemment, on d\'eduit de cette $(\tilde{G},\tilde{M})$-famille une fonction $r_{\tilde{M}}^{\tilde{G}}(\gamma,a;\lambda)$ et on pose $r_{\tilde{M}}^{\tilde{G}}(\gamma,a)=r_{\tilde{M}}^{\tilde{G}}(\gamma,a;0)$.
 
   Pour $f\in C_{c}^{\infty}(\tilde{G}(F))$, consid\'erons la fonction
 $$(1) \qquad a\mapsto \sum_{\tilde{L}\in {\cal L}(\tilde{M})}r_{\tilde{M}}^{\tilde{L}}(\gamma,a)J_{\tilde{L}}^{\tilde{G}}(a\gamma,\omega,f).$$
 
 \ass{Lemme}{La fonction (1) a une limite quand $a$ tend vers $1$ parmi les \'el\'ements en position g\'en\'erale de $A_{\tilde{M}}(F)$.}
 
 Preuve. Notons $\varphi^{Art}(a,f)$ et $\varphi(a,f)$ les fonctions d\'efinies par les formules (2) du paragraphe 1.2 et (1) ci-dessus. On peut r\'ecrire
 $$\varphi^{Art}(a,f)=D^{\tilde{G}}(a\gamma)^{1/2}\int_{M_{\gamma}(F)\backslash G(F)}\omega(g)f(g^{-1}a\gamma g)(r^{Art}v)_{\tilde{M}}(\gamma,a,g)\,dg,$$
 o\`u
 $$(r^{Art}v)_{\tilde{M}}(\gamma,a,g)=\sum_{\tilde{L}\in {\cal L}(\tilde{M})}r_{\tilde{M}}^{\tilde{L},Art}(\gamma,a)v_{\tilde{L}}(g).$$
 On reconna\^{\i}t cette expression: c'est la fonction associ\'ee comme toujours \`a la $(\tilde{G},\tilde{M})$-famille produit  $(r^{Art}_{\tilde{P}}(\gamma,a;\lambda)v_{\tilde{P}}(g;\lambda))_{\tilde{P}\in {\cal P}(\tilde{M})}$ ([A2] corollaire 6.5). On a une expression analogue pour la fonction $\varphi(a,f)$: il suffit de supprimer les exposants $Art$. Soit $(c_{\tilde{P}}(\gamma,a;\lambda))_{\tilde{P}\in {\cal P}(\tilde{M})}$ la $(\tilde{G},\tilde{M})$-famille telle que $r_{\tilde{P}}(\gamma,a;\lambda)=c_{\tilde{P}}(\gamma,a;\lambda)r_{\tilde{P}}^{Art}(\gamma,a;\lambda)$. D'apr\`es [A2] lemme 6.3, on a une \'egalit\'e
 $$(rv)_{\tilde{M}}(\gamma,a,g)=\sum_{\tilde{Q}\in {\cal F}(\tilde{M})}c'_{\tilde{Q}}(\gamma,a)(r^{Art}v)_{\tilde{M}}^{\tilde{Q}}(\gamma,a,g),$$
 o\`u $c'_{\tilde{Q}}(\gamma,a)$ est d\'efinie par [A2] 6.3. Un calcul habituel de descente des int\'egrales orbitales pond\'er\'ees conduit alors \`a l'expression
 $$\varphi(a,f)=\sum_{\tilde{Q}=\tilde{L}U_{Q}\in {\cal F}(\tilde{M})}c'_{\tilde{Q}}(\gamma,a)D^{\tilde{L}}(a\gamma)^{1/2}\int_{M_{\gamma}(F)\backslash L(F)}\omega(g)f_{\tilde{Q},\omega}(l^{-1}a\gamma l)(r^{Art}v)_{\tilde{M}}^{\tilde{L}}(\gamma,a,l)\,dl.$$
 Ou encore
 $$(2) \qquad \varphi(a,f)=\sum_{\tilde{Q}=\tilde{L}U_{Q}\in {\cal F}(\tilde{M})}c'_{\tilde{Q}}(\gamma,a)\varphi^{\tilde{L},Art}(a,f_{\tilde{Q},\omega}).$$
 Il suffit de voir que toutes les fonctions apparaissant ont une limite quand $a$ tend vers $1$. C'est le r\'esultat d'Arthur pour les fonctions $\varphi^{\tilde{L},Art}(a,f_{\tilde{Q},\omega})$. Soit $\tilde{P}\in {\cal P}(\tilde{M})$. Notons $\Sigma_{ind}(A_{\tilde{M}})$ l'ensemble des racines indivisibles de $A_{\tilde{M}}$ dans $G$. Le terme $r_{\tilde{P}}(\gamma,a;\lambda)$ est produit sur les racines $\alpha\in \Sigma_{ind}(A_{\tilde{M}})$ qui sont positives pour $P$ des expressions
 $$\prod_{n\geq1}\vert \alpha(a)^n-\alpha(a)^{-n}\vert_{F}^{<\lambda,\rho(n\alpha,\gamma)>}.$$
 La singularit\'e en $a$ de cette expression est la m\^eme que celle de
 $$\vert \alpha(a)-\alpha(a)^{-1}\vert _{F}^{<\lambda,\sum_{n\geq1}\rho(n\alpha,\gamma)>}.$$
 Un calcul analogue vaut pour $r^{Art}_{\tilde{P}}(\gamma,a;\lambda)$. La somme $\sum_{n\geq1}\rho(n\alpha,\gamma)$ y est remplac\'ee par 
 $$\sum_{n\geq1}\sum_{\beta\in \Sigma_{ind}(A_{M_{\eta}}); \beta_{\tilde{M}}=n\alpha}\rho^{Art}( \beta,u)\check{\beta}_{\tilde{M}}.$$
 Il r\'esulte de 1.4(2) que cette expression est \'egale \`a
 $$\sum_{n\geq1}\sum_{\beta\in \Sigma_{ind}(A_{M_{\eta}}); \beta_{\tilde{M}}=n\alpha}\sum_{k\geq1}\rho^{G_{\eta}}(k\beta,u)_{\tilde{M}}.$$
 Cette expression se simplifie: les racines $k\beta$ y intervenant d\'ecrivent tous les \'el\'ements de $\Sigma(A_{M_{\eta}})$ dont la restriction est un multiple positif de $\alpha$. Elle est donc \'egale \`a
 $$\sum_{n\geq1} \sum_{\beta\in \Sigma(A_{M_{\eta}}); \beta_{\tilde{M}}=n\alpha}\rho^{G_{\eta}}(\beta,u)_{\tilde{M}}.$$
 D'apr\`es nos d\'efinitions, cela est \'egal \`a  $\sum_{n\geq1}\rho(n\alpha,\gamma)$. Cela montre que les fonctions $r_{\tilde{P}}(\gamma,a;\lambda)$ et $r^{Art}_{\tilde{P}}(\gamma,a;\lambda)$ ont m\^eme singularit\'e en $a=1$, donc que
  le rapport $c_{\tilde{P}}(\gamma,a;\lambda)$ est r\'egulier en ce point. Le terme $c'_{\tilde{Q}}(\gamma,1)$ est donc d\'efini et il est facile de montrer que $c'_{\tilde{Q}}(\gamma,a)$ tend vers $c'_{\tilde{Q}}(\gamma,1)$ quand $a$ tend vers $1$. $\square$
 
 On d\'efinit $J_{\tilde{M}}^{\tilde{G}}(\gamma,\omega,f)$ comme la limite de la fonction (1) quand $a$ tend vers $1$. De nouveau, si $ \gamma$ est $\tilde{G}$-\'equisingulier, on retrouve la d\'efinition simple donn\'ee plus haut.
 
 La preuve ci-dessus, plus pr\'ecis\'ement l'\'egalit\'e (2), montre que la pr\'ecision (3) du paragraphe 1.2 vaut aussi pour nos int\'egrales orbitales pond\'er\'ees. A savoir
 
 (3)  il existe $r>0$ tel que, pour tout $\gamma\in \tilde{M}(F)$, tout $f\in C_{c}^{\infty}(\tilde{G}(F))$, il existe $C>0$ de sorte que
 $$\vert J_{\tilde{M}}^{\tilde{G}}(\gamma,\omega,f)-\sum_{\tilde{L}\in {\cal L}(\tilde{M})}r_{\tilde{M}}^{\tilde{L}}(\gamma,a)J_{\tilde{L}}^{\tilde{G}}(a\gamma,\omega,f)\vert \leq Cd(a)^{r}$$
 pour tout $a\in A_{\tilde{M}}(F)$ en position g\'en\'erale et assez proche de $1$.

 Rappelons que les donn\'ees de $\gamma$ et d'une mesure de Haar sur $M_{\gamma}(F)$ d\'efinissent un \'el\'ement $\boldsymbol{\gamma}\in D_{g\acute{e}om}(\tilde{M}(F),\omega)\otimes Mes(M(F))^*$: pour $\varphi\in C_{c}^{\infty}(\tilde{M}(F))$ et pour une mesure de Haar $dm$ sur $M(F)$, on a $I^{\tilde{M}}(\boldsymbol{\gamma},\varphi\otimes dm)= I^{\tilde{M}}(\gamma,\omega,\varphi)$, o\`u le membre de droite est calcul\'e \`a l'aide de la mesure $dm$ et de celle fix\'ee sur $M_{\gamma}(F)$.
 L'application qui, \`a $\gamma$ et \`a une mesure de Haar sur $M_{\gamma}(F)$, associe la forme lin\'eaire $f\mapsto J_{\tilde{M}}^{\tilde{G}}(\gamma,\omega,f)$ v\'erifie les propri\'et\'es  requises pour se factoriser puis s'\'etendre par lin\'earit\'e en une application d\'efinie sur $D_{g\acute{e}om}(\tilde{M}(F),\omega)\otimes Mes(M(F))^*$. C'est-\`a-dire que l'on peut d\'efinir une application lin\'eaire qui, \`a $\boldsymbol{\gamma}\in D_{g\acute{e}om}(\tilde{M},\omega)\otimes Mes(M(F))^*$, associe une forme lin\'eaire $f\mapsto J_{\tilde{M}}^{\tilde{G}}(\boldsymbol{\gamma},f)$, de sorte que si $\boldsymbol{\gamma}$ provient comme ci-dessus d'un \'el\'ement $\gamma$ et d'une mesure de Haar sur $M_{\gamma}(F)$, on ait l'\'egalit\'e $J_{\tilde{M}}^{\tilde{G}}(\boldsymbol{\gamma},f)=J_{\tilde{M}}^{\tilde{G}}(\gamma,\omega,f)$.
 
 On aura besoin plus tard de la propri\'et\'e suivante. Soit $\alpha\in \Sigma(A_{\tilde{M}})$. Introduisons le sous-groupe $G_{\alpha}$ de $G$ engendr\'e par $M$ et les sous-espaces radiciels associ\'es aux racines de la forme $k\alpha$ pour $k\in {\mathbb Z}$. Il est normalis\'e par $\tilde{M}$, c'est-\`a-dire que l'on a l'\'egalit\'e $G_{\alpha}\tilde{M}=\tilde{M}G_{\alpha}$. On note $\tilde{G}_{\alpha}$ cet espace tordu. Alors
 
 (4) on a l'\'egalit\'e $\rho^{\tilde{G}}(\alpha,\gamma)=\rho^{\tilde{G}_{\alpha}}(\alpha,\gamma)$.

Cela r\'esulte de 1.4(3) et d'un d\'evissage des d\'efinitions. 
 
   \bigskip

\subsection{D\'efinition des int\'egrales orbitales pond\'er\'ees $\omega$-\'equivariantes}
On fixe toujours $\tilde{M}$, la mesure sur ${\cal A}_{\tilde{M}}$ et le sous-groupe compact $K$. Pour simplifier, fixons  des mesures de Haar sur $G(F)$ et $M(F)$. On d\'efinit comme d'habitude des homomorphismes
$$H_{G}:G(F)\to {\cal A}_{G}\text{ et }H_{\tilde{G}}:G(F)\to {\cal A}_{\tilde{G}}$$
par
$$exp(<x^*,H_{G}(g)>)=\vert x^*(g)\vert _{F}, \text{ resp }exp(<x^*,H_{\tilde{G}}(g)>)=\vert x^*(g)\vert _{F}$$
pour tout $x^*\in X^*(G)^{\Gamma_{F}}$, resp. $x^*\in X^*(G)^{\Gamma_{F},\theta}$. Le terme $H_{\tilde{G}}(g)$  n'est autre que la projection naturelle de $H_{G}(g)$ sur ${\cal A}_{\tilde{G}}$. Notons   ${\cal A}_{\tilde{G},F}$ l'image de l'application $H_{\tilde{G}}$. C'est un r\'eseau de l'espace ${\cal A}_{\tilde{G}}$. Notons $G(F)^1$ le noyau de $H_{\tilde{G}}$ et posons 
$$\tilde{{\cal A}}_{\tilde{G},F}=G(F)^1\backslash \tilde{G}(F).$$
C'est un espace principal homog\`ene sous ${\cal A}_{\tilde{G},F}$. On note 
$$\tilde{H}_{\tilde{G}}:\tilde{G}(F)\to \tilde{{\cal A}}_{\tilde{G},F}$$
l'application naturelle. Introduisons l'espace $C^{\infty}_{ac}(\tilde{G}(F))$ form\'e des fonctions $f$ sur $\tilde{G}(F)$ telles que:

(i) il existe un sous-groupe ouvert compact $K'$ de $G(F)$ tel que $f$ soit biinvariante par $K'$;

(ii) pour tout \'el\'ement $\varphi\in C_{c}^{\infty}( \tilde{{\cal A}}_{\tilde{G},F})$ (c'est-\`a-dire que $\varphi$ est une fonction sur $\tilde{{\cal A}}_{\tilde{G},F}$ \`a support fini), le produit $f(\varphi\circ \tilde{H}_{\tilde{G}})$ appartient \`a $C_{c}^{\infty}(\tilde{G}(F))$.

Les d\'efinitions des int\'egrales orbitales ou des int\'egrales orbitales pond\'er\'ees se g\'en\'eralisent  aux  \'el\'ements de $C^{\infty}_{ac}(\tilde{G}(F)) $. En effet, soient $\gamma \in\tilde{M}(F)$ et $f\in C^{\infty}_{ac}(\tilde{G}(F))$.  La projection dans $ \tilde{{\cal A}}_{\tilde{G},F}$ de la classe de conjugaison de $\gamma$ est r\'eduite \`a un point. Choisissons $\varphi\in C_{c}^{\infty}( \tilde{{\cal A}}_{\tilde{G},F})$ valant $1$ en ce point. On pose $J_{\tilde{M}}^{\tilde{G}}(\gamma,\omega,f)=J_{\tilde{M}}^{\tilde{G}}(\gamma,\omega,f(\varphi\circ \tilde{H}_{\tilde{G}}))$. Cela ne d\'epend pas du choix de $\varphi$. On note $I_{ac}(\tilde{G}(F),\omega)$ le quotient de $C^{\infty}_{ac}(\tilde{G}(F))$ par le sous-espace des \'el\'ements $f$ tels que $I^{\tilde{G}}(\gamma,\omega,f)=0$ pour tout $\gamma\in \tilde{G}_{reg}(F)$.

A l'aide des caract\`eres pond\'er\'es, Arthur d\'efinit une application lin\'eaire
$$\phi_{\tilde{M}}:C_{c}^{\infty}(\tilde{G}(F))\to I_{ac}(\tilde{M}(F),\omega).$$
(Arthur traite le cas o\`u $\omega=1$, le cas g\'en\'eral est similaire, cf. [W1] 6.4).
Pour $\gamma\in \tilde{M}(F)$, pour une mesure fix\'ee sur $M_{\gamma}(F)$ et pour $f\in C_{c}^{\infty}(\tilde{G}(F))$, on d\'efinit $I_{\tilde{M}}^{\tilde{G}}(\gamma,\omega,f)$ par la formule de r\'ecurrence
$$I_{\tilde{M}}^{\tilde{G}}(\gamma,\omega,f)=J_{\tilde{M}}^{\tilde{G}}(\gamma,\omega,f)-\sum_{\tilde{L}\in {\cal L}(\tilde{M}), \tilde{L}\not=\tilde{G}}I_{\tilde{M}}^{\tilde{L}}(\gamma,\omega,\phi_{\tilde{L}}(f)).$$
Quand $\omega=1$, Arthur montre que cette distribution $f\mapsto I_{\tilde{M}}^{\tilde{G}}(\gamma,f)$ est invariante par conjugaison. Le cas g\'en\'eral est similaire: la distribution $f\mapsto I_{\tilde{M}}^{\tilde{G}}(\gamma,\omega,f)$ est $\omega$-\'equivariante, c'est-\`a-dire qu'elle se factorise en une application d\'efinie sur $I(\tilde{G}(F),\omega)$. Remarquons que l'on a besoin de conna\^{\i}tre cette propri\'et\'e par r\'ecurrence pour que la d\'efinition ci-dessus ait un sens. Plus exactement, on a besoin de savoir par r\'ecurrence que cette distribution s'\'etend en une application lin\'eaire d\'efinie sur $C^{\infty}_{ac}(\tilde{G}(F))$ et que celle-ci se factorise en une application lin\'eaire d\'efinie sur $I_{ac}(\tilde{G}(F),\omega)$. Cela r\'esulte des propri\'et\'es de l'application $\phi_{\tilde{M}}$. Par ailleurs, Arthur montre que la distribution $f\mapsto I_{\tilde{M}}^{\tilde{G}}(\gamma,\omega,f)$  est ind\'ependante du sous-groupe $K$ choisi, lequel peut donc dispara\^{\i}tre des donn\'ees.

{\bf Attention:} cette distribution  n'est pas, en g\'en\'eral, support\'ee par la classe de conjugaison de $\gamma$. Elle n'appartient m\^eme pas \`a $D_{g\acute{e}om}(\tilde{G}(F),\omega)$.

Encore une fois, on se d\'ebarrasse des mesures en d\'efinissant $I_{\tilde{M}}^{\tilde{G}}(\boldsymbol{\gamma},{\bf f})$ pour $\boldsymbol{\gamma}\in D_{g\acute{e}om}(\tilde{M}(F))\otimes Mes(M(F))^*$ et ${\bf f}\in C_{c}^{\infty}(\tilde{G}(F))\otimes Mes(G(F))$, ou ${\bf f}\in I(\tilde{G}(F))\otimes Mes(G(F))$.

\bigskip

\subsection{Propri\'et\'es des int\'egrales orbitales pond\'er\'ees $\omega$-\'equivariantes}
Soit toujours $\tilde{M}$ un espace de Levi de $\tilde{G}$. Pour \'enoncer les quatre premi\`eres propri\'et\'es, il est plus commode de fixer des mesures de Haar sur $M(F)$ et $G(F)$, ainsi que sur les groupes $M_{\gamma}(F)$ qui apparaissent.

Soient $\gamma\in \tilde{M}(F)$ et $f\in C_{c}^{\infty}(\tilde{G}(F))$. Pour $\varphi\in C_{c}^{\infty}( \tilde{{\cal A}}_{\tilde{G},F})$, on a l'\'egalit\'e

(1) $I_{\tilde{M}}^{\tilde{G}}(\gamma,\omega,f(\varphi\circ \tilde{H}_{\tilde{G}}))=\varphi\circ \tilde{H}_{\tilde{G}}(\gamma)I_{\tilde{M}}^{\tilde{G}}(\gamma,\omega,f)$.

A fortiori, la distribution $f\mapsto I_{\tilde{M}}^{\tilde{G}}(\gamma,\omega,f)$ est support\'ee par l'ensemble des $\gamma'\in \tilde{G}(F)$ tels que $\tilde{H}_{\tilde{G}}(\gamma')=\tilde{H}_{\tilde{G}}(\gamma)$.

On a aussi l'\'egalit\'e
$$(2) \qquad I_{\tilde{M}}^{\tilde{G}}(\gamma,\omega,f)=lim_{a\to 1}\sum_{\tilde{L}\in {\cal L}(\tilde{M})}r_{\tilde{M}}^{\tilde{L}}(\gamma,a)I_{\tilde{L}}^{\tilde{G}}(a\gamma,\omega,f),$$
la limite \'etant prise au m\^eme sens qu'en 1.2. Plus pr\'ecis\'ement

(3)   il existe $r>0$ tel que, pour tout $\gamma\in \tilde{M}(F)$, tout $f\in C_{c}^{\infty}(\tilde{G}(F))$, il existe $C>0$ de sorte que
 $$\vert I_{\tilde{M}}^{\tilde{G}}(\gamma,\omega,f)-\sum_{\tilde{L}\in {\cal L}(\tilde{M})}r_{\tilde{M}}^{\tilde{L}}(\gamma,a)I_{\tilde{L}}^{\tilde{G}}(a\gamma,\omega,f)\vert \leq Cd(a)^{r}$$
 pour tout $a\in A_{\tilde{M}}(F)$ en position g\'en\'erale et assez proche de $1$.

Supposons que $\gamma$ est $\tilde{G}$-\'equisingulier. Alors

(4) il existe $f'\in C_{c}^{\infty}(\tilde{M}(F))$ et un voisinage de $\gamma$ dans $\tilde{M}(F)$ tels que, pour $\gamma'$ dans ce voisinage, on ait l'\'egalit\'e $I_{\tilde{M}}^{\tilde{G}}(\gamma',\omega,f)=I^{\tilde{M}}(\gamma',\omega,f')$.

Pour la cinqui\`eme propri\'et\'e, il est plus simple de se d\'ebarrasser des mesures. Rappelons que pour tout $\tilde{L}\in {\cal L}(\tilde{M})$, on dispose d'applications lin\'eaires en dualit\'e
$$\begin{array}{ccc}I(\tilde{L}(F),\omega)\otimes Mes(L(F))&\to& I(\tilde{M}(F),\omega)\otimes Mes(M(F))\\ {\bf f}&\mapsto&{\bf f}_{\tilde{M}}\\ \end{array}$$
et
$$\begin{array}{ccc}D_{g\acute{e}om}(\tilde{M}(F),\omega)\otimes Mes(M(F))^*&\to&D_{g\acute{e}om}(\tilde{L}(F),\omega)\otimes Mes(L(F))^*\\ \boldsymbol{\gamma}&\mapsto&\boldsymbol{\gamma}^{\tilde{L}}.\\ \end{array}$$
Remarquons que si $\boldsymbol{\gamma}$ est l'int\'egrale orbitale dans $\tilde{M}(F)$ associ\'ee \`a un \'el\'ement $\tilde{G}$-\'equisingulier de $ \tilde{M}(F)$, $\boldsymbol{\gamma}^{\tilde{L}}$ est l'int\'egrale orbitale dans $\tilde{L}(F)$ associ\'ee au m\^eme \'el\'ement. 

 {\bf Remarque.} M\^eme si on impose que $\omega$ est trivial sur $Z(G;F)^{\theta}$, $\omega$ peut ne pas \^etre trivial sur $Z(M;F)^{\theta}$ pour un espace de Levi $\tilde{M}$. Dans ce cas, les espaces $ I(\tilde{M}(F),\omega)$ et $D_{g\acute{e}om}(\tilde{M}(F),\omega)$ sont nuls.
 
 \bigskip

Pour deux \'el\'ements $\tilde{L},\tilde{L}'\in {\cal P}(\tilde{M})$, on d\'efinit le r\'eel  $d_{\tilde{M}}^{\tilde{G}}(\tilde{L},\tilde{L}')$: il est nul sauf si ${\cal A}_{\tilde{M}}^{\tilde{G}}={\cal A}_{\tilde{M}}^{\tilde{L}}\oplus {\cal A}_{\tilde{M}}^{\tilde{L}'}$; si cette \'egalit\'e est v\'erifi\'ee, c'est le rapport entre la mesure sur le premier espace et le produit des mesures sur les deux espaces du second membre.

\ass{Lemme}{Soient $\tilde{L}\in {\cal L}(\tilde{M})$,  $\boldsymbol{\gamma}\in D_{g\acute{e}om}(\tilde{M}(F),\omega)\otimes Mes(M(F))^*$ et ${\bf f}\in I(\tilde{G}(F),\omega)\otimes Mes(G(F))$. On  a l'\'egalit\'e
$$ I_{\tilde{L}}^{\tilde{G}}(\boldsymbol{\gamma}^{\tilde{L}},{\bf f})=\sum_{\tilde{L}'\in {\cal L}(\tilde{M})}d_{\tilde{M}}^{\tilde{G}}(\tilde{L},\tilde{L}')I_{\tilde{M}}^{\tilde{L}'}(\boldsymbol{\gamma},{\bf f}_{\tilde{L}'}).$$}

Preuve. On peut fixer des mesures et supposer que $\boldsymbol{\gamma}$ est l'int\'egrale orbitale associ\'ee \`a un \'el\'ement $\gamma\in \tilde{M}(F)$. Supposons que $\gamma$ est $\tilde{G}$-\'equisingulier. Alors la preuve de la formule est essentiellement formelle \`a partir de la formule de descente des poids, cf. [A3] preuve du th\'eor\`eme 8.1. Traitons le cas g\'en\'eral.  Pour $a\in A_{\tilde{M}}(F)$ en position g\'en\'erale, posons
$$\varphi(a)=\sum_{\tilde{R}\in {\cal L}(\tilde{L})}r_{\tilde{L}}^{\tilde{R}}(\gamma,a)I_{\tilde{R}}^{\tilde{G}}(a\gamma,\omega,f).$$
Par une formule de descente, on a pour tout $\tilde{R}$ l'\'egalit\'e
$$r_{\tilde{L}}^{\tilde{R}}(\gamma,a)=\sum_{\tilde{R}'\in {\cal L}(\tilde{M}); \tilde{R}'\subset \tilde{R}}d_{\tilde{M}}^{\tilde{R}}(\tilde{L},\tilde{R}')r_{\tilde{M}}^{\tilde{R}'}(\gamma,a).$$
D'o\`u
$$\varphi(a)=\sum_{\tilde{R}'\in {\cal L}(\tilde{M})}r_{\tilde{M}}^{\tilde{R}'}(\gamma,a)\sum_{\tilde{R}\in {\cal L}(\tilde{L}); \tilde{R}'\subset \tilde{R}}d_{\tilde{M}}^{\tilde{R}}(\tilde{L},\tilde{R}')I_{\tilde{R}}^{\tilde{G}}(a\gamma,\omega,f).$$
Fixons $\tilde{R}'$ et $\tilde{R}$. Puisque $R'_{a\gamma}=G_{a\gamma}$, on peut utiliser la formule de l'\'enonc\'e pour cet \'el\'ement. D'o\`u
$$I_{\tilde{R}}^{\tilde{G}}(a\gamma,\omega,f)=\sum_{\tilde{L}'\in {\cal L}(\tilde{R}')}d_{\tilde{R}'}^{\tilde{G}}(\tilde{R},\tilde{L}')I_{\tilde{R}'}^{\tilde{L}'}(a\gamma,\omega,f_{\tilde{L}',\omega}).$$
Puis
$$\varphi(a)=\sum_{\tilde{L}'\in {\cal L}(\tilde{M})}\sum_{\tilde{R}'\in {\cal L}^{\tilde{L}'}(\tilde{M})}x(\tilde{R}',\tilde{L}')r_{\tilde{M}}^{\tilde{R}'}(\gamma,a)I_{\tilde{R}'}^{\tilde{L}'}(a\gamma,\omega,f_{\tilde{L}',\omega}),$$
o\`u $x(\tilde{R}',\tilde{L}')$ est la somme sur les $\tilde{R}\in {\cal L}(\tilde{L})$ tels que $\tilde{R}'\subset \tilde{R}$ des produits
$$d_{\tilde{M}}^{\tilde{R}}(\tilde{L},\tilde{R}') d_{\tilde{R}'}^{\tilde{G}}(\tilde{R},\tilde{L}').$$
Consid\'erons l'ensemble $A$ des couples d'espace de Levi $(\tilde{L}',\tilde{R}')$ tels que $\tilde{M}\subset \tilde{R}'\subset \tilde{L}'$ et $d_{\tilde{M}}^{\tilde{G}}(\tilde{L},\tilde{L}')\not=0$. Consid\'erons l'ensemble $B$ des triplets $(\tilde{L}',\tilde{R}',\tilde{R})$ tels que 
 $\tilde{M}\subset \tilde{R}'\subset \tilde{L}'$, $\tilde{L}\subset \tilde{R}$, $\tilde{R}'\subset \tilde{R}$ et $d_{\tilde{M}}^{\tilde{R}}(\tilde{L},\tilde{R}')d_{\tilde{R}'}^{\tilde{G}}(\tilde{R},\tilde{L}')\not=0$. Montrons que l'on a 
 
(5)   l'application $(\tilde{L}',\tilde{R}',\tilde{R})\mapsto (\tilde{L}',\tilde{R}')$ est une bijection de $B$ sur $A$;  pour $(\tilde{L}',\tilde{R}',\tilde{R})\in B$, on a l'\'egalit\'e 
 $$d_{\tilde{M}}^{\tilde{R}}(\tilde{L},\tilde{R}')d_{\tilde{R}'}^{\tilde{G}}(\tilde{R},\tilde{L}')= d_{\tilde{M}}^{\tilde{G}}(\tilde{L},\tilde{L}').$$

Soit $(\tilde{L}',\tilde{R}',\tilde{R})\in B$. La non-nullit\'e de $d_{\tilde{M}}^{\tilde{R}}(\tilde{L},\tilde{R}')d_{\tilde{R}'}^{\tilde{G}}(\tilde{R},\tilde{L}')$ \'equivaut aux relations
 $$(6) \qquad {\cal A}_{\tilde{M}}^{\tilde{R}}={\cal A}_{\tilde{M}}^{\tilde{L}}\oplus {\cal A}_{\tilde{M}}^{\tilde{R}'}$$
et
$$(7) \qquad {\cal A}_{\tilde{R}'}^{\tilde{G}}={\cal A}_{\tilde{R}'}^{\tilde{R}}\oplus {\cal A}_{\tilde{R}'}^{\tilde{L}'}.$$
D'o\`u
$$ {\cal A}_{\tilde{M}}^{\tilde{G}}={\cal A}_{\tilde{M}}^{\tilde{R}'}\oplus {\cal A}_{\tilde{R}'}^{\tilde{G}}=
{\cal A}_{\tilde{M}}^{\tilde{R}'}\oplus {\cal A}_{\tilde{R}'}^{\tilde{R}}\oplus {\cal A}_{\tilde{R}'}^{\tilde{L}'}={\cal A}_{\tilde{M}}^{\tilde{R}}\oplus {\cal A}_{\tilde{R}'}^{\tilde{L}'}={\cal A}_{\tilde{M}}^{\tilde{L}}\oplus {\cal A}_{\tilde{M}}^{\tilde{R}'}\oplus {\cal A}_{\tilde{R}'}^{\tilde{L}'}={\cal A}_{\tilde{M}}^{\tilde{L}}\oplus {\cal A}_{\tilde{M}}^{\tilde{L}'}.$$
 D'o\`u l'\'egalit\'e
$$(8) \qquad {\cal A}_{\tilde{M}}^{\tilde{G}}={\cal A}_{\tilde{M}}^{\tilde{L}}\oplus {\cal A}_{\tilde{M}}^{\tilde{L}'}$$
des termes extr\^emes, qui \'equivaut \`a la non-nullit\'e de $d_{\tilde{M}}^{\tilde{G}}(\tilde{L},\tilde{L}')$. Cela prouve que $(\tilde{L}',\tilde{R}')\in A$. Inversement, soit $(\tilde{L}',\tilde{R}')\in A$. On doit prouver qu'il existe exactement un espace de Levi $\tilde{R}$ tel que $(\tilde{L}',\tilde{R}',\tilde{R})\in B$. Il y en a au plus un:  il est d\'etermin\'e par la relation (6). Montrons que ce $\tilde{R}$ existe. On le d\'efinit comme le commutant dans $\tilde{G}$ du tore  $(A_{\tilde{L}}\cap A_{\tilde{R}'})^0$. D'apr\`es [I] 3.1(11), cet ensemble est un espace de Levi pourvu qu'il ne soit pas vide. Or il n'est pas vide puisque sa d\'efinition implique qu'il contient $\tilde{L}$ et $\tilde{R}'$. Ces deux inclusions impliquent $A_{\tilde{R}}\subset A_{\tilde{L}}$ et $A_{\tilde{R}}\subset A_{\tilde{R}'}$. On a aussi par d\'efinition
 $(A_{\tilde{L}}\cap A_{\tilde{R}'})^0\subset A_{\tilde{R}}$. On obtient donc
 $$(A_{\tilde{L}}\cap A_{\tilde{R}'})^0\subset A_{\tilde{R}}\subset A_{\tilde{L}}\cap A_{\tilde{R}'}$$
 d'o\`u l'\'egalit\'e $(A_{\tilde{L}}\cap A_{\tilde{R}'})^0=A_{\tilde{R}}$ puisque ce dernier ensemble est connexe. Cette \'egalit\'e est \'equivalente \`a ${\cal A}_{\tilde{R}}={\cal A}_{\tilde{L}}\cap {\cal A}_{\tilde{R}'}$. Par passage aux orthogonaux dans ${\cal A}_{\tilde{M}}$, on obtient $ {\cal A}_{\tilde{M}}^{\tilde{R}}={\cal A}_{\tilde{M}}^{\tilde{L}}+  {\cal A}_{\tilde{M}}^{\tilde{R}'}$. Mais ces deux derniers espaces sont en somme directe d'apr\`es l'inclusion $\tilde{R}'\subset \tilde{L}'$ et d'apr\`es (8). Donc (6) est v\'erifi\'e. On a montr\'e ci-dessus que (6) et (7) impliquaient (8). Le calcul est r\'eversible: (6) et (8) impliquent (7). Puisque $\tilde{R}$ v\'erifie (6) et (7), on a $(\tilde{L}',\tilde{R}',\tilde{R})\in B$. La derni\`ere assertion de (5) s'obtient facilement en pr\'ecisant le calcul qui a conduit ci-dessus \`a l'\'egalit\'e (8). Cela prouve (5). 
 
 Cette propri\'et\'e entra\^{\i}ne que, pour $\tilde{L}'$ et $\tilde{R}'$ intervenant dans l'expression de $\varphi(a)$ ci-dessus, on a $x(\tilde{R}',\tilde{L}')=d_{\tilde{M}}^{\tilde{G}}(\tilde{L},\tilde{L}')$. 
  Alors
$$\varphi(a)=\sum_{\tilde{L}'\in {\cal L}(\tilde{M})}d_{\tilde{M}}^{\tilde{G}}(\tilde{L},\tilde{L}')\sum_{\tilde{R}'\in {\cal L}^{\tilde{L}'}(\tilde{M})}r_{\tilde{M}}^{\tilde{R}'}(\gamma,a)I_{\tilde{R}'}^{\tilde{L}'}(a\gamma,\omega,f_{\tilde{L}',\omega}).$$
Pour tout $\tilde{L}'$, la relation (2) nous dit que la somme en $\tilde{R}'$ tend vers $I_{\tilde{M}}^{\tilde{L}'}(\gamma,\omega,f_{\tilde{L}',\omega})$ quand $a$ tend vers $1$. Donc, quand $a$ tend vers $1$, $\varphi(a)$ tend vers le membre de droite de l'\'egalit\'e de l'\'enonc\'e.

Fixons $b\in A_{\tilde{L}}(F)$ en position g\'en\'erale et faisons tendre $a$ vers $b$. Pour $\tilde{R}\in {\cal L}(\tilde{L})$, on applique (4): il existe $f'\in C_{c}^{\infty}(\tilde{R}(F))$ tel que $I_{\tilde{R}}^{\tilde{G}}(\gamma',\omega,f)=I^{\tilde{R}}(\gamma',\omega,f')$ pour tout $\gamma'\in \tilde{R}(F)$ assez proche de $b\gamma$.  Le deuxi\`eme terme est une int\'egrale orbitale ordinaire. En appliquant la formule de descente usuelle pour ces int\'egrales, on a
$$I^{\tilde{R}}(a\gamma,\omega,f')=I^{\tilde{M}}(a\gamma,\omega,(f')_{\tilde{M},\omega}).$$
La limite de cette expression quand $a$ tend vers $b$ est $I^{\tilde{M}}(b\gamma,\omega,(f')_{\tilde{M},\omega})$.  Par d\'efinition de l'induction, c'est $I^{\tilde{L}}(b\boldsymbol{\gamma}^{\tilde{L}},(f')_{\tilde{L},\omega})$. Ecrivons la distribution $\boldsymbol{\gamma}^{\tilde{L}}$ comme une combinaison lin\'eaire d'int\'egrales orbitales associ\'ees \`a des \'el\'ements $\gamma_{i}$ de l'orbite induite de $\gamma$:
$$I^{\tilde{L}}(\boldsymbol{\gamma}^{\tilde{L}},\psi)=\sum_{i=1,...,n}c_{i}I^{\tilde{L}}(\gamma_{i},\omega,\psi)$$
pour tout $\psi\in C_{c}^{\infty}(\tilde{L}(F))$.  Puisque $b\in A_{\tilde{L}}(F)$, on a la m\^eme \'egalit\'e si l'on remplace $\boldsymbol{\gamma}^{\tilde{L}}$ par $b\boldsymbol{\gamma}^{\tilde{L}}$ et les $\gamma_{i}$ par $b\gamma_{i}$. Donc
$$lim_{a\to b} I_{\tilde{R}}^{\tilde{G}}(a\gamma,\omega,f)=\sum_{i=1,...,n}c_{i}I^{\tilde{L}}(b\gamma_{i},\omega,(f')_{\tilde{L},\omega}).$$
Puisque $b$ est en position g\'en\'erale, le membre de droite n'est autre que
$$\sum_{i=1,...,n}c_{i}I^{\tilde{R}}(b\gamma_{i},\omega,f').$$
En revenant \`a la d\'efinition de $f'$, on obtient
$$lim_{a\to b} I_{\tilde{R}}^{\tilde{G}}(a\gamma,\omega,f)=\sum_{i=1,...,n}c_{i}I_{\tilde{R}}^{\tilde{G}}(b\gamma_{i},\omega,f).$$
On montrera plus loin que, pour tout $i=1,...,n$, on a l'\'egalit\'e
$$(9) \qquad lim_{a\to b}r_{\tilde{L}}^{\tilde{R}}(\gamma,a)=r_{\tilde{L}}^{\tilde{R}}(\gamma_{i},b).$$
Admettant cela, on obtient
$$lim_{a\to b}\varphi(a)=\sum_{i=1,...,n}c_{i}\sum_{\tilde{R}\in {\cal L}(\tilde{L})}r_{\tilde{L}}^{\tilde{R}}(\gamma_{i},b)I_{\tilde{R}}^{\tilde{G}}(b\gamma_{i},\omega,f).$$
En utilisant  la relation (2), on voit que cette expression a une limite quand $b$ tend vers $1$. Plus pr\'ecis\'ement,
$$lim_{b\to 1}lim_{a\to b}\varphi(a)=\sum_{i=1,...,n}c_{i}I_{\tilde{L}}^{\tilde{G}}(\gamma_{i},\omega,f).$$
La limite de gauche est \'egale \`a $lim_{a\to 1}\varphi(a)$, puisque cette derni\`ere limite existe.
Par d\'efinition, la somme de droite ci-dessus n'est autre que le membre de gauche de l'\'egalit\'e de l'\'enonc\'e.  Donc $lim_{a\to 1}\varphi(a)$ est \'egale au membre de gauche de cette \'egalit\'e. On a d\'ej\`a prouv\'e qu'elle \'etait \'egale au membre de droite. Cela conclut.

Il reste \`a prouver l'\'egalit\'e (9). On se ram\`ene ais\'ement \`a prouver que, pour $\tilde{Q}\in {\cal P}(\tilde{L})$ et $\lambda\in i{\cal A}_{\tilde{L}}^*$, on a l'\'egalit\'e similaire
$$(10) \qquad lim_{a\to b}r_{\tilde{Q}}(\gamma,a;\lambda)=r_{\tilde{Q}}(\gamma_{i},b;\lambda).$$
Fixons $\tilde{P}\in {\cal P}(\tilde{M})$ avec $\tilde{P}\subset \tilde{Q}$. Le terme $r_{\tilde{Q}}(\gamma,a;\lambda)$ est produit sur les $\alpha'\in \Sigma(A_{\tilde{M}})$ qui sont positifs pour $P$ de
$$\vert \alpha'(a)-\alpha'(a)^{-1}\vert _{F}^{<\lambda,\rho(\alpha',\gamma)>/2}.$$
Parce que $\lambda\in i{\cal A}_{\tilde{L}}^*$, ce terme vaut $1$ si la restriction de  $\alpha'$ \`a ${\cal A}_{\tilde{L}}$ est nulle.  Soit $\alpha'$ de restriction non nulle \`a ${\cal A}_{\tilde{L}}$. Alors cette restriction $\alpha$ appartient \`a $\Sigma(A_{\tilde{L}})$. Que $\alpha'$ soit positive pour $P$ \'equivaut \`a ce que $\alpha$ soit positive pour $Q$. D'autre part, 
$$lim_{a\to b}\vert \alpha'(a)-\alpha'(a)^{-1}\vert _{F}=\vert \alpha(b)-\alpha(b)^{-1}\vert _{F}.$$
Donc
$$lim_{a\to b}r_{\tilde{Q}}(\gamma,a;\lambda)=\prod_{\alpha\in \Sigma(A_{\tilde{L}});\alpha>_{Q}0}\vert \alpha(b)-\alpha(b)^{-1}\vert _{F}^{<\lambda,\rho(\alpha,\gamma)>/2},$$
o\`u on a pos\'e
$$\rho(\alpha,\gamma)=\sum_{\alpha'\in \Sigma(A_{\tilde{M}});\alpha'_{ \tilde{L}}=\alpha}\rho(\alpha',\gamma)_{\tilde{L}}.$$
En comparant avec la d\'efinition de $r_{\tilde{Q}}(\gamma_{i},b;\lambda)$, on voit que (10) r\'esulte de l'\'egalit\'e

$$(11) \qquad  \rho(\alpha,\gamma)=\rho(\alpha,\gamma_{i})$$ 
pour tout $\alpha\in \Sigma(A_{\tilde{L}})$. Ecrivons $\gamma=u\eta$ comme en 1.2. Alors on peut supposer que $\gamma_{i}=u_{i}\eta$, o\`u $u_{i}$ appartient \`a la classe de conjugaison dans $L_{\eta}$ induite par la classe de conjugaison de $u$ dans $M_{\eta}$. Par d\'efinition,
$$\rho(\alpha,\gamma_{i})=\sum_{\beta\in \Sigma(A_{L_{\eta}}); \beta_{\tilde{L}}=\alpha}\rho^{G_{\eta}}(\beta,u_{i})_{\tilde{L}},$$
$$\rho(\alpha,\gamma)=\sum_{\alpha'\in \Sigma(A_{\tilde{M}});\alpha'_{\tilde{L}}=\alpha}\sum_{\beta'\in \Sigma(A_{M_{\eta}}),\beta'_{\tilde{M}}=\alpha'}\rho^{G_{\eta}}(\beta',u)_{\tilde{L}}.$$
On peut r\'ecrire
$$\rho(\alpha,\gamma)=\sum_{\beta'\in \Sigma(A_{M_{\eta}});\beta'_{ \tilde{L}}=\alpha}\rho^{G_{\eta}}(\beta',u)_{\tilde{L}}=\sum_{\beta\in \Sigma(A_{L_{\eta}}); \beta_{\tilde{L}}=\alpha}\sum_{\beta'\in \Sigma(A_{M_{\eta}}); \beta'_{L_{\eta}}=\beta}\rho^{G_{\eta}}(\beta',u_{i})_{\tilde{L}}.$$
Mais alors l'\'egalit\'e (11) r\'esulte de 1.4(4). 
Cela ach\`eve la d\'emonstration. $\square$

On peut pr\'eciser la preuve ci-dessus: pour montrer que $\varphi(a)$ tend vers le membre de droite de l'\'egalit\'e de l'\'enonc\'e, utilisons la relation (3) au lieu de (2). On obtient (en notant $\boldsymbol{\gamma}^{\tilde{L}}$ la distribution induite de l'int\'egrale orbitale associ\'ee \`a $\gamma$):

(12) existe $r>0$ tel que, pour tout $ \gamma\in \tilde{M}(F)$ et tout $f\in C_{c}^{\infty}(\tilde{G}(F))$, il existe $C>0$ de sorte que
 $$\vert I_{\tilde{L}}^{\tilde{G}}(\boldsymbol{\gamma}^{\tilde{L}},f)-\sum_{\tilde{R}\in {\cal L}(\tilde{L})}r_{\tilde{L}}^{\tilde{R}}(\gamma,a)I_{\tilde{R}}^{\tilde{G}}(a\gamma,\omega,f)\vert \leq Cd(a)^{r}$$
 pour tout $a\in A_{\tilde{M}}(F)$ en position g\'en\'erale et assez proche de $1$.
\bigskip

\subsection{Variantes des  termes $\rho(\beta,u)$}
 On suppose dans ce paragraphe $G=\tilde{G}$ et $\omega=1$. On consid\`ere un Levi $M$ de $G$ et un \'el\'ement unipotent $u\in M(F)$. On fixe une paire de Borel \'epingl\'ee ${\cal E}=(B,T,(E_{\alpha})_{\alpha\in \Delta})$ de $G$ d\'efinie sur $\bar{F}$ de sorte que $M$ soit standard relativement \`a ${\cal E}$.  On introduit  une extension finie $F'$ de $F$ telle que ${\cal E}$ soit d\'efinie sur $F'$ et  $G$ soit d\'eploy\'e sur $F'$. On introduit aussi l'action galoisienne quasi-d\'eploy\'ee $\sigma\mapsto \sigma_{G^*}$ qui conserve ${\cal E}$, cf. [I] 1.2.  
 Notons $\Sigma(T)$ l'ensemble des racines de $T$ dans $G$. On a fix\'e en 1.2 une forme quadratique d\'efinie positive sur $X_{*}(T)\otimes {\mathbb R}$, d'o\`u, par dualit\'e, une telle forme sur $X^*(T)\otimes {\mathbb R}$, que l'on note $(.,.)$.

 On fixe une fonction $B$ sur $\Sigma(T)$ \`a valeurs dans l'ensemble ${\mathbb Q}_{>0}$ des rationnels strictement positifs. On  lui impose les conditions suivantes:

 $\bullet$ $B(-\beta)=B(\beta)$, $B(\sigma_{G^*}(\beta))=B(\beta)$ et $B(w\beta)=B(\beta)$ pour tout $\beta\in \Sigma(T)$, tout $\sigma\in \Gamma_{F}$  et tout $w\in W$;
 
 $\bullet$ pour toute composante irr\'eductible $\Sigma'$ du syst\`eme de racines $\Sigma(T)$, ou bien $B$ est constante sur $\Sigma'$, ou bien la fonction $\beta\mapsto \frac{B(\beta)}{(\beta,\beta)}$ est constante sur $\Sigma'$. 

On pose $V^*=X^*(T)\otimes {\mathbb R}$, $V_{*}=X_{*}(T)\otimes {\mathbb R}$. D\'efinissons les sous-ensembles $\Sigma(T,B)=\{B(\beta)^{-1}\beta; \beta\in \Sigma(T)\}$ de $V^*$ et $\check{\Sigma}(T,B)=\{B(\beta)\check{\beta}; \beta\in \Sigma(T)\}$ de $V_{*}$.  On va voir que $\Sigma(T,B)$ est un syst\`eme de racines dont $\check{\Sigma}(T,B)$ est l'ensemble associ\'e de coracines.   On note  $Z^*$ l'annulateur dans $V^*$   de l'ensemble de coracines $\check{\Sigma}(T)$ et $Z_{*}$ l'annulateur dans $V_{*}$ de l'ensemble $\Sigma(T)$. Pour tout sous-syst\`eme irr\'eductible $\Sigma'$ de $\Sigma(T)$, on note $V^*(\Sigma')$, resp. $V_{*}(\check{\Sigma}')$, le sous-espace de $V^*$, resp. $V_{*}$, engendr\'e par $\Sigma'$, resp. par l'ensemble correspondant $\check{\Sigma}'$ de coracines. On a
$$V^*=Z^*\oplus(\oplus_{\Sigma'\in Irr}V^*(\Sigma')),\,\,V_{*}=Z_{*}\oplus(\oplus_{\Sigma'\in Irr}V_{*}(\check{\Sigma}')),$$
o\`u $Irr$ est l'ensemble des composantes irr\'eductibles. Notons $Irr_{+}$ l'ensemble des composantes sur lesquelles $B$ est constante et notons $Irr_{-}$ le compl\'ementaire. Posons
$$V(B)^*=Z^*\oplus(\oplus_{\Sigma'\in Irr_{+}}V^*(\Sigma'))\oplus(\oplus_{\Sigma'\in Irr_{-}}V_{*}(\check{\Sigma}')),$$
$$V(B)_{*}=Z_{*}\oplus(\oplus_{\Sigma'\in Irr_{+}}V_{*}(\check{\Sigma}'))\oplus(\oplus_{\Sigma'\in Irr_{-}}V^{*}(\Sigma')).$$
Notons $\iota^*_{Z^*}$ l'identit\'e de $Z^*$. Pour $\Sigma'\in Irr_{+}$, notons $\iota^*_{\Sigma'}$ l'homoth\'etie de $V^*(\Sigma')$ de rapport la valeur constante de $B$ sur $\Sigma'$. Pour $\Sigma'\in Irr_{-}$, notons $\iota^*_{\Sigma'}:V^*(\Sigma')\to V_{*}(\check{\Sigma}')$ la compos\'ee de l'isomorphisme d\'eduit de la forme quadratique fix\'ee (cet isomorphisme envoie une racine $\beta$ sur $(\beta,\beta)\check{\beta}/2$) et de l'homoth\'etie de rapport la valeur constante sur $\Sigma'$ de la fonction $\beta\mapsto 2B(\beta)(\beta,\beta)^{-1}$. On note $\iota^*:V^*\to V(B)^*$ la somme directe de $\iota^*_{Z^*}$ et des $\iota^*_{\Sigma'}$. On note $\iota_{*}:V_{*}\to V(B)_{*}$ l'inverse du transpos\'e de $\iota^*$. On v\'erifie que $\iota_{*}$ envoie $\Sigma(T,B)$ sur
$$(1) \qquad (\sqcup_{\Sigma'\in Irr_{+}}\Sigma')\sqcup(\sqcup_{\Sigma'\in Irr_{-}}\check{\Sigma}'),$$
tandis que $\iota_{*}$ envoie $\check{\Sigma}(T,B)$ sur
$$(2)\qquad
(\sqcup_{\Sigma'\in Irr_{+}}\check{\Sigma}')\sqcup(\sqcup_{\Sigma'\in Irr_{-}}\Sigma').$$
Il est clair que l'ensemble (1) est un ensemble de racines dont l'ensemble (2) est l'ensemble associ\'e de coracines.

Notons $j:\Sigma(T)\to \Sigma(T,B)$ l'application $\beta\mapsto B(\beta)^{-1}\beta$. L'application compos\'ee $\iota^*\circ j$ est la somme, composante par composante, soit de l'identit\'e, soit de l'\'echange $\beta\mapsto \check{\beta}$. Il en r\'esulte que, pour un sous-ensemble $\Sigma_{0}\subset \Sigma(T)$, $\Sigma_{0}$ est un sous-syst\`eme de racines, resp. un sous-ensemble de Levi, de $\Sigma(T)$ si et seulement si $j(\Sigma_{0})$ est un sous-syst\`eme de racines, resp. un sous-ensemble de Levi, de $\Sigma(T,B)$. Notons que $j$ est \'equivariante pous l'action de $W$ et pour l'action galoisienne quasi-d\'eploy\'ee.

On note $\Sigma(Z(M)^0,B)$ l'ensemble des restrictions non nulles \`a $X_{*}(Z(M)^0)\times {\mathbb R}$ des $B(\beta')^{-1}\beta$ pour $\beta'\in \Sigma(T)$, ou encore des $\beta'$ pour $\beta'\in \Sigma(T,B)$. L'interpr\'etation ci-dessus montre que cet ensemble a beaucoup de  propri\'et\'es communes avec celui des racines $\Sigma(Z(M)^0)$. D'abord, pour $\beta\in \Sigma(Z(M)^0,B)$, on peut d\'efinir le groupe de Levi $M_{\beta}$ (sur $F'$) tel que $ X_{*}(Z(M_{\beta})^0)\otimes {\mathbb R}$  soit l'annulateur de $\beta$ dans $X_{*}(Z(M)^0)\times {\mathbb R}$. On a aussi

(3) supposons que $ M$ soit maximal parmi les  Levi propres de $G$; alors
 il existe il unique entier $n\geq1$ et, au signe pr\`es, un unique $\beta\in \Sigma(Z(M)^0,B)$ tel que  $\Sigma(Z(M)^0,B)=\{\pm k\beta; k=1,...,n\}$;

(4) pour $\beta\in \Sigma(Z(M)^0,B)$, l'ensemble des $\beta'\in \Sigma(T)$ tels que la restriction de $B(\beta')^{-1}\beta'$ \`a $X_{*}(Z(M)^0)\times {\mathbb R}$ soit de la forme $k\beta$ avec $k\in {\mathbb Z}$ (y compris $k=0$) est un sous-syst\`eme de racines de $\Sigma(T)$ qui contient $\Sigma^M(T)$;

plus g\'en\'eralement

(5) soient $\beta_{1},...,\beta_{n}$ des \'el\'ements lin\'eairement ind\'ependants de $\Sigma(Z(M)^0,B)$; alors l'ensemble des $\beta'\in \Sigma(T)$ tels que la restriction de $B(\beta')^{-1}\beta'$ \`a $X_{*}(Z(M)^0)\times {\mathbb R}$ appartienne au ${\mathbb Z}$-module engendr\'e par les $\beta_{i}$ est un sous-syst\`eme de racines de $\Sigma(T)$ qui contient $\Sigma^M(T)$.

 Dans la situation (4), il existe \`a isomorphisme pr\`es un unique groupe r\'eductif connexe d\'efini et d\'eploy\'e sur $F'$ qui poss\`ede un tore maximal isomorphe \`a $T$ et dont le syst\`eme de racines est l'ensemble d\'ecrit de $\beta'$. On le note $G_{\beta}$. Il poss\`ede un groupe de Levi isomorphe \`a $M$ et on identifie ce sous-groupe \`a $M$.  On note encore $B$ la restriction de $B$ \`a l'ensemble de racines de $G_{\beta}$; cette restriction v\'erifie les m\^emes conditions que la fonction de d\'epart.
 
 {\bf Attention.} Le groupe $G_{\beta}$ n'est pas, en g\'en\'eral, un sous-groupe de $G$.  Par exemple,  consid\'erons $G=SO(5)$, $M=GL(1)\times SO(3)$ et une fonction $B$ proportionnelle au carr\'e de la longueur. L'ensemble $\Sigma(Z(M)^0,B)$ a deux \'el\'ements $\alpha$ et $2\alpha$. On v\'erifie que $G_{2\alpha}=SO(3)\times SO(3)$.
 \bigskip
 
Soit $u$ un \'el\'ement unipotent de $M(F')$. Pour $\beta\in \Sigma(Z(M)^0,B)$, on d\'efinit un terme $\rho^G(\beta,u,B)\in X_{*}(Z(M)^0)\otimes {\mathbb R}$ par r\'ecurrence sur $dim(G_{SC})$. Si $dim(G_{\beta,SC})<dim(G_{SC})$, on suppose d\'efini $\rho^{G_{\beta}}(\beta,u,B)$ et on pose $\rho^G(\beta,u,B)=\rho^{G_{\beta}}(\beta,u,B)$. Si $dim(G_{\beta,SC})=dim(G_{SC})$, on a $G=G_{\beta}$ et on est dans la situation (3). On note $\beta'$ l'unique \'el\'ement de $\Sigma_{ind}(Z(M)^0)$ qui est de la forme $q\beta$ avec $q\in {\mathbb Q}_{>0}$. Les termes $\rho^G(k\beta,u,B)$ pour $k\geq2$ ont d\'ej\`a \'et\'e d\'efinis et on pose
 $$(6) \qquad \rho^G(\beta,u,B)=(\sum_{k\geq1}\rho^G(k\beta',u))-(\sum_{k\geq2}\rho^{G}(k\beta,u,B)),$$
 o\`u $\rho^G(k\beta',u)$ est le terme d\'efini en 1.4. 
 
 On redescend maintenant \`a la situation d\'efinie sur $F$. On note $\Sigma(A_{M},B)$ l'ensemble des restrictions non nulles \`a $\mathfrak{a}_{M}$ d'\'el\'ements de $\Sigma(Z(M)^0,B)$. On note cette restriction $\beta\mapsto \beta_{{\cal A}_{M}}$.  Soit $u$ un \'el\'ement unipotent de $M(F)$. Pour $\alpha\in \Sigma(A_{M},B)$, on pose
 $$\rho^G(\alpha,u,B)=\sum_{\beta\in \Sigma(Z(M)^0,B); \beta_{{\cal A}_{M}}=\alpha}\rho^G(\beta,u,B)_{{\cal A}_{M}}.$$
 
 Pour $\alpha\in \Sigma(A_{M},B)$, introduisons l'ensemble des $\beta'\in \Sigma(T)$ tels que la restriction de $B(\beta')^{-1}\beta'$ \`a ${\cal A}_{M}$ est de la forme $k\alpha$ avec $k\in {\mathbb Z}$. Comme pr\'ec\'edemment, c'est le syst\`eme de racines d'un groupe connexe $G_{\alpha}$ et $M$ s'identifie \`a un groupe de Levi de $G_{\alpha}$. On peut munir ce groupe d'une structure sur $F$ de la fa\c{c}on suivante. Tout d'abord, parce que l'action galoisienne quasi-d\'eploy\'ee (comme l'action naturelle) est triviale sur $A_{M}$, la propri\'et\'e d'invariance de $B$ implique que le syst\`eme de racines  de $G_{\alpha}$ est conserv\'e par cette action. On peut munir conform\'ement le groupe $G_{\alpha}$ d'une action galoisienne $\sigma\mapsto \sigma_{G_{\alpha}^*}$ quasi-d\'eploy\'ee sur $F$. On peut imposer que cette action co\"{\i}ncide sur $M$ avec l'action $\sigma\mapsto \sigma_{G^*}$.  On a une \'egalit\'e $\sigma_{G}=ad_{u_{{\cal E}}(\sigma)^{-1}}\circ \sigma_{G^*}$ o\`u $\sigma_{G}$ est l'action naturelle, cf. [I]1.2. Puisque $M$ est un Levi standard, on a $u_{{\cal E}}(\sigma)\in M$ et $u_{{\cal E}}^{-1}$ est un cocycle \`a valeurs dans $M/Z(G)$ (muni de l'action quasi-d\'eploy\'ee). Il est clair que $Z(G)\subset Z(G_{\alpha})$, donc $u_{{\cal E}}$ se pousse en un cocycle \`a valeurs dans $M/Z(G_{\alpha})$. On voit alors que la formule $\sigma_{G_{\alpha}}=ad_{u_{{\cal E}}(\sigma)}^{-1}\circ \sigma_{G_{\alpha}^*}$ munit $G_{\alpha}$ d'une action galoisienne qui co\"{\i}ncide sur $M$ avec l'action naturelle. En d\'evissant les d\'efinitions, on v\'erifie l'\'egalit\'e
 
 (7) $\rho^G(\alpha,u,B)=\rho^{G_{\alpha}}(\alpha,u,B)$.
 
 Il y a une bijection entre les sous-ensembles d'\'el\'ements indivisibles $\Sigma_{ind}(A_{M})$ et $\Sigma_{ind}(A_{M},B)$. A un \'el\'ement $\alpha'\in\Sigma_{ind}(A_{M})$, on associe l'unique \'el\'ement  $\alpha$ de $\Sigma_{ind}(A_{M},B)$ tel que $\alpha=q\alpha'$ avec $q\in {\mathbb Q}_{>0}$. Pour $\alpha'$ et $\alpha$ se correspondant ainsi, on a
 
 $$(8)\qquad \sum_{n\geq1}\rho^G(n\alpha',u)=\sum_{n\geq1}\rho^G(n\alpha,u,B).$$
 
Preuve.  Le membre de gauche est la somme des $\rho^G(\beta',u)_{{\cal A}_{M}}$ pour $\beta'\in \Sigma(Z(M)^0)$ se restreignant en un multiple positif de $\alpha'$. Ou encore
 $$\sum_{\beta'\in \Sigma_{ind}(Z(M)^0); \beta'_{A_{M}}\in {\mathbb N}_{>0}\alpha'}\sum_{n\geq1}\rho^G(n\beta',u)_{{\cal A}_{M}}.$$
 De m\^eme, le membre de droite est
 $$\sum_{\beta\in \Sigma_{ind}(Z(M)^0,B); \beta_{{\cal A}_{M}}\in {\mathbb N}_{>0}\alpha}\sum_{n\geq1}\rho^G(n\beta,u,B)_{{\cal A}_{M}}.$$
 Il y a une bijection similaire \`a la pr\'ec\'edente entre $\Sigma_{ind}(Z(M)^0)$ et $\Sigma_{ind}(Z(M)^0,B)$. Il est clair que si $\beta'\mapsto \beta$ par cette bijection, $\beta'$ se restreint en un multiple positif de $\alpha'$ si et seulement si $\beta$ se restreint en un multiple positif de $\alpha$. Il suffit de fixer $\beta'$ et $\beta$ indivibles et se correspondant et de prouver l'\'egalit\'e
$$\sum_{n\geq1}\rho^G(n\beta',u)= \sum_{n\geq1}\rho^G(n\beta,u,B).$$
On introduit le Levi $M'$ engendr\'e par $M$ et les espaces radiciels associ\'es aux $n\beta'$ pour $n\in {\mathbb Z}$. Il r\'esulte des d\'efinitions que les termes ci-dessus ne changent pas si l'on remplace $G$ par $M'$. Cela nous ram\`ene au cas o\`u $M$ est propre maximal. Mais alors, la relation cherch\'ee r\'esulte de la d\'efinition (6). $\square$

 Nos termes $\rho^G(\alpha,u,B)$ v\'erifient une propri\'et\'e analogue \`a 1.4(4). Pr\'ecis\'ement, dans la situation de cette relation, pour $\alpha'\in \Sigma(A_{L},B)$, on a l'\'egalit\'e
 
 (9) $\rho^G(\alpha',u',B)=\sum_{\alpha\in \Sigma(A_{M},B); \alpha_{L}=\alpha'}\rho^G(\alpha,u,B)_{L}$.
 
Preuve.  On introduit le groupe $G_{\alpha'}$ comme ci-dessus, relatif au Levi $L$. On voit  facilement, en utilisant la relation (7) ci-dessus que les deux membres de (9) ne changent pas si l'on remplace $G$ par $G_{\alpha'}$. Cela nous ram\`ene par r\'ecurrence au cas o\`u $G_{\alpha'}=G$. Dans ce cas, $L$ est un Levi propre maximal et $\alpha'\in \Sigma_{ind}(A_{L},B)$. Pour $n\geq2$, on a d\'ej\`a d\'emontr\'e la relation (9) relative \`a $n\alpha'$. La relation \`a d\'emontrer est donc \'equivalente \`a
 $$(10) \qquad \sum_{n\geq1}\rho^G(n\alpha',u',B)=\sum_{n\geq1}\sum_{\alpha\in \Sigma(A_{M},B); \alpha_{L}=n\alpha'}\rho^G(\alpha,u,B)_{L}.$$
 On peut \'ecrire le membre de droite comme
 $$\sum_{\alpha\in \Sigma_{ind}(A_{M},B);\alpha_{L}\in {\mathbb N}_{>0}\alpha'}\sum_{n\geq1}\rho^G(n\alpha,u,B)_{L}.$$
Notons ici $\beta'$ l'\'el\'ement de $\Sigma_{ind}(A_{L})$ qui correspond \`a $\alpha'$.  En utilisant (8), c'est \'egal \`a
 $$\sum_{\beta\in \Sigma_{ind}(A_{M});\beta_{L}\in {\mathbb N}_{>0}\beta'}\sum_{n\geq1}\rho^G(n\beta,u)_{L}.$$
 Ou encore \`a
 $$\sum_{n\geq1}\sum_{\beta\in \Sigma(A_{M}); \beta_{L}=n\beta'}\rho^G(\beta,u)_{L}. $$
 D'apr\`es 1.4(4),  c'est aussi
 $$\sum_{n\geq1}\rho^G(n\beta',u).$$ 
 Toujours d'apr\`es (8), c'est aussi le membre de gauche de (10). Cela prouve (9). $\square$
 
 A l'aide des  constructions ci-dessus, on peut d\'efinir des variantes $J_{M}^G(u,B,f)$ et $I_{M}^G(u,B,f)$ des int\'egrales orbitales de 1.5 et 1.6. Nous ne donnons pas les preuves n\'ecessaires car elles sont identiques \`a celles que nous ferons dans le paragraphe suivant. 
 Un \'el\'ement  $\alpha\in \Sigma(A_{M},B)$ est une forme lin\'eaire sur $\mathfrak{a}_{M}$ et n'appartient pas forc\'ement \`a $X^*(A_{M})$. Toutefois, il existe un entier $n\geq1$ tel que $n\alpha\in X^*(A_{M})$. Pour $a\in A_{M}(F)$ assez proche de $1$, on peut d\'efinir $\alpha(a)$ de la fa\c{c}on suivante. On choisit un entier $n\geq1$ tel que $n\alpha\in X^*(A_{M})$, on \'ecrit $a=exp(H)$ avec $H\in \mathfrak{a}_{M}(F)$ proche de $0$ et on pose $\alpha(a)=exp(\frac{(n\alpha)(H)}{n})$. 
 Pour  un \'el\'ement unipotent $u\in M(F)$, pour $a\in A_{M}(F)$ en position g\'en\'erale et proche de $1$, pour $P\in {\cal P}(M)$ et pour $\lambda\in i{\cal A}_{M}$, posons
 $$r_{P}(u,a,B;\lambda)=\prod_{\alpha\in \Sigma(A_{M},B); \alpha>_{P}0}\vert \alpha(a)-\alpha(a)^{-1}\vert _{F}^{<\lambda,\rho^G(\alpha,u,B)>/2}.$$
 La collection $(r_{P}(u,a,B;\lambda))_{P\in {\cal P}(M)}$ est une $(G,M)$-famille dont on d\'eduit un terme $r_{M}^G(u,a,B)$ comme en 1.5. Supposons fix\'e un sous-groupe compact sp\'ecial en bonne position relativement \`a $M$. Pour $f\in C_{c}^{\infty}(G(F))$, la fonction
 $$a\mapsto \sum_{L\in {\cal L}(M)}r_{M}^L(u,a,B)J_{L}^G(au,f)$$
 a une limite quand $a$ tend vers $1$ parmi les \'el\'ements de $A_{M}(F)$ en position g\'en\'erale. 
   On note $J_{M}^G(u,B,f)$ cette limite.  L'int\'egrale invariante $I_{M}^G(u,B,f)$ s'en d\'eduit comme en 1.6. Ces termes v\'erifient des propri\'et\'es analogues aux termes $J^G_{M}(u,f)$ et $I^G_{M}(u,f)$. Plus canoniquement, on d\'efinit $I_{M}^G(\boldsymbol{\gamma},B,{\bf f})$ pour $\boldsymbol{\gamma}\in D_{unip}(M(F))\otimes Mes(M(F))^*$ et ${\bf f}\in I(G(F))\otimes Mes(G(F))^*$, o\`u $D_{unip}(M(F))$ est le sous-espace des \'el\'ements de $D_{g\acute{e}om}(M(F))$ \`a support unipotent. 
   
\bigskip

\subsection{ Variantes des int\'egrales orbitales pond\'er\'ees  dans le cas quasi-d\'eploy\'e \`a torsion int\'erieure }
On suppose dans ce paragraphe $(G,\tilde{G},{\bf a})$ quasi-d\'eploy\'e et \`a torsion int\'erieure. Rappelons que l'on note $\tilde{G}_{ss}$ l'ensemble des \'el\'ements semi-simples de $\tilde{G}$. On appelle syst\`eme de fonctions $B$ la donn\'ee pour tout $\eta\in \tilde{G}_{ss}(F)$ d'une fonction $B_{\eta}$ sur le syst\`eme de racines de $G_{\eta}$ de sorte que les conditions (1) et (2) suivantes soient v\'erifi\'ees.

(1) Pour tout $\eta\in \tilde{G}_{ss}(F)$, $B_{\eta}$ v\'erifie les conditions de 1.8.

Soient $\eta\in \tilde{G}_{ss}(F)$ et $x\in G$ tel que $x\sigma(x)^{-1}\in G_{\eta}$ pour tout $\sigma\in \Gamma_{F}$. Posons $\eta'=ad_{x^{-1}}(\eta)$. Alors $\eta'\in \tilde{G}_{ss}(F)$ et $ad_{x}$ est un torseur int\'erieur de $G_{\eta'}$ sur $G_{\eta}$ qui permet d'identifier les syst\`emes de racines de ces deux groupes, munis de leurs actions quasi-d\'eploy\'ees. On demande

(2) par cette identification, $B_{\eta'}$ s'identifie \`a $B_{\eta}$.

Fixons un tel syst\`eme de fonctions. Soit $\tilde{M}$ un espace de Levi de $\tilde{G}$. On fixe des mesures de Haar sur tous les groupes intervenant.
Pour $\gamma\in \tilde{M}(F)$, on va d\'efinir une distribution $f\mapsto I_{\tilde{M}}^{\tilde{G}}(\gamma,B,f)$ sur $C_{c}^{\infty}(\tilde{G}(F))$. Si  $\gamma$ est $\tilde{G}$-\'equisingulier, on pose $I_{\tilde{M}}^{\tilde{G}}(\gamma,B,f)=I_{\tilde{M}}^{\tilde{G}}(\gamma,f)$. Passons au cas g\'en\'eral. On \'ecrit $\gamma=u\eta$, avec $\eta\in \tilde{M}_{ss}(F)$ et $u\in M_{\eta}(F)$.   On note $\Sigma(A_{M},B_{\eta})$ l'ensemble des restrictions non nulles \`a $\mathfrak{a}_{M}$ d'\'el\'ements de $\Sigma(A_{M_{\eta}},B_{\eta})$. Pour $\alpha\in \Sigma(A_{M},B_{\eta})$, on pose
$$(3) \qquad \rho(\alpha,\gamma,B)=\sum_{\beta\in \Sigma(A_{M_{\eta}},B_{\eta}); \beta_{{\cal A}_{M}}=\alpha}\rho^{G_{\eta}}(\beta,u,B_{\eta})_{M}.$$

 On d\'efinit comme \`a la fin de 1.8 une fonction $a\mapsto \alpha(a)$ sur un voisinage assez petit de $1$ dans $A_{\tilde{M}}(F)$.  Pour $a$ en position g\'en\'erale et assez proche de $1$, on d\'efinit ensuite une $(\tilde{G},\tilde{M})$-famille $(r_{\tilde{P}}(\gamma,a,B;\lambda))_{\tilde{P}\in {\cal P}(\tilde{M})}$ par
$$r_{\tilde{P}}(\gamma,a,B;\lambda)=\prod_{\alpha\in \Sigma(A_{\tilde{M}},B_{\eta}); \alpha>_{P}0} \vert \alpha(a)-\alpha(a)^{-1}\vert _{F}^{<\lambda,\rho(\alpha,\gamma,B)>/2}.$$
On en d\'eduit comme toujours un nombre $r_{\tilde{M}}^{\tilde{G}}(\gamma,a,B)$.
On a

(4) la fonction
$$a\mapsto \sum_{\tilde{L}\in {\cal L}(\tilde{M})}r_{\tilde{M}}^{\tilde{L}}(\gamma,a,B)I_{\tilde{L}}^{\tilde{G}}(a\gamma,B,f)$$
a une limite quand $a$ tend vers $1$ parmi les \'el\'ements de $A_{\tilde{M}}(F)$ en position g\'en\'erale.

Preuve. On introduit la $(\tilde{G},\tilde{M})$-famille $(c_{\tilde{P}}(\gamma,a,B;\lambda))_{\tilde{P}\in {\cal P}(\tilde{M})}$ telle que $r_{\tilde{P}}(\gamma,a,B;\lambda)=c_{\tilde{P}}(\gamma,a,B;\lambda)r_{\tilde{P}}(\gamma,a;\lambda)$. Pour tout $\tilde{L}\in {\cal L}(\tilde{M})$, on a l'\'egalit\'e
$$r_{\tilde{M}}^{\tilde{L}}(\gamma,a,B)=\sum_{\tilde{R}\in {\cal L}^{\tilde{L}}(\tilde{M})}c_{\tilde{M}}^{\tilde{R}}(\gamma,a,B)r_{\tilde{R}}^{\tilde{L}}(\gamma,a).$$
La fonction que l'on consid\`ere est donc
$$a\mapsto \sum_{\tilde{R}\in {\cal L}(\tilde{M})}c_{\tilde{M}}^{\tilde{R}}(\gamma,a,B)\sum_{\tilde{L}\in {\cal L}(\tilde{R})}r_{\tilde{R}}^{\tilde{L}}(\gamma,a)I_{\tilde{L}}^{\tilde{G}}(a\gamma,f).$$
La relation 1.7(12) montre que, pour tout $\tilde{R}$, la somme int\'erieure a une limite quand $a$ tend vers $1$. Il suffit de prouver que $c_{\tilde{M}}^{\tilde{R}}(\gamma,a,B)$ a aussi une limite et il suffit encore de prouver que, pour tout $\tilde{P}\in {\cal P}(\tilde{M})$, $c_{\tilde{P}}(\gamma,a,B;\lambda)$ a une limite. D'apr\`es les d\'efinitions, on a
$$r_{\tilde{P}}(\gamma,a;\lambda)=\prod_{\beta'\in \Sigma^{G_{\eta}}(A_{M_{\eta}});\beta'_{{\cal A}_{M}}>_{P}0}\vert \beta'(a)-\beta'(a)^{-1}\vert _{F}^{<\lambda,\rho^{G_{\eta}}(\beta',u)>/2}$$
et
$$r_{\tilde{P}}(\gamma,a,B;\lambda)=\prod_{\beta\in \Sigma^{G_{\eta}}(A_{M_{\eta}},B_{\eta});\beta_{{\cal A}_{M}}>_{P}0}\vert \beta(a)-\beta(a)^{-1}\vert _{F}^{<\lambda,\rho^{G_{\eta}}(\beta,u,B_{\eta})>/2}.$$
On peut r\'ecrire ces formules
$$r_{\tilde{P}}(\gamma,a;\lambda)=\prod_{\beta'\in \Sigma_{ind}^{G_{\eta}}(A_{M_{\eta}});\beta'_{{\cal A}_{M}}>_{P}0}\prod_{n\geq1}\vert \beta'(a)^n-\beta'(a)^{-n}\vert _{F}^{<\lambda,\rho^{G_{\eta}}(n\beta',u)>/2}$$
et
$$r_{\tilde{P}}(\gamma,a,B;\lambda)=\prod_{\beta\in \Sigma_{ind}^{G_{\eta}}(A_{M_{\eta}},B_{\eta});\beta_{{\cal A}_{M}}>_{P}0}\prod_{n\geq1}\vert \beta(a)^n-\beta(a)^{-n}\vert _{F}^{<\lambda,\rho^{G_{\eta}}(n\beta,u,B_{\eta})>/2}.$$
Ces fonctions ont les m\^emes singularit\'es en $a=1$ que les fonctions
$$\prod_{\beta'\in \Sigma_{ind}^{G_{\eta}}(A_{M_{\eta}});\beta'_{{\cal A}_{M}}>_{P}0}\vert \beta'(a)-\beta'(a)^{-1}\vert _{F}^{<\lambda,\sum_{n\geq1}\rho^{G_{\eta}}(n\beta',u)>/2}$$
et 
$$\prod_{\beta\in \Sigma_{ind}^{G_{\eta}}(A_{M_{\eta}},B_{\eta});\beta_{{\cal A}_{M}}>_{P}0}\vert \beta(a)-\beta(a)^{-1}\vert _{F}^{<\lambda,\sum_{n\geq1}\rho^{G_{\eta}}(n\beta,u,B_{\eta})>/2}.$$
Il y a une bijection entre les ensembles d'\'el\'ements indivisibles $\Sigma_{ind}^{G_{\eta}}(A_{M_{\eta}})$ et $\Sigma_{ind}^{G_{\eta}}(A_{M_{\eta}},B_{\eta})$: \`a un \'el\'ement  $\beta'\in \Sigma_{ind}^{G_{\eta}}(A_{M_{\eta}})$, on associe l'unique \'el\'ement indivisible $\beta$ de $\Sigma^{G_{\eta}}(A_{M_{\eta}},B_{\eta})$ qui soit de la forme $q\beta'$ avec $q\in {\mathbb Q}$, $q>0$. Cette bijection pr\'eserve la positivit\'e pour $P$.
Si $\beta$ correspond \`a $\beta'$ par la bijection ci-dessus, les fonctions $\vert \beta'(a)-\beta'(a)^{-1}\vert _{F}$ et $\vert \beta(a)-\beta(a)^{-1}\vert _{F}$ ont m\^eme singularit\'e. De plus, d'apr\`es 1.8(8), on a l'\'egalit\'e
$$\sum_{k\geq1}\rho(k\beta',u)=\sum_{k\geq1}\rho(k\beta,u,B_{\eta}).$$
 On conclut que le rapport $c_{\tilde{P}}(\gamma,a,B,\lambda)$ est r\'egulier en $a=1$. $\square$

On d\'efinit $I_{\tilde{M}}^{\tilde{G}}(\gamma,B,f)$ comme la limite de la fonction (4).

En reprenant la preuve ci-dessus et en utilisant le lemme 1.7, on obtient l'\'egalit\'e
$$(5) \qquad I_{\tilde{M}}^{\tilde{G}}(\gamma,B,f)=\sum_{\tilde{R}\in {\cal L}(\tilde{M})}c_{\tilde{M}}^{\tilde{R}}(\gamma,1,B)I_{\tilde{R}}^{\tilde{G}}(\boldsymbol{\gamma}^{\tilde{R}},f),$$
o\`u $\boldsymbol{\gamma}^{\tilde{R}}$ est la distribution induite par l'int\'egrale orbitale dans $\tilde{M}$ associ\'ee \`a $\gamma$.

En utilisant davantage 1.7(12), on obtient aussi

(6) il existe $r>0$ tel que, pour tout $ \gamma\in \tilde{M}(F)$ et tout $ f\in C_{c}^{\infty}(\tilde{G}(F))$, il existe $C>0$ de sorte que
 $$\vert I_{\tilde{M}}^{\tilde{G}}(\gamma,B,f)-\sum_{\tilde{L}\in {\cal L}(\tilde{M})}r_{\tilde{M}}^{\tilde{L}}(\gamma,a,B)I_{\tilde{L}}^{\tilde{G}}(a\gamma,B,f)\vert \leq Cd(a)^{r}$$
 pour tout $a\in A_{\tilde{M}}(F)$ en position g\'en\'erale et assez proche de $1$.

Comme en 1.4,  on se d\'ebarrasse des mesures en d\'efinissant $I_{\tilde{M}}^{\tilde{G}}(\boldsymbol{\gamma},B,{\bf f})$ pour $\boldsymbol{\gamma}\in D_{g\acute{e}om}(\tilde{M}(F))\otimes Mes(M(F))^*$ et ${\bf f}\in C_{c}^{\infty}(\tilde{G}(F))\otimes Mes(G(F))$, ou ${\bf f}\in I(\tilde{G}(F))\otimes Mes(G(F))$.

Les distributions $I_{\tilde{M}}^{\tilde{G}}(\boldsymbol{\gamma},B,{\bf f})$ ont les m\^emes propri\'et\'es  (1), (2), (3) et (4) de 1.7. Elles v\'erifient \'egalement le lemme 1.7 et le raffinement 1.7(12). La preuve est la m\^eme, en rempla\c{c}ant l'utilisation de 1.4(4) par celle de 1.8(9).

{\bf Remarque.} Si $B_{\eta}$ est la fonction constante de valeur $1$ pour tout $\eta\in \tilde{M}_{ss}(F)$, on a \'evidemment  $I_{\tilde{M}}^{\tilde{G}}(\boldsymbol{\gamma},B,{\bf f})= I_{\tilde{M}}^{\tilde{G}}(\boldsymbol{\gamma},{\bf f})$.

\bigskip

En fait, ces deux distributions sont souvent \'egales, gr\^ace au lemme suivant. On rappelle que l'on note $p$ la caract\'eristique r\'esiduelle de $F$.

\ass{Lemme}{Supposons $p$ diff\'erent de $2$, $3$ et $5$ et supposons que, pour tout $\eta\in \tilde{M}_{ss}(F)$,  les valeurs de la fonction $B_{\eta}$ soient premi\`eres \`a $p$. Alors on a l'\'egalit\'e 
$$I_{\tilde{M}}^{\tilde{G}}(\boldsymbol{\gamma},B,{\bf f})= I_{\tilde{M}}^{\tilde{G}}(\boldsymbol{\gamma},{\bf f})$$
pour tout $\boldsymbol{\gamma}\in D_{g\acute{e}om}(\tilde{M}(F))\otimes Mes(M(F))^*$ et tout ${\bf f}\in I(\tilde{G}(F))\otimes Mes(G(F))$.}

Preuve. Il suffit de prouver que, pour tout $\gamma=u\eta\in \tilde{M}(F)$, tout $\tilde{P}\in {\cal P}(\tilde{M})$, tout $a\in A_{\tilde{M}}(F)$ en position g\'en\'erale et pour tout $\lambda\in i{\cal A}_{\tilde{M}}$, on a l'\'egalit\'e
$$r_{\tilde{P}}(\gamma,a,B;\lambda)=r_{\tilde{P}}(\gamma,a;\lambda)$$
pourvu que $a$ soit assez proche de $1$.  On reprend la preuve de (4) ci-dessus. On y a utilis\'e les trois propri\'et\'es suivantes:

(7)  pour $\beta'\in \Sigma^{G_{\eta}}(A_{M_{\eta}})$ et $n\geq1$ tel que $n\beta'\in \Sigma^{G_{\eta}}(A_{M_{\eta}})$, la fonction $\frac{\vert \beta'(a)^n-\beta'(a)^{-n}\vert _{F}}{\vert \beta'(a)-\beta'(a)^{-1}\vert _{F}}$ est r\'eguli\`ere et non nulle en $a=1$;

(8) pour $\beta\in \Sigma^{G_{\eta}}(A_{M_{\eta}},B_{\eta})$ et $n\geq1$ tel que $n\beta\in \Sigma^{G_{\eta}}(A_{M_{\eta}},B_{\eta})$, la fonction $\frac{\vert \beta(a)^n-\beta(a)^{-n}\vert _{F}}{\vert \beta(a)-\beta(a)^{-1}\vert _{F}}$ est r\'eguli\`ere et non nulle en $a=1$;

(9) pour $\beta'\in \Sigma_{ind}^{G_{\eta}}(A_{M_{\eta}})$ et $\beta=q\beta'\in \Sigma_{ind}^{G_{\eta}}(A_{M_{\eta}},B_{\eta})$ se correspondant, le fonction $\frac{\vert \beta'(a)-\beta'(a)^{-1}\vert _{F}}{\vert \beta(a)-\beta(a)^{-1}\vert _{F}}$ est r\'eguli\`ere et non nulle en $a=1$.

Pour d\'emontrer l'\'egalit\'e cherch\'ee, il suffit de prouver que les valeurs en $1$ de ces fonctions sont \'egales \`a $1$. Il suffit pour cela que les entiers $n$ et les rationnels $q$ intervenant soient premiers \`a $p$. En consid\'erant tous les syst\`emes de racines possibles, on v\'erifie qu'un entier $n$ intervenant dans (7) est forc\'ement inf\'erieur ou \'egal \`a $6$. Il est donc premier \`a $p$ d'apr\`es l'hypoth\`ese. Consid\'erons un entier $n$ intervenant dans (8). Introduisons un sous-tore maximal $T^*$ de $M_{\eta}$ comme en 1.8.  D'apr\`es les d\'efinitions, il y a deux \'el\'ements $\alpha_{1}$ et $\alpha_{2}$ de $\Sigma^{G_{\eta}}(T^*)$ de sorte que
$\beta$, resp. $n\beta$, soit la restriction \`a $A_{M_{\eta}}$ de $B_{\eta}(\alpha_{1})^{-1}\alpha_{1}$, resp. $B_{\eta}(\alpha_{2})^{-1}\alpha_{2}$. Les \'el\'ements $\alpha_{1}$ et $\alpha_{2}$ se restreignent \`a $A_{M_{\eta}}$ en des multiples $n_{1}\beta'$ et $n_{2}\beta'$ d'un m\^eme \'el\'ement indivisible de $\Sigma^{G_{\eta}}(A_{M_{\eta}})$. On obtient
$$nn_{1}B_{\eta}(\alpha_{1})^{-1}=n_{2}B_{\eta}(\alpha_{2})^{-1}.$$
Comme on vient de le dire, les entiers $n_{1}$ et $n_{2}$ sont premiers \`a $p$. Par hypoth\`ese, les valeurs $B_{\eta}(\alpha_{1})$ et $B_{\eta}(\alpha_{2})$ aussi. Donc $n$ est premier \`a $p$. Une preuve analogue montre qu'un rationnel $q$ intervenant dans (9) est premier \`a $p$. $\square$

 \bigskip
\subsection{Int\'egrales orbitales pond\'er\'ees invariantes stables}
On suppose ici que $(G,\tilde{G},{\bf a})$ est quasi-d\'eploy\'e et \`a torsion int\'erieure. On fixe toujours un espace de Levi $\tilde{M}$ de $\tilde{G}$ et on fixe un syst\`eme de fonctions $B$ comme en 1.9 . 
On a d\'efini l'espace $D^{st}_{g\acute{e}om}(\tilde{M}(F))$ des distributions g\'eom\'etriques stables sur $\tilde{M}(F)$. Pour $\boldsymbol{\delta}\in D^{st}_{g\acute{e}om}(\tilde{M}(F))\otimes Mes(M(F))^*$, on va d\'efinir une forme lin\'eaire ${\bf f}\mapsto S_{\tilde{M}}^{\tilde{G}}(\boldsymbol{\delta},B,{\bf f})$ sur $C_{c}^{\infty}(\tilde{G}(F))\otimes Mes(G(F))$. Conform\'ement \`a ce que l'on a dit en 1.1, la d\'efinition se fait par r\'ecurrence sur $dim(G_{SC})$.   Pour poser cette d\'efinition, il est n\'ecessaire de conna\^{\i}tre par r\'ecurrence certaines propri\'et\'es de cette forme lin\'eaire.

La propri\'et\'e difficile est

(1) la forme lin\'eaire ${\bf f}\mapsto S_{\tilde{M}}^{\tilde{G}}(\boldsymbol{\delta},B,{\bf f})$ est stable, c'est-\`a-dire se factorise en une forme  lin\'eaire sur $SI(\tilde{G}(F))\otimes Mes(G(F))$.

Les autres propri\'et\'es sont formelles et faciles.  Pour $b\in C_{c}^{\infty}(\tilde{{\cal A}}_{\tilde{G},F})$, on a

(2) $S_{\tilde{M}}^{\tilde{G}}(\boldsymbol{\delta},B,{\bf f}(b\circ \tilde{H}_{\tilde{G}}))=S_{\tilde{M}}^{\tilde{G}}(\boldsymbol{\delta},B,{\bf f})$ pourvu que $b$ vaille $1$ sur l'image par $\tilde{H}_{\tilde{G}}$ du support de $\boldsymbol{\delta}$.

Pour simplifier, fixons des mesures de Haar sur tous les groupes intervenant. Consid\'erons des extensions compatibles
$$1\to C_{\natural}\to G_{\natural}\to G\to 1\text{ et }\tilde{G}_{\natural}\to \tilde{G}$$
o\`u $C_{\natural}$ est un tore central induit et o\`u $\tilde{G}_{\natural}$ est encore \`a torsion int\'erieure. Soit $\lambda_{\natural}$ un caract\`ere de $C_{\natural}(F)$. On d\'efinit les espaces $C_{c,\lambda_{\natural}}^{\infty}(\tilde{G}_{\natural}(F))$ et $D_{g\acute{e}om,\lambda_{\natural}}(\tilde{G}_{\natural}(F))$. Il y a un homomorphisme naturel
$$(3)\qquad \begin{array}{ccc}D_{g\acute{e}om}(\tilde{G}_{\natural}(F))&\to& D_{g\acute{e}om,\lambda_{\natural}}(\tilde{G}_{\natural}(F))\\ \dot{\boldsymbol{\gamma}}&\mapsto&\boldsymbol{\gamma}\\ \end{array}$$
que l'on peut d\'efinir par la formule suivante. Pour $\dot{\boldsymbol{\gamma}}\in D_{g\acute{e}om}(\tilde{G}_{\natural}(F))$ et $f\in C_{c,\lambda_{\natural}}^{\infty}(\tilde{G}_{\natural}(F))$, on pose $I_{\tilde{G}_{\natural}}(\boldsymbol{\gamma},f)=I_{\tilde{G}_{\natural}}(\dot{\boldsymbol{\gamma}},f(b\circ \tilde{H}_{\tilde{G}_{\natural}}))$ o\`u $b$ est n'importe quel \'el\'ement de $C_{c}^{\infty}(\tilde{{\cal A}}_{\tilde{G}_{\natural},F})$ valant $1$ sur la projection du support de $\dot{\boldsymbol{\gamma}}$. Remarquons que $f(b\circ \tilde{H}_{\tilde{G}})$ est \`a support compact dans $\tilde{G}_{\natural}(F)$. La d\'efinition ci-dessus ne d\'epend pas du choix de $b$. Il est utile de donner une autre d\'efinition. L'int\'egration d\'efinit un homomorphisme
$$(4) \qquad \begin{array}{ccc}C_{c}^{\infty}(\tilde{G}_{\natural}(F))&\to&C_{c,\lambda_{\natural}}^{\infty}(\tilde{G}_{\natural}(F))\\ \dot{f}&\mapsto&f\\ \end{array}$$
Pr\'ecis\'ement, $f(\gamma')=\int_{C_{\natural}(F)}\dot{f}^c(\gamma')\lambda_{\natural}(c)dc$ pour tout $\gamma'\in \tilde{G}_{\natural}(F)$, o\`u $\dot{f}^c(\gamma')=\dot{f}(c\gamma')$. On remarque que l'ensemble des $c$ pour lesquels $\dot{f}^c(b\circ \tilde{H}_{\tilde{G}_{\natural}}) $ n'est pas nulle est compact. Il en r\'esulte que
$$I_{\tilde{G}_{\natural}}(\dot{\boldsymbol{\gamma}},f(b\circ \tilde{H}_{\tilde{G}_{\natural}}))=I_{\tilde{G}_{\natural}}(\dot{\boldsymbol{\gamma}},\int_{C_{\natural}(F)}\dot{f}^c(b\circ \tilde{H}_{\tilde{G}_{\natural}})\lambda_{\natural}(c)dc)=\int_{C_{\natural}(F)}I_{\tilde{G}_{\natural}}(\dot{\boldsymbol{\gamma}},\dot{f}^c(b\circ \tilde{H}_{\tilde{G}_{\natural}}))\lambda_{\natural}(c)dc$$
$$=\int_{C_{\natural}(F)}I_{\tilde{G}_{\natural}}(\dot{\boldsymbol{\gamma}},\dot{f}^c)\lambda_{\natural}(c)dc.$$
Autrement dit,
$$I_{\tilde{G}_{\natural}}(\boldsymbol{\gamma},f)=\int_{C_{\natural}(F)}I_{\tilde{G}_{\natural}}(\dot{\boldsymbol{\gamma}},\dot{f}^c)\lambda_{\natural}(c)dc.$$
L'homomorphisme (3) est surjectif (les int\'egrales orbitales qui engendrent l'espace d'arriv\'ee sont clairement dans l'image). Le groupe $C_{\natural}(F)$ agit sur $C_{c}^{\infty}(\tilde{G}_{\natural}(F))$ par $(c,\dot{f})\mapsto \dot{f}^c$. On en d\'eduit une action duale sur $D_{g\acute{e}om}(\tilde{G}_{\natural}(F))$ de sorte que $I_{\tilde{G}_{\natural}}(\dot{\boldsymbol{\gamma}}^c,\dot{f}^c)=I_{\tilde{G}_{\natural}}(\dot{\boldsymbol{\gamma}},\dot{f})$. On v\'erifie en utilisant la deuxi\`eme forme de la d\'efinition que le noyau de l'homomorphisme (3) est engendr\'e par les $\dot{\boldsymbol{\gamma}}^c-\lambda_{\natural}(c)\dot{\boldsymbol{\gamma}}$ pour $\dot{\boldsymbol{\gamma}}\in D_{g\acute{e}om}(\tilde{G}_{\natural}(F))$ et $c\in C_{\natural}(F)$. Il est peut-\^etre moins clair que

(5) $D_{g\acute{e}om,\lambda_{\natural}}^{st}(\tilde{G}_{\natural}(F))$ est l'image de $ D_{g\acute{e}om}^{st}(\tilde{G}_{\natural}(F))$ par l'homomorphisme (3).

Preuve. Soient $\dot{\boldsymbol{\delta}}\in D_{g\acute{e}om}^{st}(\tilde{G}_{\natural}(F))$ et $f\in C_{c,\lambda_{\natural}}^{\infty}(\tilde{G}_{\natural}(F))$ dont toutes les int\'egrales orbitales r\'eguli\`eres stables sont nulles. On a $I_{\tilde{G}_{\natural}}(\boldsymbol{\delta},f)=I_{\tilde{G}_{\natural}}(\dot{\boldsymbol{\delta}},f(b\circ \tilde{H}_{\tilde{G}}))$ o\`u $b$ est comme ci-dessus. Il  est clair que   toutes les int\'egrales orbitales r\'eguli\`eres stables de $f(b\circ \tilde{H}_{\tilde{G}})$ sont elles-aussi nulles. Donc $I_{\tilde{G}_{\natural}}(\dot{\boldsymbol{\delta}},f(b\circ \tilde{H}_{\tilde{G}}))=0$, donc $I_{\tilde{G}_{\natural}}(\boldsymbol{\delta},f)=0$  donc $\boldsymbol{\delta}$ est stable. Inversement, soit $\dot{\boldsymbol{\delta}}\in D_{g\acute{e}om}(\tilde{G}_{\natural}(F))$  tel que $\boldsymbol{\delta}$ soit stable. Introduisons le groupe d\'eriv\'e $G_{\natural,der}$. Le groupe $C_{\natural}(F)$ agit sur $G_{\natural,der}(F)\backslash\tilde{G}_{\natural}(F)$. En ajoutant \`a $\dot{\boldsymbol{\delta}}$ un \'el\'ement du noyau de l'homomorphisme (3), on peut supposer que l'image du support de $\dot{\boldsymbol{\delta}}$ dans  $G_{\natural,der}(F)\backslash\tilde{G}_{\natural}(F)$ est de la forme $\{x_{1},...,x_{n}\}$ o\`u, pour $i\not=j$, $x_{i}$ et $x_{j}$ ne sont pas dans la m\^eme orbite pour l'action de $C_{\natural}(F)$. L'intersection  $\Delta=C_{\natural}(F)\cap G_{\natural,der}(F)$ est finie. Quitte \`a moyenner sur ce groupe, ce qui ne change pas $\boldsymbol{\delta}$, on peut supposer $\dot{\boldsymbol{\delta}}^c=\lambda_{\natural}(c)\dot{\boldsymbol{\delta}}$ pour $c\in \Delta$. On va montrer qu'alors, $\dot{\boldsymbol{\delta}}$ est stable. Soit $\dot{f}\in C_{c}^{\infty}(\tilde{G}_{\natural}(F))$ dont toutes les int\'egrales orbitales r\'eguli\`eres stables sont nulles. Introduisons sa moyenne $\dot{f}_{0}$ sur le groupe $\Delta$ de sorte que $\dot{f}_{0}^c=\lambda_{\natural}(c)^{-1}\dot{f}_{0}$ pour $c\in \Delta$. Ses int\'egrales orbitales r\'eguli\`eres stables sont nulles elles-aussi.   Fixons un suppl\'ementaire $\mathfrak{s}$ de $\mathfrak{c}_{\natural}(F)$ dans $\mathfrak{z}_{G_{\natural}}(F)$ et un voisinage ouvert $\mathfrak{U}$ de $0$ dans $\mathfrak{s}$. Posons $U=exp(\mathfrak{U})$. Pour tout $i=1,...,n$, choisissons $\delta_{i}\in \tilde{G}_{\natural}(F)$ se projetant sur $x_{i}$ et consid\'erons l'application
$$\begin{array}{cccc}p_{i}:&U\times C_{\natural}(F)\times G_{\natural,der}(F)&\to&\tilde{G}_{\natural}(F)\\&(u,c,g)&\mapsto&ucg\delta_{i}.\\ \end{array}$$
Elle est continue et ouverte. En choisissant $\mathfrak{U}$ assez petit, on peut supposer que $p_{i}$ se quotiente en un isomorphisme de $\Delta\backslash (U\times C_{\natural}(F)\times G_{\natural,der}(F))$ sur son image, o\`u $\Delta$ agit sur $U\times C_{\natural}(F)\times G_{\natural,der}(F)$ via son plongement antidiagonal dans $C_{\natural}(F)\times G_{\natural,der}(F)$. On peut aussi supposer que, pour $i\not=j$, les images de $p_{i}$ et $p_{j}$ sont d'adh\'erences disjointes. Fixons une fonction $\varphi\in C_{c}^{\infty}(U)$ telle $\varphi(1)=1$. D\'efinissons une fonction $f_{i}$ sur $U\times C_{\natural}(F)\times G_{\natural,der}(F)$ par
$$f_{i}(u,c,g)=\lambda_{\natural}(c)^{-1}\varphi(u)\dot{f}_{0}(g\delta_{i}).$$
A cause de la propri\'et\'e de transformation de $\dot{f}_{0}$ par $\Delta$, $f_{i}$ se factorise par l'application $p_{i}$. On peut donc d\'efinir une fonction $f$ sur $\tilde{G}_{\natural}(F)$ qui est nulle hors de la r\'eunion des images des $p_{i}$ et qui v\'erifie $f\circ p_{i}=f_{i}$ pour tout $i$. Cette fonction  appartient \`a $C_{c,\lambda_{\natural}}^{\infty}(\tilde{G}_{\natural}(F))$. Les int\'egrales orbitales r\'eguli\`eres stables de $f$ sont nulles. En effet, cette condition se lit sur les fibres de l'application
$$\tilde{G}_{\natural}(F)\to G_{\natural,der}(F)\backslash \tilde{G}_{\natural}(F).$$
Or, sur une telle fibre, $f$ co\"{\i}ncide \`a une translation pr\`es avec un multiple de la restriction de $\dot{f}_{0}$ \`a une fibre au-dessus de l'un des $x_{i}$. Par ailleurs, sur une telle fibre au-dessus de l'un des $x_{i}$, $f$ co\"{\i}ncide exactement avec $\dot{f}_{0}$. En appliquant la d\'efinition, on en d\'eduit que $I_{\tilde{G}_{\natural}}(\boldsymbol{\delta},f)=I^{\tilde{G}}(\dot{\boldsymbol{\delta}},\dot{f}_{0})$. Puisque $\boldsymbol{\delta}$ est stable, le premier membre de cette \'egalit\'e est nul, donc aussi le deuxi\`eme. Puisque l'on a pris soin de moyenner $\dot{\boldsymbol{\delta}}$, on a aussi $ I^{\tilde{G}}(\dot{\boldsymbol{\delta}},\dot{f}_{0})=I^{\tilde{G}}(\dot{\boldsymbol{\delta}},\dot{f})$. Donc $I^{\tilde{G}}(\dot{\boldsymbol{\delta}},\dot{f})=0$, ce qu'il fallait d\'emontrer. $\square$

Soient $\dot{\eta}\in \tilde{G}_{\natural,ss}(F)$ et $\eta$ son image dans $\tilde{G}_{\natural}(F)$. Les groupes $G_{\natural,\dot{\eta}}$ et $G_{\eta}$ ont le m\^eme syst\`eme de racines. On peut identifier la fonction $B_{\eta}$  \`a une fonction sur le syst\`eme de racines du premier groupe. Cela munit $\tilde{G}_{\natural}$ d'un syst\`eme de fonctions que nous notons encore $B$.

Soit   $\tilde{M}_{\natural}$ l'espace de Levi   de $\tilde{G}_{\natural}$ associ\'e \`a $\tilde{M}$.  Pour $\dot{\boldsymbol{\delta}}\in D_{g\acute{e}om}^{st}(\tilde{M}_{\natural}(F))$ et $f\in C_{c,\lambda_{\natural}}^{\infty}(\tilde{G}_{\natural}(F))$, posons
$$S_{\tilde{M}_{\natural}}^{\tilde{G}_{\natural}}(\dot{\boldsymbol{\delta}},B,f)=S_{\tilde{M}_{\natural}}^{\tilde{G}_{\natural}}(\dot{\boldsymbol{\delta}},B,f(b\circ \tilde{H}_{\tilde{G}_{\natural}})),$$
o\`u $b\in C_{c}^{\infty}(\tilde{{\cal A}}_{\tilde{G}_{\natural},F})$ vaut $1$ sur la projection du support de $\boldsymbol{\delta}$. Cette d\'efinition est loisible d'apr\`es la propri\'et\'e (2). Remarquons que le m\^eme raisonnement conduisant \`a la deuxi\`eme d\'efinition de l'homomorphisme (3) conduit \`a une deuxi\`eme forme de la d\'efinition ci-dessus:
$$S_{\tilde{M}_{\natural}}^{\tilde{G}_{\natural}}(\dot{\boldsymbol{\delta}},B,f)=\int_{C_{\natural}(F)}S_{\tilde{M}_{\natural}}^{\tilde{G}_{\natural}}(\dot{\boldsymbol{\delta}},B,\dot{f}^c)\lambda_{\natural}(c)dc,$$
o\`u $\dot{f}$ est reli\'e \`a $f$ par (4).  L'une des propri\'et\'es requises est:

(6) $S_{\tilde{M}_{\natural}}^{\tilde{G}_{\natural}}(\dot{\boldsymbol{\delta}},B,f)$ ne d\'epend que de l'image $\boldsymbol{\delta}$ de $\dot{\boldsymbol{\delta}}$ dans $  D_{g\acute{e}om,\lambda_{\natural}}^{st}(\tilde{M}_{\natural}(F))$. 

En utilisant (5) et (6), on peut d\'efinir $S_{\tilde{M}_{\natural},\lambda_{\natural}}^{\tilde{G}_{\natural}}(\boldsymbol{\delta},B,f)$ pour $\boldsymbol{\delta}\in D_{g\acute{e}om,\lambda_{\natural}}^{st}(\tilde{M}_{\natural}(F))$ et $f\in  C_{c,\lambda_{\natural}}^{\infty}(\tilde{G}_{\natural}(F))$, par
$$S_{\tilde{M}_{\natural},\lambda_{\natural}}^{\tilde{G}_{\natural}}(\boldsymbol{\delta},B,f)=S_{\tilde{M}_{\natural}}^{\tilde{G}_{\natural}}(\dot{\boldsymbol{\delta}},B,f)$$
o\`u $\dot{\boldsymbol{\delta}}$ est n'importe quel \'el\'ement de $D_{g\acute{e}om}^{st}(\tilde{M}_{\natural}(F))$  d'image $\boldsymbol{\delta}$. La propri\'et\'e (1) reste v\'erifi\'ee dans cette situation plus g\'en\'erale.

Consid\'erons maintenant deux  couples d'extensions
$$1\to C_{\natural}\to G_{\natural}\to G\to 1\text{ et }\tilde{G}_{\natural}\to \tilde{G}$$
$$1\to C_{\flat}\to G_{\flat}\to G\to 1\text{ et }\tilde{G}_{\flat}\to \tilde{G}$$
et deux caract\`eres $\lambda_{\natural}$ et $\lambda_{\flat}$ v\'erifiant les hypoth\`eses pr\'ec\'edentes. Introduisons le produit fibr\'e $G_{\natural,\flat}$ de $G_{\natural} $ et $G_{\flat}$ au-dessus de $G$ et le produit fibr\'e $\tilde{G}_{\natural,\flat}$ de $\tilde{G}_{\natural}$ et $\tilde{G}_{\flat}$ au-dessus de $\tilde{G}$. Supposons donn\'e un caract\`ere $\lambda_{\natural,\flat}$ de $G_{\natural,\flat}(F)$ dont la restriction \`a $C_{\natural}(F)\times C_{\flat}(F)$ soit $\lambda_{\natural}\times \lambda_{\flat}^{-1}$. Supposons donn\'ee une fonction non nulle $\tilde{\lambda}_{\natural,\flat}$ sur $\tilde{G}_{\natural,\flat}(F)$ qui se transforme selon le caract\`ere $\lambda_{\natural,\flat}$, cf. [I] 2.5(i). On d\'efinit un isomorphisme
$$\begin{array}{ccc}C_{c,\lambda_{\natural}}^{\infty}(\tilde{G}_{\natural}(F))&\to&C_{c,\lambda_{\flat}}^{\infty}(\tilde{G}_{\flat}(F))\\ f_{\natural}&\mapsto&f_{\flat}\\ \end{array}$$
par $f_{\flat}(\gamma_{\flat})=\tilde{\lambda}_{\natural,\flat}(\gamma_{\natural},\gamma_{\flat})f_{\natural}(\gamma_{\natural})$, o\`u $\gamma_{\natural}$ est n'importe quel \'el\'ement de $\tilde{G}_{\natural} (F)$ tel que $(\gamma_{\natural},\gamma_{\flat})\in \tilde{G}_{\natural,\flat}(F)$. Par restriction \`a $\tilde{M}$ puis dualit\'e, on a aussi un isomorphisme de $D_{g\acute{e}om,\lambda_{\natural}}(\tilde{M}_{\natural}(F))$ sur $D_{g\acute{e}om,\lambda_{\flat}}(\tilde{M}_{\flat}(F))$, qui se restreint en un isomorphisme entre espace de distributions stables. Pour $f_{\natural}$ et $f_{\flat}$, resp. $\boldsymbol{\delta}_{\natural}$ et $\boldsymbol{\delta}_{\flat}$, se correspondant par ces isomorphismes, on veut que

(7) $S_{M_{\natural},\lambda_{\natural}}^{G_{\natural}}(\boldsymbol{\delta}_{\natural},B,f_{\natural})=S_{M_{\flat},\lambda_{\flat}}^{G_{\flat}}(\boldsymbol{\delta}_{\flat},B,f_{\flat})$.

Soient ${\bf G}'=(G',\tilde{G}',s)$ une donn\'ee endoscopique de $(G,\tilde{G})$ (on oublie ${\bf a}$ qui est trivial) et un Levi de $G'$ associ\'e \`a $\tilde{M}$ et \`a sa donn\'ee endoscopique maximale ${\bf M}$ (cf. [I] 1.7).  On note encore $M$ ce Levi. Soit $\epsilon\in \tilde{G}'(F)$ un \'el\'ement semi-simple. Alors il lui correspond un \'el\'ement semi-simple $\eta\in \tilde{G}(F)$ (la preuve est la m\^eme que celle du lemme 1.10  de [I]).  Le syst\`eme de racines de $G'_{\epsilon}$ est un sous-syst\`eme de celui de $G_{\eta}$. On le munit de la restriction de la fonction $B_{\eta}$. On obtient ainsi un syst\`eme de fonctions sur $G'(F)$ v\'erifiant encore les hypoth\`eses de 1.9. On note encore $B$ ce syst\`eme de fonctions. Dans la d\'efinition des int\'egrales orbitales pond\'er\'ees pour $\tilde{G}$ intervient une mesure que l'on a d\'eduite en 1.2 d'une forme quadratique sur $X_{*}(T^*)\otimes {\mathbb R}$. Il convient de faire un choix analogue pour $\tilde{G}'$. Si ${\bf G}'$ n'est pas elliptique, ce choix n'importe pas. Si ${\bf G}'$ est elliptique,  on remarque que, sur $\bar{F}$, on peut identifier un tore maximal de $G'$ \`a $T^*$. On choisit alors les mesures qui se d\'eduisent de la m\^eme forme quadratique sur $X_{*}(T^*)\otimes {\mathbb R}$.
Soient $\boldsymbol{\delta}\in D_{g\acute{e}om}^{st}({\bf M})$ et ${\bf f}\in SI({\bf G}')$. Fixons des donn\'ees auxiliaires $G'_{1},...,\Delta_{1}$. Alors $\boldsymbol{\delta}$ s'identifie \`a un \'el\'ement de $D_{g\acute{e}om,\lambda_{1}}(\tilde{G}'_{1}(F))$ et ${\bf f}$ s'identifie \`a un \'el\'ement de $SI_{\lambda_{1}}(\tilde{G}'_{1}(F))$. En vertu de (1) et (6), le terme $S_{\tilde{M}_{1},\lambda_{1}}^{\tilde{G}'_{1}}(\boldsymbol{\delta},B,{\bf f})$ est d\'efini. En vertu de (7), il ne d\'epend pas du choix des donn\'ees auxiliaires. On pose
$$S_{{\bf M}}^{{\bf G}}(\boldsymbol{\delta},B,{\bf f})=S_{\tilde{M}_{1},\lambda_{1}}^{\tilde{G}'_{1}}(\boldsymbol{\delta},B,{\bf f}).$$

On r\'etablit maintenant les espaces de mesures pour donner des d\'efinitions plus canoniques. Consid\'erons  une paire de Borel \'epingl\'ee $\hat{{\cal E}}=(\hat{B},\hat{T},(\hat{E}_{\alpha})_{\alpha\in \Delta})$ de $\hat{G}$. Comme on l'a dit en [I] 1.4, on peut modifier l'action galoisienne de $\Gamma_{F}$ sur $\hat{G}$ de sorte qu'elle pr\'eserve cette paire et on peut introduire l'\'el\'ement $\hat{\theta}\in \hat{G}\boldsymbol{\hat{\theta}}$  tel que $ad_{\hat{\theta}}$ conserve cette paire. En choisissant convenablement celle-ci, on peut supposer  que $\hat{M}$ est un Levi standard et que le $L$-espace $^L\tilde{M}$ est \'egal \`a $(\hat{M}\ltimes W_{F})\hat{\theta}$. Pour $s\in Z(\hat{M})^{\Gamma_{F}}/Z(\hat{G})^{\Gamma_{F}}$, on a d\'efini en [I] 3.3 la donn\'ee endoscopique ${\bf G}'(s)$ qui v\'erifie les hypoth\`eses ci-dessus. On pose
 $$i_{\tilde{M}}(\tilde{G},\tilde{G}'(s))=\left\lbrace\begin{array}{cc}[Z(\hat{G}'(s))^{\Gamma_{F}}:Z(\hat{G})^{\Gamma_{F}}]^{-1},& \text{ si }{\bf G}'(s)\text{ est elliptique,}\\ 0,& \text{ sinon.}\\ \end{array}\right.$$
   Une d\'efinition plus g\'en\'erale sera donn\'ee en 1.12.  
   
   Apr\`es tous ces pr\'eliminaires, on peut d\'efinir $S_{\tilde{M}}^{\tilde{G}}(\boldsymbol{\delta},B,{\bf f})$ pour $\boldsymbol{\delta}\in D^{st}_{g\acute{e}om}(\tilde{M}(F))\otimes Mes(M(F))^*$ et ${\bf f}\in 
C_{c}^{\infty}(\tilde{G}(F))\otimes Mes(G(F))$ par l'\'egalit\'e
$$(8) \qquad S_{\tilde{M}}^{\tilde{G}}(\boldsymbol{\delta},B,{\bf f})=I_{\tilde{M}}^{\tilde{G}}(\boldsymbol{\delta},B,{\bf f})-\sum_{s\in Z(\hat{M})^{\Gamma_{F}}/Z(\hat{G})^{\Gamma_{F}}; s\not=1}i_{\tilde{M}}(\tilde{G},\tilde{G}'(s))S_{{\bf M}}^{{\bf G}'(s)}(\boldsymbol{\delta},B,{\bf f}^{{\bf G}'(s)}).$$
Tous les termes du membre de droite ont \'et\'e d\'efinis gr\^ace aux hypoth\`eses de r\'ecurrence.

On doit montrer que le terme ainsi d\'efini v\'erifie lui-m\^eme ces hypoth\`eses.  On va le faire ci-dessous en ce qui concerne les propri\'et\'es formelles. La propri\'et\'e (1) est \'evidemment plus difficile. Formulons-la provisoirement sous la forme d'un th\'eor\`eme \`a prouver.

\ass{Th\'eor\`eme  (\`a prouver)}{Pour $\boldsymbol{\delta}\in D^{st}_{g\acute{e}om}(\tilde{M}(F))\otimes Mes(M(F))^*$, la distribution ${\bf f}\mapsto S_{\tilde{M}}^{\tilde{G}}(\boldsymbol{\delta},B,{\bf f})$ est stable.}

{\bf Remarques.} (i) Si $\tilde{M}=\tilde{G}$, on a simplement $S_{\tilde{G}}^{\tilde{G}}(\boldsymbol{\delta},B,{\bf f})=I^{\tilde{G}}(\boldsymbol{\delta},{\bf f})$ et l'assertion du th\'eor\`eme est tautologique.

(ii) S'il n'y a pas du tout de torsion, c'est-\`a-dire si $\tilde{G}=G$,   et si de plus $B_{\eta}$ est constante de valeur $1$ pour tout $\eta\in \tilde{G}_{ss}(F)$, le th\'eor\`eme a \'et\'e prouv\'e par Arthur pour les \'el\'ements $\boldsymbol{\delta}$ dont le support est form\'e d'\'el\'ements fortement $\tilde{G}$-r\'eguliers ([A6] local theorem 1(b)). Nous prouverons dans l'article suivant  que le th\'eor\`eme ci-dessus se d\'eduit de celui d'Arthur. 
\bigskip

La v\'erification des propri\'et\'es formelles est fastidieuse mais il est peut-\^etre bon de la faire tout-de-m\^eme. Dans la suite, on ne fera plus de telles v\'erifications. 

 V\'erifions (2). Soit ${\bf G}'$ une donn\'ee endoscopique relevante et elliptique de $(G,\tilde{G})$. Appliquons [I] 1.12  en se rappelant que le groupe $G_{0}$ de ce paragraphe est \'egal \`a $G$ puisque $(G,\tilde{G},{\bf a})$ est quasi-d\'eploy\'e et sans torsion. On obtient un homomorphisme $N^{G',G}:G'_{ab}(F)\to G_{ab}(F)$ et une application $N^{\tilde{G}',\tilde{G}}:\tilde{G}'_{ab}(F)\to \tilde{G}_{ab}(F)$ compatible \`a cet homomorphisme.  Les applications  $H_{\tilde{G}}$  et $\tilde{H}_{\tilde{G}}$ d\'efinies sur $G(F)$ et $\tilde{G}(F)$ se factorisent par $G_{ab}(F)$ et $\tilde{G}_{ab}(F)$ et il y a bien s\^ur une assertion analogue pour les applications $H_{\tilde{G}'}$ et $\tilde{H}_{\tilde{G}'}$. Par ailleurs, il y a un isomorphisme ${\cal A}_{\tilde{G}'}\simeq {\cal A}_{\tilde{G}}$ puisque ${\bf G}'$ est elliptique. En reprenant les d\'efinitions, on voit qu'il y a un diagramme commutatif
 $$\begin{array}{ccc}G'_{ab}(F)&\stackrel{N^{G',G}}{\to}&G_{ab}(F)\\ H_{\tilde{G}'}\downarrow\,\,&& \,\, \downarrow H_{\tilde{G}}\\ {\cal A}_{\tilde{G}',F}&\to&{\cal A}_{\tilde{G},F}\\ \end{array}$$
 o\`u l'homomorphisme horizontal du bas est la restriction de l'isomorphisme ${\cal A}_{\tilde{G}'}\simeq {\cal A}_{\tilde{G}}$. On en d\'eduit qu'il y a un diagramme commutatif similaire
$$\begin{array}{ccc}\tilde{G}'_{ab}(F)&\stackrel{N^{\tilde{G}',\tilde{G}}}{\to}&\tilde{G}_{ab}(F)\\ \tilde{H}_{\tilde{G}'}\downarrow\,\,&& \,\, \downarrow \tilde{H}_{\tilde{G}}\\ \tilde{{\cal A}}_{\tilde{G}',F}&\to&\tilde{{\cal A}}_{\tilde{G},F}\\ \end{array}$$
o\`u la fl\`eche horizontale du bas est compatible \`a l'homomorphisme du diagramme pr\'ec\'edent. En particulier, elle est injective. Soient $\boldsymbol{\delta}$, ${\bf f}$ et $b$ comme dans la relation (2). Pour $s\in  Z(\hat{M})^{\Gamma_{F}}/Z(\hat{G})^{\Gamma_{F}}$, on v\'erifie que le transfert $({\bf f}(b\circ \tilde{H}_{\tilde{G}}))^{{\bf G}'(s)}$ est \'egal \`a ${\bf f}^{{\bf G}'(s)}(b\circ \tilde{H}_{\tilde{G}'(s)})$. Pour $s\not=1$, les hypoth\`eses de r\'ecurrence assurent que
$$S_{{\bf M}}^{{\bf G}'(s)}(\boldsymbol{\delta},B,{\bf f}^{{\bf G}'(s)}(b\circ \tilde{H}_{\tilde{G}'(s)}))=S_{{\bf M}}^{{\bf G}'(s)}(\boldsymbol{\delta},B,{\bf f}^{{\bf G}'(s)}).$$
D'apr\`es 1.7(1), on a aussi
$$I_{\tilde{M}}^{\tilde{G}}(\boldsymbol{\delta},B,{\bf f}(b\circ \tilde{H}_{\tilde{G}}))=I_{\tilde{M}}^{\tilde{G}}(\boldsymbol{\delta},B,{\bf f}).$$
Il suffit d'appliquer la relation (8) \`a ${\bf f}$ et ${\bf f}(b\circ H_{\tilde{G}})$ pour obtenir la relation (2).

V\'erifions (6). Soient $\dot{\boldsymbol{\delta}}$ et $\dot{\boldsymbol{\delta}}'$ deux \'el\'ements de $D^{st}_{g\acute{e}om}(\tilde{M}_{\natural}(F))$ ayant m\^eme image dans  $D^{st}_{g\acute{e}om,\lambda_{\natural}}(\tilde{M}_{\natural}(F))$. Soit $f\in C_{c,\lambda_{\natural}}^{\infty}(\tilde{G}_{\natural}(F))$. On veut montrer que $S_{\tilde{M}_{\natural}}^{\tilde{G}_{\natural}}(\dot{\boldsymbol{\delta}},B,f)=S_{\tilde{M}_{\natural}}^{\tilde{G}_{\natural}}(\dot{\boldsymbol{\delta}}',B,f)$. 
On choisit $\dot{f}$ reli\'e \`a $f$ par (4). Montrons que
$$(9) \qquad \int_{C_{\natural}(F)}I_{\tilde{M}_{\natural}}^{\tilde{G}_{\natural}}(\dot{\boldsymbol{\delta}},B,\dot{f}^c)\lambda_{\natural}(c)dc= \int_{C_{\natural}(F)}I_{\tilde{M}_{\natural}}^{\tilde{G}_{\natural}}(\dot{\boldsymbol{\delta}}',B,\dot{f}^c)\lambda_{\natural}(c)dc.$$
D'apr\`es la description du noyau de l'homomorphisme (3), $\dot{\boldsymbol{\delta}}-\dot{\boldsymbol{\delta}}'$ est une somme de termes $\dot{\boldsymbol{\gamma}}^c-\lambda_{\natural}(c)\dot{\boldsymbol{\gamma}}$, avec  $\dot{\boldsymbol{\gamma}}\in D_{g\acute{e}om}(\tilde{M}_{\natural}(F))$ et $c\in C_{\natural}(F)$. Alors (9) r\'esulte de l'\'egalit\'e
$$I_{\tilde{M}_{\natural}}^{\tilde{G}_{\natural}}(\dot{\boldsymbol{\gamma}}^c,B,\dot{f}^c)=I_{\tilde{M}_{\natural}}^{\tilde{G}_{\natural}}(\dot{\boldsymbol{\gamma}},B,\dot{f})$$
pour $\dot{\boldsymbol{\gamma}}$ et $c$ comme ci-dessus. On peut supposer que $\dot{\boldsymbol{\gamma}}$ est une int\'egrale orbitale. La relation 1.9(5) nous ram\`ene alors \`a prouver l'\'egalit\'e ci-dessus pour le syst\`eme de fonctions $B$ dont toutes les valeurs sont \'egales \`a $1$. Dans ce cas, l'\'egalit\'e  r\'esulte de la m\^eme \'egalit\'e pour les int\'egrales orbitales pond\'er\'ees non invariantes (qui est triviale) et de la relation $\phi_{\tilde{L}_{\natural}}(\dot{f}^c)=(\phi_{\tilde{L}_{\natural}}(\dot{f}))^c$ pour tout $\tilde{L}_{\natural}\in {\cal L}(\tilde{M}_{\natural})$.  Cette propri\'et\'e   r\'esulte imm\'ediatement de la d\'efinition de l'application $\phi_{\tilde{L}_{\natural}}$.

On a la suite exacte
$$1\to Z(\hat{G})\to Z(\hat{G}_{\natural})\to Z(\hat{C}_{\natural})\to 1$$
Le groupe $Z(\hat{C}_{\natural})^{\Gamma_{F}}$ est connexe puisque $C_{\natural}$ est induit. La suite d'invariants
$$1\to Z(\hat{G})^{\Gamma_{F}}\to Z(\hat{G}_{\natural})^{\Gamma_{F}}\to Z(\hat{C}_{\natural})^{\Gamma_{F}}\to 1$$
est donc encore exacte. On a une suite analogue en rempla\c{c}ant $\hat{G}$ par $\hat{M}$ et $\hat{G}_{\natural}$ par $\hat{M}_{\natural}$. Puisque $Z(\hat{G}_{\natural})^{\Gamma_{F}}$ se projette surjectivement sur $ Z(\hat{C}_{\natural})^{\Gamma_{F}}$, on en d\'eduit l'\'egalit\'e
$$Z(\hat{M}_{\natural})^{\Gamma_{F}}=Z(\hat{M})^{\Gamma_{F}}Z(\hat{G}_{\natural})^{\Gamma_{F}}.$$
Autrement dit, l'homomorphisme
$$Z(\hat{M})^{\Gamma_{F}}/Z(\hat{G})^{\Gamma_{F}}\to Z(\hat{M}_{\natural})^{\Gamma_{F}}/Z(\hat{G}_{\natural})^{\Gamma_{F}}$$
est surjectif. Il est aussi injectif, donc bijectif. Un \'el\'ement $s\in Z(\hat{M})^{\Gamma_{F}}/Z(\hat{G})^{\Gamma_{F}}$ d\'efinit donc \`a la fois une donn\'ee endoscopique ${\bf G}'(s)$ de $(G,\tilde{G})$ et une donn\'ee endoscopique ${\bf G}'_{\natural}(s)$ de $(G_{\natural},\tilde{G}_{\natural})$. On a une suite exacte
$$1\to C_{\natural}\to G'_{\natural}(s)\to G'(s)\to 1$$
et une application compatible
$$\tilde{G}'_{\natural}(s)\to \tilde{G}'(s).$$
Par un calcul similaire \`a celui ci-dessus, on montre que $i_{\tilde{M}_{\natural}}(\tilde{G}_{\natural},\tilde{G}'_{\natural}(s))=i_{\tilde{M}}(\tilde{G},\tilde{G}'(s))$. On simplifie les calculs ult\'erieurs en remarquant que, pour la donn\'ee ${\bf G}'(s)$, on peut choisir des donn\'ees auxiliaires $G'(s)_{1},...,\Delta(s)_{1}$ telles que $G'(s)_{1}=G'(s)$, $\tilde{G}'(s)_{1}=\tilde{G}'(s)$, $C(s)_{1}=\{1\}$. Pour le prouver, il suffit de montrer que 

(10) ${\cal G}'(s)$ est isomorphe \`a $^LG'(s)$. 

{\bf Remarque.} A premi\`ere vue, cela para\^{\i}t \'evident puisque ${\cal G}'(s)$ est \'egal au sous-groupe $\hat{G}'(s)\rtimes W_{F}$ de $^LG$. Mais  l'action de $W_{F}$ sur $\hat{G}'(s)$ est ici la restriction de l'action sur $\hat{G}$. Elle n'est pas \'equivalente, en g\'en\'eral, \`a l'action sur $\hat{G}'(s)$ consid\'er\'e comme $L$-groupe de $G'(s)$. Plus exactement, elle ne conserve pas, en g\'en\'eral, un \'epinglage de $\hat{G}'(s)$ (contre-exemple: $G=U(3)$, $s$ tel que $G'(s)=U(2)\times U(1)$). 

\bigskip

Preuve de (10).   Soit $\hat{P}$ le sous-groupe parabolique standard de Levi $\hat{M}$. Alors $\hat{P}\cap \hat{G}'(s)$ est un sous-groupe parabolique de $\hat{G}'(s)$, de Levi $\hat{M}$, et il est conserv\'e par l'action de $W_{F}$ (la restriction de celle sur $\hat{G}$). On prend  pour $\hat{B}'$ l'unique Borel contenu dans $\hat{P}\cap \hat{G}'(s)$ qui a m\^eme intersection avec $\hat{M}$ que $\hat{B}$. On prend $\hat{T}'=\hat{T}$. On prend pour \'epinglage $(\hat{E}'_{\alpha})_{\alpha\in \hat{\Delta}'(s)}$ un \'epinglage quelconque contenant $(\hat{E}_{\alpha})_{\alpha\in \hat{\Delta}^M}$, o\`u $\hat{\Delta}^M$ est le sous-ensemble de $\hat{\Delta}$ associ\'e \`a $\hat{M}$. L'action de $W_{F}$ conserve $\hat{B}'$, $\hat{T}'$ et le sous-ensemble $(\hat{E}'_{\alpha})_{\alpha\in \hat{\Delta}^M}$. Elle ne conserve pas, en g\'en\'eral, le compl\'ementaire $(\hat{E}'_{\alpha})_{\alpha\in \hat{\Delta}'(s)-\hat{\Delta}^M}$. Mais  il existe un unique cocycle $\chi_{ad}:W_{F}\to Z(\hat{M})/Z(\hat{G}'(s))$ tel que l'action $w\mapsto ad_{ \chi_{ad}(w)}w_{G}$ conserve cet \'epinglage. On peut  supposer que l'action de $W_{F}$ sur $\hat{G}'(s)$ consid\'er\'e comme le $L$-groupe de $G'(s)$ est $w\mapsto w_{G'(s)}=ad_{\chi_{ad}(w)}w_{G}$. Supposons prouv\'e que $\chi_{ad}$ se rel\`eve en un cocycle $\chi:W_{F}\to Z(\hat{M})$. On d\'efinit alors une application
$$\begin{array}{cccc}\hat{\xi}(s):&{\cal G}'(s)\simeq \hat{G}'(s)\rtimes W_{F}&\to &{^LG}'(s)\simeq \hat{G}'(s)\rtimes W_{F}\\ &(x,w)&\mapsto (x\chi(w)^{-1},w)\\ \end{array}$$
(les deux produits semi-directs sont relatifs aux deux actions  de $W_{F}$).  C'est un isomorphisme, ce qui prouve (10). Il reste \`a prouver l'assertion de rel\`evement. On va en fait prouver que $\chi_{ad}$ se rel\`eve en un cocycle $\chi_{sc}:W_{F}\to Z(\hat{M}_{sc})$. Supposons que la condition suivante soit satisfaite:

(11) pour toute racine $\alpha$ de $\hat{T}$ dans $\hat{G}$ et pour tout $w\in W_{F}$ fixant $\alpha$, l'action $w_{G}$ sur l'espace radiciel associ\'e \`a $\alpha$ soit l'identit\'e.

Dans ce cas, on modifie l'ensemble $(\hat{E}'_{\alpha})_{\alpha\in \hat{\Delta}'(s)-\hat{\Delta}^M}$ de la fa\c{c}on suivante. On fixe un ensemble de repr\'esentants $\hat{\Delta}'_{1}$ des orbites de $W_{F}$ dans $\hat{\Delta}'(s)-\hat{\Delta}^M$. On fixe arbitrairement $\hat{E}'_{\alpha}$ pour $\alpha\in \hat{\Delta}'_{1}$. Pour $\alpha\in \hat{\Delta}'(s)-\hat{\Delta}^M$ quelconque, on \'ecrit $\alpha=w_{G}\alpha_{1}$ avec $w\in W_{F}$ et $\alpha_{1}\in \hat{\Delta}'_{1}$ et on pose $\hat{E}'_{\alpha}=w_{G}(\hat{E}'_{\alpha_{1}})$. L'hypoth\`ese (11) assure que cette d\'efinition est loisible. L'\'epinglage obtenu est conserv\'e par $W_{F}$. Notons que changer d'\'epinglage ne change pas la classe du cocycle $\chi_{ad}$. Donc cette classe est triviale. Cela entra\^{\i}ne que $\chi_{ad}$ se rel\`eve bien en un cocycle $\chi_{sc}$. Revenons maintenant au cas g\'en\'eral. L'action galoisienne permute les composantes simples du groupe $\hat{G}_{AD}$ et notre probl\`eme se ram\`ene  imm\'ediatement au cas o\`u cette action est transitive, donc toutes les composantes simples sont du m\^eme type. L'action galoisienne se fait par automorphismes pr\'eservant une paire de Borel \'epingl\'ee. On sait qu'\`a une exception pr\`es, l'hypoth\`ese (11) est satisfaite (auquel cas le probl\`eme est r\'esolu), cf. [KS] 1.3. D\'ecrivons l'exception.
On consid\`ere une tour  d'extensions $F_{2}/F_{1}/F$, avec $F_{2}/F_{1}$ quadratique,  et un groupe lin\'eaire adjoint $\hat{G}_{1,AD}=PGL(2n+1,{\mathbb C})$ muni de l'action de $\Gamma_{F_{1}}$ pour laquelle un \'el\'ement de $\Gamma_{F_{1}}-\Gamma_{F_{2}}$ agit par un automorphisme ext\'erieur non trivial.  Le groupe $\hat{G}_{AD}$ est d\'eduit de $\hat{G}_{1,AD}$ par changement de base de $F_{1}$ \`a $F$. Dans ce cas, il y a des couples $(\alpha,w)$ v\'erifiant les hypoth\`eses de (11) tels que l'action $w_{G}$ sur l'espace radiciel associ\'e \`a $\alpha$ soit moins l'identit\'e. En tout cas, la classe de $\chi_{ad}$ est d'ordre au plus $2$. Puisque $\hat{G}$ n'intervient pr\'esentement que via $\hat{G}_{AD}$ et $\hat{G}_{SC}$, on peut supposer que $\hat{G}$ est d\'eduit par changement de base d'un groupe $\hat{G}=GL(2n+1,{\mathbb C})$ muni de l'action similaire \`a celle ci-dessus. Le groupe $\hat{G}'(s)$ est alors aussi un produit de groupes $GL(k,{\mathbb C})$ et son centre est connexe. De cette connexit\'e r\'esulte que l'homorphisme
$$H^1(W_{F},Z(\hat{M}))\to H^1(W_{F},Z(\hat{M})/Z(\hat{G}'(s)))$$
est surjectif ([Lan] p. 719 (1)). Relevons $\chi_{ad}$ en un cocycle $\chi:W_{F}\to Z(\hat{M})$.
De l'homomorphisme
$$GL(2n+1,{\mathbb C})\stackrel{det}{\to}GL(1,{\mathbb C}) \simeq Z(GL(2n+1,{\mathbb C})) $$
se d\'eduit un homomorphisme
$$\hat{G}\stackrel{det}{\to} Z(\hat{G}) .$$
Posons $\chi_{sc}=( det\circ \chi)^{-1}\chi^{2n+1}$. Parce que $\chi_{ad}$ est au plus d'ordre $2$, ce cocycle $\chi_{sc}$ rel\`eve encore la classe de $\chi_{ad}$. Mais $\chi_{sc}$ prend ses valeurs dans $Z(\hat{M}_{sc})$. Cela ach\`eve la preuve de (10). $\square$
 
 Choisissons donc des donn\'ees auxiliaires simples  $G'(s)_{1}=G'(s)$, $\tilde{G}'(s)_{1}=\tilde{G}'(s)$, $C(s)_{1}=\{1\}$, $\hat{\xi}(s)_{1}$, $\Delta(s)_{1}$. La preuve ci-dessus montre qu'il y a vraiment un choix, $\hat{\xi}(s)_{1}$ n'est pas canonique. Parce qu'on a d\^u tordre l'action galoisienne par un cocycle \`a valeurs dans $Z(\hat{M})$, on ne peut pas en g\'en\'eral choisir $\Delta(s)_{1}$ \'egal \`a $1$ sur la diagonale dans $\tilde{M}(F)\times \tilde{M}(F)$. Le cocycle d\'efinit un caract\`ere $\chi_{F}$ de $M(F)$. En fixant un point base $\gamma\in \tilde{M}(F)$, on a une relation $\Delta(s)_{1}(x\gamma,x\gamma)=\chi_{F}(x)\Delta(s)_{1}(\gamma,\gamma)$ pour tout $x\in M(F)$. Quoi qu'il en soit, toute distribution $\boldsymbol{\delta}$ sur $\tilde{M}(F)$, vu comme sous-groupe de $\tilde{G}(F)$, est le transfert d'une distribution $\boldsymbol{\delta}(s)$ sur $\tilde{M}(F)$, vu comme sous-groupe de $\tilde{G}'(s)(F)$.

  Pour donn\'ees auxiliaires de ${\bf G}'_{\natural}(s)$, on choisit $G'_{\natural}(s)_{1}=G'_{\natural}(s)$, $\tilde{G}'_{\natural}(s)_{1}=\tilde{G}'_{\natural}(s)$, $C_{\natural}(s)_{1}=\{1\}$.  L'homomorphisme $\hat{\xi}(s)_{1}$ s'\'etend en un isomorphisme $\hat{\xi}(s)_{\natural,1}:{\cal G}'_{\natural}(s) \to {^LG}'_{\natural}(s)$. Pour simplifier, abandonnons les indices $1$ superflus. Notons ${\cal D}\subset G'(s)(F)\times G(F)$ et ${\cal D}_{\natural}\subset G'_{\natural}(s)(F)\times G_{\natural}(F)$ les ensembles de couples d'\'el\'ements r\'eguliers se correspondant.  Le groupe $C_{\natural}(F)$ se plonge diagonalement dans $G'_{\natural}(s)(F)\times G_{\natural}(F)$. On voit que l'action de ce groupe diagonal pr\'eserve ${\cal D}_{\natural}$ et que le quotient ${\cal D}_{\natural}/diag(C_{\natural}(F))$ s'identifie par projection \`a ${\cal D}$. Notons $\Delta_{\natural}(s)$ l'image r\'eciproque de $\Delta(s)$ par cette projection. On v\'erifie que c'est un facteur de transfert pour les donn\'ees auxiliaires d\'efinies ci-dessus. Ce sont ces donn\'ees que l'on utilise pour r\'ealiser explicitement les termes  $S_{{\bf M}_{\natural}}^{{\bf G}'_{\natural}(s)}( .,.)$.
 
 Pour $c\in C_{\natural}(F)$, on a l'\'egalit\'e $(\dot{f}^c)^{\tilde{G}'_{\natural}(s)}=(\dot{f}^{\tilde{G}'_{\natural}(s)})^c$: cela r\'esulte de la d\'efinition du facteur de transfert. Posons 
 $$f^{\tilde{G}'_{\natural}(s)}=\int_{C_{\natural}(F)}(\dot{f}^{\tilde{G}'_{\natural}(s)})^c\lambda_{\natural}(c)dc.$$
 Cette fonction  appartient \`a $C_{c,\lambda_{\natural}}^{\infty}(\tilde{G}'_{\natural}(s)(F))$. 
 Alors
 $$\int_{C_{\natural}(F)}S_{\tilde{M}_{\natural}}^{\tilde{G}'_{\natural}(s)}(\dot{\boldsymbol{\delta}}(s),B,(\dot{f}^c)^{\tilde{G}'_{\natural}(s)})\lambda_{\natural}(c)dc=S_{\tilde{M}_{\natural}}^{\tilde{G}'_{\natural}(s)}(\dot{\boldsymbol{\delta}},B,f^{\tilde{G}'_{\natural}(s)}).$$
 Si $s\not=1$, on en d\'eduit l'\'egalit\'e
 $$(12) \qquad \int_{C_{\natural}(F)}S_{\tilde{M}_{\natural}}^{\tilde{G}'_{\natural}(s)}(\dot{\boldsymbol{\delta}}(s),B,(\dot{f}^c)^{\tilde{G}'_{\natural}(s)})\lambda_{\natural}(c)dc=\int_{C_{\natural}(F)}S_{\tilde{M}_{\natural}}^{\tilde{G}'_{\natural}(s)}(\dot{\boldsymbol{\delta}}'(s),B,(\dot{f}^c)^{\tilde{G}'_{\natural}(s)})\lambda_{\natural}(c)dc.$$
 En effet, c'est l'\'egalit\'e (6)   o\`u l'on remplace $\tilde{G}_{\natural}$ par $\tilde{G}'_{\natural}(s)$ et $f$ par $f^{\tilde{G}'_{\natural}(s)}$ (puisque $s\not=1$, on peut appliquer (6) par hypoth\`ese de r\'ecurrence).
 
La relation (6) r\'esulte de (9), (12) et des d\'efinitions.  
 
  Remarquons que l'on peut d\'efinir $I_{\tilde{M}_{\natural},\lambda_{\natural}}^{\tilde{G}_{\natural}}(\boldsymbol{\gamma},B,f)$  pour $\boldsymbol{\gamma}\in D_{g\acute{e}om,\lambda_{\natural}}(\tilde{M}_{\natural}(F))$ et $f\in C_{c,\lambda_{\natural}}^{\infty}(\tilde{G}_{\natural}(F))$, de m\^eme que l'on a d\'efini $S_{\tilde{M}_{\natural},\lambda_{\natural}}^{\tilde{G}_{\natural}}(\boldsymbol{\gamma},B,f)$. La relation (9) affirme que cette d\'efinition est loisible. On a aussi d\'efini ci-dessus un transfert entre $C_{c,\lambda_{\natural}}^{\infty}(\tilde{G}_{\natural}(F))$ et $C_{c,\lambda_{\natural}}^{\infty}(\tilde{G}'_{\natural}(s;F))$. Avec ces d\'efinitions, l'\'egalit\'e (8) se g\'en\'eralise \`a $\boldsymbol{\delta}\in D_{g\acute{e}om,\lambda_{\natural}}(\tilde{M}_{\natural}(F))$ et $f\in  C_{c,\lambda_{\natural}}^{\infty}(\tilde{G}_{\natural}(F))$. Consid\'erons le cas particulier o\`u $\lambda_{\natural}$ est le caract\`ere trivial ${\bf 1}$. Dans ce cas, on a des isomorphismes $C_{c,{\bf 1}}^{\infty}(\tilde{G}_{\natural}(F))\simeq C_{c}^{\infty}(\tilde{G}(F))$, $D_{g\acute{e}om,{\bf 1}}(\tilde{M}_{\natural}(F))\simeq D_{g\acute{e}om}(\tilde{M}(F))$ et $D_{g\acute{e}om,{\bf 1}}^{st}(\tilde{M}_{\natural}(F))\simeq D_{g\acute{e}om}^{st}(\tilde{M}(F))$. Notons ici $f_{\natural}\mapsto f$ et $\boldsymbol{\gamma}_{\natural}\mapsto \boldsymbol{\gamma}$ ces isomorphismes. Alors:
  
  (13)  soient $\boldsymbol{\delta}_{\natural}\in D_{g\acute{e}om,{\bf 1}}^{st}(\tilde{M}_{\natural}(F))$ et $f_{\natural}\in C_{c,{\bf 1}}^{\infty}(\tilde{G}_{\natural}(F))$; on a l'\'egalit\'e
  $$S_{\tilde{M}_{\natural},{\bf 1}}^{\tilde{G}_{\natural}}(\boldsymbol{\delta}_{\natural},B,f_{\natural})=S_{\tilde{M}}^{\tilde{G}}(\boldsymbol{\delta},B,f).$$
  
  Cette assertion se d\'ecompose en deux:
  
  (14) pour $s\in Z(\hat{M})^{\Gamma_{F}}/Z(\hat{G})^{\Gamma_{F}}$, $s\not=1$,  
  $$S_{\tilde{M}_{\natural},{\bf 1}}^{\tilde{G}'_{\natural}(s)}(\boldsymbol{\delta}_{\natural}(s),B,(f_{\natural})^{\tilde{G}'_{\natural}(s)})=S_{\tilde{M}}^{\tilde{G}'(s)}(\boldsymbol{\delta}(s),B,f^{\tilde{G}'(s)});$$
  
  (15) pour $\boldsymbol{\gamma}_{\natural}\in D_{g\acute{e}om,{\bf 1}}(\tilde{M}_{\natural}(F))$, 
   $$I_{\tilde{M}_{\natural},{\bf 1}}^{\tilde{G}_{\natural}}(\boldsymbol{\gamma}_{\natural},B,f_{\natural})=I_{\tilde{M}}^{\tilde{G}}(\boldsymbol{\gamma},B,f).$$
   
   On v\'erifie sur la d\'efinition ci-dessus du transfert que $(f_{\natural})^{\tilde{G}'_{\natural}(s)}=(f^{\tilde{G}'(s)})_{\natural}$. Donc (14) est la m\^eme assertion que (13) avec $\tilde{G}$ remplac\'e par $\tilde{G}'(s)$. On peut l'admettre par r\'ecurrence. Pour (15), on peut supposer que $\boldsymbol{\gamma}$ est une int\'egrale orbitale associ\'ee \`a un \'el\'ement $\gamma\in \tilde{M}(F)$. La relation 1.9(5) nous ram\`ene au cas o\`u le syst\`eme de fonctions $B$ a toutes ses valeurs \'egales \`a  $1$. Fixons $\dot{\gamma}_{\natural}\in \tilde{M}_{\natural}(F)$ se projetant sur $\gamma$. On v\'erifie que $\boldsymbol{\gamma}_{\natural}$ est l'image par l'homomorphisme $D_{g\acute{e}om}(\tilde{M}_{\natural}(F))\to D_{g\acute{e}om,{\bf 1}}(\tilde{M}_{\natural}(F))$ de l'int\'egrale orbitale associ\'ee \`a $\dot{\gamma}$. Donc
    $$I_{\tilde{M}_{\natural},{\bf 1}}^{\tilde{G}_{\natural}}(\boldsymbol{\gamma}_{\natural},f_{\natural})= I_{\tilde{M}_{\natural}}^{\tilde{G}_{\natural}}(\dot{\gamma}_{\natural},f_{\natural}).$$
    Chosissons $\dot{f}_{\natural}\in C_{c}^{\infty}(\tilde{G}_{\natural}(F))$ tel que
    $$f_{\natural}(\gamma')=\int_{C_{\natural}(F)}\dot{f}_{\natural}^c(\gamma')dc$$
    pour tout $\gamma'\in \tilde{G}_{\natural}(F)$. Alors
  $$ I_{\tilde{M}_{\natural}}^{\tilde{G}_{\natural}}(\dot{\gamma}_{\natural},f_{\natural}) =\int_{C_{\natural}(F)}I_{\tilde{M}_{\natural}}^{\tilde{G}_{\natural}}(\dot{\gamma}_{\natural},\dot{f}^c_{\natural})dc.$$
    On fixe un sous-groupe compact sp\'ecial $K$ de $G(F)$ en bonne position relativement \`a $M$.  Il lui correspond un tel sous-groupe $K_{\natural}$ de $G_{\natural}(F)$ (par la bijection entre facettes sp\'eciales des immeubles de $G$ et $G_{\natural}$). On utilise ces sous-groupes pour d\'efinir les int\'egrales pond\'er\'ees suivantes. On montre d'abord que
    $$(16) \qquad \int_{C_{\natural}(F)}J_{\tilde{M}_{\natural}}^{\tilde{G}_{\natural}}(\dot{\gamma}_{\natural},\dot{f}^c_{\natural})dc=J_{\tilde{M}}^{\tilde{G}}(\gamma,f).$$
    Si $M_{\gamma}=G_{\gamma}$, il suffit d'appliquer les d\'efinitions: pour $x_{\natural}\in G_{\natural}(F)$ se projetant sur $x\in G(F)$, on a
    $$\int_{C_{\natural}(F)}\dot{f}_{\natural}^c(x_{\natural}^{-1}\dot{\gamma}_{\natural}x_{\natural})dc=f(x^{-1}\gamma x)$$
    et $v_{\tilde{M}_{\natural}}^{\tilde{G}_{\natural}}(x_{\natural})=v_{\tilde{M}}^{\tilde{G}}(x)$. Pour $\gamma$ quelconque, on v\'erifie que pour $\tilde{P}\in {\cal P}(\tilde{M})$, $a_{\natural}\in  A_{\tilde{M}_{\natural}}(F)$ se projetant en $a\in A_{\tilde{M}}(F)$ et pour $\lambda\in {\cal A}_{\tilde{M}_{\natural},{\mathbb C}}^{\tilde{G}_{\natural},*}\simeq {\cal A}_{\tilde{M},{\mathbb C}}^{\tilde{G},*}$, on a $r_{\tilde{P}_{\natural}}^{\tilde{G}_{\natural}}(\dot{\gamma}_{\natural},a;\lambda)=r_{\tilde{P}}^{\tilde{G}}(\gamma,a;\lambda)$. L'\'egalit\'e (16) se d\'eduit alors pour $\gamma$ par passage \`a la limite \`a partir du cas o\`u $M_{\gamma}=G_{\gamma}$. Il faut ensuite montrer que pour tout $\tilde{L}\in {\cal L}(\tilde{M})$ avec $\tilde{L}\not=\tilde{G}$, on a
 $$(17) \qquad \int_{C_{\natural}(F)}I_{\tilde{M}_{\natural}}^{\tilde{L}_{\natural}}(\dot{\gamma}_{\natural},\phi_{\tilde{L}_{\natural}}(\dot{f}^c_{\natural}))dc=I_{\tilde{M}}^{\tilde{L}}(\gamma,\phi_{\tilde{L}}(f)).$$
 On a besoin pour cela de propri\'et\'es des applications $\phi_{\tilde{L}}$ et $\phi_{\tilde{L}_{\natural}}$, qui sont essentiellement formelles. A savoir que $\phi_{\tilde{L}_{\natural}}(\dot{f}^c_{\natural})=(\phi_{\tilde{L}_{\natural}}(\dot{f}_{\natural})) ^c$ comme on l'a d\'ej\`a dit et que $\phi_{\tilde{L}_{\natural}}(\dot{f}_{\natural})$ et $\phi_{\tilde{L}}(f)$ sont reli\'ees de la m\^eme fa\c{c}on que $\dot{f}_{\natural}$ et $f$ (\`a ceci pr\`es qu'elles ne sont pas \`a support compact mais appartiennent \`a des espaces $C^{\infty}_{ac}$; le passage \`a ces espaces ne pose pas de probl\`eme). Alors l'\'egalit\'e (17) n'est autre que (15) o\`u l'on change $\tilde{G}$ en $\tilde{L}$ et $f$ en $\phi_{\tilde{L}}(f)$. On peut l'admettre par r\'ecurence. L'assertion (15) r\'esulte de (16), (17) et des d\'efinitions. Cela ach\`eve la preuve de (13).
 
 Il nous reste \`a prouver la relation (7). Consid\'erons les extensions
 $$1\to C_{\natural}\times C_{\flat}\to G_{_{\natural,\flat}}\to G\to 1 \text{ et }\tilde{G}_{\natural,\flat}\to \tilde{G}$$
 ainsi que le caract\`ere $\lambda_{\natural}\times {\bf 1}_{\flat}$ de $C_{\natural}(F)\times C_{\flat}(F)$, o\`u ${\bf 1}_{\flat}$ est le caract\`ere trivial de $C_{\flat}(F)$. On a de nouveau des isomorphismes $C_{c,\lambda_{\natural}\times{\bf 1}_{\flat}}^{\infty}(\tilde{G}_{\natural,\flat}(F))\simeq C_{c,\lambda_{\natural}}^{\infty}(\tilde{G}_{\natural})$, $D_{g\acute{e}om,\lambda_{\natural}\times{\bf 1}_{\flat}}^{st}(\tilde{M}_{\natural,\flat}(F))\simeq D_{g\acute{e}om,\lambda_{\natural}}(\tilde{M}_{\natural})$. Notons $f_{\natural,\flat}\in C_{c,\lambda_{\natural}\times{\bf 1}_{\flat}}^{\infty}(\tilde{G}_{\natural,\flat}(F))$ et $\boldsymbol{\delta}_{\natural,\flat}\in D_{g\acute{e}om,\lambda_{\natural}\times{\bf 1}_{\flat}}^{st}(\tilde{M}_{\natural,\flat}(F))$ les \'el\'ements auxquels s'identifient $f_{\natural}$ et $\boldsymbol{\delta}_{\natural}$. La relation (13) se g\'en\'eralise en
 
 (18) $S_{\tilde{M}_{\natural,\flat},\lambda_{\natural}\times {\bf 1}_{\flat}}^{\tilde{G}_{\natural,\flat}}(\boldsymbol{\delta}_{\natural,\flat},B,f_{\natural,\flat})=S_{\tilde{M}_{\natural},\lambda_{\natural}}^{\tilde{G}_{\natural}}(\boldsymbol{\delta}_{\natural},B,f_{\natural})$.

Pour le prouver, on choisit $\dot{\boldsymbol{\delta}}_{\natural}\in D_{g\acute{e}om}(\tilde{M}_{\natural}(F))$ se projetant sur $\boldsymbol{\gamma}_{\natural}$ par l'homomorphisme (3) et $\dot{f}_{\natural}\in C_{c}^{\infty}(\tilde{G}_{\natural}(F))$ reli\'ee \`a $f_{\natural}$ par (4). Alors
 $$S_{\tilde{M}_{\natural},\lambda_{\natural}}^{\tilde{G}_{\natural}}(\boldsymbol{\delta}_{\natural},B,f_{\natural})=\int_{C_{\natural}(F)}S_{\tilde{M}_{\natural}}^{\tilde{G}_{\natural}}(\dot{\boldsymbol{\delta}}_{\natural},B,\dot{f}_{\natural})\lambda_{\natural}(c)dc.$$ 
 On identifie $\dot{f}_{\natural}$ \`a un \'el\'ement de $C_{c,{\bf 1}_{\flat}}^{\infty}(\tilde{G}_{\natural,\flat}(F))$ et on choisit $\dot{f}_{\natural,\flat}\in C_{c}^{\infty}(\tilde{G}_{\natural,\flat}(F))$ tel que cet \'el\'ement soit $\int_{C_{\flat}(F)}\dot{f}_{\natural,\flat}^{c'}\,dc'$. On identifie $\dot{\boldsymbol{\delta}}_{\natural}$ \`a un \'el\'ement de $D_{g\acute{e}om,{\bf 1}_{\flat}}^{st}(\tilde{M}_{\natural,\flat}(F))$ et on choisit $\dot{\delta}_{\natural,\flat}\in D_{g\acute{e}om}^{st}(\tilde{M}_{\natural,\flat}(F))$ se projetant sur cet \'el\'ement. Pour $c\in C_{\natural}(F)$, (13) implique que
 $$  S_{\tilde{M}_{\natural}}^{\tilde{G}_{\natural}}(\dot{\boldsymbol{\delta}}_{\natural},B,\dot{f}_{\natural})=\int_{C_{\flat}(F)}S_{\tilde{M}_{\natural,\flat}}^{\tilde{G}_{\natural,\flat}}(\dot{\boldsymbol{\delta}}_{\natural,\flat},B,\dot{f}_{\natural,\flat}^{c'})dc'.$$
 Donc
  $$S_{\tilde{M}_{\natural},\lambda_{\natural}}^{\tilde{G}_{\natural}}(\boldsymbol{\delta}_{\natural},B,f_{\natural})=\int_{C_{\natural}(F)\times C_{\flat}(F)} S_{\tilde{M}_{\natural,\flat}}^{\tilde{G}_{\natural,\flat}}(\dot{\boldsymbol{\delta}}_{\natural,\flat},B,\dot{f}_{\natural,\flat}^{cc'})\lambda_{\natural}(c)\,dc'\,dc.$$
 Mais $\dot{\boldsymbol{\delta}}_{\natural,\flat}$ se projette sur $\boldsymbol{\delta}_{\natural,\flat}$ par (3) et $\dot{f}_{\natural,\flat}$ est reli\'e \`a $f_{\natural,\flat}$ par (4). Le membre de droite ci-dessus est donc \'egal \`a $S_{\tilde{M}_{\natural},\lambda_{\natural}\times {\bf 1}_{\flat}}^{\tilde{G}_{\natural}}(\boldsymbol{\delta}_{\natural},B,f_{\natural})$, ce qui prouve (18).
 
 On effectue la m\^eme construction \`a partir de $f_{\flat}$ et $\boldsymbol{\delta}_{\flat}$. On obtient des \'el\'ements disons $f_{\flat,\natural}\in C_{c,{\bf 1}_{\natural}\times \lambda_{\flat}}^{\infty}(\tilde{G}_{\natural,\flat}(F))$ et $\boldsymbol{\delta}_{\flat,\natural}\in D_{g\acute{e}om,{\bf 1}_{\natural}\times \lambda_{\flat}}^{st}(\tilde{M}_{\natural,\flat}(F))$ et l'\'egalit\'e
 $$(19)\qquad S_{\tilde{M}_{\natural,\flat},{\bf 1}_{\natural}\times \lambda_{\flat}}^{\tilde{G}_{
\natural,\flat}}(\boldsymbol{\delta}_{\flat,\natural},B,f_{\flat,\natural})=S_{\tilde{M}_{\flat},\lambda_{\flat}}^{\tilde{G}_{\flat}}(\boldsymbol{\delta}_{\flat},B,f_{\flat}).$$
La multiplication par $\tilde{\lambda}_{\natural,\flat}$ d\'efinit un isomorphisme de $C_{c,\lambda_{\natural}\times{\bf 1}_{\flat}}^{\infty}(\tilde{G}_{\natural,\flat}(F))$ sur $C_{c,{\bf 1}_{\natural}\times \lambda_{\flat}}^{\infty}(\tilde{G}_{\natural,\flat}(F))$. Par restriction puis dualit\'e, on obtient un isomorphisme de $D_{g\acute{e}om,\lambda_{\natural}\times{\bf 1}_{\flat}}^{st}(\tilde{M}_{\natural,\flat}(F))$ sur $D_{g\acute{e}om,{\bf 1}_{\natural}\times \lambda_{\flat}}^{st}(\tilde{M}_{\natural,\flat}(F))$. Les \'el\'ements $f_{\natural,\flat}$ et $f_{\flat,\natural}$, resp. $\boldsymbol{\delta}_{\natural,\flat}$ et $\boldsymbol{\delta}_{\flat,\natural}$, se correspondent par ces isomorphismes. En \'ecrivant les d\'efinitions des membres de gauche de (18) et (19), on est ramen\'e au probl\`eme suivant, o\`u les extensions ont disparu. On consid\`ere un caract\`ere $\lambda$ de $G(F)$ et une fonction non nulle $\tilde{\lambda}$ sur $\tilde{G}(F)$ se transformant selon $\lambda$. On consid\`ere $\boldsymbol{\delta}\in D_{g\acute{e}om}(\tilde{M}(F))$ et $f\in C_{c}^{\infty}(\tilde{G}(F))$. On veut prouver

(20) $S_{\tilde{M}}^{\tilde{G}}(\tilde{\lambda}\boldsymbol{\delta},B,\tilde{\lambda}f)=S_{\tilde{M}}^{\tilde{G}}(\boldsymbol{\delta},B,f)$,

\noindent o\`u $\boldsymbol{\delta}\mapsto \tilde{\lambda}\boldsymbol{\delta}$ est l'analogue de l'isomorphisme ci-dessus. Cette assertion se d\'ecompose en deux:

(21)  $I_{\tilde{M}}^{\tilde{G}}(\tilde{\lambda}\boldsymbol{\delta},B,\tilde{\lambda}f)=I_{\tilde{M}}^{\tilde{G}}(\boldsymbol{\delta},B,f)$;

(22) pour $s\in Z(\hat{M})^{\Gamma_{F}}/Z(\hat{G})^{\Gamma_{F}}$, $s\not=1$, $S_{\tilde{M}}^{\tilde{G}'(s)}((\tilde{\lambda}\boldsymbol{\delta})(s),B,(\tilde{\lambda}f)^{\tilde{G}'(s)})=S_{\tilde{M}}^{\tilde{G}'(s)}(\boldsymbol{\delta}(s),B,f^{\tilde{G}'(s)})$.

Pour (21), on se ram\`ene encore une fois au cas o\`u le syst\`eme de fonctions $B$ a toutes ses valeurs \'egales \`a $1$. Dans ce cas, l'assertion similaire pour les int\'egrales orbitales pond\'er\'ees non invariantes est imm\'ediate. On doit encore utiliser une propri\'et\'e formelle des applications $\phi_{\tilde{L}}$, \`a savoir que $\phi_{\tilde{L}}(\tilde{\lambda}f)=\tilde{\lambda}\phi_{\tilde{L}}(f)$. Pour (22), on utilise l'homomorphisme $G'(s)_{ab}(F)\to G_{ab}(F)$ et l'application  compatible $\tilde{G}'(s)_{ab}(F)\to \tilde{G}_{ab}(F)$.  Gr\^ace \`a ces applications, $\lambda$ se restreint en un caract\`ere de $G'(s)(F)$ et $\tilde{\lambda}$ se restreint en une fonction sur $\tilde{G}'(s)(F)$. On v\'erifie que $(\tilde{\lambda}f)^{\tilde{G}'(s)}=\tilde{\lambda}(f^{\tilde{G}'(s)})$ et, dualement, que $(\tilde{\lambda}\boldsymbol{\delta})(s)=\tilde{\lambda}\boldsymbol{\delta}(s)$.  Alors (22) n'est autre que (20) o\`u l'on a remplac\'e $\tilde{G}$ par $\tilde{G}'(s)$ et $f$ par $f^{\tilde{G}'(s)}$. Par r\'ecurrence, on peut admettre (22). Cela ach\`eve la preuve de (20) et la v\'erification des propri\'et\'es formelles. $\square$

Dans le cas o\`u le syst\`eme de fonctions $B$ a toutes ses valeurs \'egales \`a $1$, on note simplement $S_{\tilde{M}}^{\tilde{G}}(\boldsymbol{\delta},{\bf f})=S_{\tilde{M}}^{\tilde{G}}(\boldsymbol{\delta},B,{\bf f})$.

Pour $B$ quelconque, on a

(23) supposons que le support de $\boldsymbol{\delta}$ soit form\'e d'\'el\'ements de $ \tilde{M}(F)$ qui sont $\tilde{G}$-\'equisinguliers; alors on a l'\'egalit\'e  $S_{\tilde{M}}^{\tilde{G}}(\boldsymbol{\delta},B,{\bf f})=S_{\tilde{M}}^{\tilde{G}}(\boldsymbol{\delta},{\bf f})$.

En effet, on a encore $M_{\gamma}=G'(s)_{\gamma}$ pour tout \'el\'ement $\gamma$ du support
et pour tout $s\in Z(\hat{M})^{\Gamma_{F}}$. En raisonnant par r\'ecurrence, l'assertion r\'esulte de la relation  $I_{\tilde{M}}^{\tilde{G}}(\boldsymbol{\delta},B,{\bf f})=I_{\tilde{M}}^{\tilde{G}}(\boldsymbol{\delta},{\bf f})$, laquelle r\'esulte des d\'efinitions.

\bigskip

 \subsection{D\'efinition d'un syst\`eme de fonctions $B^{\tilde{G}}$}
 
 On revient au cas o\`u $(G,\tilde{G},{\bf a})$ est quelconque. Soit ${\bf G}'=(G',{\cal G}',\tilde{s})$ une donn\'ee endoscopique de $(G,\tilde{G},{\bf a})$. On peut fixer une paire de Borel \'epingl\'ee $\hat{{\cal E}}=(\hat{B},\hat{T},(\hat{E}_{\hat{\alpha}})_{\hat{\alpha}\in \hat{\Delta}})$ de $\hat{G}$ de sorte que $\tilde{s}=s\hat{\theta}$ avec $s\in \hat{T}$. Introduisons "les" paires de Borel \'epingl\'ees ${\cal E}=(B,T,(E_{\alpha})_{\alpha\in \Delta})$ et ${\cal E}^{'}=(B^{'},T^{'},(E^{'}_{\alpha})_{\alpha\in \Delta'})$ de $G$ et $G'$. Le choix de $\hat{\cal E}$ permet d'identifier $T^{'}$ \`a $T/(1-\theta)(T)$, cf. [I] 1.5. On fixe $e\in {\cal Z}(\tilde{G})$ et on note $e'$ son image dans ${\cal Z}(\tilde{G}')$. Soit $\epsilon\in \tilde{G}'_{ss}(F)$. On peut identifier $(B^{'},T^{'})$ \`a une paire de Borel conserv\'ee par $ad_{\epsilon}$. Alors $\epsilon$ s'\'ecrit $\mu e'$, avec $\mu \in T^{'}$ et on rel\`eve $\mu$ en un \'el\'ement $\nu\in T$. On utilise les notations de [I] 1.6. En particulier, on note $\Sigma(T)$ l'ensemble des racines de $T$ dans $G$ et $\Sigma^{G'_{\epsilon}}(T^{'})$ celui des racines de $T'$ dans $G'_{\epsilon}$.  D'apr\`es [W2] 3.3, l'ensemble $\Sigma^{G'_{\epsilon}}(T^{'})$ est alors la r\'eunion des ensembles suivants
 
 (a) les $N\alpha$, pour $\alpha\in \Sigma(T)$ de type 1 tels que $N\alpha(\nu)=1$ et $N\hat{\alpha}(s)=1$;
 
 (b) les $2N\alpha$ pour  $\alpha\in \Sigma(T)$ de type 2 tels que $N\alpha(\nu)=1$ et $N\hat{\alpha}(s)=1$

(c) les $2N\alpha$ pour $\alpha\in \Sigma(T)$ de type 2 tels que $N\alpha(\nu)=-1$ et $N\hat{\alpha}(s)=1$;

 (d) les $N\alpha$ pour $\alpha\in \Sigma(T)$ de type 3 tels que $N\alpha(\nu)=1$ et $N\hat{\alpha}(s)=-1$. 

On d\'efinit une fonction $B^{\tilde{G}}_{\epsilon}$ sur cet ensemble de la fa\c{c}on suivante. Dans le cas (a), $B^{\tilde{G}}_{\epsilon}(N\alpha)=n_{\alpha}$ (rappelons que $n_{\alpha}$ est le plus petit entier $n\geq1$ tel que $\theta^n(\alpha)=\alpha$). Dans le cas (b), $B^{\tilde{G}}_{\epsilon}(2N\alpha)=2n_{\alpha}$. Dans le cas (c), $B^{\tilde{G}}_{\epsilon}(2N\alpha)=n_{\alpha}$. Dans le cas (d), $B^{\tilde{G}}_{\epsilon}(N\alpha)=2n_{\alpha}$.

\ass{Lemme}{La fonction $B^{\tilde{G}}_{\epsilon}$ ne d\'epend pas des choix faits dans sa construction. Elle v\'erifie les conditions de 1.8. }

Preuve. La fonction $B^{\tilde{G}}_{\epsilon}$ v\'erifie $B^{\tilde{G}}_{\epsilon}(-\beta)=B^{\tilde{G}}_{\epsilon}(\beta)$ pour tout $\beta\in \Sigma^{G'_{\epsilon}}(T^{'})$. Introduisons les groupes de Weyl $W^{G'_{\epsilon}}$ de $G'_{\epsilon}$ et $W^{G'}$ de $G'$ tous deux relatifs \`a $T^{'}$ et le groupe de Weyl $W$ de $G$ relatif \`a $T$. On a

(1) $W^{G'_{\epsilon}}\subset \{w\in W^{G'}; w(\mu)=\mu\}$.

En effet, puisque $e'$ commute \`a tout \'el\'ement de $G'$, un \'el\'ement $x\in G'$ v\'erifie $ad_{x}(\epsilon)=\epsilon$ si et seulement s'il v\'erifie $ad_{x}(\mu)=\mu$.

(2) $W^{G'}$ s'identifie \`a un sous-groupe de l'ensemble des $w\in W^{\theta}$ qui fixent l'image de $s$ dans $\hat{T}/(1-\hat{\theta})(\hat{T})$.

Le groupe $W^{G'}$ s'identifie au groupe de Weyl de $\hat{G}'$ relatif \`a $\hat{T}^{\hat{\theta},0}$. Un \'el\'ement $x\in \hat{G}'$ qui normalise ce tore  normalise aussi $\hat{T}$. Il v\'erifie en outre $ s\hat{\theta}(x)s^{-1}=x$. Cela entra\^{\i}ne que son image $w$ dans $W$ est fixe par $\hat{\theta}$. On peut alors l'\'ecrire $x=tn$, o\`u $t\in \hat{T}$ et $n$ est fixe par $\hat{\theta}$. Alors on a l'\'egalit\'e $nsn^{-1}=t^{-1}s\hat{\theta}(t)$, d'o\`u $w(s)\in s(1-\hat{\theta})(\hat{T})$. $\square$

Il  r\'esulte de (1) et (2) que l'action sur $T'$ d'un \'el\'ement de $W^{G'_{\epsilon}}$ co\"{\i}ncide avec celle d'un \'el\'ement de $W^{\theta}$ qui conserve l'image de $\nu$ dans $T/(1-\theta)(T)$ et celle de $s$ dans $\hat{T}/(1-\hat{\theta})(\hat{T})$.
 L'action d'un tel \'el\'ement ne change ni le type d'une racine $\alpha\in \Sigma(T)$, ni le nombre $n_{\alpha}$, ni les valeurs $N\alpha(\nu)$ et $N\hat{\alpha}(s)$. D'o\`u l'\'egalit\'e $B^{\tilde{G}}_{\epsilon}(w\beta)=B^{\tilde{G}}_{\epsilon}(\beta)$ pour tout $\beta\in \Sigma^{G'_{\epsilon}}(T^{'})$ et tout $w\in W^{G'_{\epsilon}}$. 
 
 On compl\`ete la paire de Borel $(B'\cap G'_{\epsilon},T')$ de $G'_{\epsilon}$ en une paire de Borel \'epingl\'ee ${\cal E}'_{\epsilon}$. On introduit les actions galoisiennes quasi-d\'eploy\'ees $\sigma\mapsto \sigma_{G^{'*}_{\epsilon}}$, resp. $\sigma\mapsto \sigma_{G^{'*}}$, $\sigma\mapsto \sigma_{G^*}$, de $\Gamma_{F}$ sur $G'_{\epsilon}$, resp. $G'$, $G$, relatives aux paires de Borel \'epingl\'ees ${\cal E}'_{\epsilon}$, resp. ${\cal E}'$, ${\cal E}$. En notant $\sigma\mapsto \sigma_{G'}$ et $\sigma\mapsto \sigma_{G}$ les actions naturelles sur $G'$ et $G$, ou $\tilde{G}' $ et $\tilde{G}$, on a des \'egalit\'es $\sigma_{G^{'*}_{\epsilon}}=ad_{u'_{\epsilon}(\sigma)}\circ \sigma_{G'}$, $\sigma_{G^{'*}}=ad_{u'(\sigma)}\circ \sigma_{G'}$ et $\sigma_{G^*}=ad_{u(\sigma)}\circ \sigma_{G}$, o\`u $u'_{\epsilon}(\sigma)\in G'_{\epsilon}$, $u'(\sigma)\in G'$ et $u(\sigma)\in G$. On introduit le cocycle $z:\Gamma_{F}\to Z(G)$ tel que $ad_{u(\sigma)}\circ \sigma_{G}(e)=z(\sigma)^{-1}e$ pour tout $\sigma\in \Gamma_{F}$. On a
 
 (3) $\sigma_{G^{'*}_{\epsilon}}(\mu)=z(\sigma)\mu$ pour tout $\sigma\in \Gamma_{F}$.
 
  Comme on l'a dit dans la preuve de (1), $u'_{\epsilon}(\sigma)$ commute \`a $\mu$ puisqu'il commute \`a $\epsilon$. Donc $\sigma_{G^{'*}_{\epsilon}}(\mu)=\sigma_{G'}(\mu)$. Par d\'efinition de l'action galoisienne sur ${\cal Z}(\tilde{G}')$, on a $\sigma_{G'}(e')=z(\sigma)^{-1}e'$. Puisque $\epsilon\in G'(F)$, on a $\sigma_{G'}(\epsilon)=\epsilon$, c'est-\`a-dire $\sigma_{G'}(\mu e')=\mu e'$. D'o\`u $\sigma_{G'}(\mu)=z(\sigma)\mu$. $\square$
 
On a

(4) pour tout $\sigma\in \Gamma_{F}$, il existe $w'(\sigma)\in W^{G'}$ tel que l'on ait l'\'egalit\'e $\sigma_{G^{'*}_{\epsilon}}=w'(\sigma)\circ \sigma_{G^{'*}}$ sur $T'$.

En effet, $\sigma_{G^{'*}_{\epsilon}}=ad_{u'_{\epsilon}(\sigma)u'(\sigma)^{-1}}\circ \sigma_{G^{'*}}$. Puisque ces deux actions conservent $T'$, l'\'el\'ement $u'_{\epsilon}(\sigma)u'(\sigma)^{-1}$ normalise ce tore et d\'efinit l'\'el\'ement $w'(\sigma)$ cherch\'e.

Pour simplifier, on  conserve la notation $\sigma\mapsto \sigma_{G^*}$ pour l'action galoisienne sur $\hat{G}$. On a

(5) pour tout $\sigma\in \Gamma_{F}$, il existe $w(\sigma)\in W^{\theta}$ tel que  l'on ait l'\'egalit\'e $\sigma_{G^{'*}}=w(\sigma)\circ \sigma_{G^*}$ et que, sur $\hat{T}$, $w(\sigma)\circ \sigma_{G^*} $ conserve l'image de $s$ dans $\hat{T}/Z(\hat{G})(1-\hat{\theta})(\hat{T})$. 

Pour $\sigma\in \Gamma_{F}$, on rel\`eve $\sigma$ en $u\in W_{F}$ et on choisit $g_{u}=(g(u),u)\in {\cal G}'$ tel que $ad_{g_{u}}$ co\"{\i}ncide avec $u_{G^{'*}}$ sur $\hat{G}'$. L'\'el\'ement $g(u) $ normalise $\hat{T}^{\hat{\theta},0}$ donc aussi $\hat{T}$.  On a $s\hat{\theta}(g(u))\sigma_{G^*}(s)^{-1}=a(u)g(u)$.  Cela entra\^{\i}ne que l'image $w(\sigma)$ de $g(u)$ dans  $W$ est fixe par $\hat{\theta}$. On peut \'ecrire $g(u)=tn$ avec $t\in \hat{T}$ et $n$ fixe par $\hat{\theta}$.  Alors $n\sigma_{G^*}(s)n^{-1}=a(u)^{-1}t^{-1}s\hat{\theta}(t)$, d'o\`u $w(\sigma)\circ \sigma_{G^*}(s)\in sZ(\hat{G})(1-\hat{\theta})(\hat{T})$. L'assertion (5) en r\'esulte. $\square$

Il r\'esulte de  (4) et (5) que, pour $\sigma\in \Gamma_{F}$, $\sigma_{G^{'*}_{\epsilon}}$ co\"{\i}ncide sur $T'$ avec  $w'(\sigma)w(\sigma)\sigma_{G^*}$. Il r\'esulte de (2), (3) et (5) que cette action conserve l'image de $\nu$ dans $T/Z(G)(1-\theta)(T)$ et l'image de $s$ dans $\hat{T}/Z(\hat{G})(1-\hat{\theta})(\hat{T})$. Une telle action ne change ni le type d'une racine $\alpha\in \Sigma(T)$, ni le nombre $n_{\alpha}$, ni les valeurs $N\alpha(\nu)$ et $N\hat{\alpha}(s)$. D'o\`u l'\'egalit\'e $B^{\tilde{G}}_{\epsilon}(\sigma_{G^{'*}_{\epsilon}}(\beta))=B^{\tilde{G}}_{\epsilon}(\beta)$ pour tout $\beta\in \Sigma^{G'_{\epsilon}}(T^{'})$ et tout $\sigma\in \Gamma_{F}$.  

On a ainsi v\'erifi\'e la premi\`ere condition de 1.8. V\'erifions la seconde. Posons $\Sigma_{1}=\Sigma^{G'_{\epsilon}}(T')$ et notons $\check{\Sigma}_{1}$ l'ensemble associ\'e de coracines. On consid\`ere $\Sigma_{1}$, resp. $\check{\Sigma}_{1}$, comme un sous-ensemble de 
$X^{*}(T'_{\epsilon,SC})\otimes {\mathbb R}$, resp.  $X_{*}(T'_{\epsilon,SC})\otimes {\mathbb R}$, o\`u $T'_{\epsilon,sc}$ est l'image r\'eciproque de $T'$ dans $G'_{\epsilon,SC}$. Posons $b=B^{\tilde{G}}_{\epsilon}$, 
$$\Sigma_{2}=\{\alpha/b(\alpha); \alpha\in \Sigma_{1}\},$$
$$\check{\Sigma}_{2}=\{b(\alpha)\check{\alpha}; \alpha\in \Sigma_{1}\},$$
o\`u $\check{\alpha}$ est la coracine associ\'ee \`a $\alpha$. On a

(6) $\Sigma_{2}$ est un syst\`eme de racines dont $\check{\Sigma}_{2}$ est l'ensemble associ\'e de coracines.

En effet, posons $\eta=\nu e \in \tilde{G}$. Si on oublie les actions galoisiennes qui ne comptent pas pour ce que l'on veut prouver, on a construit en [W2] 3.5 un groupe $\bar{H}$ qui est un groupe endoscopique de $G_{\eta,SC}$ et tel que $G'_{\epsilon,SC}$ et $\bar{H}_{SC}$ sont en situation d'endoscopie non standard (cf. [W2] 1.7). Alors $\Sigma_{2}$ est l'ensemble de racines de ce groupe $\bar{H}$ et $\check{\Sigma}_{2}$ est l'ensemble de coracines associ\'e ([W2] 3.3(2)).

{\bf Remarque.} La fonction $B^{\tilde{G}}_{\epsilon}$ a  \'evidemment \'et\'e d\'efinie pour que (6) soit v\'erifi\'ee.

\bigskip

Il est facile de classifier les triplets $(\Sigma_{1},\Sigma_{2},b)$ v\'erifiant la condition (6), cf. [W2] 1.7. Il sont produits de triplets analogues tels que $\Sigma_{1}$ et $\Sigma_{2}$ sont irr\'eductibles. Dans le cas irr\'eductible, \`a homoth\'etie pr\`es (c'est-\`a-dire quitte \`a multiplier $b$ par un rationnel strictement positif), les possibilit\'es sont les suivantes:

- $\Sigma_{1}$ et $\Sigma_{2}$ sont de m\^eme type et $b$ est constante;

- $\Sigma_{1}$ est de type $B_{n}$, $C_{n}$, $F_{4}$ ou $G_{2}$, $\Sigma_{2}$ est respectivement de type $C_{n}$, $B_{n}$, $F_{4}$ ou $G_{2}$ et $b$ est le carr\'e de la fonction longueur.

Cela v\'erifie exactement la seconde condition de 1.8.

On doit montrer que la d\'efinition ne d\'epend pas des choix. On peut voir que changer de choix revient \`a remplacer composer l'identification de $T'$ \`a $T/(1-\theta)(T)$ par l'action d'un \'el\'ement $w$ de $W^{\theta}$, remplacer $\nu$ par un \'el\'ement de $w(\nu)Z(G)(1-\theta)(T)$ et $s$ par un \'el\'ement de $w(s)Z(\hat{G})(1-\hat{\theta})(\hat{T})$. On laisse la v\'erification fastidieuse de ce fait au lecteur. Il est clair qu'une telle modification laisse $B^{\tilde{G}}_{\epsilon}$ inchang\'ee. $\square$ 

Pour tout $\epsilon\in \tilde{G}'_{ss}(F)$, on vient de d\'efinir une fonction $B^{\tilde{G}}_{\epsilon}$ sur l'ensemble de racines de $G'_{\epsilon}$. Il r\'esulte de la d\'efinition que ces fonctions v\'erifient la condition (2) de 1.9. Elles se regroupent donc en un syst\`eme de fonctions $B^{\tilde{G}}$ sur $\tilde{G}'(F)$ au sens de ce paragraphe. 

Soit ${\bf G}'_{0}=(G'_{0},{\cal G}'_{0},\tilde{s}_{0})$ une donn\'ee endoscopique \'equivalente \`a ${\bf G}'$. Il y a alors un isomorphisme $\tilde{\alpha}:\tilde{G}'\to \tilde{G}'_{0}$ d\'efini sur $F$, unique modulo composition avec un automorphisme int\'erieur (cf. [I] 1.5). Les d\'efinitions entra\^{\i}nent que cet isomorphisme est compatible aux  syst\`emes de fonctions $B^{\tilde{G}}$ d\'efinis  sur $\tilde{G}'(F)$  et $\tilde{G}'_{0}(F)$.

\bigskip

\subsection{Int\'egrales orbitales pond\'er\'ees $\omega$-\'equivariantes et endoscopie}
Soit $(G,\tilde{G},{\bf a})$ un triplet quelconque. Soient $\tilde{M}$ un espace de Levi de $\tilde{G}$ et ${\bf M}'=(M',{\cal M}',\tilde{\zeta})$ une donn\'ee endoscopique elliptique et relevante de $\tilde{M}$. Comme en [I] 3.2, on r\'ealise $^LM$ comme espace de Levi standard de $^LG$ et on impose que le cocycle $a_{M}$ associ\'e \`a cette donn\'ee prend ses valeurs dans $Z(\hat{G})$ et que sa classe dans $H^1(W_{F};Z(\hat{G}))$ est ${\bf a}$. Soit $\tilde{s}\in \tilde{\zeta}Z(\hat{M})^{\Gamma_{F},\hat{\theta}}$. On construit la donn\'ee endoscopique ${\bf G}'(\tilde{s})=(G'(\tilde{s}),{\cal G}'(\tilde{s}),\tilde{s})$, cf. [I] 3.3. On introduit le syst\`eme de fonctions $B^{\tilde{G}}$ sur $\tilde{G}'(s)(F)$.

Pour $\boldsymbol{\delta}\in D^{st}_{g\acute{e}om}({\bf M}')\otimes Mes(M'(F))^*$ et ${\bf f}\in C_{c}^{\infty}(\tilde{G}(F))\otimes Mes(G(F))$, on pose
$$ I_{\tilde{M}}^{\tilde{G},{\cal E}}({\bf M}',\boldsymbol{\delta},{\bf f})=\sum_{\tilde{s}\in \tilde{\zeta} Z(\hat{M})^{\Gamma_F,\hat{\theta}}/Z(\hat{G})^{\Gamma_F,\hat{\theta}}}i_{\tilde{M}'}(\tilde{G},\tilde{G}'(\tilde{s}))
S_{{\bf M}'}^{{\bf G}'(\tilde{s})}(\boldsymbol{\delta},B^{\tilde{G}},{\bf f}^{{\bf G}'(\tilde{s})}).$$
Expliquons cette formule. Le coefficient $i_{M'}(G,G'(\tilde{s}))$ g\'en\'eralise celui du paragraphe 1.10. Il est nul si ${\bf G}'(\tilde{s})$ n'est pas elliptique. Si ${\bf G}'(\tilde{s})$ est elliptique, on pose
$$i_{\tilde{M}'}(\tilde{G},\tilde{G}'(\tilde{s}))=[Z(\hat{M}')^{\Gamma_F}:(Z(\hat{M}')^{\Gamma_F}\cap Z(\hat{M}))][Z(\hat{G}'(\tilde{s}))^{\Gamma_F}:(Z(\hat{G}'(\tilde{s}))^{\Gamma_F}\cap Z(\hat{G}))]^{-1}.$$
Pour d\'efinir les termes du membre de droite, on a besoin de choisir des mesures sur des espaces analogues \`a ${\cal A}_{\tilde{M}}^{\tilde{G}}$, cf. 1.1. Pour cela, on a fix\'e une forme quadratique sur $X_{*}(T^*)\otimes {\mathbb R}$. On remarque que pour chaque groupe $G'(\tilde{s})$ intervenant ci-dessus, un tore maximal $T'(\tilde{s})$ de ce groupe s'identifie sur $\bar{F}$ \`a $T^*/(1-\theta^*)(T^*)$. Il s'en d\'eduit un isomorphisme $X_{*}(T'(\tilde{s}))\otimes {\mathbb R}\simeq (X_{*}(T^*)\otimes {\mathbb R})^{\theta^*}$ On choisit pour forme quadratique sur le premier espace  la restriction au second de la forme que l'on a fix\'ee. Si $(G,\tilde{G},{\bf a})$ n'est pas quasi-d\'eploy\'e et \`a torsion int\'erieure ou si $(G,\tilde{G},{\bf a})$ est  quasi-d\'eploy\'e et \`a torsion int\'erieure et ${\bf M}'\not={\bf M}$, tous les termes $S_{{\bf M}'}^{{\bf G}'(\tilde{s})}(\boldsymbol{\delta},B^{\tilde{G}},{\bf f}^{{\bf G}'(\tilde{s})})$ sont bien d\'efinis gr\^ace aux hypoth\`eses de r\'ecurrence pos\'ees en 1.1. Si $(G,\tilde{G},{\bf a})$ est quasi-d\'eploy\'e et \`a torsion int\'erieure et si ${\bf M}'={\bf M}$, un seul terme ne l'est pas, \`a savoir le terme $S_{{\bf M}}^{{\bf G}}(\boldsymbol{\delta},B^{\tilde{G}},{\bf f}^{{\bf G}})$ associ\'e \`a $\tilde{s}=\tilde{\zeta}=1$. On le remplace dans ce cas par $S_{\tilde{M}}^{\tilde{G}}(\boldsymbol{\delta},B^{\tilde{G}},{\bf f})$ qui est bien d\'efini. On a dans ce cas la simple \'egalit\'e $I_{\tilde{M}}^{\tilde{G},{\cal E}}({\bf M},\boldsymbol{\delta},{\bf f})=I_{\tilde{M}}^{\tilde{G}}(\boldsymbol{\delta},{\bf f})$.

L'application ${\bf f}\to I_{\tilde{M}}^{\tilde{G},{\cal E}}({\bf M}',\boldsymbol{\delta},{\bf f})$ se factorise en une application d\'efinie sur $I(\tilde{G}(F),\omega)\otimes Mes(G(F))$.

{\bf Remarque.} Conform\'ement \`a [I] 3.3(2), en supposant $\hat{M}$ standard, on pourrait aussi sommer sur 
$$\tilde{s}\in \tilde{\zeta} Z(\hat{M})^{\Gamma_F}/(Z(\hat{G})^{\Gamma_F}(1-\hat{\theta})(Z(\hat{M})^{\Gamma_F})),$$
(ou plus canoniquement sur $\tilde{\zeta}Z(\hat{M})^{\Gamma_F}/Z(\hat{G})^{\Gamma_F}$ \`a conjugaison pr\`es par $Z(\hat{M})^{\Gamma_{F}}$) \`a condition de multiplier les coefficients par $\vert det((1-\theta^{\tilde{M}})_{{\cal A}_{M}/({\cal A}_{\tilde{M}}+{\cal A}_{G}})\vert $.

\bigskip
Donnons une autre d\'efinition du coefficient $i_{\tilde{M}'}(\tilde{G},\tilde{G}'(\tilde{s}))$. On suppose ${\bf G}'(\tilde{s})$ elliptique. Parce que
$$Z(\hat{M})^{\Gamma_{F},\hat{\theta}}=Z(\hat{M})^{\Gamma_{F},\hat{\theta},0}Z(\hat{G})^{\Gamma_{F},\hat{\theta}}$$
et 
$$Z(\hat{M})^{\Gamma_{F},\hat{\theta},0}\subset Z(\hat{M})^{\Gamma_{F}}\cap \hat{T}^{\hat{\theta},0}\subset Z(\hat{M})^{\Gamma_{F},\hat{\theta}},$$
on voit que l'homomorphisme naturel
$$(Z(\hat{M})^{\Gamma_{F}}\cap \hat{T}^{\hat{\theta},0})/(Z(\hat{G})^{\Gamma_{F}}\cap \hat{T}^{\hat{\theta},0})\to Z(\hat{M})^{\Gamma_{F},\hat{\theta}}/Z(\hat{G})^{\Gamma_{F},\hat{\theta}}$$
est un isomorphisme. D'autre part, on a un homomorphisme naturel
$$(Z(\hat{M})^{\Gamma_{F}}\cap \hat{T}^{\hat{\theta},0})/(Z(\hat{G})^{\Gamma_{F}}\cap \hat{T}^{\hat{\theta},0})\to Z(\hat{M}')^{\Gamma_{F}}/Z(\hat{G}'(\tilde{s}))^{\Gamma_{F}}.$$
On peut donc l'interpr\'eter comme un homomorphisme
$$(1) \qquad  Z(\hat{M})^{\Gamma_{F},\hat{\theta}}/Z(\hat{G})^{\Gamma_{F},\hat{\theta}}\to Z(\hat{M}')^{\Gamma_{F}}/Z(\hat{G}'(\tilde{s}))^{\Gamma_{F}}.$$
Le groupe d'arriv\'ee n'est autre que $Z(\hat{M}'_{ad})^{\Gamma_{F}}$, o\`u $\hat{M}'_{ad}$ est 'image de $\hat{M}'$ dans $\hat{G}'(\tilde{s})_{AD}$. Donc ce groupe est connexe et l'homomorphisme (1) est surjectif. On a

(2) $i_{\tilde{M}'}(\tilde{G},\tilde{G}'(\tilde{s}))$ est l'inverse du nombre d'\'el\'ements du noyau de (1).

Preuve. On a un diagramme commutatif
$$\begin{array}{ccccccccc}&&&&&&1&&\\ &&&&&&\downarrow&&\\ &&1&&1&&A&&\\ &&\downarrow&&\downarrow&&\downarrow\\1&\to&Z(\hat{G})^{\Gamma_{F}}\cap \hat{T}^{\hat{\theta},0}&\to&Z(\hat{M})^{\Gamma_{F}}\cap \hat{T}^{\hat{\theta},0}&\to& Z(\hat{M})^{\Gamma_{F},\hat{\theta}}/Z(\hat{G})^{\Gamma_{F},\hat{\theta}}&\to&1\\&&\downarrow&&\downarrow&&\downarrow&&\\ 1&\to&Z(\hat{G}'(\tilde{s})^{\Gamma_{F}}&\to&Z(\hat{M}')^{\Gamma_{F}}&\to&Z(\hat{M}')^{\Gamma_{F}}/Z(\hat{G}'(\tilde{s}))^{\Gamma_{F}}&\to&1\\ &&\downarrow&&\downarrow&&\downarrow&&\\ &&B&&C&&1&&\\ &&\downarrow&&\downarrow&&&&\\&&1&&1&&&&\\ \end{array}$$ 
Les groupes $A$, $B$, $C$ sont d\'efinis de sorte que les colonnes soient exactes. Par un raisonnement d'alg\`ebre \'el\'ementaire, on en d\'eduit une suite exacte
$$1\to A\to B\to C\to 1$$
On a $Z(\hat{M})^{\Gamma_{F}}\cap \hat{T}^{\hat{\theta},0}=Z(\hat{M})\cap Z(\hat{M}')^{\Gamma_{F}}$ et $Z(\hat{G})^{\Gamma_{F}}\cap \hat{T}^{\hat{\theta},0}=Z(\hat{G})\cap Z(\hat{G}'(\tilde{s}))^{\Gamma_{F}}$. Donc 
$$i_{\tilde{M}'}(\tilde{G},\tilde{G}'(\tilde{s}))=\vert B\vert^{-1} \vert C\vert=\vert A\vert ^{-1}.$$
Mais $A$ est le noyau de (1). $\square$

\bigskip

{\bf Variante.} Supposons $(G,\tilde{G},{\bf a})$ quasi-d\'eploy\'e et \`a torsion int\'erieure (auquel cas on note simplement $ \zeta=\tilde{\zeta}$). Fixons un syst\`eme de fonctions $B$ comme en 1.9. Comme en 1.10, pour tout $s\in \zeta Z(\hat{M})^{\Gamma_F}/Z(\hat{G})^{\Gamma_F}$, ce syst\`eme de fonctions en d\'etermine un sur $\tilde{G}'(s)(F)$, que l'on note encore $B$.  Pour $\boldsymbol{\delta}\in D^{st}_{g\acute{e}om}({\bf M}')\otimes Mes(M'(F))^*$ et ${\bf f}\in C_{c}^{\infty}(\tilde{G}(F))\otimes Mes(G(F))$, on pose
$$ I_{\tilde{M}}^{\tilde{G},{\cal E}}({\bf M}',\boldsymbol{\delta},B,{\bf f})=\sum_{s\in \zeta Z(\hat{M})^{\Gamma_F}/Z(\hat{G})^{\Gamma_F}}i_{\tilde{M}'}(\tilde{G},\tilde{G}'(s))
S_{{\bf M}'}^{{\bf G}'(s)}(\boldsymbol{\delta},B,{\bf f}^{{\bf G}'(s)}).$$

{\bf Variante.} Supposons $G=\tilde{G}$ et ${\bf a}=1$. Fixons une fonction $B$ comme en 1.8. Comme dans la variante pr\'ec\'edente, pour tout $s\in \zeta Z(\hat{M})^{\Gamma_F}/Z(\hat{G})^{\Gamma_F}$, cette fonction se restreint en une fonction pour $G'(s)(F)$,  a fortiori comme un syst\`eme de fonctions  comme en 1.9 (pour tout $\epsilon\in G'(s)_{ss}(F)$, $B_{\epsilon}$ est la restriction de $B$ au syst\`eme de racines de $G'(s)_{\epsilon}$). Pour $\boldsymbol{\delta}\in D^{st}_{unip}({\bf M}')\otimes Mes(M'(F))^*$ et ${\bf f}\in C_{c}^{\infty}(G(F))\otimes Mes(G(F))$, on pose
$$ I_{M}^{G,{\cal E}}({\bf M}',\boldsymbol{\delta},B,{\bf f})=\sum_{s\in \zeta Z(\hat{M})^{\Gamma_F}/Z(\hat{G})^{\Gamma_F}}i_{M'}(G,G'(s))
S_{{\bf M}'}^{{\bf G}'(s)}(\boldsymbol{\delta},B,{\bf f}^{{\bf G}'(s)}).$$

 \bigskip

 \subsection{Action d'un groupe d'automorphismes}
 Les donn\'ees sont les m\^emes que dans le paragraphe pr\'ec\'edent. On a introduit en [I] 3.2 le groupe d'automorphismes $Aut(\tilde{M},{\bf M}')$. Il agit sur $D_{g\acute{e}om}^{st}({\bf M}')$ ou $D_{g\acute{e}om}^{st}({\bf M}')\otimes Mes(M'(F))^*$. On note cette action $(x,\boldsymbol{\delta})\mapsto x(\boldsymbol{\delta})$.
 
 \ass{Lemme}{Soient $x\in Aut(\tilde{M},{\bf M}')$ et $\boldsymbol{\delta}\in D_{g\acute{e}om}^{st}({\bf M}')\otimes Mes(M'(F))^*$. Pour tout ${\bf f}\in I(\tilde{G}(F),\omega)\otimes Mes(G(F))$, on a l'\'egalit\'e
 $$I_{\tilde{M}}^{\tilde{G},{\cal E}}({\bf M}',x(\boldsymbol{\delta}),{\bf f})=I_{\tilde{M}}^{\tilde{G},{\cal E}}({\bf M}',\boldsymbol{\delta},{\bf f}).$$}
 
 Preuve.  On r\'ealise $^LM$ comme espace de Levi standard de $^LG$. D'apr\`es la d\'efinition de [I] 3.2, $x$ est un \'el\'ement de $\hat{G}$ tel que $ad_{x}(\hat{M})=\hat{M}$, $ad_{x}({\cal M}')={\cal M}'$ et $ad_{x}(\tilde{\zeta})\in Z(\hat{M})\tilde{\zeta}$. Montrons que
 
 (1) $ad_{x}(\tilde{\zeta}) \in Z(\hat{M})^{\Gamma_{F}}\tilde{\zeta}$. 
 
  Rappelons que $\tilde{\zeta}m'=a(w)m'w(\tilde{\zeta})$ pour tout $(m',w)\in {\cal M}'$, o\`u $a$ est \`a valeurs dans $Z(\hat{G})$. Ecrivons $ad_{x}(\tilde{\zeta})=z\tilde{\zeta} $, avec $z\in Z(\hat{M})$. Pour $(m',w)\in {\cal M}'$, posons $m''=x^{-1}m'w(x)$. Alors $(m'',w)\in {\cal M}'$, donc $\tilde{\zeta}m''=a(w)m''w(\tilde{\zeta})$, c'est-\`a-dire
 $$\tilde{\zeta}x^{-1}m'w(x)=a(w)x^{-1}m'w(x)w(\tilde{\zeta}).$$
 Cela \'equivaut \`a $ad_{x}(\tilde{\zeta})m'=a(w)m'w(ad_{x}(\tilde{\zeta}))$, ou encore \`a $z\tilde{\zeta}m'=a(w)m'w(z)w(\tilde{\zeta})$. En comparant avec la premi\`ere \'egalit\'e de la preuve, on obtient $w(z)=z$. D'o\`u (1).

 D'autre part, $x$ normalise $\hat{M}$ et la classe $x\hat{M}$ est conserv\'ee par l'action galoisienne  et par $\hat{\theta}$ (parce que $ad_{x}$ normalise $^LM\hat{\theta}$, cf. [I] 3.2). Il en r\'esulte que la restriction de $ad_{x}$ \`a $Z(\hat{M})$ conserve aussi ces actions. Alors l'application $\tilde{s}\mapsto x\tilde{s}x^{-1}$ d\'efinit une bijection de $ \tilde{\zeta}Z(\hat{M})^{\Gamma_{F}}/(Z(\hat{G})^{\Gamma_{F}}(1-\hat{\theta})(Z(\hat{M})^{\Gamma_{F}}))$ sur lui-m\^eme. Les donn\'ees endoscopiques ${\bf G}'(\tilde{s})$ et ${\bf G}'(x\tilde{s}x^{-1})$ sont \'equivalentes, l'\'equivalence \'etant d\'efinie par $x$. Cette \'equivalence \'echange  les  syst\`emes de fonctions $B^{\tilde{G}}$ relatives aux deux  groupes, ainsi qu'il r\'esulte de leur d\'efinition. L'\'equivalence  d\'efinit un isomorphisme de $SI({\bf G}'(\tilde{s}))$ sur $SI({\bf G}'(x\tilde{s}x^{-1}))$. Cet isomorphisme envoie ${\bf f}^{{\bf G}'(\tilde{s})}$ sur ${\bf f}^{{\bf G}'(x\tilde{s}x^{-1})}$. Par restriction \`a ${\bf M}'$ puis par dualit\'e, il s'en d\'eduit un automorphisme de $D_{g\acute{e}om}^{st}({\bf M}')$ qui n'est autre que celui introduit avant l'\'enonc\'e. Il en r\'esulte que
 $$S_{{\bf M}'}^{{\bf G}'(x\tilde{s}x^{-1})}(x\boldsymbol{\delta},B^{\tilde{G}},{\bf f}^{{\bf G}'(x\tilde{s}x^{-1})})=S_{{\bf M}'}^{{\bf G}'(\tilde{s})}(\boldsymbol{\delta},B^{\tilde{G}},{\bf f}^{{\bf G}'(\tilde{s})}).$$
 On a aussi \'evidemment l'\'egalit\'e $i_{\tilde{M}'}(\tilde{G},\tilde{G}'(x\tilde{s}x^{-1}))=i_{\tilde{M}'}(\tilde{G},\tilde{G}'(\tilde{s}))$. L'\'enonc\'e r\'esulte alors du simple changement de variables $\tilde{s}\mapsto x\tilde{s}x^{-1}$ dans la d\'efinition de $I_{\tilde{M}}^{\tilde{G},{\cal E}}({\bf M}',x(\boldsymbol{\delta}),{\bf f})$. $\square$
 
 \ass{Corollaire}{ Soit  $\boldsymbol{\delta}\in D_{g\acute{e}om}^{st}({\bf M}')\otimes Mes(M'(F))^*$. Supposons que l'une des conditions suivantes soit v\'erifi\'ee:
 
 (i) la projection de $\boldsymbol{\delta}$ sur le sous-espace des \'el\'ements de $D_{g\acute{e}om}^{st}({\bf M}')\otimes Mes(M'(F))^*$ invariants par l'action de $Aut(\tilde{M},{\bf M}')$ est nulle;
 
 (ii) le support de $\boldsymbol{\delta}$ ne coupe pas l'ensemble des $\delta\in \tilde{M}'(F)$ tels que $N^{\tilde{M}',\tilde{M}}(\delta)\in N^{\tilde{M}}(\tilde{M}_{ab}(F))$. 
 
 Alors on a l'\'egalit\'e $I_{\tilde{M}}^{\tilde{G},{\cal E}}({\bf M}',\boldsymbol{\delta},{\bf f})=0$ pour  tout ${\bf f}\in I(\tilde{G}(F))\otimes Mes(G(F))$.}

Preuve.   D'apr\`es le lemme, on peut remplacer $\boldsymbol{\delta}$ par sa projection sur le sous-espace des \'el\'ements de $D_{g\acute{e}om}^{st}({\bf M}')\otimes Mes(M'(F))^*$ invariants par l'action de $Aut(\tilde{M},{\bf M}')$. La conclusion s'ensuit sous l'hypoth\`ese (i). Par ailleurs,  gr\^ace \`a [I] lemme 2.6, on voit que  (ii) entra\^{\i}ne (i). $\square$

{\bf Variante.} Supposons $(G,\tilde{G},{\bf a})$ quasi-d\'eploy\'e et \`a torsion int\'erieure. Fixons un syst\`eme de fonctions $B$ comme en 1.9. Les r\'esultats ci-dessus valent aussi pour les distributions $I_{\tilde{M}}^{\tilde{G},{\cal E}}({\bf M}',\boldsymbol{\delta},B,{\bf f})$.

{\bf Variante.} Supposons $G=\tilde{G}$ et ${\bf a}=1$. Fixons une fonction $B$ comme en 1.8. Les r\'esultats ci-dessus valent aussi pour les distributions $I_{M}^{G,{\cal E}}({\bf M}',\boldsymbol{\delta},B,{\bf f})$. 

\bigskip

\subsection{Formules de descente}
  Les donn\'ees sont les m\^emes qu'en 1.12. On va consid\'erer trois situations dans lesquelles on a des formules de descente pour les distributions introduites en 1.10 et 1.12.
  
  (a) Soit $R'$ un groupe de Levi de $M'$ qui est relevant. Modulo certains choix, on construit comme en [I] 3.4 un sous-espace de Levi $\tilde{R}$ de $\tilde{M}$ et une donn\'ee endoscopique elliptique et relevante ${\bf R}'$ de $\tilde{R}$. On dispose d'un homomorphisme
  $$\begin{array}{ccc}I({\bf M}')\otimes Mes(M'(F))&\to&I({\bf R}')\otimes Mes(R'(F))\\ \boldsymbol{\varphi}&\mapsto &\boldsymbol{\varphi}_{{\bf R}'}\\ \end{array}$$
  et d'un homomorphisme dual
  $$\begin{array}{ccc}D_{g\acute{e}om}({\bf R}')\otimes Mes(R'(F))^*&\to&D_{g\acute{e}om}({\bf M}')\otimes Mes(M'(F))^*\\ \boldsymbol{\delta}&\mapsto &\boldsymbol{\delta}^{\bf M'}\end{array}$$
  qui pr\'eserve la stabilit\'e.
  
  (b)  On suppose  que $(G,\tilde{G},{\bf a})$ est quasi-d\'eploy\'e et \`a torsion int\'erieure et on fixe un syst\`eme de fonctions $B$ comme en 1.9. Soit $R$ un groupe de Levi de $M$. On a les m\^emes homomorphismes que ci-dessus pour ${\bf M}'={\bf M}$ et ${\bf R}'={\bf R}$.  On r\'ealise $\hat{R}$ et $\hat{M}$ comme des groupes de Levi standard de $\hat{G}$. Posons dans ce cas la d\'efinition suivante. Soit $\tilde{L}\in {\cal L}(\tilde{R})$. Il lui correspond un \'el\'ement $\hat{L}\in {\cal L}(\hat{R})$. Alors
  $$e_{\tilde{R}}^{\tilde{G}}(\tilde{M},\tilde{L})=\left\lbrace\begin{array}{cc}d_{\tilde{R}}^{\tilde{G}}(\tilde{M},\tilde{L})[(Z(\hat{M})^{\Gamma_{F}}\cap Z(\hat{L})^{\Gamma_{F}}):Z(\hat{G})^{\Gamma_{F}}]^{-1},&\text{ si } {\cal A}_{M}^G\oplus {\cal A}_{L}^G={\cal A}_{R}^G,\\ 0,&\text{ sinon.}\\ \end{array}\right.$$
  Remarquons que l'hypoth\`ese ${\cal A}_{M}^G\oplus {\cal A}_{L}^G={\cal A}_{R}^G$ entra\^{\i}ne que le quotient  $(Z(\hat{M})^{\Gamma_{F}}\cap Z(\hat{L})^{\Gamma_{F}})/Z(\hat{G})^{\Gamma_{F}}$ est fini. Remarquons aussi que le terme $[(Z(\hat{M})^{\Gamma_{F}}\cap Z(\hat{L})^{\Gamma_{F}}):Z(\hat{G})^{\Gamma_{F}}]^{-1}$ peut s'interpr\'eter comme l'inverse du nombre d'\'el\'ements du noyau de l'homomorphisme naturel
  $$Z(\hat{M})^{\Gamma_{F}}/Z(\hat{G})^{\Gamma_{F}}  \to Z(\hat{R})^{\Gamma_{F}}/Z(\hat{L})^{\Gamma_{F}}.$$
  
  (c) Soit $R'$ un groupe de Levi de $M'$ qui n'est pas relevant. L'espace $D_{g\acute{e}om}({\bf R}')$ n'est pas d\'efini. N\'eanmoins, fixons des donn\'ees suppl\'ementaires $M'_{1},...,\Delta_{1}$ pour ${\bf M}'$. On a alors un homomorphisme d'induction
  $$\begin{array}{ccc}D_{g\acute{e}om,\lambda_{1}}(\tilde{R}'_{1}(F))\otimes Mes(R'(F))^*&\to &D_{g\acute{e}om,\lambda_{1}}(\tilde{M}'_{1}(F))\otimes Mes(M'(F))^*\\ && \parallel\\&& D_{g\acute{e}om}({\bf M}')\otimes Mes(M'(F))^*\\\boldsymbol{\delta}&\mapsto&\boldsymbol{\delta}^{\bf M'}\\ \end{array}$$

  \ass{Proposition}{(i) Dans la situation (a), soient $\boldsymbol{\delta}\in D_{g\acute{e}om}^{st}({\bf R}')\otimes Mes(R'(F))^*$ et ${\bf f}\in I(\tilde{G}(F),\omega)\otimes Mes(G(F))$. On a l'\'egalit\'e
  $$I_{\tilde{M}}^{\tilde{G},{\cal E}}({\bf M}',\boldsymbol{\delta}^{\bf M'},{\bf f})=\sum_{\tilde{L}\in {\cal L}(\tilde{R})}d_{\tilde{R}}^{\tilde{G}}(\tilde{M},\tilde{L})I_{\tilde{R}}^{\tilde{L},{\cal E}}({\bf R}',\boldsymbol{\delta},{\bf f}_{\tilde{L},\omega}).$$

  (ii) Dans la situation (b), soient $\boldsymbol{\delta}\in D_{g\acute{e}om}^{st}({\bf R})\otimes Mes(R(F))^*$ et ${\bf f}\in I(\tilde{G}(F))\otimes Mes(G(F))$. On a l'\'egalit\'e
  $$S_{\tilde{M}}^{\tilde{G}}(\boldsymbol{\delta}^{M},B,{\bf f})=\sum_{\tilde{L}\in {\cal L}(\tilde{R})}e_{\tilde{R}}^{\tilde{G}}(\tilde{M},\tilde{L})S_{\tilde{R}}^{\tilde{L}}(\boldsymbol{\delta},B,{\bf f}_{\tilde{L}}).$$
  
  (iii) Dans la situation (c), soient $\boldsymbol{\delta}\in D_{g\acute{e}om,\lambda_{1}}(\tilde{R}'_{1}(F))\otimes Mes(R'(F))^*$ et ${\bf f}\in I(\tilde{G}(F),\omega)\otimes Mes(G(F))$. On a l'\'egalit\'e $I_{\tilde{M}}^{\tilde{G},{\cal E}}({\bf M}',\boldsymbol{\delta}^{\bf M'},{\bf f})=0$.}

  Preuve. On choisit une paire de Borel \'epingl\'ee $\hat{{\cal E}}=(\hat{B},\hat{T}, (\hat{E}_{\alpha})_{\alpha\in \Delta})$ de $\hat{G}$ comme en [I] 1.5.   Si $(G,\tilde{G},{\bf a})$ est quasi-d\'eploy\'e et \`a torsion int\'erieure et si ${\bf M}'={\bf M}$, la formule (i) n'est autre que celle du lemme 1.7. On exclut ce cas. On peut supposer que $\hat{R}\subset \hat{M}$ sont des Levi standard de $\hat{M}$. On \'ecrit ${\bf R}'=(R',{\cal R}',\tilde{\zeta})$. On peut supposer que ${\bf M}'=(M',{\cal M}',\tilde{\zeta})$, avec ${\cal M}'=\hat{M}'{\cal R}'$, cf. [I] 3.4. Rappelons la d\'efinition
 $$(1) \qquad I_{\tilde{M}}^{\tilde{G},{\cal E}}({\bf M}',\boldsymbol{\delta}^{\bf M'},{\bf f})=\sum_{\tilde{s}\in \tilde{\zeta} Z(\hat{M})^{\Gamma_F,\hat{\theta}}/Z(\hat{G})^{\Gamma_F,\hat{\theta}}}i_{\tilde{M}'}(\tilde{G},\tilde{G}'(\tilde{s}))
S_{{\bf M}'}^{{\bf G}'(\tilde{s})}(\boldsymbol{\delta}^{\bf M'},B^{\tilde{G}},{\bf f}^{{\bf G}'(\tilde{s})}).$$
Pour chaque $\tilde{s}$, on fixe des donn\'ees auxiliaires $G'(\tilde{s})_{1},...,\Delta(\tilde{s})_{1}$.
On note $\lambda(\tilde{s})_{1}$ le caract\`ere  associ\'e de $C(\tilde{s})_{1}(F)$ et $\tilde{M}'(\tilde{s})_{1}$  l'image r\'eciproque de $\tilde{M}'$ dans $\tilde{G}'(\tilde{s})_{1}$. On peut remplacer $S_{{\bf M}'}^{{\bf G}'(\tilde{s})}(\boldsymbol{\delta}^{\bf M'},{\bf f}^{{\bf G}'(\tilde{s})})$ par $S_{\tilde{M}'(\tilde{s})_{1},\lambda(\tilde{s})_{1}}^{\tilde{G}'(\tilde{s})_{1}}(\boldsymbol{\delta}(\tilde{s})_{1}^{M'(\tilde{s})_{1}},{\bf f}^{\tilde{G}'(\tilde{s})_{1}})$. L'assertion (ii) se g\'en\'eralise au cas o\`u les fonctions et distributions se transforment selon un caract\`ere d'un tore central. La preuve est formelle. On applique cette assertion par r\'ecurrence \`a chacun des termes du second membre.   On obtient
  $$(2) \qquad I_{\tilde{M}}^{\tilde{G},{\cal E}}({\bf M}',\boldsymbol{\delta}^{\bf M'},{\bf f})=\sum_{\tilde{s}\in \tilde{\zeta} Z(\hat{M})^{\Gamma_F,\hat{\theta}}/Z(\hat{G})^{\Gamma_F,\hat{\theta}}}i_{\tilde{M}'}(\tilde{G},\tilde{G}'(\tilde{s}))$$
  $$\sum_{\tilde{L}'_{\tilde{s}}\in{\cal L}^{\tilde{G}'(\tilde{s})}(\tilde{R}')}e_{\tilde{R}'(\tilde{s})_{1}}^{\tilde{G}'(\tilde{s})_{1}}(\tilde{M}'(\tilde{s})_{1},\tilde{L}'_{\tilde{s},1})S_{\tilde{R'}(\tilde{s})_{1},\lambda(\tilde{s})_{1}}^{\tilde{L}'_{\tilde{s},1}}(\boldsymbol{\delta}(\tilde{s})_{1},B^{\tilde{G}},({\bf f}^{\tilde{G}'(\tilde{s})_{1}})_{\tilde{ L}'_{\tilde{s},1}})$$
  ($\tilde{L}'_{\tilde{s},1}$ est l'image r\'eciproque de $\tilde{L}'_{\tilde{s}}$ dans $\tilde{G}'(\tilde{s})_{1}$).  Rappelons que l'on peut identifier ${\cal L}(\tilde{R})$ \`a un sous-ensemble de ${\cal L}(\hat{R})$ et de m\^eme, pour tout $\tilde{s}$, ${\cal L}^{\tilde{G}'(\tilde{s})}(\tilde{R}')$ \`a un sous-ensemble de ${\cal L}^{\hat{G}'(\tilde{s})}(\hat{R}')$. Soient $\tilde{s}\in \tilde{\zeta} Z(\hat{M})^{\Gamma_F,\hat{\theta}}/Z(\hat{G})^{\Gamma_F,\hat{\theta}}$ et $\tilde{L}'_{\tilde{s}}\in{\cal L}^{\tilde{G}'(\tilde{s})}(\tilde{R}')$. L'espace ${\cal A}_{L'}$ est inclus dans ${\cal A}_{R'}$ qui s'identifie \`a ${\cal A}_{\tilde{R}}$ puisque ${\bf R}'$ est une donn\'ee elliptique de $(R,\tilde{R},{\bf a})$. Un raisonnement standard montre qu'il existe un unique $\tilde{L}\in {\cal L}(\tilde{R})$ de sorte que ${\cal A}_{L'}$ s'identifie \`a ${\cal A}_{\tilde{L}}$. Alors $\hat{L}'_{\tilde{s}}$ est \'egal \`a l'intersection de $\hat{G}'(\tilde{s})$ avec $\hat{L}$ et aussi \`a la composante neutre du commutant de $\tilde{s}$ dans $\hat{M}$. On introduit le groupe ${\cal L}'(\tilde{s})=\hat{L}'_{\tilde{s}}{\cal R}'$. Alors $(L'_{\tilde{s}},{\cal L}'(\tilde{s}),\tilde{s})$ n'est autre que la donn\'ee endoscopique ${\bf L}'(\tilde{s})$ de $(L,\tilde{L},{\bf a})$. Cette donn\'ee est elliptique par construction et est relevante puisqu'elle "contient" ${\bf R}'$ qui l'est par hypoth\`ese. Les donn\'ees $L'_{\tilde{s},1},...$ obtenues par restriction de celles fix\'ees pour ${\bf G}'(\tilde{s})$ sont des donn\'ees auxiliaires pour ${\bf L}'(\tilde{s})$. Enfin on a l'\'egalit\'e $({\bf f}^{\tilde{G}'(\tilde{s})_{1}})_{\tilde{ L}'_{\tilde{s},1}}=({\bf f}_{\tilde{L},\omega})^{\tilde{L}'_{\tilde{s},1}}$. Tout cela montre que l'on a
$$S_{\tilde{R'}(\tilde{s})_{1},\lambda(\tilde{s})_{1}}^{\tilde{L}'_{\tilde{s},1}}(\boldsymbol{\delta}(\tilde{s})_{1},B^{\tilde{G}},({\bf f}^{\tilde{G}'(\tilde{s})_{1}})_{\tilde{ L}'_{\tilde{s},1}})=S_{{\bf R}'}^{{\bf L}'(\tilde{s})}(\boldsymbol{\delta},B^{\tilde{G}},({\bf f}_{\tilde{L},\omega})^{{\bf L}'(\tilde{s})}).$$
Il est clair que $d_{\tilde{R}'(\tilde{s})_{1}}^{\tilde{G}'(\tilde{s})_{1}}(\tilde{M}'(\tilde{s})_{1},\tilde{L}'_{\tilde{s},1})=d_{\tilde{R}'}^{\tilde{G}'(\tilde{s})}(\tilde{M}',\tilde{L}'(\tilde{s}))$ (par exemple, ${\cal A}_{R'(\tilde{s})_{1}}^{G'(\tilde{s})_{1}}={\cal A}_{R'}^{G'(\tilde{s})}$). On v\'erifie que l'homomorphisme naturel
$$(Z(\hat{M}')^{\Gamma_{F}}\cap Z(\hat{L}'(\tilde{s}))^{\Gamma_{F}})/Z(\hat{G}'(\tilde{s}))^{\Gamma_{F}}\to (Z(\hat{M}'(\tilde{s})_{1})^{\Gamma_{F}}\cap Z(\hat{L}'_{\tilde{s},1})^{\Gamma_{F}})/Z(\hat{G}'(\tilde{s})_{1})^{\Gamma_{F}}$$
est bijectif (cf. la preuve de la relation (6) de 1.10). On en d\'eduit que $e_{\tilde{R}'(\tilde{s})_{1}}^{\tilde{G}'(\tilde{s})_{1}}(\tilde{M}'(\tilde{s})_{1},\tilde{L}'_{\tilde{s},1})=e_{\tilde{R}'}^{\tilde{G}'(\tilde{s})}(\tilde{M}',\tilde{L}'(\tilde{s}))$. 

On est parti d'un couple $(\tilde{s},\tilde{L}'_{\tilde{s}})$ et on lui a associ\'e $\tilde{L}\in {\cal L}(\tilde{R})$. A fortiori, on peut lui associer le couple $(\tilde{s},\tilde{L})$. On voit que l'on obtient une bijection de notre ensemble de couples $(\tilde{s},\tilde{L}'_{\tilde{s}})$ sur celui des couples $(\tilde{s},\tilde{L})$ pour lequel la donn\'ee endoscopique ${\bf L}'(\tilde{s})$ est elliptique.

 Utilisons les relations ci-dessus et regroupons les $\tilde{L}'_{\tilde{s}}$ qui interviennent dans la formule (2) selon l'espace de Levi $\tilde{L}$ que l'on vient de leur associer. On obtient
$$(3) \qquad I_{\tilde{M}}^{\tilde{G},{\cal E}}({\bf M}',\boldsymbol{\delta}^{\bf M'},{\bf f})=\sum_{\tilde{L}\in {\cal L}(\tilde{R})}\sum_{\tilde{s}\in \tilde{\zeta }Z(\hat{M})^{\Gamma_F,\hat{\theta}}/Z(\hat{G})^{\Gamma_F,\hat{\theta}}; {\bf L}'(\tilde{s})\text{ elliptique}}i_{\tilde{M}'}(\tilde{G},\tilde{G}'(\tilde{s}))$$
$$e_{\tilde{R}'}^{\tilde{G}'(\tilde{s})}(\tilde{M}',\tilde{L}'(\tilde{s}))S_{{\bf R}'}^{{\bf L}'(\tilde{s})}(\boldsymbol{\delta},B^{\tilde{G}},({\bf f}_{\tilde{L},\omega})^{{\bf L}'(\tilde{s})}).$$
Fixons $\tilde{L}\in {\cal L}(\tilde{R})$.Pour chaque espace $\tilde{L}'(\tilde{s})$ apparaissant ci-dessus, les syst\`emes de fonctions $B^{\tilde{G}}$ et $B^{\tilde{L}}$ sont les m\^emes, ce qui nous autorise \`a remplacer le premier  par le second. Un \'el\'ement $\tilde{s}\in \tilde{\zeta} Z(\hat{M})^{\Gamma_F,\hat{\theta}}/Z(\hat{G})^{\Gamma_F,\hat{\theta}}$ n'intervient effectivement dans la formule ci-dessus que si ${\bf L}'(\tilde{s})$ est elliptique, ${\bf G}'(\tilde{s})$ l'est aussi (d'apr\`es la d\'efinition de $i_{\tilde{M}'}(\tilde{G},\tilde{G}'(\tilde{s}))$) et ${\cal A}_{R'}^{G'(\tilde{s})}={\cal A}_{R'}^{M'}\oplus {\cal A}_{R'}^{L'(\tilde{s})}$ (d'apr\`es la d\'efinition de $e_{\tilde{R}'}^{\tilde{G}'(\tilde{s})}(\tilde{M}',\tilde{L}'(\tilde{s}))$). Les deux premi\`eres conditions plus les hypoth\`eses que ${\bf M}'$ et ${\bf R}'$ sont elliptiques entra\^{\i}nent les \'egalit\'es
${\cal A}_{R'}^{G'(\tilde{s})}={\cal A}_{\tilde{R}}^{\tilde{G}}$, ${\cal A}_{R'}^{ M'}={\cal A}_{\tilde{R}}^{\tilde{M}}$, ${\cal A}_{R'}^{L'(\tilde{s})}={\cal A}_{\tilde{R}}^{\tilde{L}}$. L'\'egalit\'e pr\'ec\'edente devient ${\cal A}_{\tilde{R}}^{\tilde{G}}={\cal A}_{\tilde{R}}^{\tilde{M}}\oplus {\cal A}_{\tilde{R}}^{\tilde{L}}$. Plus pr\'ecis\'ement les rapports de mesures sont les m\^emes, c'est-\`a-dire que $d_{\tilde{R}'}^{\tilde{G}'(\tilde{s})}(\tilde{M}',\tilde{L}'(\tilde{s}))=d_{\tilde{R}}^{\tilde{G}}(\tilde{M},\tilde{L})$. Inversement, si ce dernier nombre n'est pas nul et si ${\bf L}'(\tilde{s})$ est elliptique, ${\bf G}'(\tilde{s})$ l'est aussi. En effet, l'espace ${\cal A}_{R'}^{G'(\tilde{s})}$ contient  ${\cal A}_{R'}^{ M'}$ et ${\cal A}_{R'}^{L'(\tilde{s})}$ puisque $M'$ et $L'(\tilde{s})$ sont des Levi de $G'(\tilde{s})$. Il contient donc leur somme, laquelle est ${\cal A}_{\tilde{R}}^{\tilde{G}}$, ce qui assure l'ellipticit\'e.  

On suppose d\'esormais $d_{\tilde{R}}^{\tilde{G}}(\tilde{M},\tilde{L})\not=0$. En se rappelant la d\'efinition des diff\'erents coefficients, on peut donc r\'ecrire la sous-somme de (3) index\'ee par $\tilde{L}$ sous la forme
$$(4) \qquad d_{\tilde{R}}^{\tilde{G}}(\tilde{M},\tilde{L})\sum_{\tilde{s}\in \tilde{\zeta} Z(\hat{M})^{\Gamma_F,\hat{\theta}}/Z(\hat{G})^{\Gamma_F,\hat{\theta}}; {\bf L}'(\tilde{s})\text{ elliptique}}X(\tilde{s})S_{{\bf R}'}^{{\bf L}'(\tilde{s})}(\boldsymbol{\delta},B^{\tilde{L}},({\bf f}_{\tilde{L},\omega})^{{\bf L}'(\tilde{s})}),$$
o\`u
$$X(\tilde{s})= i_{\tilde{M}'}(\tilde{G},\tilde{G}'(\tilde{s}))
[(Z(\hat{M}')^{\Gamma_{F}}\cap Z(\hat{L}'(\tilde{s}))^{\Gamma_{F}}):Z(\hat{G}'(\tilde{s}))^{\Gamma_{F}}]^{-1}.$$
Les deux facteurs composant $X(\tilde{s})$ sont les inverses des nombres d'\'el\'ements des noyaux des homomorphismes
$$Z(\hat{M})^{\Gamma_{F},\hat{\theta}}/Z(\hat{G})^{\Gamma_{F},\hat{\theta}}\to Z(\hat{M}')^{\Gamma_{F}}/Z(\hat{G}'(\tilde{s}))^{\Gamma_{F}}$$
et 
$$Z(\hat{M}')^{\Gamma_{F}}/Z(\hat{G}'(\tilde{s}))^{\Gamma_{F}}\to Z(\hat{R}')^{\Gamma_{F}}/Z(\hat{L}'(\tilde{s}))^{\Gamma_{F}}.$$
Ces deux homomorphismes \'etant surjectifs, $X(\tilde{s})$ est l'inverse du nombre d'\'el\'ements du noyau de l'homomorphisme compos\'e
$$p_{1}(\tilde{s}):Z(\hat{M})^{\Gamma_{F},\hat{\theta}}/Z(\hat{G})^{\Gamma_{F},\hat{\theta}}\to  Z(\hat{R}')^{\Gamma_{F}}/Z(\hat{L}'(\tilde{s}))^{\Gamma_{F}}.$$
La donn\'ee ${\bf G}'(\tilde{s})$ a maintenant disparu et tous les termes ne d\'ependent que de la classe $\tilde{s}Z(\hat{L})^{\Gamma_{F},\hat{\theta}}$. Rappelons que l'hypoth\`ese ${\cal A}_{\tilde{R}}^{\tilde{G}}={\cal A}_{\tilde{R}}^{\tilde{M}}\oplus {\cal A}_{\tilde{R}}^{\tilde{L}}$ entra\^{\i}ne dualement que l'homomorphisme
$$(Z(\hat{M})^{\Gamma_F,\hat{\theta},0}\times Z(\hat{L})^{\Gamma_F,\hat{\theta},0})/diag_{-}(Z(\hat{G})^{\Gamma_F,\hat{\theta},0})\to Z(\hat{R})^{\Gamma_F,\hat{\theta},0}$$
est surjectif de noyau fini ($diag_{-}$ est le plongement antidiagonal). Rappelons aussi que
  $$Z(\hat{R})^{\Gamma_F,\hat{\theta}}=Z(\hat{R})^{\Gamma_F,\hat{\theta},0}Z(\hat{L})^{\Gamma_F,\hat{\theta}}.$$
  On en d\'eduit que l'homomorphisme naturel.
$$p_{2}:Z(\hat{M})^{\Gamma_F,\hat{\theta}}/Z(\hat{G})^{\Gamma_F,\hat{\theta}}\to Z(\hat{R})^{\Gamma_F,\hat{\theta}}/Z(\hat{L})^{\Gamma_F,\hat{\theta}}$$
 est aussi surjectif de noyau fini.  On peut r\'ecrire la formule (4) en sommant sur l'espace d'arriv\'ee de cet homomorphisme plut\^ot que sur son espace de d\'epart. On obtient
 $$(5) \qquad d_{\tilde{R}}^{\tilde{G}}(\tilde{M},\tilde{L})\sum_{\tilde{s}\in \tilde{\zeta} Z(\hat{R})^{\Gamma_F,\hat{\theta}}/Z(\hat{L})^{\Gamma_F,\hat{\theta}}; {\bf L}'(\tilde{s})\text{ elliptique}}\vert Ker(p_{2})\vert X(\tilde{s})S_{{\bf R}'}^{{\bf L}'(\tilde{s})}(\boldsymbol{\delta},B^{\tilde{L}},({\bf f}_{\tilde{L},\omega})^{{\bf L}'(\tilde{s})}).$$
 On a l'\'egalit\'e
 $$(6) \qquad  \vert Ker(p_{2})\vert X(\tilde{s})=i_{\tilde{R}'}(\tilde{L},\tilde{L}'(\tilde{s})).$$
 En effet, l'homomorphisme $p_{1}(\tilde{s})$ se factorise en 
 $$Z(\hat{M})^{\Gamma_{F},\hat{\theta}}/Z(\hat{G})^{\Gamma_{F},\hat{\theta}}\stackrel{p_{2}}{\to}Z(\hat{R})^{\Gamma_F,\hat{\theta}}/Z(\hat{L})^{\Gamma_F,\hat{\theta}} \to  Z(\hat{R}')^{\Gamma_{F}}/Z(\hat{L}'(\tilde{s}))^{\Gamma_{F}}.$$
 Ces deux homomorphismes \'etant surjectifs et $X(\tilde{s})$ \'etant l'inverse du nombre d'\'el\'ements de leur compos\'e, le produit $\vert Ker(p_{2})\vert X(\tilde{s})$ est l'inverse du nombre d'\'el\'ements du noyau du second homomorphisme. C'est $i_{\tilde{R}'}(\tilde{L},\tilde{L}'(\tilde{s}))$ par d\'efinition de ce terme.

Rempla\c{c}ons $\vert Ker(p_{2})\vert X(\tilde{s})$ par  $i_{\tilde{R}'}(\tilde{L},\tilde{L}'(\tilde{s}))$ dans la formule (5).  Cela nous permet de supprimer la condition ${\bf L}'(\tilde{s})$ elliptique puisque ce terme est nul si cette condition n'est pas v\'erifi\'ee.   Alors (5) co\"{\i}ncide avec le produit de $d_{\tilde{R}}^{\tilde{G}}(\tilde{M},\tilde{L})$ avec la formule qui d\'efinit $I_{\tilde{R}}^{\tilde{L},{\cal E}}({\bf R}',\boldsymbol{\delta},{\bf f}_{\tilde{L},\omega})$. Reportons ensuite (5) dans l'\'egalit\'e (3). On obtient
$$I_{\tilde{M}}^{\tilde{G},{\cal E}}({\bf M}',\boldsymbol{\delta}^{\bf M'},{\bf f})=\sum_{\tilde{L}\in {\cal L}(\tilde{R})}d_{\tilde{R}}^{\tilde{G}}(\tilde{M},\tilde{L}) I_{\tilde{R}}^{\tilde{L},{\cal E}}({\bf R}',\boldsymbol{\delta},{\bf f}_{\tilde{L},\omega}).$$
C'est l'\'egalit\'e du (i) de l'\'enonc\'e.

La preuve de (ii) est similaire, \`a ceci pr\`es que l'on raisonne par r\'ecurrence. On part de l'\'egalit\'e analogue \`a (1)
$$I_{\tilde{M}}^{\tilde{G}}(\boldsymbol{\delta}^{M},B,{\bf f})=S_{\tilde{M}}^{\tilde{G}}(\boldsymbol{\delta}^M,B,{\bf f})+\sum_{s\in Z(\hat{M})^{\Gamma_{F}}/Z(\hat{G})^{\Gamma_{F}}, s\not=1}i_{\tilde{M}}(\tilde{G},\tilde{G}'(s))S_{{\bf M}}^{{\bf G}'(s)}(\boldsymbol{\delta}^{M},B,{\bf f}^{{\bf G}'(s)}).$$
Pour le terme associ\'e \`a $s\not=1$, on peut appliquer par r\'ecurrence la relation (ii) au terme index\'e par $s$. Pour le premier terme du membre de droite, on l'applique aussi mais, puisqu'on ne sait pas encore qu'elle est vraie, on doit ajouter la diff\'erence $X$ entre le membre de gauche et celui de droite de l'\'egalit\'e du (ii). Le calcul se poursuit (c'en est un cas particulier) et on obtient finalement
$$I_{\tilde{M}}^{\tilde{G}}(\boldsymbol{\delta}^{M},B,{\bf f})=X+\sum_{\tilde{L}\in {\cal L}(\tilde{R})}d_{\tilde{R}}^{\tilde{G}}(\tilde{M},\tilde{L}) I_{\tilde{R}}^{\tilde{L}}(\boldsymbol{\delta},B,{\bf f}_{\tilde{L}}).$$
Il reste \`a appliquer le lemme 1.7 (dont on a dit qu'il se g\'en\'eralisait aux distributions relatives au syst\`eme de  fonctions $B$) pour conclure $X=0$, ce que l'on voulait prouver. 

La preuve de (iii) est plus d\'elicate.  Montrons d'abord que l'on peut imposer des hypoth\`eses suppl\'ementaires aux donn\'ees $R'$ et $\boldsymbol{\delta}$.  Du c\^ot\'e des groupes duaux, la situation est la m\^eme que dans le cas (i). On peut d\'efinir $\hat{R}$ comme la composante neutre du commutant de $Z(\hat{R}')^{\Gamma_{F},0}$ dans $\hat{G}$.  On a d\'ej\`a suppos\'e $\hat{M}$ standard et, par  un proc\'ed\'e analogue, on peut supposer que $\hat{R}$ est lui-aussi  standard. Les deux Levi $\hat{M}$ et $\hat{R}$ sont invariants par $\Gamma_{F}$ et par $\hat{\theta}$. On pose ${\cal R}'={\cal M}'\cap(\hat{R}\rtimes W_{F})$.  On a encore ${\cal M}'=\hat{M}'{\cal R}'$. Deux cas sont possibles. Le  premier est

(7) $\hat{R}$ ne correspond \`a aucun Levi de $G$.

Supposons au contraire que $\hat{R}$ corresponde \`a un Levi $R$ de $G$. Dans ce cas, $R$ s'\'etend naturellement en un espace de Levi $\tilde{R}$ de $\tilde{G}$ et $(R',{\cal R}',\tilde{\zeta})$ est une donn\'ee endoscopique de $(R,\tilde{R},{\bf a})$. Alors l'hypoth\`ese que $R'$ n'est pas relevant signifie que cette donn\'ee endoscopique n'est pas relevante.  On dispose des applications
$$\begin{array}{ccc}\tilde{R}(F)&&\\&\searrow N^{\tilde{R}}&\\ &&\tilde{R}_{0,ab}(F)\\&\nearrow N^{\tilde{R}',\tilde{R}}&\\ \tilde{R}'(F)&&\\ \end{array}.$$
Notons $\tilde{R}'(F)^{in}$, resp. $\tilde{R}'(F)^{out}$, l'ensemble des $\gamma\in \tilde{R}'(F)$ tels que $N^{\tilde{R}',\tilde{R}}(\gamma)$ appartient \`a l'image de $N^{\tilde{R}}$, resp. n'appartient pas \`a cette image. L'ensemble $\tilde{R}'(F)$ est union disjointe de $\tilde{R}'(F)^{in}$ et $\tilde{R}'(F)^{out}$. Ces deux ensembles sont ouverts, ferm\'es et invariants par conjugaison stable. Parce que $R'$ n'est pas relevant, $\tilde{R}'(F)^{in}$ ne contient aucun \'el\'ement elliptique et fortement $\tilde{R}$-r\'egulier ([I] proposition 1.14). Par lin\'earite, on peut supposer que le support de $\boldsymbol{\delta}$ est form\'e d'\'el\'ements dont la partie semi-simple appartient \`a une classe de conjugaison stable fix\'ee. Supposons que cette classe soit contenue dans $\tilde{R}'(F)^{in}$.  Fixons $\epsilon$ dans cette classe. Parce que $\epsilon$ n'appartient pas \`a un sous-tore tordu elliptique de $\tilde{R}$, l'inclusion $A_{R'}\subset A_{R'_{\epsilon}}$ est stricte. On introduit le Levi $S'$ de $R'$ tel que $A_{S'}=A_{R'_{\epsilon}}$. C'est un Levi propre. D'apr\`es [I] lemme 5.12,  il existe   un \'el\'ement $\boldsymbol{\sigma}\in D^{st}_{g\acute{e}om,\lambda_{1}}(\tilde{S}'_{1}(F))\otimes Mes(S'(F))^*$  tels que $\boldsymbol{\delta}=\boldsymbol{\sigma}^{R'}$.   Evidemment, le Levi $S'$ est encore moins relevant que $R'$. En raisonnant par r\'ecurrence sur la dimension de $R'$, on peut supposer $I_{\tilde{M}}^{\tilde{G},{\cal E}}({\bf M}',\boldsymbol{\sigma}^{\bf M'},{\bf f})=0$.  Mais  $I_{\tilde{M}}^{\tilde{G},{\cal E}}({\bf M}',\boldsymbol{\delta}^{\bf M'},{\bf f})=I_{\tilde{M}}^{\tilde{G},{\cal E}}({\bf M}',\boldsymbol{\sigma}^{\bf M'},{\bf f})$ et la conclusion cherch\'ee s'ensuit. On est donc ramen\'e au cas

(8) $\hat{R}$ correspond \`a un espace de Levi $\tilde{R}$ de $\tilde{G}$ et le support de $\boldsymbol{\delta}$ est contenu dans l'ensemble $\tilde{R}'(F)^{out}$ ci-dessus.

  Partons de la formule (1) et  introduisons pour chaque $\tilde{s}$ des donn\'ees auxiliaires $G'(\tilde{s})_{1},...,\Delta(\tilde{s})_{1}$.   On dispose de l'isomorphisme de transition
$$C_{c,\lambda_{1}}^{\infty}(\tilde{M}'_{1}(F))\to C_{c,\lambda(\tilde{s})_{1}}^{\infty}(\tilde{M}'(s)_{1}),$$
qui se restreint en un isomorphisme
$$C_{c,\lambda_{1}}^{\infty}(\tilde{R}'_{1}(F))\to C_{c,\lambda(\tilde{s})_{1}}^{\infty}(\tilde{R}'(s)_{1}(F)).$$
Par dualit\'e, on a aussi un isomorphisme
$$D_{g\acute{e}om,\lambda_{1}}^{st}(\tilde{R}'_{1}(F))\to D_{g\acute{e}om,\lambda(\tilde{s})_{1}}^{st}(\tilde{R}'(s)_{1}(F)).$$
Par cet isomorphisme, $\boldsymbol{\delta}$ s'identifie \`a un \'el\'ement $\boldsymbol{\delta}(\tilde{s})_{1}$ de l'espace d'arriv\'ee. Alors $\boldsymbol{\delta}^{\bf M'}$ s'identifie \`a $(\boldsymbol{\delta}(\tilde{s})_{1})^{M'(\tilde{s})_{1}}$. Avec cette d\'efinition, la formule (2) reste valable. Soient $\tilde{s}$ et $\tilde{L}'_{\tilde{s}}$  intervenant dans cette formule. La d\'efinition de $\tilde{L}$ n'a plus de sens puisque $\tilde{R}$ n'existe plus. Mais on peut d\'efinir $\hat{L}\in {\cal L}(\hat{R})$ comme le commutant de $Z(\hat{L}'(\tilde{s}))^{\Gamma_{F},0}$ dans $\hat{G}$. C'est un Levi de $\hat{G}$ et il existe un sous-groupe parabolique $\hat{Q}\in {\cal P}(\hat{L})$ qui est  invariant par $\Gamma_{F}$ et  $\hat{\theta}$. Si $G$ \'etait quasi-d\'eploy\'e, il correspondrait \`a $\hat{L}$ un espace de Levi de $\tilde{G}$. Mais $G$ n'est pas suppos\'e quasi-d\'eploy\'e. On construit comme pr\'ec\'edemment le triplet ${\bf L}'(\tilde{s})=(L'(\tilde{s})=L'_{\tilde{s}},{\cal L}'(\tilde{s}),\tilde{s})$. Il est elliptique pour $\hat{L}$ au sens o\`u $Z(\hat{L}'_{\tilde{s}})^{\Gamma_{F},0}=Z(\hat{L})^{\Gamma_{F},\hat{\theta},0}$. On regroupe les termes de (2) selon le Levi $\hat{L}$ et on obtient une formule parall\`ele \`a (3):
$$I_{\tilde{M}}^{\tilde{G},{\cal E}}({\bf M}',\boldsymbol{\delta}^{\bf M'},{\bf f})=\sum_{\hat{L}}\sum_{\tilde{s}\in \tilde{\zeta} Z(\hat{M})^{\Gamma_{F},\hat{\theta}}/Z(\hat{G})^{\Gamma_{F},\hat{\theta}}; {\bf L}'(\tilde{s}) \text{ elliptique}}i_{\tilde{M}'}(\tilde{G},\tilde{G}'(\tilde{s}))$$
$$e_{\tilde{R}'}^{\tilde{G}'(\tilde{s})}(\tilde{M}',\tilde{L}'(\tilde{s}))S_{\tilde{R}'(\tilde{s})_{1},\lambda(\tilde{s})_{1}}^{\tilde{L}'(\tilde{s})_{1}}(\boldsymbol{\delta}(\tilde{s})_{1},B^{\tilde{G}},({\bf f}^{\tilde{G}'(\tilde{s})_{1}})_{\tilde{L}'(\tilde{s})_{1}}).$$
Ici $\hat{L}$ parcourt l'ensemble des \'el\'ements de ${\cal L}(\hat{R})$ qui v\'erifient la condition ci-dessus: il existe $\hat{Q}\in {\cal P}(\hat{L})$ qui est invariant par $\Gamma_{F}$ et $\hat{\theta}$.   Fixons $\hat{L}$. On va montrer que la sous-somme index\'ee par $\hat{L}$ dans l'expression ci-dessus est nulle.  Si elle est non nulle, il y a un $\tilde{s}$ pour lequel $({\bf f}^{\tilde{G}'(\tilde{s})_{1}})_{\tilde{L}'(\tilde{s})_{1}}$ est non nulle. Cela entra\^{\i}ne que $L'(\tilde{s})$ est relevant. A fortiori, $\hat{L}$ correspond \`a un espace de Levi de $\tilde{G}$, ou plus exactement \`a une classe de conjugaison de tels Levi. On peut donc fixer un espace de Levi $\tilde{L}$ de $\tilde{G}$ et supposer que $\hat{L}$ est le groupe dual de $L$. Alors ${\bf L}'(\tilde{s})$ est une donn\'ee endoscopique elliptique de $(L,\tilde{L},{\bf a})$ et on a l'\'egalit\'e
$({\bf f}^{\tilde{G}'(\tilde{s})_{1}})_{\tilde{L}'(\tilde{s})_{1}}=({\bf f}_{\tilde{L},\omega})^{\tilde{L}'(\tilde{s})_{1}}$.  On peut aussi imposer que pour un $\tilde{s}$, le produit des coefficients soit non nul. Cela impose que l'homomorphisme
$$(Z(\hat{M})^{\Gamma_{F},\hat{\theta},0}\times Z(\hat{L})^{\Gamma_{F},\hat{\theta},0})/diag_{-}(Z(\hat{G})^{\Gamma_{F},\hat{\theta},0})\to Z(\hat{R}')^{\Gamma_{F},0}$$ est surjectif et de noyau fini.
La condition d'ellipticit\'e impos\'ee \`a ${\bf L}'(\tilde{s})$ entra\^{\i}ne que l'espace ${\cal A}_{R'}^{L'(\tilde{s})}$ ne d\'epend pas de $\tilde{s}$. Le coefficient $d_{\tilde{R}'}^{\tilde{G}'(\tilde{s})}(\tilde{M}',\tilde{L}'(\tilde{s}))$ n'en d\'epend pas non plus. En notant $d$ sa valeur constante, on obtient une formule parall\`ele \`a (4)
$$\sum_{\tilde{s}\in \tilde{\zeta} Z(\hat{M})^{\Gamma_{F},\hat{\theta}}/Z(\hat{G})^{\Gamma_{F},\hat{\theta}}; {\bf L}'(\tilde{s}) \text{ elliptique}}dX(\tilde{s})S_{\tilde{R}'(\tilde{s})_{1},\lambda(\tilde{s})_{1}}^{\tilde{L}'(\tilde{s})_{1}}(\boldsymbol{\delta}(\tilde{s})_{1},B^{\tilde{L}},({\bf f}_{\tilde{L},\omega})^{\tilde{L}'(\tilde{s})_{1}}),$$
o\`u $X(\tilde{s})$ est comme pr\'ec\'edemment. Introduisons le groupe $Z(\hat{R})_{*}$ image r\'eciproque dans $Z(\hat{R})$ de $(Z(\hat{R})/(Z(\hat{R})\cap \hat{T}^{\hat{\theta},0}))^{\Gamma_{F}}$. L'ensemble
$$(Z(\hat{M})^{\Gamma_{F},\hat{\theta}}\cap (Z(\hat{L})^{\Gamma_{F}}(1-\hat{\theta})(Z(\hat{R})_{*})))/Z(\hat{G})^{\Gamma_{F},\hat{\theta}}$$
est un sous-groupe fini de $Z(\hat{M})^{\Gamma_{F},\hat{\theta}}/Z(\hat{G})^{\Gamma_{F},\hat{\theta}}$.    Il nous suffit de trouver un sous-groupe ${\cal Z}$ de ce groupe tel que, pour tout  $\tilde{s}_{0}\in \tilde{\zeta} Z(\hat{M})^{\Gamma_{F},\hat{\theta}}/Z(\hat{G})^{\Gamma_{F},\hat{\theta}}$ la sous-somme sur $\tilde{s}\in {\cal Z}\tilde{s}_{0}$ de l'expression ci-dessus soit  nulle. Fixons-donc un tel sous-groupe ${\cal Z}$ que nous pr\'eciserons plus tard. Pour prouver la nullit\'e ci-dessus, on ne perd pas grand'chose \`a supposer que $\tilde{s}_{0}=\tilde{\zeta}$, ce que nous ferons pour simplifier. On a

(9) pour $\tilde{s}\in {\cal Z}\tilde{\zeta}$, les donn\'ees endoscopiques ${\bf L}'(\tilde{s})$ et ${\bf L}'(\tilde{\zeta})$ sont \'equivalentes; si ces donn\'ees sont  elliptiques, on a $X(\tilde{s})=X(\tilde{\zeta})$. 

Preuve. Soit $z\in {\cal Z}$ (ou plus exactement un repr\'esentant dans $\hat{G}$, les \'el\'ements de ${\cal Z}$ \'etant des classes modulo $Z(\hat{G})^{\Gamma_{F},\hat{\theta}}$). Ecrivons $z=\tau(1-\hat{\theta})(\rho)$, avec $\tau\in Z(\hat{L})^{\Gamma_{F}}$ et $\rho\in Z(\hat{R})_{*}$.   L'automorphisme $ad_{\rho}$ conserve $\hat{L}$ puisque $\rho\in \hat{R}\subset \hat{L}$.  On a l'\'egalit\'e $ z\tilde{\zeta}=\rho \tau\tilde{\zeta}\rho^{-1}$. Donc $ad_{\rho}$ envoie  $\hat{L}'(\tau\tilde{\zeta}  )$ sur $\hat{L}'(z\tilde{\zeta} )$. Puisque $\tau\in Z(\hat{L})$, on a $\hat{L}'(\tau\tilde{\zeta})=\hat{L}'(\tilde{\zeta})$ donc $ad_{\rho}$ envoie $\hat{L}'(\tilde{\zeta})$ sur $\hat{L}'(z\tilde{\zeta} )$. Puisque $\rho\in Z(\hat{R})_{*}$, $ad_{\rho}$ conserve ${\cal R}'$. Donc $ad_{\rho}$ envoie ${\cal L}'(\tilde{\zeta})$ sur ${\cal L}'(z\tilde{\zeta })$. Autrement dit $\rho$ d\'efinit une \'equivalence entre les donn\'ees ${\bf L}'(\tilde{\zeta})$ et ${\bf L}'(z\tilde{\zeta} )$.

Les calculs conduisant \`a l'\'egalit\'e (6) restent valables: ils se placent enti\`erement dans les groupes duaux et dans ces groupes, la situation n'a pas chang\'e. Cette \'egalit\'e  montre que $X(\tilde{s})$ ne d\'epend que de la classe d'\'equivalence de la donn\'ee ${\bf L}'(\tilde{s})$. D'o\`u l'assertion (9).

Si ${\bf L}'(\tilde{\zeta})$ n'est pas elliptique, la sous-somme sur $\tilde{s}\in {\cal Z}\tilde{\zeta}$ est nulle. De m\^eme, si ${\bf L}'(\tilde{\zeta})$ n'est pas relevant, les fonctions $({\bf f}_{\tilde{L},\omega})^{\tilde{L}'(\tilde{s})_{1}}$ sont nulles. Supposons ${\bf L}'(\tilde{\zeta})$ elliptique et relevant. Gr\^ace \`a (9), l'assertion \`a prouver se r\'eduit \`a
$$(10) \qquad \sum_{z\in {\cal Z}}S_{\tilde{R}'(z\tilde{\zeta} )_{1},\lambda(z\tilde{\zeta})_{1}}^{\tilde{L}'(z\tilde{\zeta} )_{1}}(\boldsymbol{\delta}(z\tilde{\zeta} )_{1},B^{\tilde{L}},({\bf f}_{\tilde{L},\omega})^{\tilde{L}'(z\tilde{\zeta} )_{1}})=0.$$

Fixons $z\in {\cal Z}$. On va calculer  $S_{\tilde{R}'(z\tilde{\zeta} )_{1},\lambda(z\tilde{\zeta}_{1}}^{\tilde{L}'(z\tilde{\zeta} )_{1}}(\boldsymbol{\delta}(z\tilde{\zeta} )_{1},B^{\tilde{L}},({\bf f}_{\tilde{L},\omega})^{\tilde{L}'(z\tilde{\zeta} )_{1}})$. Pour cela, on a besoin de fixer une d\'ecomposition $z=\tau(1-\hat{\theta})(\rho)$ comme dans la preuve de (9). On a deux donn\'ees auxiliaires pour ${\bf M}'$: les donn\'ees $M'(\tilde{\zeta})_{1}$,... et les donn\'ees $M'(z\tilde{\zeta} )_{1}$,... D'o\`u une fonction de recollement $\tilde{\lambda}(z)^M$ d\'efinie sur le produit fibr\'e de $\tilde{M}'(\tilde{\zeta})_{1}(F)$ et $\tilde{M}'(z\tilde{\zeta} )_{1}(F)$ au-dessus de $\tilde{M}'(F)$. On note un tel produit fibr\'e $\tilde{M}'(\tilde{\zeta})_{1}(F)\times_{\tilde{M}'(F)}\tilde{M}'(z\tilde{\zeta })_{1}(F)$. Par restriction, cette fonction d\'efinit un isomorphisme
$$\iota(z)^M:C_{c,\lambda(\tilde{\zeta})_{1}}^{\infty}(\tilde{R}'(\tilde{\zeta})_{1}(F))\simeq C_{c,\lambda(z\tilde{\zeta})_{1}}^{\infty}(\tilde{R}'(z\tilde{\zeta} )_{1}(F)).$$
On en d\'eduit un isomorphisme dual
$$\iota(z)^{M,*}:D_{g\acute{e}om,\lambda(z\tilde{\zeta})_{1}}^{st}(\tilde{R}'(z\tilde{\zeta} )_{1}(F))\simeq D_{g\acute{e}om,\lambda(\tilde{\zeta})_{1}}^{st}(\tilde{R}'(\tilde{\zeta} )_{1}(F)).$$
Par construction, on a $\boldsymbol{\delta}(\tilde{\zeta})_{1}=\iota(z)^{M,*}(\boldsymbol{\delta}(z\tilde{\zeta} )_{1})$. 

L'action galoisienne sur $\hat{L}'(\tilde{\zeta})$ est h\'erit\'ee de celle sur $\hat{G}'(\tilde{\zeta})$. L'intersection de $(\hat{B},\hat{T})$ avec $\hat{L}'(\tilde{\zeta})$ est une paire de Borel de ce groupe,  invariante pour cette action et pour laquelle $\hat{R}'$ est standard.  On peut compl\'eter cette paire en une paire de Borel \'epingl\'ee invariante par $\Gamma_{F}$. On effectue  la m\^eme construction pour $\hat{L}'(z\tilde{\zeta })$. Il est loisible de supposer  que la restriction de l'\'epinglage de ce groupe \`a $\hat{R}'$ co\"{\i}ncide avec celle de l'\'epinglage de $\hat{L}'(\tilde{\zeta})$. Ecrivons $z=\tau(1-\hat{\theta})(\rho)$ comme dans la preuve de (9).  L'automorphisme $ad_{\rho}$ transporte la paire de Borel de $\hat{L}'(\tilde{\zeta})$ sur celle de $\hat{L}'(z\tilde{\zeta} )$. Quitte \`a multiplier $\rho$ par un \'el\'ement de $Z(\hat{R})\cap \hat{T}^{\hat{\theta},0}$,   ce qui ne change pas $(1-\hat{\theta})(\rho)$, on peut supposer que $ad_{\rho}$ transporte aussi les \'epinglages. Alors $ad_{\rho}$ est \'equivariant pour les actions galoisiennes. On peut identifier les groupes $L'(\tilde{\zeta})$ et $L'(z\tilde{\zeta })$, ainsi que les espaces $\tilde{L}'(\tilde{\zeta})$ et $\tilde{L}'(z\tilde{\zeta} )$. Comme en [I] 2.6, les donn\'ees auxiliaires $L'(z\tilde{\zeta })_{1}$,...  pour ${\bf L}'(z\tilde{\zeta} )$ se transportent en des donn\'ees auxiliaires pour ${\bf L}'(\tilde{\zeta})$. C'est-\`a-dire que, via les isomorphismes pr\'ec\'edents, on consid\`ere $L'(z\tilde{\zeta} )_{1}$ comme une extension de $L'(\tilde{\zeta})$ et $\tilde{L}'(z\tilde{\zeta})_{1}$ comme un espace au-dessus de $\tilde{L}'(\tilde{\zeta})$. On compl\`ete ces donn\'ees par le plongement
$${\cal L}'(\tilde{\zeta})\stackrel{ad_{\rho}}{\to}{\cal L}'(z\tilde{\zeta})\stackrel{\hat{\xi}(z\tilde{\zeta })_{1}}{\to}{^LL}'(z\tilde{\zeta })_{1}$$
et par le facteur de transfert $\Delta(z\tilde{\zeta })_{1}$. Il y a une fonction de recollement $\tilde{\lambda}(z,\rho)^L$ d\'efinie sur $\tilde{L}'(\tilde{\zeta} )_{1}(F)\times_{\tilde{L}'(\tilde{\zeta})(F)}\tilde{L}'(z\tilde{\zeta })_{1}(F)$ qui fait passer  des donn\'ees choisies pour ${\bf L}'(\tilde{\zeta})$ \`a ces nouvelles donn\'ees. Comme ci-dessus, il s'en d\'eduit un isomorphisme
$$\iota(z,\rho)^{L,*}:D_{g\acute{e}om,\lambda(z\tilde{\zeta})_{1}}^{st}(\tilde{R}'(z\tilde{\zeta} )_{1}(F))\simeq D_{g\acute{e}om,\lambda(\tilde{\zeta})_{1}}^{st}(\tilde{R}'(\tilde{\zeta} )_{1}(F)).$$
Le transfert comute au recollement  donc celui-ci envoie $({\bf f}_{\tilde{L},\omega})^{\tilde{L}'(z\tilde{\zeta} )_{1}}$ sur $({\bf f}_{\tilde{L},\omega})^{\tilde{L}'(\tilde{\zeta} )_{1}}$. Il envoie $\boldsymbol{\delta}(z\tilde{\zeta} )_{1}$ sur $\iota(z,\rho)^{L,*}(\boldsymbol{\delta}(z\tilde{\zeta })_{1})=\iota(z,\rho)^{L,*}\circ(\iota(z)^{M,*})^{-1}(\boldsymbol{\delta}(\tilde{\zeta})_{1})$. On a donc
$$(11) \qquad S_{\tilde{R}'(z\tilde{\zeta} )_{1},\lambda(z\tilde{\zeta})_{1}}^{\tilde{L}'(z\tilde{\zeta} )_{1}}(\boldsymbol{\delta}(z\tilde{\zeta} )_{1},B^{\tilde{L}},({\bf f}_{\tilde{L},\omega})^{\tilde{L}'(z\tilde{\zeta} )_{1}})=$$
$$S_{\tilde{R}'(\tilde{\zeta} )_{1},\lambda(\tilde{\zeta})_{1}}^{\tilde{L}'(\tilde{\zeta} )_{1}}(\iota(z,\rho)^{L,*}\circ(\iota(z)^{M,*})^{-1}(\boldsymbol{\delta}(\tilde{\zeta} )_{1}),B^{\tilde{L}},({\bf f}_{\tilde{L},\omega})^{\tilde{L}'(\tilde{\zeta} )_{1}}).$$
 Le compos\'e $(\iota(z)^M)^{-1}\circ \iota(z,\rho)^L$ est un automorphisme de $C_{c,\lambda(\tilde{\zeta})_{1}}^{\infty}(\tilde{R}'(\tilde{\zeta})_{1}(F))$. Il est de la forme $\varphi\mapsto \tilde{\lambda}_{z,\rho}\varphi$, o\`u $\tilde{\lambda}_{z,\rho}$ est une certaine fonction continue sur $\tilde{R}'(\tilde{\zeta})_{1}(F)$ que nous allons calculer. On simplifie les notations en supprimant autant que possible $\tilde{\zeta}$ et $z$ des notations. On pose $\tilde{\lambda}=\tilde{\lambda}_{z,\rho}$. On conserve  les indices $1$ pour  les termes associ\'es aux donn\'ees en $\tilde{\zeta}$ et on les convertit en indices $2$ pour ceux associ\'es aux donn\'ees en $z\tilde{\zeta }$. Par exemple, on note $\Delta_{1}$ et $\Delta_{2}$ les termes not\'es pr\'ecedemment $\Delta(\tilde{\zeta})_{1}$ et $\Delta(z\tilde{\zeta })_{1}$. Soit $r'_{1}\in \tilde{R}'_{1}(F)$. Notons $r'$ sa projection dans $\tilde{R}'(F)$ et choisissons $r'_{2}\in \tilde{R}'_{2}(F)$ se projetant sur $r'$. Par d\'efinition
 $$\tilde{\lambda}(r'_{1})=\tilde{\lambda}^M(r'_{1},r'_{2})^{-1}\tilde{\lambda}^L(r'_{1},r'_{2}).$$
 Fixons $l\in \tilde{L}(F)$ semi-simple et assez r\'egulier et $l'\in \tilde{L}'(F)$ de sorte que leurs classes de conjugaison stable se correspondent. C'est possible puisqu'on a suppos\'e $L'$ relevant. Fixons des \'el\'ements $l'_{1}\in \tilde{L}'_{1}(F)$ et $l'_{2}\in \tilde{L}'_{2}(F)$ se projetant sur $l'$. Notons $(a_{1},a_{2})$ l'\'el\'ement de $L'_{1}(F)\times_{ L'(F)}L'_{2}(F)$ tel que $(r'_{1},r'_{2})=(a_{1}l'_{1},a_{2}l'_{2})$. On a l'\'egalit\'e
 $$\tilde{\lambda}^L(r'_{1},r'_{2})=\lambda^L(a_{1},a_{2})\tilde{\lambda}^L(l'_{1},l'_{2})=\lambda^L(a_{1},a_{2})\Delta_{2}(l'_{2},l)\Delta_{1}(l'_{1},l)^{-1}.$$
 On introduit  $m\in \tilde{M}(F)$, $m'_{1}\in \tilde{M}'_{1}(F)$ et $m'_{2}\in \tilde{M}'_{2}(F)$ v\'erifiant des conditions analogues et $(b_{1},b_{2})\in M'_{1}(F)\times_{M'(F)}M'_{2}(F)$ tel que $(r'_{1},r'_{2})=(b_{1}m'_{1},b_{2}m'_{2})$. On a une relation analogue \`a celle ci-dessus. On se rappelle que $\Delta_{i}(m'_{i},m)\Delta_{i}(l'_{i},l)^{-1}=\boldsymbol{\Delta}_{i}(m'_{i},m';l'_{i},l')$ pour $i=1,2$. On obtient
 $$(12) \qquad \tilde{\lambda}(r'_{1})=\lambda^M(b_{1},b_{2})^{-1}\lambda^L(a_{1},a_{2})\boldsymbol{\Delta}_{1}(m'_{1},m;l'_{1},l)\boldsymbol{\Delta}_{2}(m'_{2},m;l'_{2},l)^{-1}.$$
 
 On  calcule les facteurs de transfert ci-dessus en utilisant les d\'efinitions de [I] 2.2. Pour rendre les calculs plus clairs, on modifie les notations de cette r\'ef\'erence: on y avait deux s\'eries d'objets, la deuxi\`eme \'etant soulign\'ee ($T$,  $\underline{T}$, etc...); on affecte maintenant la premi\`ere s\'erie d'un exposant $M$ et la deuxi\`eme s\'erie d'un exposant $L$; d'autre part, on conserve les indices $1$ pour le facteur $\boldsymbol{\Delta}_{1}$ et on les transforme en indices $2$ pour le facteur $\boldsymbol{\Delta}_{2}$. On fixe des diagrammes $(m',B^{ M'},T^{M'},B^M,T^M,m)$ et $(l',B^{L'},T^{L'},B^L,T^L,l)$   et on utilise ces diagrammes pour calculer les deux facteurs. De m\^eme, on utilise les m\^emes $a$-data et $\chi$-data. On suppose que ces $\chi$-data sont triviales sur les orbites galoisiennes asym\'etriques. Montrons que les facteurs $\Delta_{II}$ sont les m\^emes pour les deux facteurs. Consid\'erons par exemple les facteurs $\Delta_{II,1}(l',l)$ et $\Delta_{II,2}(l',l)$. Ce sont des produits sur les orbites pour l'action de $\Gamma_{F}$ dans $\Sigma(T^L)_{res,ind}$. Le Levi $\hat{L}$ d\'etermine un sous-ensemble $\Sigma^L(T^L)_{res,ind}\subset \Sigma(T^L)_{res,ind}$. Les contributions de ce sous-ensemble aux deux facteurs sont les m\^emes par d\'efinition. Pour $\alpha_{res}\in \Sigma(T^L)_{res,ind}-\Sigma^L(T^L)_{res,ind}$, l'orbite galoisienne de $\alpha_{res}$ est asym\'etrique  (un sous-groupe parabolique $\hat{P}\in {\cal P}(\hat{L})$ invariant par $\Gamma_{F}$ d\'etermine une  partition de   $ \Sigma(T^L)_{res,ind}-\Sigma^L(T^L)_{res,ind}$ en deux ensembles oppos\'es  invariants par $\Gamma_{F}$). Puisqu'on a suppos\'e les $\chi$-data triviales sur les orbites asym\'etriques, la contribution de $\Sigma(T^L)_{res,ind}-\Sigma^L(T^L)_{res,ind}$ est \'egale \`a $1$. Cela d\'emontre l'assertion. Donc
 $$(13) \qquad \boldsymbol{\Delta}_{1}(m'_{1},m;l'_{1},l)\boldsymbol{\Delta}_{2}(m'_{2},m;l'_{2},l)^{-1}=\Delta_{1,imp}(m'_{1},m;l'_{1},l)\Delta_{2,imp}(m'_{2},m;l'_{2},l)^{-1}.$$
 Le tore $U=(T^M_{sc}\times T^L_{sc})/diag_{-}(Z(G_{SC}))$ qui intervient dans les d\'efinitions est le m\^eme pour les deux facteurs. Le cocycle $V$ \`a valeurs dans ce tore est aussi le m\^eme. Les termes \`a valeurs dans $\hat{U}$ sont $\boldsymbol{\zeta}=(\zeta_{sc},\zeta_{sc})$ pour le facteur $\Delta_{1,imp}$ et ${\bf z}\boldsymbol{\zeta}=(z_{sc}\zeta_{sc},z_{sc}\zeta_{sc})$ pour le facteur $\Delta_{2,imp}$, o\`u on a \'ecrit $\tilde{\zeta}=\zeta\hat{\theta}$. Les tores $S_{1}$ et $S_{2}$ sont diff\'erents. Introduisons le tore $\mathfrak{ T}_{12}^M$ produit fibr\'e de $T_{1}^M$, $T_{2}^M$ et $T^M$ au-dessus de $T^{\bf M'}$   et le tore analogue $\mathfrak{ T}_{12}^L$. Notons $\mathfrak{Z}_{12}$ le produit fibr\'e de $\mathfrak{Z}_{1}$ et $\mathfrak{Z}_{2}$ au-dessus de $Z(G)$. Posons $S_{12}=(\mathfrak{T}_{12}^M\times\mathfrak{T}_{12}^L)/diag_{-}(\mathfrak{Z}_{12})$. Il y a des homomorphismes naturels d'oubli d'une s\'erie de variables
 $$\begin{array}{ccccc}&&S_{12}&&\\ &\,\,\swarrow p_{1}&&\,\,\searrow p_{2}&\\ S_{1}&&&&S_{2}\\ \end{array}$$
 Posons $\nu_{12}^M=(\mu_{1}^M,\mu_{2}^M,\nu^M)\in \mathfrak{T}_{12}^M$, d\'efinissons de m\^eme $\nu_{12}^L$, notons $\boldsymbol{\nu}_{12}$ l'image de $(\nu_{12}^M,(\nu_{12}^L)^{-1})$ dans $S_{12}$. Alors $p_{i}(\boldsymbol{\nu}_{12})=\boldsymbol{\nu}_{i}$ pour $i=1,2$. De plus, le couple $(V,\boldsymbol{\nu}_{12})$ d\'efinit un \'el\'ement de $H^1(\Gamma_{F};U\stackrel{1-\theta}{\to }S_{12})$. Pour $i=1,2$, on a donc $(V,\boldsymbol{\nu}_{i})=p_{i}(V,\boldsymbol{\nu}_{12})$ en notant encore selon notre habitude $p_{i}:H^1(\Gamma_{F};U\stackrel{1-\theta}{\to }S_{12})\to H^1(\Gamma_{F};U\stackrel{1-\theta}{\to }S_{i})$ l'homomorphisme d\'eduit fonctoriellement du $p_{i}$ pr\'ec\'edent. Par une propri\'et\'e de compatibilit\'e, on obtient
 $$\Delta_{i,imp}(m'_{i},m;l'_{i},l)=\left\lbrace\begin{array}{cc}<(V,\boldsymbol{\nu}_{12}),(\hat{p}_{1}(\hat{V}_{1}),\boldsymbol{\zeta})>^{-1},&\text{ si }i=1,\\ <(V,\boldsymbol{\nu}_{12}),(\hat{p}_{2}(\hat{V}_{2}),{\bf z}\boldsymbol{\zeta})>^{-1},& \text{ si }i=2,\\ \end{array}\right.$$
 o\`u $\hat{p}_{i}:H^1(W_{F};\hat{S}_{i}\stackrel{1-\hat{\theta}}{\to}\hat{U})\to H^1(W_{F};\hat{S}_{12}\stackrel{1-\hat{\theta}}{\to}\hat{U})$ est dual de $p_{i}$. Donc
$$(14)\qquad \Delta_{1,imp}(m'_{1},m;l'_{1},l)\Delta_{2,imp}(m'_{2},m;l'_{2},l)^{-1}=<(V,\boldsymbol{\nu}_{12}),(\hat{V}_{12},{\bf z})>,$$
 o\`u  $\hat{V}_{12}=\hat{p}_{1}(\hat{V}_{1})^{-1}\hat{p}_{2}(\hat{V}_{2})$.  Le tore dual $\hat{\mathfrak{T}}_{12}^M$ est le quotient de $\hat{T}_{1}^M\times \hat{T}_{2}^M\times \hat{T}^M$ par le groupe $\{(t_{1},t_{2},t)\in \hat{T}^{ M'}; t_{1}t_{2}t=1\}$ plong\'e par $(t_{1},t_{2},t)\mapsto (\hat{\xi}_{1}(t_{1}),\hat{\xi}_{2}(t_{2}),t)$. Le tore $\hat{S}_{12}$ est le sous-tore de $\hat{\mathfrak{T}}_{12}^M\times\hat{\mathfrak{T}}_{12}^L\times \hat{T}_{sc}$ form\'e des $(t^M,t^L,t_{sc})$ tels que $j(t_{sc})=t^M(t^L)^{-1}$ (les notations sont adapt\'ees de [I] 2.2;  il faut aussi mettre sur ce tore une action galoisienne d\'efinie comme dans ce paragraphe). Notons $\hat{V}_{12}^M$, $\hat{V}_{12}^L$ et $\hat{V}_{12,sc}$ les trois composantes de $\hat{V}_{12}$.
 
  Calculons $\hat{V}_{12}^L(w)$ pour $w\in W_{F}$. Il convient de choisir $(g(w),w)\in {\cal G}'(\tilde{\zeta})$ tel que $ad_{g(w)}\circ w_{G}=w_{G'(\tilde{\zeta})}$ et $(g_{z}(w),w)\in {\cal G}'(z\tilde{\zeta} )$ tel que $ad_{g_{z}(w)}\circ w_{G}=w_{G'(z\tilde{\zeta} )}$. Puisque ${\cal G}'(\tilde{\zeta})=\hat{G}'(\tilde{\zeta}){\cal R}'$, on peut certainement supposer que $(g(w),w)$ appartient \`a ${\cal R}'$. On pose plut\^ot $g(w)=r(w)\in \hat{R}$. Puisque $\hat{ G}'(\tilde{\zeta})$ et $\hat{ G}'(z\tilde{\zeta} )$ ont en commun le Levi standard $\hat{M}'$, on peut supposer que leurs actions galoisiennes co\"{\i}ncident sur ce Levi. Donc $g_{z}(w)=m'(w)r(w)$, avec $m'(w)\in Z(\hat{M}')$. Puisque ce groupe est produit de $Z(\hat{M}')^0$ et de $Z(\hat{G}'(z\tilde{\zeta} ))$ et que l'on peut modifier $g_{z}(w)$ par un \'el\'ement de ce dernier groupe, on peut m\^eme supposer $m'(w)\in Z(\hat{M}')^0$. On pose $\hat{\xi}_{1}(r(w),w)=(\zeta_{1}(w),w)$, $\hat{\xi}_{2}(m'(w)r(w),w)=(\zeta_{2}(w),w)$. On a $\hat{V}_{12}^L(w)=(\zeta_{1}(w)^{-1},\zeta_{2}(w),t_{T^L,1}(w)^{-1}t_{T^L,2}(w))$ (ici encore, on affecte d'un indice $1$, resp. $2$, le terme provenant des donn\'ees en $\tilde{\zeta}$, resp. $z\tilde{\zeta} $). Dans les termes $t_{T^L,i}(w)$ interviennent des termes $\hat{n}(\omega_{T^L}(w))$ et $\hat{r}_{T^L}(w)$. Ce sont les m\^emes pour $i=1,2$ et ils disparaissent dans le quotient  ci-dessus. De m\^eme, les termes $r(w)$ disparaissent. Donc
 $$(15)\qquad t_{T^L,1}(w)^{-1}t_{T^L,2}(w)= \hat{r}_{T^L,G'(\tilde{\zeta})}(w)\hat{n}_{G'(\tilde{\zeta})}(\omega_{T^L,G'(\tilde{\zeta})}(w))m'(w)^{-1}$$
 $$\hat{n}_{G'(z\tilde{\zeta} )}(\omega_{T^L,G'(z\tilde{\zeta} )}(w))^{-1}\hat{r}_{T^L,G'(z\tilde{\zeta} )}(w)^{-1}.$$
 
 Consid\'erons les termes $\hat{r}_{T^L,G'(\tilde{\zeta})}(w)$ et $\hat{n}_{G'(\tilde{\zeta})}(\omega_{T^L,G'(\tilde{\zeta})}(w))$. En rempla\c{c}ant dans leurs d\'efinitions $\hat{G}'(\tilde{\zeta})$ par $\hat{L}'(\tilde{\zeta})$ et en utilisant la paire de Borel \'epingl\'ee de ce groupe que l'on a fix\'ee plus haut, on  obtient des termes $\hat{r}_{T^L,L'(\tilde{\zeta})}(w)$ et $\hat{n}_{L'(\tilde{\zeta})}(\omega_{T^L,L'(\tilde{\zeta})}(w))$. Montrons que:
 
 (16) il existe $x\in \hat{T}^{L'}_{sc}$ tel que 
 $$ \hat{r}_{T^L,G'(\tilde{\zeta})}(w)\hat{n}_{G'(\tilde{\zeta})}(\omega_{T^L,G'(\tilde{\zeta})}(w))=xw_{T^L}(x)^{-1} \hat{r}_{T^L,L'(\tilde{\zeta})}(w)\hat{n}_{L'(\tilde{\zeta})}(\omega_{T^L,L'(\tilde{\zeta})}(w))$$
 pour tout $w\in W_{F}$.
 
On supprime les $\tilde{\zeta}$ le temps de cette preuve. Si $\hat{L}'$ \'etait standard dans $\hat{G}'$, les termes pour $\hat{G}'$ et $\hat{L}'$ seraient \'egaux et l'assertion serait claire. En tout cas, on a l'\'egalit\'e $\omega_{T^L,G'}(w)=\omega_{T^L,L'}(w)$ par d\'efinition. Notons simplement $\omega(w)=\omega_{T^L,G'}(w)$  et $\hat{n}=\hat{n}_{G'}$. La section de Springer $\hat{n}$ est relative \`a la paire de Borel \'epingl\'ee fix\'ee  $\hat{{\cal E}}'$ de $\hat{G}'$. Notons $(\hat{ B}',\hat{T}')$ la paire de Borel sous-jacente \`a $\hat{{\cal E}}'$. On rappelle que $\hat{B}'=\hat{B}\cap \hat{G}'$ et $\hat{T}'=\hat{T}^{\hat{\theta},0}$. Fixons un sous-groupe parabolique $\hat{P}'\in {\cal P}^{\hat{G}'}(\hat{L}')$ invariant par $\Gamma_{F}$. Notons $\hat{B}'_{\sharp}$ le sous-groupe de Borel de $\hat{G}'$ contenu dans $\hat{P}'$ et dont l'intersection avec $\hat{L}'$ co\"{\i}ncide avec  celle de $\hat{B}'$. On peut compl\'eter la paire $(\hat{B}'_{\sharp},\hat{T}')$ en une paire de Borel \'epingl\'ee $\hat{{\cal E}}'_{\sharp}$ conserv\'ee par $\Gamma_{F}$, de sorte que sa restriction \`a $\hat{L}'$  soit la paire de Borel \'epingl\'ee de $\hat{L}'$. De cette paire se d\'eduit une autre section de Springer pour $\hat{G}'$ dont $\hat{n}_{L'}$ est la restriction puisque $\hat{L}'$ est standard pour cette paire. Soit $g_{ad}\in \hat{G}'_{AD}$ tel que $ad_{g_{ad}}$ envoie $\hat{{\cal E}}'$ sur $\hat{{\cal E}}'_{\sharp}$. C'est un \'el\'ement de $\hat{G}_{AD}^{_{'}\Gamma_{F}}$. D'apr\`es [K] lemme 1.6, on peut le relever en un \'el\'ement $g\in \hat{G}_{SC}^{_{'}\Gamma_{F}}$. Cet \'el\'ement normalise $\hat{T}'$ donc d\'efinit un \'el\'ement $u\in W^{G'}$. Par transport de structure, on a $\hat{n}_{\hat{L}'}(\omega(w))=g\hat{n}(u^{-1}\omega(w)u)g^{-1}$. Pour $u_{1},u_{2}\in W^{G'}$, on a l'\'egalit\'e $\hat{n}(u_{1}u_{2})=t(u_{1},u_{2})\hat{n}(u_{1})\hat{n}(u_{2})$ (cf. [LS] lemme 2.1.A) o\`u
 $$t(u_{1},u_{2})=\prod_{\alpha>0, u_{1}^{-1}(\alpha)<0,u_{2}^{-1}u_{1}^{-1}(\alpha)>0}\check{\alpha}(-1).$$
 Ici les $\alpha$ parcourent les racines de $\hat{T}'$ dans $\hat{G}'$ et la positivit\'e est relative \`a $\hat{B}'$. On calcule
 $$\hat{n}(u^{-1}\omega(w)u)=t(u^{-1},\omega(w)u)\hat{n}(u^{-1})t(\omega(w),u)\hat{n}(\omega(w))\hat{n}(u),$$
 $$\hat{n}(u^{-1})=t(u^{-1},u)\hat{n}(u)^{-1},$$
 d'o\`u
 $$\hat{n}(u^{-1}\omega(w)u)=t(u^{-1},\omega(w)u)t(u^{-1},u)u^{-1}(t(\omega(w),u))\hat{n}(u)^{-1}\hat{n}(\omega(w))\hat{n}(u).$$
 On calcule
 $$t(u^{-1},\omega(w)u)t(u^{-1},u)=\prod_{\alpha>0, u(\alpha)<0,u^{-1}\omega(w)^{-1}u(\alpha)<0}\check{\alpha}(-1).$$
 Mais ce produit est vide. En effet $\alpha>0$ \'equivaut \`a $u(\alpha)>_{\sharp}0$ o\`u $>_{\sharp}$ est l'ordre d\'efini par $\hat{B}'_{\sharp}$. Par d\'efinition de ce Borel, les conditions $u(\alpha)<0$ et $u(\alpha)>_{\sharp}0$ interdisent \`a $u(\alpha)$ d'\^etre dans $\hat{L}'$. Pour les racines hors de ce Levi, $\omega(w)$ conserve l'ordre $>_{\sharp}$. Donc $\omega(w)^{-1}u(\alpha)>_{\sharp}0$ puis $u^{-1}\omega(w)^{-1}u(\alpha)>0$. Cela d\'emontre l'assertion. Donc $t(u^{-1},\omega(w)u)t(u^{-1},u)=1$. On en d\'eduit
 $$\hat{n}_{\hat{L}'}(\omega(w))= t(\omega(w),u)x^{-1}\hat{n}(\omega(w))x,$$
 o\`u $x=g\hat{n}(u)^{-1}$. Cet \'el\'ement $x$ appartient \`a $\hat{T}'_{sc}\simeq \hat{T}^{L'}_{sc}$. Puisque $g$ est fixe par $\Gamma_{F}$, $u$ l'est aussi. Puisque la section de Springer est \'equivariante par $\Gamma_{F}$, $\hat{n}(u) $ est fixe par $\Gamma_{F}$, donc $x$ aussi. Alors $\hat{n}(\omega(w))x=\omega(x)\circ w_{G'}(x) \hat{n}(\omega(w))=w_{T^L}(x)\hat{n}(\omega(w))$. D'o\`u
 $$\hat{n}_{\hat{L}'}(\omega(w))= t(\omega(w),u)x^{-1}w_{T^L}(x)\hat{n}(\omega(w)).$$
 Les termes $ \hat{r}_{T^L,G'}(w)$ et $\hat{r}_{T^L,L'}(w)$ sont des produits sur les racines $\alpha$. En fait, $\hat{r}_{T^L,L'}(w)$ est exactement \'egal \`a la contribution \`a $ \hat{r}_{T^L,G'}(w)$ des racines dans $\hat{L}'$. Comme on l'a remarqu\'e, les racines hors de ce Levi appartiennent \`a des orbites asym\'etriques pour lesquelles les $\chi$-data sont triviales. Pour un couple $({\cal O},-{\cal O})$ de telles orbites, on peut supposer que ${\cal O}$ est form\'e d'\'el\'ements positifs pour l'ordre $>_{\sharp}$. On voit alors que
 $$\hat{r}_{T^L,G'}(w)\hat{r}_{T^L,L'}(w)^{-1}=\prod_{\alpha\text{ hors de }\hat{L}',\alpha>0, \alpha>_{\sharp}0, w_{T^L}^{-1}(\alpha)<0}\check{\alpha}(-1).$$
 Puisque $w_{T^L}=\omega(w)w_{G'}$ et que l'action $w_{G'}$ pr\'eserve $\hat{B}'$, la condition $w_{T^L}^{-1}(\alpha)<0$ \'equivaut \`a $\omega(w)^{-1}(\alpha)<0$. Pour $\alpha$ hors de $\hat{L}'$, $\omega(w)$ pr\'eserve l'ordre $>_{\sharp}$. La condition $\alpha>_{\sharp}0$ \'equivaut \`a $\omega(w)^{-1}(\alpha)>_{\sharp}0$. Le produit ci-dessus est donc sur l'ensemble des $\alpha$ hors de $\hat{L}'$ tels que $\alpha>0$, $\omega(w)^{-1}(\alpha)<0$ et 
$\omega(w)^{-1}(\alpha)>_{\sharp}0$. Mais alors, la condition $\alpha$ hors de $\hat{L}'$ devient superflue car les deux derni\`eres relations interdisent \`a $\alpha$ d'\^etre dans ce Levi. En rempla\c{c}ant la condition $\omega(w)^{-1}(\alpha)>_{\sharp}0$ par la condition \'equivalente $u^{-1}\omega(w)^{-1}(\alpha)>0$, on obtient l'\'egalit\'e
    $$\hat{r}_{T^L,G'}(w)\hat{r}_{T^L,L'}(w)^{-1}=t(\omega(w),u).$$
    En rassemblant ces calculs, on obtient (16).
    
    {\bf Remarque.} On n'a pas pris soin du sens des indices $sc$. L'\'el\'ement $x$ construit appartient \`a $\hat{G}'_{SC}$. On n'aura besoin que de son image naturelle dans $\hat{G}_{SC}$.

On effectue les m\^emes constructions  pour les donn\'ees relatives \`a $z\tilde{\zeta} $. Cela nous fournit deux \'el\'ements $x_{1}$ et $x_{2}$.  On a \'evidemment l'\'egalit\'e $\hat{r}_{T^L,L'(\tilde{\zeta})}(w)=\hat{r}_{T^L,L'(z\tilde{\zeta} )}(w)$.   On a ajust\'e $\rho$ de sorte que $ad_{\rho}$ transporte la paire de Borel \'epingl\'ee de $\hat{L}'(\tilde{\zeta})$ sur celle de $\hat{L}'(z\tilde{\zeta })$. Il en r\'esulte que
$$\hat{n}_{L'(z\tilde{\zeta} )}(\omega_{T^L,L'(z\tilde{\zeta })}(w)) =\rho \hat{n}_{L'(\tilde{\zeta}) }(\omega_{T^L,L'(\tilde{\zeta})}(w))\rho^{-1}.$$
A ce point, la formule (15) se r\'ecrit
$$t_{T^L,1}(w)^{-1}t_{T^L,2}(w)=xw_{T^L}(x)^{-1}\hat{n}_{L'(\tilde{\zeta})}(\omega_{T^L,L'(\tilde{\zeta})}(w))m'(w)^{-1}\rho\hat{n}_{L'(\tilde{\zeta})}(\omega_{T^L,L'(\tilde{\zeta})}(w))^{-1}\rho^{-1},$$
o\`u $x=x_{1}x_{2}^{-1}$. Utilisons encore que $ad_{\rho}:\hat{L}'(\tilde{\zeta})\to \hat{L}'(z\tilde{\zeta} )$ est \'equivariante pour les actions galoisiennes. On a $w_{L'(\tilde{\zeta})}=ad_{r(w)}\circ w_{G}$ et $w_{L'(z\tilde{\zeta} )}=ad_{m'(w)r(w)}\circ w_{G}$. De plus $\rho\in Z(\hat{R})$ donc commute \`a $r(w)$. Alors l'\'egalit\'e $w_{L'(z\tilde{\zeta} )}\circ ad_{\rho}=ad_{\rho}\circ w_{L'(\tilde{\zeta})}$ entra\^{\i}ne que $\rho^{-1}m'(w)w_{G}(\rho)\in Z(\hat{L}'(\tilde{\zeta}))$. Donc cet \'el\'ement commute \`a $\hat{n}_{L'(\tilde{\zeta})}(\omega_{T^L,L'(\tilde{\zeta})}(w))$. Alors
$$t_{T^L,1}(w)^{-1}t_{T^L,2}(w)=xw_{T^L}(x)^{-1}\hat{n}_{L'(\tilde{\zeta})}(\omega_{T^L,L'(\tilde{\zeta})}(w))w_{G}(\rho)w_{G}(\rho)^{-1}$$
$$\qquad m'(w)^{-1}\rho \hat{n}_{L'(\tilde{\zeta})}(\omega_{T^L,L'(\tilde{\zeta})}(w))^{-1}\rho^{-1}$$
$$=xw_{T^L}(x)^{-1}\hat{n}_{L'(\tilde{\zeta})}(\omega_{T^L,L'(\tilde{\zeta})}(w))w_{G}(\rho)\hat{n}_{L'(\tilde{\zeta})}(\omega_{T^L,L'(\tilde{\zeta})}(w))^{-1}w_{G}(\rho)^{-1}m'(w)^{-1}$$
$$=xw_{T^L}(x)^{-1}w_{T^L}(\rho)w_{G}(\rho)^{-1}m'(w)^{-1},$$
puisque $ad_{\hat{n}_{L'(\tilde{\zeta})}(\omega_{T^L,L'(\tilde{\zeta})}(w))}\circ w_{G}=w_{T^L}$. On obtient
$$\hat{V}_{12}^L(w)=(\zeta_{1}(w)^{-1},\zeta_{2}(w),xw_{T^L}(x)^{-1}w_{T^L}(\rho)w_{G}(\rho)^{-1}m'(w)^{-1}).$$
On a effectu\'e les calculs dans le tore $\hat{T}^L_{1}\times \hat{T}^L_{2}\times \hat{T}^L$. En fait $\hat{V}_{12}^L(w)$ appartient au tore $\hat{\mathfrak{T}}_{12}^L$ qui est un quotient du pr\'ec\'edent. En particulier, puisque $\hat{m'}(w)\in Z(\hat{M}')\subset \hat{T}^{\hat{\theta},0}$, on peut multiplier l'expression pr\'ec\'edente par $(1,\hat{\xi}_{2}(m'(w))^{-1},m'(w))$ qui appartient au noyau de la projection. On obtient finalement
$$\hat{V}_{12}^L(w)=(\zeta_{1}(w)^{-1},\zeta'_{2}(w),xw_{T^L}(x)^{-1}w_{T^L}(\rho)w_{G}(\rho)^{-1}),$$
o\`u $\zeta'_{2}(w)=\hat{\xi}_{2}(m'(w))^{-1}\zeta_{2}(w)$. Notons que $\hat{\xi}_{2}(r(w),w)=(\zeta'_{2}(w),w)$.

On calcule de m\^eme les composantes $\hat{V}_{12}^M(w)$ et $\hat{V}_{12,sc}(w)$. Le calcul est beaucoup plus simple pour $\hat{V}_{12}^M(w)$ puisque le groupe $\hat{M}'$ est le m\^eme pour les deux donn\'ees et est standard. La conjugaison par $\rho$ n'intervient plus. On obtient
$$\hat{V}_{12}^M(w)=(\zeta_{1}(w)^{-1},\zeta'_{2}(w),1),$$
$$\hat{V}_{12,sc}(w)=x^{-1}w_{T^L}(x)w_{G}(\rho_{sc})w_{T^L}(\rho_{sc})^{-1},$$
o\`u, comme toujours, $\rho_{sc}$ est un \'el\'ement de $\hat{G}_{SC}$ qui a m\^eme image que $\rho$ dans $\hat{G}_{AD}$. A ce point, on voit que l'on peut supprimer les $x$ des formules ci-dessus en multipliant $(\hat{V}_{12},{\bf z})$ par le cobord de l'\'el\'ement $(1,x^{-1},1)\in \hat{S}_{12}$.

Introduisons le tore $\hat{T}_{12}^M$ quotient de $\hat{T}^M_{1}\times \hat{T}^M_{2}$ par $\hat{T}^{ M'}$ plong\'e par $\hat{\xi}_{1}\times \hat{\xi}^{-1}_{2}$. Notons $\hat{\Sigma}_{ML}$ le sous-groupe des $(t^M,t^L,t_{sc})\in \hat{T}_{12}^M\times \hat{\cal T}^L_{12}\times \hat{T}_{sc}^L$ tels que $j(t_{sc})=t^M(t^L)^{-1}$. On a un diagramme commutatif
$$\begin{array}{ccc}\hat{\Sigma}_{ML}&\stackrel{1-\hat{\theta}}{\to}&\hat{T}^L_{sc}\\ \downarrow&&\downarrow\\ \hat{S}_{12}&\stackrel{1-\hat{\theta}}{\to}&\hat{U}\\ \end{array}$$
La fl\`eche $1-\hat{\theta}$ du haut est $(t^M,t^L,t_{sc})\mapsto (1-\hat{\theta})(t_{sc})^{-1}$. Le tore $\hat{T}_{12}^M$ s'envoie naturellement dans  $\hat{{\cal T}}_{12}^M$ donc $\hat{\Sigma}^{ML}$ s'envoie naturellement dans $\hat{S}_{12}$. C'est la fl\`eche verticale de gauche. La fl\`eche verticale de droite est $t_{sc}\mapsto (1,t_{sc})$. On v\'erifie que tous ces homomorphismes sont \'equivariants pour les actions galoisiennes. De ce diagramme se d\'eduisent des homomorphismes duaux
$$\begin{array}{ccc}H^1(\Gamma_{F};U\stackrel{1-\theta}{\to}S_{12})&\times& H^1(W_{F};\hat{S}_{12}\stackrel{1-\hat{\theta}}{\to}\hat{U})\\ \downarrow &&\uparrow\\ H^1(\Gamma_{F};T^L_{ad}\stackrel{1-\theta}{\to} \Sigma_{ML})&\times&H^1(W_{F};\hat{\Sigma}_{ML}\stackrel{1-\hat{\theta}}{\to}\hat{T}_{sc}^L)\\ \end{array}$$
o\`u bien s\^ur, $\Sigma_{ML}$ est le tore dual de $\hat{\Sigma}_{ML}$. Puisque $z\in Z(\hat{M})^{\Gamma_{F},\hat{\theta}}$ et puisque l'image de ce groupe dans $\hat{G}_{AD}$ est connexe, on peut supposer $z_{sc}\in Z(\hat{M}_{sc})^{\Gamma_{F},\hat{\theta},0}\subset (\hat{T}_{sc}^M)^{\Gamma_{F},0}$. Donc l'\'el\'ement $(z_{sc},1)$ appartient \`a $\hat{U}^{\Gamma_{F},0}$. Or ce groupe est le noyau de l'accouplement avec $H^1(\Gamma_{F};U\stackrel{1-\theta}{\to }S_{12})$ ([KS], lemme A3B). Dans la formule (14), on peut donc remplacer le terme ${\bf z}$ par $(1,z_{sc})$. Pour $w\in W_{F}$, posons
$$X_{ML}^M(w)=(\zeta_{1}(w)^{-1},\zeta'_{2}(w))\in \hat{T}_{12}^M,\,\,X_{ML}^L(w)=\hat{V}_{12}^L(w)\in \hat{\cal T}_{12}^L$$
$$X_{ML,sc}(w)=\hat{V}_{12,sc}(w)\in \hat{T}_{sc}\simeq \hat{T}^L_{sc}$$
les deux derniers termes \'etant d\'ebarrass\'es des $x$ comme on l'a dit ci-dessus. On note $X_{ML}(w)=(X_{ML}^M(w),X_{ML}^L(w),X_{ML,sc}(w))$. C'est un \'el\'ement de $\hat{\Sigma}_{ML}$. On v\'erifie que le couple $(X_{ML},z_{sc})$ est un cocycle et d\'efinit donc un \'el\'ement de $H^1(W_{F};\hat{\Sigma}_{ML}\stackrel{1-\hat{\theta}}{\to}\hat{T}_{sc}^L)$. Le cocycle $(\hat{V}_{12},(1,z_{sc}))$ est l'image de cet \'el\'ement par la fl\`eche de droite du diagramme ci-dessus. 
Notons $\boldsymbol{\nu}_{ML}$ l'image naturelle de $\boldsymbol{\nu}_{12}$ dans $\Sigma_{ML}$. L'image par la fl\`eche de gauche de $(V,\boldsymbol{\nu}_{12})$ est $(V_{T^L,ad}^{-1},\boldsymbol{\nu}_{ML})$. Par compatibilit\'e des produits, (12) se r\'ecrit
$$ \Delta_{1,imp}(m'_{1},m;l'_{1},l)\Delta_{2,imp}(m'_{2},m;l'_{2},l)^{-1}=<(V_{T^L,ad}^{-1},\boldsymbol{\nu}_{ML}),(X_{ML},z_{sc})>.$$
Le tore  $T_{12}^M$ dual de $\hat{T}_{12}^M$ est le produit fibr\'e de $T_{1}^M$ et $T_{2}^M$ au-dessus de $T^{M'}$. Alors $\Sigma_{ML}$ est le quotient de $T_{12}^M\times {\cal T}_{12}^L$ par l'image antidiagonale de $\mathfrak{Z}_{12}$ (ce groupe est un sous-groupe de ${\cal T}_{12}^L$ et il s'envoie naturellement dans $T_{12}^M$). On  note abusivement ce quotient $(T_{12}^M\times {\cal T}_{12}^L)/diag_{-}(\mathfrak{Z}_{12})$. On a introduit des \'el\'ements $r'_{1}\in \tilde{R}'_{1}(F)$, $r'_{2}\in \tilde{R}'_{2}(F)$ et $r'\in \tilde{R}'(F)$. Supposons-les assez r\'eguliers. On note leurs commutants $T_{1}^R$,$T_{2}^R$ et $T^{R'}$ et on introduit leur produit fibr\'e $T_{12}^R$ au-dessus de de $T^{R'}$ (l'exposant $R$ est ici formel, il n'y a pas de groupe $R$). On introduit les tores
$$\Sigma_{MRL}=(T_{12}^M\times T_{12}^R\times {\cal T}_{12}^L)/\{(z^M,z^R,z^L)\in (\mathfrak{Z}_{12})^3; z^Mz^Rz^L=1\},$$
$$\Sigma_{MR}=(T_{12}^M\times T_{12}^R)/diag_{-}(\mathfrak{Z}_{12}),$$
$$\Sigma_{RL}=(T_{12}^R\times {\cal T}_{12}^L)/diag_{-}(\mathfrak{Z}_{12}),$$
avec les m\^emes abus d'\'ecriture que ci-dessus. Il y a des homomorphismes
$$(17) \qquad \begin{array}{ccccc}\Sigma_{MR}\times \Sigma_{RL}&&&&\Sigma_{ML}\\ &\searrow&&\swarrow\\ &&\Sigma_{MRL}&&\\ \end{array}$$
Celui de gauche est $(t^M,t^R),(u^R,t^L)\mapsto (t^M,t^Ru^R,t^L)$, celui de droite est $(t^M,t^L)\mapsto (t^M,1,t^L)$. On en d\'eduit ais\'ement un diagramme d'homomorphismes duaux
$$(18)\begin{array}{ccc}H^1(\Gamma_{F};T^L_{ad}\stackrel{1-\theta}{\to}\Sigma_{ML})&\times&H^1(W_{F};\hat{\Sigma}_{ML}\stackrel{1-\hat{\theta}}{\to}\hat{T}^L_{sc})\\ \downarrow&&\uparrow\\H^1(\Gamma_{F};T^L_{ad}\stackrel{1-\theta}{\to}\Sigma_{MRL})&\times&H^1(W_{F};\hat{\Sigma}_{MRL}\stackrel{1-\hat{\theta}}{\to}\hat{T}^L_{sc})\\ \uparrow&&\downarrow\\ (H^0(\Gamma_{F};\Sigma_{MR})\times H^1(\Gamma_{F};T^L_{ad}\stackrel{1-\theta}{\to}\Sigma_{RL}))&\times&(H^1(W_{F};\hat{\Sigma}_{MR})\times H^1(W_{F};\hat{\Sigma}_{RL}\stackrel{1-\hat{\theta}}{\to}\hat{T}^L_{sc}))\\ \end{array}$$
Le tore $\hat{\Sigma}_{MRL}$ dual de $\Sigma_{MRL}$ est le groupe des $(t^M,t^R,t^L,t^{MR}_{sc},t^{RL}_{sc})\in \hat{T}_{12}^M\times \hat{T}_{12}^{R}\times \hat{\cal T}_{12}^L\times \hat{T}^{\hat{\theta}}_{sc}\times \hat{T}_{sc}$ tels que $j(t_{sc}^{MR})=t^M(t^R)^{-1}$ et $j(t_{sc}^{RL})=t^R(t^L)^{-1}$, muni d'une action galoisienne convenable. Pour $w\in W_{F}$, posons
$$X_{MRL}^M(w)=X_{ML}^M(w)\in \hat{T}_{12}^M,\,\,X_{MRL}^R(w)=(\zeta_{1}(w)^{-1},\zeta'_{2}(w))\in \hat{T}_{12}^R,$$
$$X_{MRL}^L(w)=X_{ML}^L(w)\in \hat{\cal T}_{12}^L,\,\,X_{MRL,sc}^{MR}(w)=1,\,\,X_{MRL,sc}^{RL}(w)=X_{ML,sc}(w)\in \hat{T}_{sc}.$$
On note $X_{MRL}(w)$ l'\'el\'ement $(X_{MRL}^M(w),X_{MRL}^R(w),X_{MRL}^L(w),X_{MRL,sc}^{MR}(w),X_{MRL,sc}^{RL}(w))$ de $\hat{\Sigma}_{MRL}$. On v\'erifie que le couple $(X_{MRL},z_{sc})$ est un cocycle et d\'efinit un \'el\'ement de $H^1(W_{F};\hat{\Sigma}_{MRL}\stackrel{1-\hat{\theta}}{\to}\hat{T}^L_{sc})$. Le cocycle $(X_{MR},z_{sc})$ en est l'image par la fl\`eche de droite sup\'erieure du diagramme ci-dessus. Par compatibilit\'e des produits, on en d\'eduit
$$ \Delta_{1,imp}(m'_{1},m;l'_{1},l)\Delta_{2,imp}(m'_{2},m;l'_{2},l)^{-1}=<(V_{T^L,ad}^{-1},\boldsymbol{\nu}_{MRL}),(X_{MRL},z_{sc})>,$$
o\`u $\boldsymbol{\nu}_{MRL}$ est l'image de $\boldsymbol{\nu}_{ML}$ dans $\Sigma_{MRL}$. Notons $\boldsymbol{\mu}_{MR}$ l'image naturelle de $((\mu_{1}^M,\mu_{2}^M),(\mu_{1}^R,\mu_{2}^R)^{-1})$ dans $\Sigma_{MR}$ et $\boldsymbol{\nu}_{RL}$ l'image naturelle de $((\mu_{1}^R,\mu_{2}^R),(\nu_{12}^L)^{-1})$ dans $\Sigma_{RL}$. On v\'erifie que $\boldsymbol{\nu}_{MRL}$ est l'image de $(\boldsymbol{\mu}_{MR},\boldsymbol{\nu}_{RL})$ par la fl\`eche de gauche de (17). On v\'erifie aussi que $\boldsymbol{\mu}_{MR}\in \Sigma_{MR}^{\Gamma_{F}}=H^0(\Gamma_{F};\Sigma_{MR})$ et que $(V_{T^L,ad}^{-1},\boldsymbol{\nu}_{RL})$ d\'efinit un \'el\'ement de $H^1(\Gamma_{F};T^L_{ad}\stackrel{1-\theta}{\to}\Sigma_{RL}))$. Donc $(V_{T^L,ad}^{-1},\boldsymbol{\nu}_{MRL})$ est l'image de $(\boldsymbol{\mu}_{MR},(V_{T^L,ad}^{-1},\boldsymbol{\nu}_{RL}))$ par la fl\`eche de gauche inf\'erieure de (18). Par compatibilit\'e des produits,
 $$(19)\qquad  \Delta_{1,imp}(m'_{1},m;l'_{1},l)\Delta_{2,imp}(m'_{2},m;l'_{2},l)^{-1}=<\boldsymbol{\mu}_{MR},X_{MR}>$$
 $$<(V_{T^L,ad}^{-1},\boldsymbol{\nu}_{RL}),(X_{RL},z_{sc})>,$$
 o\`u $(X_{MR},X_{RL})$ est l'image de $X_{MRL}$ par l'homomorphisme dual de celui de gauche de (17). 
 
 Pour $w\in W_{F}$, $X_{MR}(w)$ est l'image de $(X_{MRL}^M(w),X_{MRL}^R(w),X_{MRL,sc}^{MR}(w))$ dans $\hat{\Sigma}_{MR}$.  Notons $\hat{M}'_{12}$ le groupe dual du produit fibr\'e $M'_{1}\times_{M'}M'_{2}$. On a  l'inclusion naturelle diagonale $Z(\hat{M}'_{12})\to \hat{\Sigma}_{MR}$. On voit que $X_{MR}$ est l'image par cette inclusion du cocycle $w\mapsto (\zeta_{1}(w)^{-1},\zeta'_{2}(w))\in Z(\hat{M}'_{12})$. On se rappelle que $\hat{\xi}_{1}(r(w),w)=(\zeta_{1}(w),w)$, $\hat{\xi}_{2}(r(w),w)=(\zeta'_{2}(w),w)$.  Alors le cocycle pr\'ec\'edent est l'inverse du cocycle qui d\'efinit le caract\`ere $\lambda^M$, cf. [I] 2.5. En reprenant la preuve du lemme [I] 2.5, on calcule
 $$(20)\qquad <\boldsymbol{\mu}_{MR},X_{MR}>=\lambda^M(b_{1},b_{2}).$$
 
 Notons $\hat{L}'_{12}$ le groupe dual du produit fibr\'e $L'_{1}\times_{L'}L'_{2}$. On se rappelle que le recollement est ici relatif non pas aux homomorphismes $\hat{\xi}_{1}$ et $\hat{\xi}_{2}$, mais aux homomorphismes $\hat{\xi}_{1}$ et $\hat{\xi}_{2}\circ ad_{\rho}$. On a $\hat{\xi}_{2}\circ ad_{\rho}(r(w),w)=(\hat{\xi}_{2}(\rho w_{G}(\rho)^{-1})\zeta'_{2}(w),w)$. Le caract\`ere $\lambda^L$ est donc d\'efini par le cocycle $w\mapsto (\zeta_{1}(w),\hat{\xi}_{2}(w_{G}(\rho)\rho^{-1})\zeta'_{2}(w)^{-1})\in Z(\hat{L}'_{12})$. Ce groupe s'envoie naturellement dans $\hat{\Sigma}_{RL}$. Notons $D$ le cocycle de $W_{F}$ \`a valeurs dans $\hat{\Sigma}_{RL}$ qui est l'image de l'inverse du pr\'ec\'edent. Introduisons le tore $\hat{\cal Y}_{RL}$ form\'e des $(t^R,t^L,t_{sc})\in \hat{T}^{R'}\times \hat{T}^L\times \hat{T}_{sc}$ tels que $j(t_{sc})=t^R(t^L)^{-1}$, muni d'une action galoisienne similaire \`a celle sur $\Sigma_{RL}$. Il y a un homomorphisme naturel $\hat{{\cal Y}}_{RL}\to \hat{\Sigma}_{RL}$.  On se rappelle que $\rho\in Z(\hat{R})_{*}$. Donc, pour $w\in W_{F}$,    $w_{G}(\rho)\rho^{-1}$ appartient \`a $Z(\hat{R})\cap \hat{T}^{\hat{\theta},0}\subset \hat{T}^{R'}$. Posons
 $$Y_{RL}^R(w)=w_{G}(\rho)\rho^{-1}\in \hat{T}^{R'},\,\,Y_{RL}^L(w)=w_{T^L}(\rho)\rho^{-1}\in \hat{T}^L,\,\, Y_{RL,sc}=w_{G}(\rho_{sc})w_{T^L}(\rho_{sc})^{-1}\in \hat{T}_{sc},$$
 puis $Y_{RL}(w)=(Y^R_{RL}(w),Y^L_{RL}(w),Y_{RL,sc}(w))\in \hat{\cal Y}_{RL}$. Il y a un homomorphisme naturel $\hat{p}:\hat{{\cal Y}}_{RL}\to \Sigma_{RL}$. Pour $w\in W_{F}$, $X_{RL}(w)$ est l'image de $(X_{MRL}^R(w),X_{MRL}^L(w),X_{MRL,sc}^{RL}(w))$ dans $\hat{\Sigma}_{RL}$.  On v\'erifie que $X_{RL}$ est le produit de $D$ et de $\hat{p}(Y_{RL})$. Plus pr\'ecisement le cocycle $(X_{RL},z_{sc})$ est le produit de $(D,1)$ et de $(\hat{p}(Y_{RL}),z_{sc})$.  Le produit $<(V_{T^L,ad}^{-1},\boldsymbol{\nu}_{RL}),(D,1)>$ se calcule comme on a calcul\'e $<\boldsymbol{\mu}_{MR},X_{MR}>$. Il vaut $\lambda^L(a_{1},a_{2})^{-1}$. Notons ${\cal Y}_{RL}$ le tore dual de $\hat{{\cal Y}}_{RL}$ et $p:\Sigma_{RL}\to {\cal Y}_{RL}$ l'homomorphisme dual de $\hat{p}$.  Par compatibilit\'e des produits, le second produit est \'egal \`a $<(V_{T^L,ad}^{-1},y_{RL}),(Y_{RL},z_{sc})>$ o\`u $y_{RL}=p(\boldsymbol{\nu}_{RL})$. D'o\`u
 $$<(V_{T^L,ad}^{-1},\boldsymbol{\nu}_{RL}),(X_{RL},z_{sc})>=\lambda^L(a_{1},a_{2})^{-1}<(V_{T^L,ad}^{-1},y_{RL}),(Y_{RL},z_{sc})>.$$
  Revenons \`a l'\'egalit\'e (12) o\`u on r\'etablit les indices $z$ et $\rho$ de $\tilde{\lambda}$, utilisons (19), (20) et l'\'egalit\'e pr\'ec\'edente. On obtient
  $$(21) \qquad \tilde{\lambda}_{z,\rho}(r'_{1})=<(V_{T^L,ad}^{-1},y_{RL}),(Y_{RL},z_{sc})>.$$
  
  On se rappelle qu'au cours de la d\'emonstration, on a d\^u ajuster $\rho$ de sorte que $ad_{\rho}$ envoie l'\'epinglage fix\'e de $\hat{L}'(\tilde{\zeta})$ sur celui de $\hat{L}'(z\tilde{\zeta})$. Pour cela, on a multipli\'e $\rho$ par un \'el\'ement de $Z(\hat{R})\cap \hat{T}^{\hat{\theta},0}$. On peut maintenant oublier cette modification, car la formule (21) y est insensible. En effet, si on multiplie $\rho$ par $\rho' \in Z(\hat{R})\cap \hat{T}^{\hat{\theta},0}$  cela ne change que $Y_{RL}$, qui est multipli\'e par le cocycle $Y'_{RL}$ d\'efini par
  $$Y'_{RL}(w)= (w_{G}(\rho')(\rho')^{-1},w_{T^L}(\rho')(\rho')^{-1},w_{G}(\rho'_{sc})w_{T^L}(\rho'_{sc})^{-1}).$$
  Or le couple $(Y'_{RL},1)$ est le cobord de l'\'el\'ement $(\rho',\rho',1)\in \hat{{\cal Y}}_{RL}$, donc dispara\^{\i}t par passage aux groupes de cohomologie. 
  
  Consid\'erons  l'ensemble des couples $(z,\rho)\in {\cal Z}\times Z(\hat{R})_{*} $ tels que $z(1-\hat{\theta})(\rho)^{-1}\in Z(\hat{L})^{\Gamma_{F}}$. C'est un groupe qui se projette sur ${\cal Z}$. Fixons-en un sous-ensemble fini $\underline{{\cal Z}}$ tel que la projection $\underline{{\cal Z}}\to {\cal Z}$ soit surjective et que ses fibres aient toutes le m\^eme nombre d'\'el\'ements. La formule (21) est valable pour tout $(z,\rho)\in \underline{{\cal Z}}$. La formule (10) \`a prouver est \'equivalente \`a
 $$  \sum_{(z,\rho)\in \underline{{\cal Z}}}S_{\tilde{R}'(z\tilde{\zeta} )_{1},\lambda(z\tilde{\zeta})_{1}}^{\tilde{L}'(z\tilde{\zeta} )_{1}}(\boldsymbol{\delta}(z\tilde{\zeta} )_{1},B^{\tilde{L}},({\bf f}_{\tilde{L},\omega})^{\tilde{L}'(z\tilde{\zeta} )_{1}})=0,$$
 ou encore
 $$ \sum_{(z,\rho)\in \underline{{\cal Z}}}S_{\tilde{R}'(\tilde{\zeta} )_{1},\lambda(\tilde{\zeta})_{1}}^{\tilde{L}'(\tilde{\zeta} )_{1}}(\iota(z,\rho)^{L,*}\circ(\iota(z)^{M,*})^{-1}(\boldsymbol{\delta}(\tilde{\zeta} )_{1}),B^{\tilde{L}},({\bf f}_{\tilde{L},\omega})^{\tilde{L}'(z\tilde{\zeta} )_{1}})=0.$$
 Il suffit pour cela de prouver que
 $$\sum_{(z,\rho)\in \underline{{\cal Z}}}\iota(z,\rho)^{L,*}\circ(\iota(z)^{M,*})^{-1}(\boldsymbol{\delta}(\tilde{\zeta} )_{1})=0$$
 et il suffit encore de prouver

$$(22) \qquad \sum_{(z,\rho)\in \underline{{\cal Z}}}\tilde{\lambda}_{z,\rho}(r'_{1})=0$$
 pour tout $r'_{1}$ dans un voisinage du support de $\boldsymbol{\delta}(\tilde{\zeta})_{1}$ dans $R'(\tilde{\zeta})_{1}(F)$.

 Supposons que l'on soit dans le cas (7), c'est-\`a-dire que $\hat{R}$ ne corresponde pas \`a un Levi  de $G$.  On choisit 
 $${\cal Z}=(Z(\hat{M})^{\Gamma_{F},\hat{\theta}}\cap (Z(\hat{L})^{\Gamma_{F}}(1-\hat{\theta})\circ\pi(Z(\hat{R}_{sc})^{\Gamma_{F}})))/Z(\hat{G})^{\Gamma_{F},\hat{\theta}}.$$
 Notons $\underline{{\cal Z}}_{0}$ l'ensemble des $(z,\rho)\in {\cal Z}\times \pi(Z(\hat{R}_{sc})^{\Gamma_{F}})$  tels que $z(1-\hat{\theta})(\rho)^{-1}\in Z(\hat{L})^{\Gamma_{F}}$. 
Appliquons les calculs pr\'ec\'edents \`a un couple $(z,\rho)\in \underline{{\cal Z}}_{0}$. On peut supposer $\rho_{sc}\in Z(\hat{R}_{sc})^{\Gamma_{F}}$ et $\rho=\pi(\rho_{sc})$. Le terme $(1,\rho,\rho_{sc}^{-1})$ appartient \`a $\hat{\cal Y}_{RL}$. On peut remplacer $(Y_{RL},z_{sc})$ par son produit avec le cobord associ\'e cet \'el\'ement. Ce produit n'est autre que $(1,\tau_{sc})$, o\`u $\tau_{sc}=z_{sc}(1-\hat{\theta})(\rho_{sc})^{-1}\in \hat{T}^{L,\Gamma_{F}}_{sc}$. On a \'ecrit $z=\tau(1-\hat{\theta})(\rho)$. On voit que $\tau_{sc}$ a m\^eme image que $\tau$ dans $\hat{G}_{AD}$. Donc $\tau_{sc}\in Z(\hat{L}_{sc})^{\Gamma_{F}}$. Notons encore $\tau_{sc}$ son image dans $Z(\hat{L}_{sc})^{\Gamma_{F}}/Z(\hat{L}_{sc})^{\Gamma_{F},0}$ et notons $u^L$ l'image de $V_{T^L,ad}^{-1}$ dans $H^1(\Gamma_{F},L_{ad})$. On obtient que $\tilde{\lambda}_{z,\rho}$ est la fonction constante de valeur $<u^L,\tau_{sc}>$, o\`u il s'agit du produit   sur
  $$H^1(\Gamma_{F};L_{ad})\times   Z(\hat{L}_{sc})^{\Gamma_{F}}/Z(\hat{L}_{sc})^{\Gamma_{F},0}.$$
  L'homomorphisme
  $$Z(\hat{G}_{SC})^{\Gamma_{F}}/Z(\hat{G}_{SC})^{\Gamma_{F},0}\to Z(\hat{L}_{sc})^{\Gamma_{F}}/Z(\hat{L}_{sc})^{\Gamma_{F},0}$$
  est surjectif. Choisissons $v\in Z(\hat{G}_{SC})^{\Gamma_{F}}/Z(\hat{G}_{SC})^{\Gamma_{F},0}$ qui s'envoie sur $\tau_{sc}$. Notons $u$ l'image de $u^L$ dans $H^1(\Gamma_{F},G_{AD})$. Alors la valeur de $\tilde{\lambda}_{z,\rho}$ est aussi \'egale \`a $<u,v>$, o\`u il s'agit du produit sur
  $$H^1(\Gamma_{F};G_{AD})\times Z(\hat{G}_{SC})^{\Gamma_{F}}/Z(\hat{G}_{SC})^{\Gamma_{F},0}.$$
  Pour construire $V_{T^L}$, cf. [I] 2.2, on a fix\'e une paire de Borel \'epingl\'ee ${\cal E}$ de $G$ et on a \'enonc\'e l'\'egalit\'e $dV_{T^L}=du_{{\cal E}}$ dans $H^2(\Gamma_{F};Z(G_{SC}))$. L'application $\sigma\mapsto u_{{\cal E}}(\sigma)_{ad}$ est un cocycle \`a valeurs dans $G_{AD}$ dont la classe de cohomologie ne d\'epend pas de ${\cal E}$. Notons-la $u_{G}$. Parce que l'application $H^1(\Gamma_{F};G_{AD})\to H^2(\Gamma_{F};Z(G_{SC}))$ est injective, l'\'egalit\'e rappel\'ee ci-dessus montre que $u=u_{G}$. Par l'accouplement ci-dessus, $u_{G}$ d\'efinit un caract\`ere de $Z(\hat{G}_{SC})^{\Gamma_{F}} $. Notons $Ann(u_{G})\subset Z(\hat{G}_{SC})^{\Gamma_{F}} $ le noyau de ce caract\`ere et $v_{z,\rho}$ l'image de $v$ dans le quotient $Z(\hat{G}_{SC})^{\Gamma_{F}}/Ann(u_{G})$. On a effectu\'e divers choix pour construire cet \'el\'ement. Mais le r\'esultat de notre calcul montre que celui-ci ne d\'epend pas de ces choix. On a donc une application
 $$\begin{array}{ccc}\underline{{\cal Z}}_{0}&\to&Z(\hat{G}_{SC})^{\Gamma_{F}}/Ann(u_{G})\\( z,\rho)&\mapsto&v_{z,\rho}. \\ \end{array}$$
 C'est un homomorphisme \`a valeurs dans un groupe fini. Notons $J$ son image. Fixons un sous-ensemble $\underline{{\cal Z}}_{1}$ de l'ensemble de d\'epart se projetant bijectivement sur $J$.  Fixons aussi un sous-ensemble $\underline{{\cal Z}}_{2}$ se projetant bijectivement sur ${\cal Z}$. Notons $\underline{{\cal Z}}$ l'ensemble des produits $(z_{1},\rho_{1})(z_{2},\rho_{2})$, pour $(z_{1},\rho_{1})\in \underline{{\cal Z}}_{1}$ et $(z_{2},\rho_{2})\in \underline{{\cal Z}}_{2}$. On v\'erifie que les deux projections
 $$\begin{array}{ccccc}&&\underline{{\cal Z}}&&\\ &\swarrow&&\searrow&\\ {\cal Z}&&&&J\\ \end{array}$$
 sont surjectives et toutes leurs fibres ont m\^eme nombre d'\'el\'ements. La somme (22) est donc proportionnelle \`a
 $$\sum_{v\in J}<u,v>.$$
 Pour d\'emontrer la relation (22), il suffit de prouver que $J\not=\{1\}$. On utilise le lemme 2.1 de [A4]: puisque $\hat{R}$ ne correspond pas \`a un Levi de $G$, l'image dans  $Z(\hat{G}_{SC})^{\Gamma_{F}}/Ann(u_{G})$ du groupe $Z(\hat{G}_{SC})^{\Gamma_{F}}\cap Z(\hat{R}_{sc})^{\Gamma_{F},0}$ n'est pas r\'eduite \`a l'identit\'e (Arthur \'enonce ce lemme apr\`es passage \`a une forme quasi-d\'eploy\'ee, mais c'est \'equivalent \`a notre assertion). Il  nous suffit de prouver que, pour tout $v\in Z(\hat{G}_{SC})^{\Gamma_{F}}\cap Z(\hat{R}_{sc})^{\Gamma_{F},0}$, on peut trouver $(z,\rho)\in \underline{{\cal Z}}_{0}$  et effectuer les divers choix n\'ecessaires  de sorte que $v$ soit l'\'el\'ement associ\'e ci-dessus \`a $(z,\rho)$. Rappelons les deux \'egalit\'es
 $$Z(\hat{R}_{sc})^{\Gamma_{F},0}=Z(\hat{R}_{sc})^{\Gamma_{F},\hat{\theta},0}(1-\hat{\theta})(Z(\hat{R}_{sc})^{\Gamma_{F},0}),$$
$$Z(\hat{R}_{sc})^{\Gamma_{F},\hat{\theta},0}=Z(\hat{M}_{sc})^{\Gamma_{F},\hat{\theta},0}Z(\hat{L}_{sc})^{\Gamma_{F},\hat{\theta},0}.$$
Un \'el\'ement $v\in Z(\hat{G}_{SC})^{\Gamma_{F}}\cap Z(\hat{R}_{sc})^{\Gamma_{F},0}$ peut donc s'\'ecrire $v=(\tau'_{sc})^{-1}z_{sc}(1-\hat{\theta})(\rho_{sc}^{-1})$, avec $\tau'_{sc}\in Z(\hat{L}_{sc})^{\Gamma_{F},\hat{\theta},0}$, $z_{sc}\in Z(\hat{M}_{sc})^{\Gamma_{F},\hat{\theta},0}$ et $\rho_{sc}\in Z(\hat{R}_{sc})^{\Gamma_{F},0}$. Posons $z=\pi(z_{sc})$. On a $z_{sc}=\tau_{sc}(1-\hat{\theta})(\rho_{sc})$, o\`u $\tau_{sc}=\tau'_{sc}v\in Z(\hat{L}_{sc})^{\Gamma_{F}}$. Donc $(z,\rho)\in \underline{{\cal Z}}_{0}$. On peut choisir les \'el\'ements $\rho_{sc}$ et $z_{sc}$ pour effectuer nos calculs. L'\'el\'ement $\tau_{sc}$ apparaissant plus haut est celui que l'on vient de d\'efinir et on peut ensuite choisir $v$ comme rel\`evement dans $Z(\hat{G}_{SC})^{\Gamma_{F}}$ de l'image de $\tau_{sc}$ dans $Z(\hat{L}_{sc})^{\Gamma_{F}}/Z(\hat{L}_{sc})^{\Gamma_{F},0}$. Cela prouve l'assertion et ach\`eve de prouver l'assertion (iii) de la proposition dans le cas o\`u $\hat{R}$ ne correspond pas \`a un Levi de $G$.

Supposons maintenant que l'on soit dans le cas (8), c'est-\`a-dire que $\hat{R}$ corresponde \`a un Levi de $G$ donc aussi \`a un espace de Levi  $\tilde{R}$ de $\tilde{G}$, et que le  support de $\boldsymbol{\delta}$ soit contenu dans l'ensemble not\'e  $\tilde{R}'(F)^{out}$ en (8).   On prend 
${\cal Z}$  maximal, c'est-\`a-dire
$${\cal Z}=(Z(\hat{M})^{\Gamma_{F},\hat{\theta}}\cap (Z(\hat{L})^{\Gamma_{F}}(1-\hat{\theta})(Z(\hat{R})_{*})))/Z(\hat{G})^{\Gamma_{F},\hat{\theta}}$$
et on note $\underline{{\cal Z}}_{0}$ l'ensemble des $(z,\rho)\in {\cal Z}\times Z(\hat{R})_{*}$ tels que $z(1-\hat{\theta})(\rho)^{-1}\in Z(\hat{L})^{\Gamma_{F}}$. Soit $(z,\rho)$ un \'el\'ement de cet ensemble.
On  peut supposer que $\tilde{R}$ est inclus dans $\tilde{M}$ et $\tilde{L}$. Le triplet ${\bf R}'=(R',{\cal R}',\zeta)$ est une donn\'ee endoscopique elliptique de $\tilde{R}$ mais elle n'est pas relevante. Fixons un \'el\'ement assez r\'egulier $r_{\sharp}\in \tilde{R}(F)$. Bien qu'il ne corresponde \`a aucun \'el\'ement de $\tilde{R}'_{1}(F)$, on peut lui associer une partie des donn\'ees que l'on a associ\'ees \`a $l$ ou $m$: le tore $T^R_{\sharp}$ (commutant de $R_{r_{\sharp}}$), le tore $T_{\sharp}^{R'}=T^{R}_{\sharp}/(1-\theta)(T^R_{\sharp})$ (ou ici $\theta=ad_{r_{\sharp}}$),  un \'el\'ement $\nu^R_{\sharp}\in T^R_{\sharp}$ et son image $\mu^R_{\sharp}$ dans $T^{R'}_{\sharp}$, une cocha\^{\i}ne $V_{T^R_{\sharp}}$. Posons
$${\cal Y}_{R\sharp L}=(T^{R'}\times T^{R'}_{\sharp}\times T^L)/\{(z^{R'},z^{R'}_{\sharp},z^L)\in Z(G)^3; z^{R'}z^{R'}_{\sharp}z^L=1\}$$
(par abus d'\'ecriture, on ne distingue pas un \'el\'ement de $Z(G)$ de ses images naturelles dans diff\'erents quotients);
$${\cal Y}^{R\sharp}=(T^{R'}\times T^{R'}_{\sharp})/diag_{-}(Z(G)),$$
$${\cal Y}^{\sharp L}=(T^{R'}_{\sharp}\times T^L)/diag_{-}(Z(G)).$$
Leurs tores duaux se d\'ecrivent de fa\c{c}on similaire aux pr\'ec\'edents, par exemple $\hat{{\cal Y}}^{R\sharp L}$ est le groupe des $(t^{R'},t^{R'}_{\sharp},t^L,t_{sc}^{R\sharp},t_{sc}^{\sharp L})\in \hat{T}^{R'}\times \hat{T}^{R'}_{\sharp}\times \hat{T}_{sc}^{\hat{\theta}}\times \hat{T}_{sc}$ tels que $j(t_{sc}^{R\sharp})=t^{R'}(t^{R'}_{\sharp})^{-1}$, $j(t_{sc}^{\sharp L})=t^{R'}_{\sharp}(t^L)^{-1}$. On a un diagramme similaire \`a (18):
$$\begin{array}{ccc}H^1(\Gamma_{F};T^L_{ad}\stackrel{1-\theta}{\to}{\cal Y}_{RL})&\times&H^1(W_{F};\hat{{\cal Y}}_{RL}\stackrel{1-\hat{\theta}}{\to}\hat{T}^L_{sc})\\ \downarrow&&\uparrow\\H^1(\Gamma_{F};T^L_{ad}\stackrel{1-\theta}{\to}{\cal Y}_{R\sharp L})&\times&H^1(W_{F};\hat{{\cal Y}}_{R\sharp L}\stackrel{1-\hat{\theta}}{\to}\hat{T}^L_{sc})\\ \uparrow&&\downarrow\\ H^0(\Gamma_{F};{\cal Y}_{R\sharp})\times H^1(\Gamma_{F};T^L_{ad}\stackrel{1-\theta}{\to}{\cal Y}_{\sharp L})&\times&H^1(W_{F};\hat{{\cal Y}}_{R\sharp})\times H^1(W_{F};\hat{{\cal Y}}_{\sharp L}\stackrel{1-\hat{\theta}}{\to}\hat{T}^L_{sc})\\ \end{array}$$
Pour $w\in W_{F}$, posons
$$Y_{R\sharp L}(w)=(w_{G}(\rho)\rho^{-1} ,w_{G}(\rho)\rho^{-1} ,w_{T^L}(\rho)\rho^{-1},1,w_{G}(\rho_{sc})w_{T^L}(\rho_{sc})^{-1})\in  \hat{{\cal Y}}_{R\sharp L}.$$
 Le couple $(Y_{R\sharp L},z_{sc})$ d\'efinit un \'el\'ement de $H^1(W_{F};\hat{{\cal Y}}_{R\sharp L}\stackrel{1-\hat{\theta}}{\to}\hat{T}^L_{sc})$ qui s'envoie sur $(Y_{RL},z_{sc})$ par l'homomorphisme en haut \`a droite du diagramme ci-dessus. Notons $y_{R\sharp L}$ l'image de $(\mu^R,1,(\nu^L)^{-1})$ dans ${\cal Y}_{R\sharp L}$. Alors $(V_{T^L,ad}^{-1},y_{R\sharp L})$ est l'image de $(V_{T^L,ad}^{-1},y_{RL})$ par l'homomorphisme en haut \`a gauche du diagramme. Par compatibilit\'e des produits,
$$<(V_{T^L,ad}^{-1},y_{RL}),(Y_{RL},z_{sc})>=<(V_{T^L,ad}^{-1},y_{R\sharp L}),(Y_{R\sharp L},z_{sc})>.$$
Notons $y_{R\sharp}$ l'image de $(\mu^R,(\mu^{R}_{\sharp})^{-1})$ dans ${\cal Y}_{R\sharp}$ et $y_{\sharp L}$ celle de $(\mu^R_{\sharp}, (\nu^L)^{-1})$ dans ${\cal Y}_{\sharp L}$. Pour $w\in W_{F}$, posons
$$Y_{R\sharp}(w)=(w_{G}(\rho)\rho^{-1} ,w_{G}(\rho)\rho^{-1},1 )\in \hat{{\cal Y}}_{R\sharp},$$
$$Y_{\sharp L}(w)=(w_{G}(\rho)\rho^{-1} ,w_{T^L}(\rho)\rho^{-1},w_{G}(\rho_{sc})w_{T^L}(\rho_{sc})^{-1})\in \hat{\cal Y}_{\sharp L}.$$
Le cocycle $(V_{T^L,ad}^{-1},y_{R\sharp L})$ est l'image de la paire de cocycles $(y_{R\sharp},(V_{T^L,ad}^{-1},y_{\sharp L}))$ par la fl\`eche en bas \`a gauche du diagramme. La paire $(Y_{R\sharp},(Y_{\sharp L},z_{sc}))$ est l'image de $(Y_{R\sharp L},z_{sc})$ par la fl\`eche en bas \`a droite. Par compatibilit\'e des produits, on obtient
$$(23) \qquad <(V_{T^L,ad}^{-1},y_{R\sharp L}),(Y_{R\sharp L},z_{sc})>=<y_{R\sharp},Y_{R\sharp}><(V_{T^L,ad}^{-1},y_{\sharp L}),(Y_{\sharp L},z_{sc})>.$$

Montrons que 
$$(24) \qquad <(V_{T^L,ad}^{-1},y_{\sharp L}),(Y_{\sharp L},z_{sc})>=1.$$
Posons $U_{\sharp,L}=(T^R_{\sharp,sc}\times T^L_{sc})/diag_{-}(Z(G_{SC}))$ et ${\cal X}_{\sharp L}=(T^R_{\sharp}\times T^L)/diag_{-}(Z(G))$. Notons $x_{\sharp L}$ l'image de $(\nu^R_{\sharp},\mu^L)$  dans ${\cal X}_{\sharp L}$. Il y a un homomorphisme naturel
$$(25) \qquad H^1(\Gamma_{F};U_{\sharp L}\stackrel{1-\theta}{\to} {\cal X}_{\sharp L})\to H^1(\Gamma_{F}; T^L_{ad}\stackrel{1-\theta}{\to}{\cal Y}_{\sharp L}).$$
Le couple $((V_{T^R_{\sharp}},V_{T^L}^{-1}),x_{\sharp L})$ d\'efinit un \'el\'ement du premier groupe qui s'envoie sur l'\'el\'ement $(V_{T^L,ad}^{-1},y_{\sharp L})$ du second. Le tore $\hat{\cal X}_{\sharp L}$ est form\'e des $(t^R_{\sharp},t^L,t_{sc})\in \hat{T}^R_{\sharp}\times \hat{T}^L\times \hat{T}_{sc}$ tels que $j(t_{sc})=t^R_{\sharp}(t^L)^{-1}$ tandis que $\hat{U}_{\sharp L}=(\hat{T}^R_{\sharp}\times \hat{T}^L_{sc})/diag(Z(\hat{G}_{SC}))$. Pour $w\in W_{F}$, on pose
$$X_{\sharp L}(w)=(w_{G}(\rho)\rho^{-1},w_{T^L}(\rho)\rho^{-1},w_{G}(\rho_{sc})w_{T^L}(\rho_{sc})^{-1})\in \hat{\cal X}_{\sharp L}.$$
Alors $(X_{\sharp L},(1,z_{sc}))$ d\'efinit un \'el\'ement de $H^1(W_{F};\hat{{\cal X}}_{\sharp L}\stackrel{1-\hat{\theta}}{\to}\hat{U}_{\sharp L})$. C'est l'image de $(Y_{\sharp L},z_{sc})$ par l'homomorphisme dual de (25). D'o\`u l'\'egalit\'e
$$<(V_{T^L,ad}^{-1},y_{\sharp L}),(Y_{\sharp L},z_{sc})>=<((V_{T^R_{\sharp}},V_{T^L}^{-1}),x_{\sharp L}),(X_{\sharp L},(1,z_{sc}))>.$$
Le triplet $(\rho,\rho,1)$ appartient \`a $\hat{{\cal X}}_{\sharp L}$. On peut multiplier $(X_{\sharp L},(1,z_{sc}))$ par le cobord de cet \'el\'ement. On obtient un cocycle qui est l'image par l'homomorphisme naturel
$$\hat{U}_{\sharp L}^{\Gamma_{F}}=H^0(W_{F};\hat{U}_{\sharp L})\to H^1(W_{F};\hat{{\cal X}}_{\sharp L}\stackrel{1-\hat{\theta}}{\to}\hat{U}_{\sharp L})$$
de l'\'el\'ement $u_{\sharp L}=((1-\hat{\theta})(\rho_{sc})^{-1},z_{sc}(1-\hat{\theta})(\rho_{sc})^{-1})$ de $\hat{U}_{\sharp L}^{\Gamma_{F}}$.   Donc
$$<((V_{T^R_{\sharp}},V_{T^L}^{-1}),x_{\sharp L}),(X_{\sharp L},(1,z_{sc}))>=<(V_{T^R_{\sharp}},V_{T^L}^{-1}),u_{\sharp L}>.$$
On sait que $\hat{U}_{\sharp L}^{\Gamma_{F},0}$ est contenu dans le noyau de l'accouplement intervenant ici. On va montrer que $u_{\sharp L}$ appartient \`a ce sous-groupe, ce qui prouvera (24). On \'ecrit
$$u_{\sharp L}=(z_{sc}^{-1},1)(\tau_{sc},\tau_{sc})$$
o\`u $\tau_{sc}=z_{sc}(1-\hat{\theta})(\rho_{sc})^{-1}$. On a suppos\'e $z_{sc}\in Z(\hat{M}_{sc})^{\Gamma_{F},\hat{\theta},0}$, a fortiori $z_{sc}\in Z(\hat{R}_{sc})^{\Gamma_{F},0}$.   Le couple $(z_{sc},1)$ (ou plut\^ot son image dans $\hat{U}_{\sharp L}$) appartient donc \`a $\hat{U}_{\sharp L}^{\Gamma_{F},0}$. L'\'el\'ement $\tau_{sc}$ a m\^eme image que $\tau$ dans $\hat{G}_{AD}$, cf. (9) pour la d\'efinition de $\tau$. D'apr\`es la d\'efinition de $\hat{U}_{\sharp L}$ (qui est un quotient par $Z(\hat{G}_{SC})$), on peut aussi bien remplacer $\tau_{sc}$ par un \'el\'ement quelconque de $\hat{G}_{SC}$  qui a m\^eme image que $\tau$ dans $\hat{G}_{AD}$. Puisque $\tau\in Z(\hat{L})^{\Gamma_{F}}$, on peut supposer $\tau_{sc}\in Z(\hat{L}_{sc})^{\Gamma_{F},0}$. Mais alors $(\tau_{sc},\tau_{sc})\in \hat{U}_{\sharp L}^{\Gamma_{F},0}$. Cela ach\`eve la preuve de (24).

Introduisons le groupe $R_{0}$ quasi-d\'eploy\'e et dual de $\hat{R}^{\hat{\theta},0}$, c'est-\`a-dire l'analogue du groupe $G_{0}$ de [I] 1.12 quand on remplace $\tilde{G}$ par $\tilde{R}$. Les tores $T^{R'}$ et $T^{R'}_{\sharp}$ se r\'ealisent naturellement comme sous-tores de $R_{0}$. Pour $\rho\in Z(\hat{R})_{*}$, le cocycle $w\mapsto  w_{G}(\rho)\rho^{-1}$ est \`a valeurs dans $Z(\hat{R}_{0})$ et d\'efinit un caract\`ere   de $R_{0}(F)$, cf. [I] 1.13, que l'on note ici $\chi_{\rho}$. Ce caract\`ere se factorise par $R_{0,ab}(F)$ et est alors trivial sur $N^{R}(R_{ab}(F))$.  Par un calcul d\'ej\`a fait plusieurs fois, on a
$$<y_{R\sharp},Y_{R\sharp}>=\chi_{\rho}(\mu^R(\mu^R_{\sharp})^{-1}).$$
Notons $(R_{0,ab}(F)/N^{R}(R_{ab}(F)))^{\vee}$ le groupe dual du groupe fini $R_{0,ab}(F)/N^{R}(R_{ab}(F))$. On obtient un homomorphisme
$$(26)\qquad \begin{array}{ccc}\underline{{\cal Z}}_{0}&\to &(R_{0,ab}(F)/N^{R}(R_{ab}(F)))^{\vee}\\ (z,\rho)&\mapsto &\chi_{\rho}.\\ \end{array}$$
Notons $J$ son image. Comme dans la preuve du cas o\`u (7) est v\'erifi\'ee, on peut d\'efinir $\underline{{\cal Z}}$ de sorte que la somme (22) soit proportionnelle \`a
$$\sum_{\chi\in J}\chi(\mu^R(\mu^R_{\sharp})^{-1}).$$
Il reste \`a montrer que cette somme est nulle sous l'hypoth\`ese de (22).

Comme en [I] 1.13, on d\'eduit de $\chi\in (R_{0,ab}(F)/N^{R}(R_{ab}(F)))^{\vee}$ une application $\tilde{\chi}$ sur $\tilde{R}_{0,ab}(F)$ qui vaut $1$ sur l'image de $\tilde{R}_{ab}(F)$ par $N^{\tilde{R}}$ et qui v\'erifie $\tilde{\chi}(x\gamma)=\chi(x)\tilde{\chi}(\gamma)$ pour tous $x\in R_{0,ab}(F)$ et $\gamma\in \tilde{R}_{0,ab}(F)$. En reprenant les constructions, on v\'erifie que l'on a l'\'egalit\'e $N^{\tilde{R}',\tilde{R}}(r')=\mu^R(\mu^R_{\sharp})^{-1}N^{\tilde{R}}(r_{\sharp})$. Donc
$$\chi(\mu^R(\mu^R_{\sharp})^{-1})=\tilde{\chi}(N^{\tilde{R}',\tilde{R}}(r')).$$
Reportons-nous aux hypoth\`eses de (22) et (8): on peut supposer que $r'$ appartient \`a l'ensemble $\tilde{R}'(F)^{out}$ d\'efini en (8). D'apr\`es sa d\'efinition, on a
$$\sum_{\chi\in (R_{0,ab}(F)/N^{R}(R_{ab}(F)))^{\vee}}\tilde{\chi}(N^{\tilde{R}',\tilde{R}}(r'))=0.$$
Pour achever de prouver (22), il suffit de prouver que $J=  (R_{0,ab}(F)/N^{R}(R_{ab}(F)))^{\vee}$, autrement dit que l'homomorphisme (26) est surjectif. Puisque $(R_{0,ab}(F)/N^{R}(R_{ab}(F)))^{\vee}$ est l'image par $\rho\mapsto \chi_{\rho}$ de $Z(\hat{R})_{*}/(Z(\hat{R})\cap \hat{T}^{\hat{\theta},0})Z(\hat{R})^{\Gamma_{F}}$, il suffit de prouver que l'homomorphisme
$$\begin{array}{ccc}\underline{{\cal Z}}_{0}&\to &Z(\hat{R})_{*}/(Z(\hat{R})\cap \hat{T}^{\hat{\theta},0})Z(\hat{R})^{\Gamma_{F}} \\ (z,\rho)&\mapsto &\rho (Z(\hat{R})\cap \hat{T}^{\hat{\theta},0})Z(\hat{R})^{\Gamma_{F}}\\ \end{array}$$
est surjectif. On a les relations
$$(1-\hat{\theta})(Z(\hat{R})_{*})\subset Z(\hat{R})^{\Gamma}=Z(\hat{G})^{\Gamma}Z(\hat{R})^{\Gamma,\hat{\theta},0}(1-\hat{\theta})(Z(\hat{R})^{\Gamma})$$
$$=Z(\hat{G})^{\Gamma}Z(\hat{M})^{\Gamma,\hat{\theta},0}Z(\hat{L})^{\Gamma,\hat{\theta},0}(1-\hat{\theta})(Z(\hat{R})^{\Gamma})\subset Z(\hat{M})^{\Gamma,\hat{\theta}}Z(\hat{L})^{\Gamma}(1-\hat{\theta})(Z(\hat{R})^{\Gamma}).$$
Soit $\rho\in Z(\hat{R})_{*}$. Ecrivons $(1-\hat{\theta})(\rho)=z\tau^{-1}(1-\hat{\theta})((\rho')^{-1})$, avec $z\in Z(\hat{M})^{\Gamma,\hat{\theta}}$, $\tau\in Z(\hat{L})^{\Gamma}$, $\rho'\in Z(\hat{R})^{\Gamma}$. Alors $z=\tau(1-\hat{\theta})(\rho\rho')$ appartient \`a ${\cal Z}$, $(z,\rho\rho')$ appartient \`a $\underline{{\cal Z}}_{0}$ et l'image par l'homomorphisme ci-dessus est $\rho (Z(\hat{R})\cap \hat{T}^{\hat{\theta},0})Z(\hat{R})^{\Gamma_{F}}$. Cela ach\`eve enfin la preuve. $\square$

\bigskip

{\bf Variante.} Supposons $(G,\tilde{G},{\bf a})$ quasi-d\'eploy\'e et \`a torsion int\'erieure. Fixons un syst\`eme de fonctions $B$ comme en 1.9. On a des assertions analogues \`a (i) et (iii) pour les int\'egrales $I_{\tilde{M}}^{\tilde{G},{\cal E}}({\bf M}',\boldsymbol{\delta}^{\bf M'},B,{\bf f})$. En fait, sur notre corps $F$ non-archim\'edien, les hypoth\`eses de (iii) ne sont jamais v\'erifi\'ees car, dans la situation quasi-d\'eploy\'ee et \`a torsion int\'erieure, un groupe de Levi $R'$ de $M'$ est toujours relevant. 

{\bf Variante.} Supposons $G=\tilde{G}$ et ${\bf a}=1$. Fixons une fonction $B$ comme en 1.8. On a des assertions analogues \`a (i) et (iii) pour les int\'egrales $I_{M}^{G,{\cal E}}({\bf M}',\boldsymbol{\delta}^{\bf M'},B,{\bf f})$, en supposant $\boldsymbol{\delta}^{\bf M'}$ \`a support unipotent (puisque l'on n'a d\'efini ces termes que sous cette hypoth\`ese).

  \bigskip
  
  \subsection{Int\'egrales orbitales pond\'er\'ees $\omega$-\'equivariantes endoscopiques}
  
     Soient $(G,\tilde{G},{\bf a})$ un triplet quelconque et $\tilde{M}$ un espace de Levi de $\tilde{G}$.
  
  \ass{Lemme}{Pour tout ${\bf M}'\in {\cal E}(\tilde{M},{\bf a})$, soit $\boldsymbol{\delta}_{{\bf M}'}\in D^{st}_{g\acute{e}om}({\bf M}')\otimes Mes(M'(F))^*$. Supposons
  $$\sum_{{\bf M}'\in {\cal E}(\tilde{M},{\bf a})}transfert(\boldsymbol{\delta}_{{\bf M}'})=0.$$
  Alors 
  $$\sum_{{\bf M}'\in {\cal E}(\tilde{M},{\bf a})}I^{\tilde{G},{\cal E}}_{\tilde{M}}({\bf M}',\boldsymbol{\delta}_{{\bf M}'},{\bf f})=0$$
  pour tout ${\bf f}\in I(\tilde{G}(F),\omega)\otimes Mes(G(F))$.}
  
  Preuve. Par lin\'earit\'e, on peut fixer une classe de conjugaison g\'eom\'etrique semi-simple ${\cal O}$ dans $\tilde{G}(F)$ et  supposer que, pour tout ${\bf M}'$, $\boldsymbol{\delta}_{{\bf M}'}$ appartient \`a $D^{st}_{g\acute{e}om}({\cal O}_{\tilde{M}'})\otimes Mes(M'(F))^*$, o\`u ${\cal O}_{\tilde{M}'}$ est la r\'eunion des classes de conjugaison g\'eom\'etriques dans $\tilde{M}'(F)$ correspondant \`a une classe dans ${\cal O}\cap \tilde{M}(F)$.  D'apr\`es [I] proposition 5.7, il suffit de prouver la conclusion de l'\'enonc\'e quand la famille $(\boldsymbol{\delta}_{{\bf M}'})_{{\bf M}'\in {\cal E}(\tilde{M},{\bf a})}$ appartient \`a l'un des sous-espaces d\'ecrits par chacune des conditions (3), (4), (5) de cette r\'ef\'erence. Dans le cas (3), l'assertion r\'esulte du (iii) de la proposition 1.14. Dans le cas (4), elle r\'esulte du corollaire 1.13. Dans le cas (5), elle r\'esulte du (i) de la proposition 1.14. $\square$
  
    Soit $\boldsymbol{\gamma}\in D_{g\acute{e}om}(\tilde{M}(F),\omega)\otimes Mes(M(F))^*$. Gr\^ace \`a la proposition 5.7 de [I], il existe une famille $(\boldsymbol{\delta}_{{\bf M}'})_{{\bf M}'\in {\cal E}(\tilde{M},{\bf a})}$, avec $\boldsymbol{\delta}_{{\bf M}'}\in D^{st}_{g\acute{e}om}({\bf M}')\otimes Mes(M'(F))^*$ pour tout ${\bf M}'$, de sorte que 
  $$\boldsymbol{\gamma}=\sum_{{\bf M}'\in {\cal E}(\tilde{M},{\bf a})}transfert(\boldsymbol{\delta}_{{\bf M}'}).$$
  Pour ${\bf f}\in I(\tilde{G}(F),\omega)\otimes Mes(G(F))$, on pose
  $$I_{\tilde{M}}^{\tilde{G},{\cal E}}(\boldsymbol{\gamma},{\bf f})=\sum_{{\bf M}'\in {\cal E}(\tilde{M},{\bf a})}I^{\tilde{G},{\cal E}}_{\tilde{M}}({\bf M}',\boldsymbol{\delta}_{{\bf M}'},{\bf f}).$$
  Cette d\'efinition est loisible puisque le lemme ci-dessus nous dit que le membre de droite ne d\'epend pas de la famille $(\boldsymbol{\delta}_{{\bf M}'})_{{\bf M}'\in {\cal E}(\tilde{M},{\bf a})}$ choisie. 
  
  Soit $\tilde{R}$ un espace de Levi de $\tilde{M}$ et soit $\boldsymbol{\gamma}\in D_{g\acute{e}om}(\tilde{R}(F),\omega)\otimes Mes(R(F))^*$. Alors on a l'\'egalit\'e
  
  $$(1)\qquad I_{\tilde{M}}^{\tilde{G},{\cal E}}(\boldsymbol{\gamma}^{\tilde{M}},{\bf f})=\sum_{\tilde{L}\in {\cal L}(\tilde{R})}d^{\tilde{G}}_{\tilde{R}}(\tilde{M},\tilde{L})I_{\tilde{R}}^{\tilde{L},{\cal E}}(\boldsymbol{\gamma},{\bf f}_{\tilde{L},\omega}).$$
  
  Preuve. Par lin\'earit\'e, on peut supposer qu'il existe ${\bf R}'\in {\cal E}(\tilde{R},{\bf a})$ et $\boldsymbol{\delta}\in D_{g\acute{e}om}^{st}({\bf R}')\otimes Mes(R'(F))^*$ tels que $\boldsymbol{\gamma}$ soit le transfert de $\boldsymbol{\delta}$. Ecrivons ${\bf R}'=(R',{\cal R}',\tilde{s})$. De $\tilde{M}$ se d\'eduit une donn\'ee endoscopique ${\bf M}'=(M',{\cal M}',\tilde{s})$. Quitte \`a multiplier $\tilde{s}$ par un \'el\'ement de $Z(\hat{R})^{\Gamma_{F},\hat{\theta}}$, on peut supposer ${\bf M}'$ elliptique. Par compatibilit\'e du transfert \`a l'induction, $\boldsymbol{\gamma}^{\tilde{M}}$ est le transfert de $\boldsymbol{\delta}^{\bf M'}$. Par d\'efinition, on a alors
  $$I_{\tilde{M}}^{\tilde{G},{\cal E}}(\boldsymbol{\gamma}^{\tilde{M}},{\bf f})=I_{\tilde{M}}^{\tilde{G},{\cal E}}({\bf M}',\boldsymbol{\delta}^{\bf M'},{\bf f}).$$
  On applique la proposition 1.14(i):
  $$I_{\tilde{M}}^{\tilde{G},{\cal E}}(\boldsymbol{\gamma}^{\tilde{M}},{\bf f})= \sum_{\tilde{L}\in {\cal L}(\tilde{R})}d^{\tilde{G}}_{\tilde{R}}(\tilde{M},\tilde{L})I_{\tilde{R}}^{\tilde{L},{\cal E}}({\bf R}',\boldsymbol{\delta},{\bf f}_{\tilde{L},\omega}).$$
  Toujours par d\'efinition, on a pour tout $\tilde{L}$ l'\'egalit\'e
  $$I_{\tilde{R}}^{\tilde{L},{\cal E}}({\bf R}',\boldsymbol{\delta},{\bf f}_{\tilde{L},\omega})=I_{\tilde{R}}^{\tilde{L},{\cal E}}(\boldsymbol{\gamma},{\bf f}_{\tilde{L},\omega}).$$
  L'\'egalit\'e (1) s'ensuit. $\square$
  
  {\bf Variante.} Supposons que $(G,\tilde{G},{\bf a})$ soit quasi-d\'eploy\'e et \`a torsion int\'erieure. Fixons un syst\`eme de fonctions $B$ comme en 1.9. Pour $\boldsymbol{\gamma}\in D_{g\acute{e}om}(\tilde{M}(F))\otimes Mes(M(F))^*$ et ${\bf f}\in I(\tilde{G}(F))\otimes Mes(G(F))$, on d\'efinit $I_{\tilde{M}}^{\tilde{G},{\cal E}}(\boldsymbol{\gamma},B,{\bf f})$ de la m\^eme  fa\c{c}on que ci-dessus. Ce terme v\'erifie l'analogue de la relation (1).  
  
  {\bf Variante.} Supposons $G=\tilde{G}$ et ${\bf a}=1$. Soit $B$ une fonction comme en 1.8. Pour $\boldsymbol{\gamma}\in D_{unip}(M(F))\otimes Mes(M(F))^*$ et ${\bf f}\in I(G(F))\otimes Mes(G(F))$, on d\'efinit $I_{M}^{G,{\cal E}}(\boldsymbol{\gamma},B,{\bf f})$ de la m\^eme  fa\c{c}on que ci-dessus. Ce terme v\'erifie l'analogue de la relation (1).
  
  \bigskip
  
  \subsection{Le th\'eor\`eme principal}
  
    Soient $(G,\tilde{G},{\bf a})$ un triplet quelconque et  $\tilde{M}$ un espace de Levi de $\tilde{G}$.
  
  \ass{Th\'eor\`eme (\`a prouver)}{(i) Soient $\boldsymbol{\gamma}\in D_{g\acute{e}om}(\tilde{M}(F),\omega)\otimes Mes(M(F))^*$ et ${\bf f}\in I(\tilde{G}(F),\omega)\otimes Mes(G(F))$.  Alors on a l'\'egalit\'e
  $$I_{\tilde{M}}^{\tilde{G},{\cal E}}(\boldsymbol{\gamma},{\bf f})=I_{\tilde{M}}^{\tilde{G}}(\boldsymbol{\gamma},{\bf f})$$.
  
  (ii) Supposons que $(G,\tilde{G},{\bf a})$ soit quasi-d\'eploy\'e et \`a torsion int\'erieure. Fixons un syst\`eme de fonctions $B$ comme en 1.9.. Soient $\boldsymbol{\gamma}\in D_{g\acute{e}om}(\tilde{M}(F))\otimes Mes(M(F))^*$ et ${\bf f}\in I(\tilde{G}(F))\otimes Mes(G(F))$.  Alors on a l'\'egalit\'e
  $$I_{\tilde{M}}^{\tilde{G},{\cal E}}(\boldsymbol{\gamma},B,{\bf f})=I_{\tilde{M}}^{\tilde{G}}(\boldsymbol{\gamma},B,{\bf f})$$.}

  {\bf Remarque.} Quand $G$ est quasi-d\'eploy\'e, $\tilde{G}=G$, ${\bf a}=1$ et $B$ est le syst\`eme de fonctions constant de valeur $1$, le (ii) a \'et\'e prouv\'e par Arthur pour $\boldsymbol{\gamma}$ \`a support fortement $\tilde{G}$-r\'egulier. Nous prouverons dans l'article suivant que l'assertion (ii) dans notre situation un peu plus g\'en\'erale se d\'eduit du r\'esultat d'Arthur.
  
  \bigskip

  \section{Germes de Shalika}
  
  \bigskip
  
  \subsection{Germes de Shalika ordinaires}
  
  Soit $(G,\tilde{G},{\bf a})$ un triplet quelconque comme en 1.1. Soit $\eta\in \tilde{G}_{ss}(F)$. Fixons un ensemble de repr\'esentants $\{u_{i};i,...,n\}$ des classes de conjugaison par $Z_{G}(\eta;F)$ dans l'ensemble des \'el\'ements unipotents $u\in G_{\eta}(F)$ tels que $\omega$ est trivial sur $Z_{G}(\eta u;F)$. On fixe des mesures sur tous les groupes intervenant. La th\'eorie des germes de Shalika nous dit que, pour tout $i=1,...,n$, il existe un unique germe $g_{i}(.,\omega)$ de fonction d\'efinie au voisinage de $\eta$ dans $\tilde{G}(F)$, de sorte que, pour tout $f\in C_{c}^{\infty}(\tilde{G}(F))$ et tout $\gamma\in \tilde{G}(F)$, on ait l'\'egalit\'e
  $$I^{\tilde{G}}(\gamma,\omega,f)=\sum_{i=1,...,n}g_{i}(\gamma,\omega)I^{\tilde{G}}(\eta u_{i},\omega,f)$$
  pourvu que $\gamma$ soit assez proche de $\eta$.
  
  Comme on vient de le dire, ces germes sont d\'efinis au voisinage de $\eta$ dans $\tilde{G}(F)$, mais leurs restrictions \`a $\tilde{G}_{reg}(F)$ sont d'un int\'er\^et particulier. On sait que les restrictions des $g_{i}(.,\omega)$ \`a $\tilde{G}_{reg}(F)$ sont homog\`enes, c'est-\`a-dire qu'il existe $d_{i}\in {\mathbb N}$ tel que $$g_{i}(exp(\lambda^2 X)\eta,\omega)=\vert \lambda\vert _{F}^{d_{i}}g_{i}(exp(X)\eta,\omega)$$
  pour tout $X\in \mathfrak{g}_{\eta,reg}(F)$ assez proche de $0$ et tout $\lambda\in F^{\times}$ de valuation positive ou nulle. On sait aussi que les restrictions de ces germes \`a $\tilde{G}_{reg}(F)$  s\'eparent les orbites $u_{i}$. C'est-\`a-dire que l'on peut trouver des familles $(\gamma_{i})_{i=1,...,n}$ form\'ees d'\'el\'ements de $\tilde{G}_{reg}(F)$ aussi proches de $\eta$ que l'on veut, de sorte que la matrice $(g_{i}(\gamma_{j},\omega))_{i,j=1,...,n}$ soit inversible. On peut raffiner ce r\'esultat: fixons un sous-ensemble $\tilde{V}$ ouvert et dense dans $ \tilde{G}_{reg}(F)$; alors on peut imposer aux $\gamma_{i}$ d'appartenir \`a $\tilde{V}$. 
  
  Reformulons les d\'efinitions de fa\c{c}on plus abstraite.  Soit ${\cal O}$ une classe de conjugaison (par $G(F)$) semi-simple. Consid\'erons l'ensemble ${\cal U}({\cal O})$ des voisinages ouverts et ferm\'es $\tilde{U}$ de ${\cal O}$ qui sont invariants par conjugaison et tels que, pour  tout $\gamma\in \tilde{G}(F)$, $\gamma$ appartient \`a $ \tilde{U}$ si et seulement si la partie semi-simple $\gamma_{ss}$ appartient \`a $\tilde{U}$. Pour un tel voisinage $\tilde{U}$, on note $D_{g\acute{e}om}(\tilde{U},\omega)$ le sous-espace des \'el\'ements de $D_{g\acute{e}om}(\tilde{G}(F),\omega)$  \`a support dans $\tilde{U}$. Pour une propri\'et\'e d\'ependant d'un \'el\'ement $\boldsymbol{\gamma}\in D_{g\acute{e}om}(\tilde{G}(F),\omega)$  nous dirons qu'elle est v\'erifi\'ee "pour $\boldsymbol{\gamma}$ assez proche de ${\cal O}$" si et seulement s'il existe $\tilde{U}\in {\cal U}({\cal O})$ tel que la propri\'et\'e soit v\'erifi\'ee pour $\boldsymbol{\gamma}\in D_{g\acute{e}om}(\tilde{U},\omega)$. Soit $E$ un espace vectoriel sur ${\mathbb C}$.  Consid\'erons l'ensemble  des  couples $(\tilde{U},g)$, o\`u $\tilde{U}\in {\cal U}({\cal O})$ et $g: D_{g\acute{e}om}(\tilde{U})\otimes Mes(G(F))^*\to E$ est une application lin\'eaire. Disons que deux couples $(\tilde{U},g)$ et $(\tilde{U}',g')$ sont \'equivalents si et seulement s'il existe $\tilde{U}''\in {\cal U}({\cal O})$, avec $\tilde{U}''\subset \tilde{U}\cap \tilde{U}'$, tel que les restrictions de $g$ et $g'$ \`a  $D_{g\acute{e}om}(\tilde{U}'')\otimes Mes(G(F))^*$ co\"{\i}ncident. Une classe d'\'equivalence sera appel\'ee un germe d'application lin\'eaire sur   
$D_{g\acute{e}om}(\tilde{G}(F))\otimes Mes(G(F))^*$ au voisinage de ${\cal O}$, \`a valeurs dans $E$. Un tel germe sera not\'e simplement $g$ et sera consid\'er\'e comme une application lin\'eaire $D_{g\acute{e}om}(\tilde{G}(F))\otimes Mes(G(F))^*\to E$, dont la valeur n'est bien d\'efinie que pour des \'el\'ements de $D_{g\acute{e}om}(\tilde{G}(F))\otimes Mes(G(F))^*$ assez proches de ${\cal O}$.

On peut reformuler la d\'efinition des germes de Shalika en disant qu'il existe un unique germe d'application lin\'eaire $g_{{\cal O}}$ sur   
$D_{g\acute{e}om}(\tilde{G}(F))\otimes Mes(G(F))^*$ au voisinage de ${\cal O}$, \`a valeurs dans $D_{g\acute{e}om}({\cal O},\omega)\otimes Mes(G(F))^*$, de sorte que, pour tout ${\bf f}\in I(\tilde{G}(F),\omega)\otimes Mes(G(F))$ et tout $\boldsymbol{\gamma}\in D_{g\acute{e}om}(\tilde{G}(F),\omega)\otimes Mes(G(F))^*$, on ait l'\'egalit\'e
$$I^{\tilde{G}}(\boldsymbol{\gamma},{\bf f})=I^{\tilde{G}}(g_{{\cal O}}(\boldsymbol{\gamma}),{\bf f})$$
pourvu que $\boldsymbol{\gamma}$ soit assez proche de ${\cal O}$. La propri\'et\'e de s\'eparation des orbites se traduit de la fa\c{c}on suivante. Soit $\tilde{V}$ un sous-ensemble ouvert et dense dans $\tilde{G}_{reg}(F)$, invariant par conjugaison par $G(F)$. Alors

(1) pour tout $\boldsymbol{\tau}\in D_{g\acute{e}om}({\cal O},\omega)\otimes Mes(G(F))^*$, il existe $\boldsymbol{\gamma}\in D_{g\acute{e}om}(\tilde{G}(F),\omega)\otimes Mes(G(F))^*$, \`a support dans $\tilde{V}$ et aussi proche que l'on veut de ${\cal O}$, de sorte que $g_{{\cal O}}(\boldsymbol{\gamma})=\boldsymbol{\tau}$.

\bigskip

\subsection{Germes de Shalika  et stabilit\'e}

Dans le paragraphe pr\'ec\'edent, on a d\'efini la notion de germe d'application lin\'eaire sur $D_{g\acute{e}om}(\tilde{G}(F),\omega)\otimes Mes(G(F))^*$ au voisinage de ${\cal O}$, \`a valeurs dans un espace vectoriel complexe $E$, quand ${\cal O}$ \'etait une classe de conjugaison semi-simple par $G(F)$. La d\'efinition se g\'en\'eralise au cas o\`u ${\cal O}$ est une r\'eunion finie de telles classes. Dans ce cas, d\'ecomposons ${\cal O}$ en union finie $\cup_{i=1,...,n}{\cal O}_{i}$ de classes de conjugaison. Pour tout $i=1,...,n$ consid\'erons le germe de Shalika $g_{{\cal O}_{i}}$ relatif \`a ${\cal O}_{i}$. On peut fixer pour tout $i$ un voisinage $\tilde{U}_{i}\in {\cal U}({\cal O}_{i})$ de sorte que $g_{{\cal O}_{i}}$ soit d\'efini sur $D_{g\acute{e}om}(\tilde{U}_{i},\omega)\otimes Mes(G(F))^*$. On peut supposer les $\tilde{U}_{i}$ deux \`a deux disjoints. Posons $\tilde{U}=\cup_{i=1,...,n}\tilde{U}_{i}$. Alors $g_{{\cal O}}=\oplus_{i=1,...,n}g_{{\cal O}_{i}}$ est une application lin\'eaire \`a valeurs dans $D_{g\acute{e}om}({\cal O},\omega)\otimes Mes(G(F))^*$, d\'efinie sur
$$\oplus_{i=1,...,n}D_{g\acute{e}om}(\tilde{U}_{i},\omega)\otimes Mes(G(F))^*=D_{g\acute{e}om}(\tilde{U},\omega)\otimes Mes(G(F))^*.$$
Le germe de cette application est uniquement d\'efini.

Supposons maintenant $(G,\tilde{G},{\bf a})$ quasi-d\'eploy\'e et \`a torsion int\'erieure. Soit ${\cal O}\subset \tilde{G}(F)$ une classe de conjugaison stable semi-simple. Pour \'etudier les distributions stables, il convient d'adapter les d\'efinitions en rempla\c{c}ant l'ensemble de voisinages ${\cal U}({\cal O})$ par son sous-ensemble ${\cal U}^{st}({\cal O})$ des $\tilde{U}$ qui v\'erifient la condition: pour $\gamma,\gamma'\in \tilde{G}_{reg}(F)$ stablement conjugu\'es, $\gamma$ appartient \`a $\tilde{U}$ si et seulement si $\gamma'\in \tilde{U}$. Cela ne cr\'ee pas de difficult\'es car tout \'el\'ement de ${\cal U}({\cal O})$ contient un voisinage v\'erifiant cette condition ([I] 4.6). On dispose du germe $g_{{\cal O}}$ ci-dessus.Fixons un sous-ensemble $\tilde{V}\subset \tilde{G}_{reg}(F)$, ouvert et dense et invariant par conjugaison stable.

\ass{Lemme}{Pour $\boldsymbol{\delta}\in D^{st}_{g\acute{e}om}(\tilde{G}(F))\otimes Mes(G(F))^*$, $g_{{\cal O}}(\boldsymbol{\delta})$ appartient \`a $D^{st}_{g\acute{e}om}({\cal O})\otimes Mes(G(F))^*$ pourvu que $\boldsymbol{\delta}$ soit assez proche de ${\cal O}$. Pour tout $\boldsymbol{\tau}\in D^{st}_{g\acute{e}om}({\cal O})\otimes Mes(G(F))^*$, il existe $\boldsymbol{\delta}\in D^{st}_{g\acute{e}om}(\tilde{G}(F))\otimes Mes(G(F))^*$, \`a support dans $\tilde{V}$ et aussi proche que l'on veut de ${\cal O}$, de sorte que $g_{{\cal O}}(\boldsymbol{\delta})=\boldsymbol{\tau}$.}

Preuve. On oublie les espaces de mesures. L'espace $D_{g\acute{e}om}({\cal O})$ est le dual de $I(\tilde{G}(F))_{{\cal O},loc}$, cf. [I] 5.1, tandis que $D^{st}_{g\acute{e}om}({\cal O})$ est le dual de son quotient $SI(\tilde{G}(F))_{{\cal O},loc}$. On peut donc fixer une base $(\boldsymbol{\tau}_{i})_{i=1,...,n}$ de $D_{g\acute{e}om}({\cal O})$, un entier $s\in \{0,...,n\}$ et une famille $(f_{i})_{i=1,...,n}$ d'\'el\'ements de $I(\tilde{G}(F))$ de sorte que

- $(\boldsymbol{\tau}_{i})_{i=1,...,s}$ est une base de $D^{st}_{g\acute{e}om}({\cal O})$  et $(\boldsymbol{\tau}_{i})_{i=s+1,...,n}$ est une base d'un suppl\'ementaire de $D^{st}_{g\acute{e}om}({\cal O})$ dans $D_{g\acute{e}om}({\cal O})$;

- l'image de $(f_{i})_{i=1,...,n}$ dans $I(\tilde{G}(F))_{{\cal O},loc}$ est une base de cet espace et l'image de $(f_{i})_{i=s+1,...,n}$ est une base du noyau de la projection $I(\tilde{G}(F))_{{\cal O},loc}\to SI(\tilde{G}(F))_{{\cal O},loc}$;

- $ I^{\tilde{G}}(\boldsymbol{\tau}_{i},f_{j})=\left\lbrace\begin{array}{cc}1,&\text{ si }i=j,\\ 0,&\text{ si }i\not=j.\\ \end{array}\right.$

Pour $i=s+1,...,n$, la condition que l'image de $f_{i}$ appartient au noyau de la projection dans $ SI(\tilde{G}(F))_{{\cal O},loc}$ signifie que $S^{\tilde{G}}(\boldsymbol{\delta},f_{i})=0$ pour tout $\boldsymbol{\delta}\in D^{st}_{g\acute{e}om}(\tilde{G}(F))$ assez proche de ${\cal O}$. Fixons un voisinage $\tilde{U}\in {\cal U}^{st}({\cal O})$ de sorte que

-  $g_{{\cal O}}$ soit d\'efini sur  $D_{g\acute{e}om}(\tilde{U})$;

- pour tout $i=s+1,...,n$, on ait $S(\boldsymbol{\delta},f_{i})=0$ pour tout $\boldsymbol{\delta}\in D^{st}_{g\acute{e}om}(\tilde{U})$ (cet ensemble \'etant bien s\^ur $D^{st}_{g\acute{e}om}(\tilde{G}(F))\cap D_{g\acute{e}om}(\tilde{U}))$;

 - pour tout $i=1,...,n$, on ait l'\'egalit\'e
$$(1) \qquad I^{\tilde{G}}(\boldsymbol{\gamma},f_{i})=I^{\tilde{G}}(g_{{\cal O}}(\boldsymbol{\gamma}),f_{i})$$
pour $\boldsymbol{\gamma}\in D_{g\acute{e}om}(\tilde{U})$. 

 Soient $\boldsymbol{\delta}\in D^{st}_{g\acute{e}om}(\tilde{U})$ et $i=s+1,...,n$. On a $I^{\tilde{G}}(\boldsymbol{\delta},f_{i})=S^{\tilde{G}}(\boldsymbol{\delta},f_{i})=0$. L'\'egalit\'e ci-dessus nous dit que $I^{\tilde{G}}(g_{{\cal O}}(\boldsymbol{\delta}),f_{i})=0$. La composante de $g_{{\cal O}}(\boldsymbol{\delta})$ sur l'\'el\'ement de base $\boldsymbol{\tau}_{i}$ est donc nulle. C'est la condition pour que $g_{{\cal O}}(\boldsymbol{\delta})$ appartienne \`a $D^{st}_{g\acute{e}om}({\cal O})$, ce qui prouve la premi\`ere assertion de l'\'enonc\'e.
 
 La famille $(f_{i})_{i=1,...,s}$ est lin\'eairement ind\'ependante du noyau de la projection $I(\tilde{G}(F))\to SI(\tilde{G}(F))_{{\cal O},loc}$. Cela entra\^{\i}ne que l'on peut trouver une famille $(\boldsymbol{\delta}_{i})_{i=1,...,s}$ d'\'el\'ements de $D^{st}_{g\acute{e}om,reg}(\tilde{U})$ de sorte que la matrice $(S^{\tilde{G}}(\boldsymbol{\delta}_{i},f_{j}))_{i,j=1,...,s}$ soit inversible. On peut remplacer les $\boldsymbol{\delta}_{i}$ par des \'el\'ements assez proches et \`a support dans $\tilde{V}$. En prenant des combinaisons lin\'eaires convenables de ces \'el\'ements, on obtient une nouvelle famille $(\boldsymbol{\delta}_{i})_{i=1,...,s}$ d'\'el\'ements de $D^{st}_{g\acute{e}om,reg}(\tilde{U})$, \`a supports dans $\tilde{V}$, de sorte que,
 pour $i,j=1,...,s$,
$$S^{\tilde{G}}(\boldsymbol{\delta}_{i},f_{j})=\left\lbrace\begin{array}{cc}1,&\text{ si }i=j,\\ 0,&\text{ si }i\not=j.\\ \end{array}\right.$$
D'apr\`es ce qui pr\'ec\`ede, cette \'egalit\'e vaut m\^eme pour $i=1,...,s$ et $j=1,...,n$.
Appliquons la relation (1) pour $\boldsymbol{\gamma}=\boldsymbol{\delta}_{i}$. On obtient que
$$I^{\tilde{G}}(g_{{\cal O}}(\boldsymbol{\delta}_{i}),f_{j})=\left\lbrace\begin{array}{cc}1,&\text{ si }i=j,\\ 0,&\text{ si }i\not=j\\ \end{array}\right.$$
pour $j=1,...,n$. Donc $g_{{\cal O}}(\boldsymbol{\delta}_{i})=\boldsymbol{\tau}_{i}$. Cela d\'emontre la seconde assertion. $\square$

De nouveau, les d\'efinitions et r\'esultats se g\'en\'eralisent au cas o\`u ${\cal O}$ est une r\'eunion finie de classes de conjugaison stable semi-simples.

\bigskip

\subsection{Int\'egrales orbitales pond\'er\'ees $\omega$-\'equivariantes}

Soit $\tilde{M}$ un espace de Levi de $\tilde{G}$. On note $D_{g\acute{e}om,\tilde{G}-\acute{e}qui}(\tilde{M}(F),\omega)$ le sous-espace des \'el\'ements de $D_{g\acute{e}om}(\tilde{M}(F),\omega)$ dont le support est form\'e d'\'el\'ements de $\tilde{M}(F)$ qui sont $\tilde{G}$-\'equisinguliers. Soit  ${\cal O}\subset \tilde{M}(F)$ une classe de conjugaison semi-simple par $M(F)$.  On note ${\cal O}^{\tilde{G}}\subset \tilde{G}(F)$ l'unique classe de conjugaison par $G(F)$ qui contient ${\cal O}$.  

\ass{Proposition}{ Il existe un unique germe d'application lin\'eaire $g_{\tilde{M},{\cal O}}^{\tilde{G}}$ sur   
$D_{g\acute{e}om,\tilde{G}-\acute{e}qui}(\tilde{M}(F),\omega)\otimes Mes(M(F))^*$ au voisinage de ${\cal O}$, \`a valeurs dans $D_{g\acute{e}om}({\cal O}^{\tilde{G}},\omega)\otimes Mes(G(F))^*$, de sorte que
 pour tout ${\bf f}\in I(\tilde{G}(F),\omega)\otimes Mes(G(F))$ et tout $\boldsymbol{\gamma}\in D_{g\acute{e}om,\tilde{G}-\acute{e}qui}(\tilde{M}(F),\omega)\otimes Mes(M(F))^*$, on ait l'\'egalit\'e
$$I_{\tilde{M}}^{\tilde{G}}(\boldsymbol{\gamma},{\bf f})=\sum_{\tilde{L}\in {\cal L}(\tilde{M})}I_{\tilde{L}}^{\tilde{G}}(g_{\tilde{M},{\cal O}}^{\tilde{L}}(\boldsymbol{\gamma}),{\bf f})$$
pourvu que $\boldsymbol{\gamma}$ soit assez proche de ${\cal O}$.}

C'est une reformulation de  [A3] 2.5 (voir aussi [A1] proposition 9.1). En fait, ces germes se d\'efinissent aussi bien sans imposer \`a $\boldsymbol{\gamma}$ la restriction d'\'equisingularit\'e de son support. Mais nous ne les utiliserons que pour les $\boldsymbol{\gamma}$ indiqu\'es et il est plus simple de se limiter d\`es le d\'ebut \`a de tels \'el\'ements.

En particulier,  si $\tilde{M}=\tilde{G}$, le germe $g_{\tilde{G},{\cal O}}^{\tilde{G}}$ est le germe de Shalika ordinaire   not\'e $g_{{\cal O}}$ dans les paragraphes pr\'ec\'edents. On a une propri\'et\'e suppl\'ementaire:

(1) supposons que ${\cal O}$ soit form\'e d'\'el\'ements $\tilde{G}$-\'equisinguliers de $\tilde{M}(F)$; alors $g_{\tilde{M},{\cal O}}^{\tilde{G}}=0$ si $\tilde{G}\not=\tilde{M}$.

C'est [A1] remarque page 270.

{\bf Variante.} Supposons $(G,\tilde{G},{\bf a})$ quasi-d\'eploy\'e et \`a torsion int\'erieure. Fixons un syst\`eme de fonctions $B$ comme en 1.9. Il y a une proposition similaire \`a la pr\'ec\'edente, o\`u l'on remplace les int\'egrales orbitales $I_{\tilde{L}}^{\tilde{G}}(\boldsymbol{\gamma},{\bf f})$ par leurs variantes $I_{\tilde{L}}^{\tilde{G}}(\boldsymbol{\gamma},B,{\bf f})$. Cette variante se d\'eduit de la proposition pr\'ec\'edente en utilisant 1.9(5). On note $g_{\tilde{M},{\cal O}}^{\tilde{L}}(\boldsymbol{\gamma},B)$ les germes dans cette situation.

{\bf Variante.} Supposons $G=\tilde{G}$ et ${\bf a}=1$. Fixons une fonction $B$ comme en 1.8. Supposons ${\cal O}=\{1\}$. On a d\'efini les int\'egrales orbitales $I_{M}^G(\boldsymbol{\gamma},B,{\bf f})$ pour $\boldsymbol{\gamma}$ \`a support unipotent. On les d\'efinit aussi pour $\boldsymbol{\gamma}$ \`a support $G$-\'equisingulier par la simple \'egalit\'e $I_{M}^G(\boldsymbol{\gamma},B,{\bf f})=I_{M}^G(\boldsymbol{\gamma},{\bf f})$. Alors ces int\'egrales v\'erifient une proposition similaire \`a la pr\'ec\'edente. On note $g_{M,unip}^L(\boldsymbol{\gamma},B)$ les germes dans cette situation.

\bigskip

Comme pr\'ec\'edemment, les d\'efinitions et r\'esultats se g\'en\'eralisent au cas o\`u ${\cal O}$ est une r\'eunion finie de classes de conjugaison  semi-simples.

\bigskip

\subsection{D\'efinition des germes stables}
 On suppose $(G,\tilde{G},{\bf a})$ quasi-d\'eploy\'e et \`a torsion int\'erieure. On fixe un syst\`eme de fonctions $B$ comme en 1.9. Soit $\tilde{M}$ un espace de Levi de $\tilde{G}$. On note $D^{st}_{g\acute{e}om,\tilde{G}-\acute{e}qui}(\tilde{M}(F))$ l'intersection de $D^{st}_{g\acute{e}om}(\tilde{M}(F))$ et de $D_{g\acute{e}om,\tilde{G}-\acute{e}qui}(\tilde{M}(F))$.  Soit ${\cal O}\subset \tilde{M}(F)$ une classe de conjugaison stable d'\'el\'ements semi-simples dans $\tilde{M}(F)$. On note ${\cal O}^{\tilde{G}}\subset \tilde{G}(F)$ l'unique classe de conjugaison stable dans $\tilde{G}(F)$ contenant ${\cal O}$. On va d\'efinir un germe  d'application lin\'eaire $Sg_{\tilde{M},{\cal O}}^{\tilde{G}}(B)$ sur   
$D^{st}_{g\acute{e}om,\tilde{G}-\acute{e}qui}(\tilde{M}(F))\otimes Mes(M(F))^*$ au voisinage de ${\cal O}$, \`a valeurs dans $D_{g\acute{e}om}({\cal O}^{\tilde{G}})\otimes Mes(G(F))^*$.  La proposition ci-dessous affirme qu'il prend en fait ses valeurs dans $D^{st}_{g\acute{e}om}({\cal O}^{\tilde{G}})\otimes Mes(G(F))^*$. Comme toujours, on admet cette propri\'et\'e par r\'ecurrence pour les triplets $(G',\tilde{G}',{\bf a}')$ quasi-d\'eploy\'es et \`a torsion int\'erieure tels que $dim(G'_{SC})< dim(G_{SC})$.
On peut alors poser, pour
$\boldsymbol{\delta}\in D^{st}_{g\acute{e}om,\tilde{G}-\acute{e}qui}(\tilde{M}(F)\otimes Mes(M(F))^*$,  
$$(1) \qquad Sg_{\tilde{M},{\cal O}}^{\tilde{G}}(\boldsymbol{\delta},B)=g_{\tilde{M},{\cal O}}^{\tilde{G}}(\boldsymbol{\delta},B)-\sum_{s\in Z(\hat{M})^{\Gamma_{F}}/Z(\hat{G})^{\Gamma_{F}}; s\not=1} i_{\tilde{M}}(\tilde{G},\tilde{G}'(s))transfert(Sg_{{\bf M},{\cal O}}^{{\bf G}'(s)}(\boldsymbol{\delta},B)).$$
On s'est dispens\'e du formalisme n\'ecessaire pour attacher des germes \`a des donn\'ees ${\bf G}'(s)$ plut\^ot qu'\`a des espaces $\tilde{G}'(s)$. Modulo cet oubli, tous les termes $Sg_{{\bf M},{\cal O}}^{{\bf G}'(s)}(\boldsymbol{\delta},B)$ sont d\'efinis par r\'ecurrence et sont des distributions stables pour ${\bf G}'(s)$.  Remarquons que le support de la distribution $\boldsymbol{\delta}$  est form\'e d'\'el\'ements $\tilde{G}'(s)$-\'equisinguliers: si $\gamma\in \tilde{M}(F)$, le syst\`eme de racines de $G'(s)_{\gamma}$ est contenu dans celui de $G_{\gamma}$, donc dans celui de $M_{\gamma}$ si $\gamma$ est $\tilde{G}$-\'equisingulier. 

\ass{Proposition (\`a prouver)}{ Pour tout $\boldsymbol{\delta}\in D^{st}_{g\acute{e}om,\tilde{G}-\acute{e}qui}(\tilde{M}(F))\otimes Mes(M(F))^*$ assez proche de ${\cal O}$, le terme $ Sg_{\tilde{M},{\cal O}}^{\tilde{G}}(\boldsymbol{\delta},B)$ appartient \`a $D^{st}_{g\acute{e}om}({\cal O}^{\tilde{G}})\otimes Mes(G(F))^*$.}

 Par r\'ecurence, la relation 2.3(1) et la d\'efinition (1) ci-dessus entra\^{\i}nent

(2) si ${\cal O}$ est form\'ee d'\'el\'ements  $\tilde{G}$-\'equisinguliers de $\tilde{M}(F)$, 
$$Sg_{\tilde{M},{\cal O}}^{\tilde{G}}(B)=\left\lbrace\begin{array}{cc}0,& \text{ si } \tilde{M}\not=\tilde{G},\\ g_{\tilde{M},{\cal O}}^{\tilde{M}}(B)=g_{\tilde{M},{\cal O}}^{\tilde{M}},&\text{ si }\tilde{M}=\tilde{G}.\\ \end{array}\right.$$

\bigskip

\subsection{Int\'egrales orbitales pond\'er\'ees invariantes stables}
On conserve les hypoth\`eses du paragraphe pr\'ec\'edent. 
  
   \ass{Proposition}{(i) Pour tout ${\bf f}\in I(\tilde{G}(F))\otimes Mes(G(F))$ et tout $\boldsymbol{\delta}\in D^{st}_{g\acute{e}om,\tilde{G}-\acute{e}qui}(\tilde{M}(F))\otimes Mes(M(F))^*$, on a l'\'egalit\'e
$$S_{\tilde{M}}^{\tilde{G}}(\boldsymbol{\delta},B,{\bf f})=I^{\tilde{G}}(Sg_{\tilde{M},{\cal O}}^{\tilde{G}}(\boldsymbol{\delta},B),B,{\bf f})+\sum_{\tilde{L}\in {\cal L}(\tilde{M}), \tilde{L}\not=\tilde{G}}S_{\tilde{L}}^{\tilde{G}}(Sg_{\tilde{M},{\cal O}}^{\tilde{L}}(\boldsymbol{\delta},B),B,{\bf f})$$
pourvu que $\boldsymbol{\delta}$ soit assez proche de ${\cal O}$.

(ii) Supposons  que $ Sg_{\tilde{M},{\cal O}}^{\tilde{G}}(\boldsymbol{\delta},B)$ appartienne \`a $D^{st}_{g\acute{e}om}({\cal O}^{\tilde{G}})\otimes Mes(G(F))^*$. La formule ci-dessus  devient 
$$S_{\tilde{M}}^{\tilde{G}}(\boldsymbol{\delta},B,{\bf f})=\sum_{\tilde{L}\in {\cal L}(\tilde{M})}S_{\tilde{L}}^{\tilde{G}}(Sg_{\tilde{M},{\cal O}}^{\tilde{L}}(\boldsymbol{\delta},B),B,{\bf f})$$
pourvu que $\boldsymbol{\delta}$ soit assez proche de ${\cal O}$.}

Preuve. Pour simplifier les notations, on oublie les espaces de mesures.  On note $X$ la diff\'erence entre $S_{\tilde{M}}^{\tilde{G}}(\boldsymbol{\delta},B,{\bf f})$  et le membre de droite de l'\'egalit\'e du (i) de l'\'enonc\'e. On utilise la formule de d\'efinition
$$(1) \qquad I_{\tilde{M}}^{\tilde{G}}(\boldsymbol{\delta},B,{\bf f})=S_{\tilde{M}}^{\tilde{G}}(\boldsymbol{\delta},B,{\bf f})+\sum_{s\in Z(\hat{M})^{\Gamma_{F}}/Z(\hat{G})^{\Gamma_{F}};s\not=1}i_{\tilde{M}}(\tilde{G},\tilde{G}'(s))S_{{\bf M}}^{{\bf G}'(s)}(\boldsymbol{\delta},B,{\bf f}^{{\bf G}'(s)}).$$
On d\'eveloppe le membre de gauche par la proposition 2.3.  On \'ecrit le premier terme du membre de droite comme $X$ plus le membre de droite de la formule du (i) de l'\'enonc\'e. On d\'eveloppe les autres termes en utilisant la formule du (ii) de l'\'enonc\'e, que l'on peut utiliser par r\'ecurrence.  Toutes ces formules sont des sommes sur des ensembles d'espaces de Levi. Rappelons que, pour $s\in Z(\hat{M})^{\Gamma_{F}}/Z(\hat{G})^{\Gamma_{F}}$, un espace de Levi $\tilde{L}'_{s}\in {\cal L}^{\tilde{G}'(s)}(\tilde{M})$ d\'etermine un espace de Levi $\tilde{L}$ de $\tilde{G}$ par la formule ${\cal A}_{\tilde{L}}={\cal A}_{L'_{s}}$. On a tacitement r\'ealis\'e $^LM$ comme un sous-groupe de Levi standard de $^LG$. En posant ${\cal L}'_{s}=\hat{L}'_{s}\rtimes W_{F}\subset {^LG}$, le triplet $(L'_{s},{\cal L}'_{s},s)$ est \'egal \`a la donn\'ee endoscopique ${\bf L}'(s)$ de $(L,\tilde{L})$.  On isole dans chaque formule l'espace de Levi maximal. Ainsi, pourvu que $\boldsymbol{\delta}$ soit assez proche de ${\cal O}$, $I_{\tilde{M}}^{\tilde{G}}(\boldsymbol{\delta},B,{\bf f})$ est la somme de

(2) $I^{\tilde{G}}(g_{\tilde{M},{\cal O}}^{\tilde{G}}(\boldsymbol{\delta},B),{\bf f})$

et de 

(3) $\sum_{\tilde{L}\in {\cal L}(\tilde{M}),\tilde{L}\not=\tilde{G}}I_{\tilde{L}}^{\tilde{G}}(g_{\tilde{M},{\cal O}}^{\tilde{L}}(\boldsymbol{\delta},B),B,{\bf f})$.

Le membre de droite de (1) est \'egal \`a la somme de $X$, de

 $$(4) \qquad I^{\tilde{G}}(Sg_{\tilde{M},{\cal O}}^{\tilde{G}}(\boldsymbol{\delta},B),{\bf f})+\sum_{s\in Z(\hat{M})^{\Gamma_{F}}/Z(\hat{G})^{\Gamma_{F}};s\not=1}i_{\tilde{M}}(\tilde{G},\tilde{G}'(s))S^{{\bf G}'(s)}(Sg_{{\bf M},{\cal O}}^{{\bf G}'(s)}(\boldsymbol{\delta},B),{\bf f}^{{\bf G}'(s)}),$$
 et de
 
 $$(5)\qquad \sum_{s\in Z(\hat{M})^{\Gamma_{F}}/Z(\hat{G})^{\Gamma_{F}}}i_{\tilde{M}}(\tilde{G},\tilde{G}'(s))\sum_{\tilde{L}'_{s}\in {\cal L}^{\tilde{G}'(s)}(\tilde{M}); \tilde{L}'_{s}\not=\tilde{G}'(s)}S_{{\bf L}'(s)}^{{\bf G}'(s)}(Sg_{{\bf M},{\cal O}}^{{\bf L}'(s)}(\boldsymbol{\delta},B),B,{\bf f}^{{\bf G}'(s)}),$$
 \'etant entendu que, si $s=1$,  $S_{{\bf L}'(s)}^{{\bf G}'(s)}(Sg_{{\bf M},{\cal O}}^{{\bf L}'(s)}(\boldsymbol{\delta},B),B,{\bf f}^{{\bf G}'(s)})=S_{\tilde{L}}^{\tilde{G}}(Sg_{\tilde{M},{\cal O}}
^{\tilde{L}}(\boldsymbol{\delta},B),B,{\bf f})$.
 
 On veut prouver que $X=0$. Il suffit pour cela de prouver que (2) est \'egal \`a (4) et que (3) est \'egal \`a (5). Par d\'efinition du transfert, on peut remplacer dans (4) les termes $S^{{\bf G}'(s)}(Sg_{{\bf M},{\cal O}}^{{\bf G}'(s)}(\boldsymbol{\delta},B),B,{\bf f}^{{\bf G}'(s)})$ par $I^{\tilde{G}}(transfert(Sg_{{\bf M},{\cal O}}^{{\bf G}'(s)}(\boldsymbol{\delta},B)),{\bf f})$. Alors l'\'egalit\'e de (2) et (4) r\'esulte de la d\'efinition (1) de 2.4. 
 
 Dans (5), on regroupant les couples $(s,\tilde{L}'_{s})$ selon l'espace de Levi $\tilde{L}$  associ\'e et  l'image de $s$ dans $Z(\hat{M})^{\Gamma_{F}}/Z(\hat{L})^{\Gamma_{F}}$ (le triplet ${\bf L}'(s)$ ne d\'epend que de cette image). L'expression (5) devient
 $$\sum_{\tilde{L}\in {\cal L}(\tilde{M}), \tilde{L}\not=\tilde{G}}\sum_{s\in Z(\hat{M})^{\Gamma_{F}}/Z(\hat{L})^{\Gamma_{F}}, {\bf L}'(s) \text{ elliptique}}\sum_{t\in sZ(\hat{L})^{\Gamma_{F}}/Z(\hat{G})^{\Gamma_{F}}}i_{\tilde{M}}(\tilde{G},\tilde{G}'(t))$$
 $$S_{{\bf L}'(s)}^{{\bf G}'(t)}(Sg_{{\bf M},{\cal O}}^{{\bf L}'(s)}(\boldsymbol{\delta},B),B,f^{{\bf G}'(t)}).$$
 Soient $s$ et $t$ intervenant dans cette formule. On v\'erifie l'\'egalit\'e
 $$i_{\tilde{M}}(\tilde{G},\tilde{G}'(t))=i_{\tilde{M}}(\tilde{L},\tilde{L}'(s))i_{\tilde{L}'(s)}(\tilde{G},\tilde{G}'(t)).$$
 La non nullit\'e du terme de droite implique que ${\bf L}'(s)$ est elliptique. L'expression (5) devient
 $$(6) \qquad \sum_{\tilde{L}\in {\cal L}(\tilde{M}), \tilde{L}\not=\tilde{G}}\sum_{s\in Z(\hat{M})^{\Gamma_{F}}/Z(\hat{L})^{\Gamma_{F}}}i_{\tilde{M}}(\tilde{L},\tilde{L}'(s))$$
 $$\sum_{t\in sZ(\hat{L})^{\Gamma_{F}}/Z(\hat{G})^{\Gamma_{F}}}i_{\tilde{L}'(s)}(\tilde{G},\tilde{G}'(t))S_{{\bf L}'(s)}^{{\bf G}'(t)}(Sg_{{\bf M},{\cal O}}^{{\bf L}'(s)}(\boldsymbol{\delta},B),B,{\bf f}^{{\bf G}'(t)}).$$
 Consid\'erons la contribution de $\tilde{L}=\tilde{M}$ (qui n'intervient que si $\tilde{M}\not=\tilde{G}$). C'est simplement
 $$\sum_{t\in Z(\hat{M})^{\Gamma_{F}}/Z(\hat{G})^{\Gamma_{F}}}i_{\tilde{M}}(\tilde{G},\tilde{G}'(t))S_{{\bf M}}^{{\bf G}'(t)}(Sg_{{\bf M},{\cal O}}^{{\bf M}}(\boldsymbol{\delta},B),B,{\bf f}^{{\bf G}'(t)}),$$
 avec la m\^eme convention que plus haut si $t=1$. Par d\'efinition de $S_{\tilde{M}}^{\tilde{G}}(Sg_{\tilde{M},{\cal O}}^{\tilde{M}}(\boldsymbol{\delta},B),B,{\bf f})$, ceci n'est autre que $I_{\tilde{M}}^{\tilde{G}}(Sg_{\tilde{M},{\cal O}}^{\tilde{M}}(\boldsymbol{\delta},B),B,{\bf f})$. Consid\'erons maintenant la contribution \`a (6) d'un $\tilde{L}\not=\tilde{M}$.
Par d\'efinition, la somme int\'erieure en $t$ n'est autre que $I_{\tilde{L}}^{\tilde{G},{\cal E}}({\bf L}'(s),Sg_{{\bf M},{\cal O}}^{{\bf L}'(s)}(\boldsymbol{\delta},B),B,{\bf f})$, ou encore $I_{\tilde{L}}^{\tilde{G},{\cal E}}(transfert(Sg_{{\bf M},{\cal O}}^{{\bf L}'(s)}(\boldsymbol{\delta},B)),B,{\bf f})$.   Parce que $\tilde{L}\not=\tilde{M}$, nos hypoth\`eses de r\'ecurrence nous autorisent \`a appliquer  le th\'eor\`eme 1.16. Il nous dit que le terme ci-dessus est aussi \'egal \`a 
$I_{\tilde{L}}^{\tilde{G}}(transfert(Sg_{{\bf M},{\cal O}}^{{\bf L}'(s)}(\boldsymbol{\delta},B)),B,{\bf f})$. L'expression (6) devient
$$\sum_{\tilde{L}\in {\cal L}(\tilde{M}), \tilde{L}\not=\tilde{G}}\sum_{s\in Z(\hat{M})^{\Gamma_{F}}/Z(\hat{L})^{\Gamma_{F}}}i_{\tilde{M}}(\tilde{L},\tilde{L}'(s))I_{\tilde{L}}^{\tilde{G}}(transfert(Sg_{{\bf M},{\cal O}}^{{\bf L}'(s)}(\boldsymbol{\delta},B)),B,{\bf f}).$$
En utilisant la d\'efinition (1) de 2.4 avec $\tilde{G}$ remplac\'e par $\tilde{L}$, la somme int\'erieure en $s$ devient $I_{\tilde{L}}^{\tilde{G}}(g_{\tilde{M},{\cal O}}^{\tilde{L}}(\boldsymbol{\delta},B),B,{\bf f})$ et (6) est \'egal \`a (3). Cela d\'emontre le (i) de l'\'enonc\'e. Le (ii) est imm\'ediat.

   $\square$

\bigskip

\subsection{D\'eveloppement en germes d'int\'egrales orbitales pond\'er\'ees $\omega$-\'equivariantes endoscopiques}

On revient au cas g\'en\'eral. Soient $\tilde{M}$ un espace de Levi de $\tilde{G}$ et ${\cal O}$ une classe de conjugaison stable semi-simple dans $\tilde{M}(F)$.   

\ass{Proposition}{ Il existe un unique germe d'application lin\'eaire $g_{\tilde{M},{\cal O}}^{\tilde{G},{\cal E}}$ sur   
$D_{g\acute{e}om,\tilde{G}-\acute{e}qui}(\tilde{M}(F),\omega)\otimes Mes(M(F))^*$ au voisinage de ${\cal O}$, \`a valeurs dans $D_{g\acute{e}om}({\cal O}^{\tilde{G}},\omega)\otimes Mes(G(F))^*$, de sorte que, pour tout ${\bf f}\in I(\tilde{G}(F),\omega)\otimes Mes(G(F))$ et tout $\boldsymbol{\gamma}\in D_{g\acute{e}om,\tilde{G}-\acute{e}qui}(\tilde{M}(F),\omega)\otimes Mes(M(F))^*$, on a l'\'egalit\'e
$$I_{\tilde{M}}^{\tilde{G},{\cal E}}(\boldsymbol{\gamma},{\bf f})=\sum_{\tilde{L}\in {\cal L}(\tilde{M})}I_{\tilde{L}}^{\tilde{G},{\cal E}}(g_{\tilde{M},{\cal O}}^{\tilde{L},{\cal E}}(\boldsymbol{\gamma}),{\bf f})$$
pourvu que $\boldsymbol{\gamma}$ soit assez proche de ${\cal O}$.}

Avant de d\'emontrer cette proposition, posons une d\'efinition.  Soit  ${\bf M}'=(M',{\cal M}',\tilde{\zeta})$ une donn\'ee endoscopique elliptique et relevante  de $(M,\tilde{M},{\bf a}_{M})$ et soit $\delta\in \tilde{M}'(F)$. Notons $\epsilon$ la partie semi-simple de $\delta$. Il correspond \`a $\epsilon$ une classe de conjugaison  par $M(\bar{F})$ dans $\tilde{M}(\bar{F})$. On dit que $\delta$ est $\tilde{G}$-\'equisingulier si cette classe est form\'ee d'\'el\'ements $\tilde{G}$-\'equisinguliers. On a

(1) si $\delta$ est $\tilde{G}$-\'equisingulier,  $\delta$ est $\tilde{G}'(\tilde{s})$-\'equisingulier pour tout $\tilde{s}\in \tilde{\zeta}Z(\hat{M})^{\Gamma_{F},\hat{\theta}}/Z(\hat{G})^{\Gamma_{F},\hat{\theta}}$. 

En effet, soit $\eta\in \tilde{M}(\bar{F})$ dans la classe de conjugaison correspondant \`a $\epsilon$. On sait d\'ecrire les syst\`emes de racines de $G_{\eta}$, $M_{\eta}$, $G'(\tilde{s})_{\epsilon}$ et $M'_{\epsilon}$, cf. [W2] 3.3. Les syst\`emes de racines de $G_{\eta}$ et $G'(\tilde{s})_{\epsilon}$, resp. $M_{\eta}$ et $M'_{\epsilon}$ ne sont pas \'egaux en g\'en\'eral, mais sont en bijection, les deux bijections \'etant compatibles. Il en r\'esulte que, si $G_{\eta}=M_{\eta}$, on a aussi $G'(\tilde{s})_{\epsilon}=M'_{\epsilon}$. $\square$

On note $D^{st}_{g\acute{e}om,\tilde{G}-\acute{e}qui}({\bf M}')$ le sous-espace des \'el\'ements de $D^{st}_{g\acute{e}om}({\bf M}')$ dont le support est form\'e d'\'el\'ements $\tilde{G}$-\'equisinguliers.
\bigskip

Preuve de la proposition. Comme toujours, oublions les espaces de mesures. Supposons d'abord que $(G,\tilde{G},{\bf a})$ n'est pas quasi-d\'eploy\'e et \`a torsion int\'erieure. Consid\'erons d'abord une donn\'ee endoscopique elliptique et relevante ${\bf M}'=(M',{\cal M}',\tilde{\zeta})$ de $(M,\tilde{M},{\bf a}_{M})$. Il correspond \`a ${\cal O}$ une r\'eunion ${\cal O}_{\tilde{M}'}$ de classes de conjugaison stable  semi-simple dans $\tilde{M}'(F)$.   On d\'efinit un germe d'application lin\'eaire $\boldsymbol{\delta}\mapsto g_{\tilde{M},{\cal O}}^{\tilde{G},{\cal E}}({\bf M}',\boldsymbol{\delta})$ sur $D_{g\acute{e}om,\tilde{G}-\acute{e}qui}^{st}({\bf M}')$ au voisinage de ${\cal O}_{\tilde{M}'}$, \`a valeurs dans $D_{g\acute{e}om}({\cal O}^{\tilde{G}},\omega)$, par la formule
$$(2)\qquad g_{\tilde{M},{\cal O}}^{\tilde{G},{\cal E}}({\bf M}',\boldsymbol{\delta})=\sum_{\tilde{s}\in \tilde{\zeta}Z(\hat{M})^{\Gamma_{F},\hat{\theta}}/Z(\hat{G})^{\Gamma_{F},\hat{\theta}}}i_{\tilde{M}'}(\tilde{G},\tilde{G}'(\tilde{s}))transfert(Sg_{{\bf M}',{\cal O}_{\tilde{M}'}}^{{\bf G}'(\tilde{s})}(\boldsymbol{\delta},B^{\tilde{G}})).$$
Les hypoth\`eses  de r\'ecurrence assurent que, pour tout $\tilde{s}$,  la proposition 2.4  est v\'erifi\'ee pour ${\bf G}'(\tilde{s})$. Donc les termes $Sg_{{\bf M}',{\cal O}_{\tilde{M}'}}^{{\bf G}'(\tilde{s})}(\boldsymbol{\delta},B^{\tilde{G}})$ sont bien d\'efinis et sont stables.
Avec cette d\'efinition, montrons que l'on a l'\'egalit\'e
$$(3)\qquad I_{\tilde{M}}^{\tilde{G},{\cal E}}({\bf M}',\boldsymbol{\delta},{\bf f})=\sum_{\tilde{L}\in {\cal L}(\tilde{M})}I_{\tilde{L}}^{\tilde{G},{\cal E}}(g_{\tilde{M},{\cal O}}^{\tilde{L},{\cal E}}({\bf M}',\boldsymbol{\delta}),{\bf f})$$
pourvu que $\boldsymbol{\delta}$ soit assez proche de ${\cal O}_{\tilde{M}'}$.

La preuve est similaire \`a celle de la proposition pr\'ec\'edente. Faisons-la rapidement. On a par d\'efinition
$$I_{\tilde{M}}^{\tilde{G},{\cal E}}({\bf M}',\boldsymbol{\delta},{\bf f})=\sum_{\tilde{s}\in \tilde{\zeta}Z(\hat{M})^{\Gamma_{F},\hat{\theta}}/Z(\hat{G})^{\Gamma_{F},\hat{\theta}}}i_{\tilde{M}'}(\tilde{G},\tilde{G}'(\tilde{s}))S_{{\bf M}'}^{{\bf G}'(\tilde{s})}(\boldsymbol{\delta},B^{\tilde{G}},{\bf f}^{{\bf G}'(\tilde{s})}).$$
En utilisant la proposition 2.5, on obtient
$$I_{\tilde{M}}^{\tilde{G},{\cal E}}({\bf M}',\boldsymbol{\delta},{\bf f})=\sum_{\tilde{s}\in \tilde{\zeta}Z(\hat{M})^{\Gamma_{F},\hat{\theta}}/Z(\hat{G})^{\Gamma_{F},\hat{\theta}}}i_{\tilde{M}'}(\tilde{G},\tilde{G}'(\tilde{s}))\sum_{\tilde{L}'_{\tilde{s}}\in {\cal L}^{\tilde{G}'(\tilde{s})}(\tilde{M}')}S_{{\bf L}'_{\tilde{s}}}^{{\bf G}'(\tilde{s})}(Sg_{{\bf M}',{\cal O}_{\tilde{M}'}}^{{\bf L}'_{\tilde{s}}}(\boldsymbol{\delta},B^{\tilde{G}}),B^{\tilde{G}},{\bf f}^{{\bf G}'(\tilde{s})})$$
pouvu que $\boldsymbol{\delta}$ soit assez proche de ${\cal O}_{\tilde{M}'}$. On regroupe les couples $(\tilde{s},\tilde{L}_{\tilde{s}})$ selon l'espace de Levi $\tilde{L}$ de $\tilde{G}$ d\'etermin\'e par l'\'egalit\'e ${\cal A}_{\tilde{L}}={\cal A}_{\tilde{L}'_{\tilde{s}}}$. On obtient
$$I_{\tilde{M}}^{\tilde{G},{\cal E}}({\bf M}',\boldsymbol{\delta},{\bf f})=\sum_{\tilde{L}\in {\cal L}(\tilde{M})}\sum_{\tilde{s}\in \tilde{\zeta}Z(\hat{M})^{\Gamma_{F},\hat{\theta}}/Z(\hat{L})^{\Gamma_{F},\hat{\theta}}, {\bf L}'(\tilde{s})\text{ elliptique}}$$
$$\sum_{\tilde{t}\in \tilde{s}Z(\hat{L})^{\Gamma_{F},\hat{\theta}}/Z(\hat{G})^{\Gamma_{F},\hat{\theta}}}i_{\tilde{M}'}(\tilde{G},\tilde{G}'(\tilde{t}))S_{{\bf L}'(\tilde{s})}^{{\bf G}'(\tilde{t})}(Sg_{{\bf M}',{\cal O}_{\tilde{M}'}}^{{\bf L}'(\tilde{s})}(\boldsymbol{\delta},B^{\tilde{G}}),B^{\tilde{G}},{\bf f}^{{\bf G}'(\tilde{s})}).$$
On a encore l'\'egalit\'e 
$$i_{\tilde{M}'}(\tilde{G},\tilde{G}'(\tilde{t}))=i_{\tilde{M}'}(\tilde{L},\tilde{L}'(\tilde{s}))i_{\tilde{L}'(\tilde{s})}(\tilde{G},\tilde{G}'(\tilde{t}))$$
et l'expression ci-dessus se transforme en
$$(4) \qquad I_{\tilde{M}}^{\tilde{G},{\cal E}}({\bf M}',\boldsymbol{\delta},{\bf f})=\sum_{\tilde{L}\in {\cal L}(\tilde{M})}\sum_{\tilde{s}\in \tilde{\zeta}Z(\hat{M})^{\Gamma_{F},\hat{\theta}}/Z(\hat{L})^{\Gamma_{F},\hat{\theta}}}i_{\tilde{M}'}(\tilde{L},\tilde{L}'(\tilde{s}))$$
$$\sum_{\tilde{t}\in \tilde{s}Z(\hat{L})^{\Gamma_{F},\hat{\theta}}/Z(\hat{G})^{\Gamma_{F},\hat{\theta}}}i_{\tilde{L}'(\tilde{s})}(\tilde{G},\tilde{G}'(\tilde{t}))S_{{\bf L}'(\tilde{s})}^{{\bf G}'(\tilde{t})}(Sg_{{\bf M}',{\cal O}_{\tilde{M}'}}^{{\bf L}'(\tilde{s})}(\boldsymbol{\delta},B^{\tilde{G}}),B^{\tilde{G}},{\bf f}^{{\bf G}'(\tilde{s})}).$$
La somme int\'erieure en $\tilde{t}$ n'est autre que $I_{\tilde{L}}^{\tilde{G},{\cal E}}({\bf L}'(\tilde{s}),Sg_{{\bf M}',{\cal O}_{\tilde{M}'}}^{{\bf L}'(\tilde{s})}(\boldsymbol{\delta},B^{\tilde{G}}),{\bf f})$, ou encore 
$$I_{\tilde{L}}^{\tilde{G},{\cal E}}(transfert(Sg_{{\bf M}',{\cal O}_{\tilde{M}'}}^{{\bf L}'(\tilde{s})}(\boldsymbol{\delta},B^{\tilde{G}})),{\bf f}).$$ 
En utilisant la d\'efinition (2), (4) devient l'\'egalit\'e (3).
 
 Soit maintenant $\boldsymbol{\gamma}\in D_{g\acute{e}om}(\tilde{M}(F),\omega)$ assez proche de ${\cal O}$. On peut \'ecrire 
 $$(5)\qquad \boldsymbol{\gamma}=\sum_{{\bf M}'\in {\cal E}(\tilde{M},{\bf a})}transfert(\boldsymbol{\delta}_{{\bf M}'}),$$
 avec des $\boldsymbol{\delta}_{{\bf M}'}\in D^{st}_{g\acute{e}om}({\bf M}')$, cf. [I] proposition 5.7. On peut supposer que, pour tout ${\bf M}'$ et tout \'el\'ement $\delta$ du support $\boldsymbol{\delta}_{{\bf M}'}$, il existe un \'el\'ement $\gamma$ du support de $\boldsymbol{\gamma}$ de sorte que sa partie semi-simple $\eta$ appartienne \`a la classe de conjugaison  dans $\tilde{M}(\bar{F})$ associ\'ee \`a la partie semi-simple $\epsilon$ de $\delta$. Un tel $\delta$ est $\tilde{G}$-\'equisingulier et proche de ${\cal O}^{\tilde{M}'}$ si $\boldsymbol{\gamma}$ est assez proche de ${\cal O}$. 
 Cela \'etant, on pose
 $$(6) \qquad g_{\tilde{M},{\cal O}}^{\tilde{G},{\cal E}}(\boldsymbol{\gamma})=\sum_{{\bf M}'\in {\cal E}(\tilde{M},{\bf a})}g_{\tilde{M},{\cal O}}^{\tilde{G},{\cal E}}({\bf M}',\boldsymbol{\delta}_{{\bf M}'}).$$
 Pour que cette d\'efinition soit loisible, il faut montrer que:
 
 (7)  ce terme ne d\'epend pas de la d\'ecomposition (5) choisie.
 
Fixons cette d\'ecomposition. Par d\'efinition, on a
 $$I_{\tilde{M}}^{\tilde{G},{\cal E}}(\boldsymbol{\gamma},{\bf f})=\sum_{{\bf M}'\in {\cal E}(\tilde{M},{\bf a})}I_{\tilde{M}}^{\tilde{G},{\cal E}}({\bf M}',\boldsymbol{\delta}_{{\bf M}'},{\bf f}).$$
 En utilisant (3), on obtient
 $$I_{\tilde{M}}^{\tilde{G},{\cal E}}(\boldsymbol{\gamma},{\bf f})=\sum_{{\bf M}'\in {\cal E}(\tilde{M},{\bf a})}\sum_{\tilde{L}\in {\cal L}(\tilde{M})}I_{\tilde{L}}^{\tilde{G},{\cal E}}(g_{\tilde{M},{\cal O}}^{\tilde{L},{\cal E}}({\bf M}',\boldsymbol{\delta}_{{\bf M}'}),{\bf f})$$
 $$=\sum_{L\in {\cal L}(\tilde{M})}I_{\tilde{L}}^{\tilde{G},{\cal E}}(g_{\tilde{M},{\cal O}}^{\tilde{L},{\cal E}}(\boldsymbol{\gamma}),{\bf f}).$$
 En raisonnant par r\'ecurrence, on peut supposer que, pour $\tilde{L}\not=\tilde{G}$, le germe $g_{\tilde{M},{\cal O}}^{\tilde{L},{\cal E}}(\boldsymbol{\gamma})$ ne d\'epend pas de la d\'ecomposition (5). L'unique terme restant, \`a savoir le terme pour $\tilde{L}=\tilde{G}$, n'en d\'epend donc pas non plus. Cela d\'emontre l'assertion (7) et en m\^eme temps l'\'egalit\'e de l'\'enonc\'e. 
 
 Supposons maintenant que $(G,\tilde{G},{\bf a})$ soit quasi-d\'eploy\'e et \`a torsion int\'erieure. Dans le raisonnement pr\'ec\'edent, l'unique probl\`eme qui se pose est qu'on ne conna\^{\i}t pas la stabilit\'e de l'un des germes que l'on manipule.  Il s'agit du germe $Sg_{{\bf M},{\cal O}}^{{\bf G}}$. Mais il n'intervient que dans $S^{{\bf G}}(Sg_{{\bf M},{\cal O}}^{{\bf G}}(\boldsymbol{\delta}),{\bf f}^{{\bf G}})$. Il suffit de remplacer cette expression par $I^{\tilde{G}}(Sg_{\tilde{M},{\cal O}}^{\tilde{G}}(\boldsymbol{\delta}),{\bf f})$ et la d\'emonstration s'applique.  $\square$
 
 Supposons que ${\cal O}$ soit form\'e d'\'el\'ements $\tilde{G}$-\'equisinguliers de $ \tilde{M}(F)$. Alors, pour ${\bf M}'=(M',{\cal M}',\tilde{\zeta})\in {\cal E}(\tilde{M},{\bf a})$, l'ensemble  $ {\cal O}_{\tilde{M}'}$ est form\'e d'\'el\'ements $\tilde{G}$-\'equisinguliers.   Les d\'efinitions (2) et (6) et  la relation 2.4(2) entra\^{\i}nent alors  que $g_{\tilde{M},{\cal O}}^{\tilde{G},{\cal E}}(\boldsymbol{\gamma})=0$ au voisinage de ${\cal O}$ si $\tilde{M}\not=\tilde{G}$, tandis que, si $\tilde{M}=\tilde{G}$, on a l'\'egalit\'e
 $$g_{\tilde{G},{\cal O}}^{\tilde{G},{\cal E}}(\boldsymbol{\gamma})=\sum_{{\bf G}'\in {\cal E}(\tilde{G},{\bf a})}transfert(g_{{\bf G}',{\cal O}_{\tilde{G}'}}^{{\bf G}'}(\boldsymbol{\delta}_{{\bf G}'})),$$
 avec les notations de (5) adapt\'ees au cas $\tilde{M}=\tilde{G}$.  Pour ${\bf G}'$ apparaissant ci-dessus et pour $f\in C_{c}^{\infty}(\tilde{G}(F))$, on a les \'egalit\'es
 $$I^{\tilde{G}}(transfert(g_{{\bf G}',{\cal O}_{\tilde{G}'}}^{{\bf G}'}(\boldsymbol{\delta}_{{\bf G}'})),f)=S^{{\bf G}'}(g_{{\bf G}',{\cal O}_{\tilde{G}'}}^{{\bf G}'}(\boldsymbol{\delta}_{{\bf G}'}),f^{{\bf G}'})=S^{{\bf G}'}(\boldsymbol{\delta}_{{\bf G}'},f^{{\bf G}'})$$
 $$=I^{\tilde{G}}(transfert(\boldsymbol{\delta}_{{\bf G}'}),f)=I^{\tilde{G}}(g_{\tilde{G},{\cal O}}^{\tilde{G}}(transfert(\boldsymbol{\delta}_{{\bf G}'}),f),$$
  pourvu que $\boldsymbol{\gamma}$ soit assez proche  de ${\cal O}$. Donc 
  $$transfert(g_{{\bf G}',{\cal O}_{\tilde{G}'}}^{{\bf G}'}(\boldsymbol{\delta}_{{\bf G}'}))=g_{\tilde{G},{\cal O}}^{\tilde{G}}(transfert(\boldsymbol{\delta}_{{\bf G}'})).$$
  La formule plus haut devient
 
 (8) si ${\cal O}$ est form\'e d'\'el\'ements $\tilde{G}$-\'equisinguliers de $ \tilde{M}(F)$, 
 $$g_{\tilde{M},{\cal O}}^{\tilde{G},{\cal E}}=\left\lbrace\begin{array}{cc}0,&\text{ si }\tilde{M}\not=\tilde{G},\\ g_{\tilde{M},{\cal O}}^{\tilde{M}},&\text{ si }\tilde{M}=\tilde{G}.\\ \end{array}\right.$$

 {\bf Variante.} Supposons $(G,\tilde{G},{\bf a})$ quasi-d\'eploy\'e et \`a torsion int\'erieure. Fixons un syst\`eme de fonctions $B$ comme en 1.9. Il y a une proposition similaire \`a celle ci-dessus concernant les distributions $I_{\tilde{M}}^{\tilde{G},{\cal E}}(\boldsymbol{\gamma},B,{\bf f})$. On note $g_{\tilde{M},{\cal O}}^{\tilde{G},{\cal E}}(\boldsymbol{\gamma},B)$ les germes correspondants.
 
 {\bf Variante.} Supposons $G=\tilde{G}$, ${\bf a}=1$ et ${\cal O}=\{1\}$. Fixons une fonction $B$ comme en 1.8. Il y a une proposition similaire \`a celle ci-dessus concernant les distributions $I_{M}^{G,{\cal E}}(\boldsymbol{\gamma},B,{\bf f})$. On note $g_{M,unip}^{G,{\cal E}}(\boldsymbol{\gamma},B)$ les germes correspondants.

 \bigskip
 
 \subsection{Une \'egalit\'e de germes}
 
 On conserve les m\^emes donn\'ees que dans le paragraphe pr\'ec\'edent. 
 
 \ass{Proposition (\`a prouver)}{Sous les hypoth\`eses ci-dessus, on a l'\'egalit\'e $g_{\tilde{M},{\cal O}}^{\tilde{G}}=g_{\tilde{M},{\cal O}}^{\tilde{G},{\cal E}}$.}
 
\bigskip

\subsection{Relation entre la proposition 2.7 et le th\'eor\`eme 1.16 } 

On conserve les m\^emes donn\'ees.   Toutes nos assertions sont tautologiques dans le cas $\tilde{M}=\tilde{G}$, on suppose donc ici $\tilde{M}\not=\tilde{G}$. Soient $\boldsymbol{\gamma}\in D_{g\acute{e}om,\tilde{G}-\acute{e}qui}(\tilde{M}(F),\omega)\otimes Mes(M(F))^*$ et ${\bf f}\in I(\tilde{G}(F),\omega)\otimes Mes(G(F))$. Supposons $\boldsymbol{\gamma}$ assez proche de ${\cal O}$. Consid\'erons les d\'eveloppements des propositions  2.3(ii) et 2.6(ii) et faisons leur diff\'erence. Soit $\tilde{L}\in {\cal L}(\tilde{M})$. Si $\tilde{L}\not=\tilde{G}$, nos hypoth\`eses de r\'ecurrence permettent d'appliquer la proposition ci-dessus: on a $g_{\tilde{M},{\cal O}}^{\tilde{L}}=g_{\tilde{M},{\cal O}}^{\tilde{L},{\cal E}}$. Si $\tilde{L}\not=\tilde{M} $, ces hypoth\`eses permettent d'appliquer le th\'eor\`eme 1.16: on a $I_{\tilde{L}}^{\tilde{G}}(\boldsymbol{\gamma}',{\bf f})
=I_{\tilde{L}}^{\tilde{G},{\cal E}}(\boldsymbol{\gamma}',{\bf f})$ pour tout $\boldsymbol{\gamma}'$. On obtient
$$(1) \qquad I_{\tilde{M}}^{\tilde{G}}(\boldsymbol{\gamma},{\bf f})-I_{\tilde{M}}^{\tilde{G},{\cal E}}(\boldsymbol{\gamma},{\bf f})=I_{\tilde{M}}^{\tilde{G}}(g_{\tilde{M},{\cal O}}^{\tilde{M}}(\boldsymbol{\gamma}),{\bf f})-I_{\tilde{M}}^{\tilde{G},{\cal E}}(g_{\tilde{M},{\cal O}}^{\tilde{M}}(\boldsymbol{\gamma}),{\bf f})$$
$$+I^{\tilde{G}}(g_{\tilde{M},{\cal O}}^{\tilde{G}}(\boldsymbol{\gamma})-g_{\tilde{M},{\cal O}}^{\tilde{G},{\cal E}}(\boldsymbol{\gamma}),{\bf f}).$$

Soit ${\cal D}$ un sous-ensemble de  $D_{g\acute{e}om,\tilde{G}-\acute{e}qui}(\tilde{M}(F),\omega)\otimes Mes(M(F))^*$ v\'erifiant la propri\'et\'e suivante

(2) pour tout $\boldsymbol{\tau}\in D_{g\acute{e}om}({\cal O},\omega)\otimes Mes(M(F))^*$, il existe $\boldsymbol{\gamma}\in {\cal D}$ aussi proche que l'on veut de ${\cal O}$, de sorte que $g_{\tilde{M},{\cal O}}^{\tilde{M}}(\boldsymbol{\gamma})=\boldsymbol{\tau}$.

{\bf Exemple.} L'ensemble ${\cal D}=D_{g\acute{e}om,\tilde{G}-reg}(\tilde{M}(F),\omega)\otimes Mes(M(F))^*$ des \'el\'ements de  $D_{g\acute{e}om}(\tilde{M}(F),\omega)\otimes Mes(M(F))^*$ \`a support fortement $\tilde{G}$-r\'egulier v\'erifie cette propri\'et\'e d'apr\`es 2.1(1).

\ass{Lemme}{Supposons que l'on ait l'\'egalit\'e $ I_{\tilde{M}}^{\tilde{G}}(\boldsymbol{\gamma},{\bf f})= I_{\tilde{M}}^{\tilde{G},{\cal E}}(\boldsymbol{\gamma},{\bf f})$ pour tout $\boldsymbol{\gamma}\in  {\cal D}$ et pour tout ${\bf f}\in I(\tilde{G}(F),\omega)\otimes Mes(G(F))$. Alors les deux propri\'et\'es  suivantes sont \'equivalentes:

(i) on a l'\'egalit\'e $ I_{\tilde{M}}^{\tilde{G}}(\boldsymbol{\gamma},{\bf f})= I_{\tilde{M}}^{\tilde{G},{\cal E}}(\boldsymbol{\gamma},{\bf f})$  pour tout $\boldsymbol{\gamma}\in D_{g\acute{e}om}({\cal O},\omega)\otimes Mes(M(F))^*$ et tout ${\bf f}\in I(\tilde{G}(F),\omega)\otimes Mes(G(F))$;

(ii) on a l'\'egalit\'e $g_{\tilde{M},{\cal O}}^{\tilde{G}}(\boldsymbol{\gamma})=g_{\tilde{M},{\cal O}}^{\tilde{G},{\cal E}}(\boldsymbol{\gamma})$ pour tout $\boldsymbol{\gamma}\in {\cal D}$ assez proche de ${\cal O}$.}

Preuve. L'hypoth\`ese implique que, pour $\boldsymbol{\gamma}\in {\cal D}$  , le membre de gauche de (1) est nul donc aussi celui de droite. Si (i) est v\'erifi\'e, la premi\`ere diff\'erence de ce membre de droite est nulle. La deuxi\`eme l'est donc aussi, d'o\`u la conclusion de (ii). 
 En sens inverse, (ii) implique de la m\^eme fa\c{c}on l'\'egalit\'e
$$I_{\tilde{M}}^{\tilde{G}}(g_{\tilde{M},{\cal O}}^{\tilde{M}}(\boldsymbol{\gamma}),{\bf f})=I_{\tilde{M}}^{\tilde{G},{\cal E}}(g_{\tilde{M},{\cal O}}^{\tilde{M}}(\boldsymbol{\gamma}),{\bf f})$$
pour tout ${\bf f}$ et tout $\boldsymbol{\gamma}\in {\cal D}$. En utilisant (2), cela entra\^{\i}ne 
$$I_{\tilde{M}}^{\tilde{G}}(\boldsymbol{\tau},{\bf f})=I_{\tilde{M}}^{\tilde{G},{\cal E}}(\boldsymbol{\tau},{\bf f})$$
pour tout  $\boldsymbol{\tau}\in D_{g\acute{e}om}({\cal O},\omega)\otimes Mes(M(F))^*$. C'est l'assertion (i). $\square$

De nouveau, il y a des variantes dans les deux situations suivantes: $(G,\tilde{G},{\bf a})$ est quasi-d\'eploy\'e et \`a torsion int\'erieure et on fixe un syst\`eme de fonctions $B$ comme en 1.9; ou $G=\tilde{G}$, ${\bf a}=1$ et on fixe une fonction $B$ comme en 1.8.

\bigskip

\subsection{Relation entre la proposition 2.4   et le th\'eor\`eme 1.10.}
On suppose $(G,\tilde{G},{\bf a})$ quasi-d\'eploy\'e et \`a torsion int\'erieure. On suppose donn\'e un syst\`eme de fonctions $B$ comme en 1.9. Soit $\tilde{M}$ un espace de Levi de $\tilde{G}$ et ${\cal O}$ une classe de conjugaison stable semi-simple dans $\tilde{M}(F)$. Toutes nos assertions sont tautologiques dans le cas $\tilde{M}=\tilde{G}$, on suppose donc ici $\tilde{M}\not=\tilde{G}$. Soient $\boldsymbol{\delta}\in D_{g\acute{e}om,\tilde{G}-\acute{e}qui}^{st}(\tilde{M}(F))\otimes Mes(M(F))^*$ et ${\bf f}\in I(\tilde{G}(F))\otimes Mes(G(F))$. On suppose que l'image de ${\bf f}$ dans $S I(\tilde{G}(F))\otimes Mes(G(F))$ est nulle, autrement dit que les int\'egrales orbitales stables fortement r\'eguli\`eres de ${\bf f}$ sont  nulles. Consid\'erons le d\'eveloppement de la proposition 2.5(i). Pour $\tilde{L}\in {\cal L}(\tilde{M})$ tel que $\tilde{L}\not=\tilde{M}$ et $\tilde{L}\not=\tilde{G}$, les hypoth\`eses de r\'ecurrence impliquent que $Sg_{\tilde{M},{\cal O}}^{\tilde{L}}(\boldsymbol{\delta},B)$ est stable et que la forme lin\'eaire ${\bf f}'\mapsto S_{\tilde{L}}^{\tilde{G}}(\boldsymbol{\delta}',B,{\bf f}')$ est stable pour tout $\boldsymbol{\delta}'$ stable. Par ailleurs, pour $\tilde{L}=\tilde{M}$, on a simplement $Sg_{\tilde{M},{\cal O}}^{\tilde{M}}(\boldsymbol{\delta},B)=g_{\tilde{M},{\cal O}}^{\tilde{M}}(\boldsymbol{\delta})$ et ce terme est stable d'apr\`es le lemme 2.2. En vertu de l'hypoth\`ese sur ${\bf f}$, le d\'eveloppement se r\'eduit \`a
$$(1) \qquad S_{\tilde{M}}^{\tilde{G}}(\boldsymbol{\delta},B,{\bf f})=  S_{\tilde{M}}^{\tilde{G}}(g_{\tilde{M},{\cal O}}^{\tilde{M}}(\boldsymbol{\delta}),B,{\bf f})+I^{\tilde{G}}(Sg_{\tilde{M},{\cal O}}^{\tilde{G}}(\boldsymbol{\delta},B),B,{\bf f}).$$

 Soit ${\cal D}^{st}$ un sous-ensemble de  $D_{g\acute{e}om,\tilde{G}-\acute{e}qui}^{st}(\tilde{M}(F))\otimes Mes(M(F))^*$ v\'erifiant la propri\'et\'e suivante

(2) pour tout $\boldsymbol{\tau}\in D_{g\acute{e}om}^{st}({\cal O})\otimes Mes(M(F))^*$, il existe $\boldsymbol{\delta}\in {\cal D}^{st}$ aussi proche que l'on veut de ${\cal O}$, de sorte que $g_{\tilde{M},{\cal O}}^{\tilde{M}}(\boldsymbol{\delta})=\boldsymbol{\tau}$.

{\bf Exemple.} L'ensemble ${\cal D}^{st}=D_{g\acute{e}om,\tilde{G}-reg}^{st}(\tilde{M}(F),\omega)\otimes Mes(M(F))^*$ des \'el\'ements de  $D_{g\acute{e}om}^{st}(\tilde{M}(F),\omega)\otimes Mes(M(F))^*$ \`a support fortement $\tilde{G}$-r\'egulier v\'erifie cette propri\'et\'e d'apr\`es le lemme 2.2.

\ass{Lemme}{Supposons que la distribution ${\bf f}'\mapsto S_{\tilde{M}}^{\tilde{G}}(\boldsymbol{\delta},B,{\bf f}')$ soit stable pour tout $\boldsymbol{\delta}\in {\cal D}^{st}$. Alors les propri\'et\'es suivantes sont \'equivalentes:

(i)  la distribution ${\bf f}'\mapsto S_{\tilde{M}}^{\tilde{G}}(\boldsymbol{\delta},B,{\bf f}')$ est stable pour tout $\boldsymbol{\delta}\in D_{g\acute{e}om}^{st}({\cal O})\otimes Mes(M(F))^*$;

(ii) $Sg_{\tilde{M},{\cal O}}^{\tilde{G}}(\boldsymbol{\delta},B)$ est stable pour tout $\boldsymbol{\delta}\in {\cal D}^{st}$ assez proche de ${\cal O}$.}

La preuve est similaire \`a celle du lemme pr\'ec\'edent.

\bigskip

 \subsection{Premi\`eres cons\'equences}
 Soient $(G,\tilde{G},{\bf a})$ un triplet quelconque et $\tilde{M}$ un espace de Levi de $\tilde{G}$.
 
 \ass{Proposition}{(i) Soit ${\cal O}$ une classe de conjugaison semi-simple dans $\tilde{M}(F)$ form\'ee d'\'el\'ements $\tilde{G}$-\'equisinguliers. Supposons que l'on ait l'\'egalit\'e $I_{\tilde{M}}^{\tilde{G},{\cal E}}(\boldsymbol{\gamma},{\bf f})=I_{\tilde{M}}^{\tilde{G}}(\boldsymbol{\gamma},{\bf f})$ pour tout $\boldsymbol{\gamma}\in D_{g\acute{e}om,\tilde{G}-reg}(\tilde{M}(F),\omega)\otimes Mes(M(F))^*$ et tout ${\bf f}\in C_{c}^{\infty}(\tilde{G}(F))\otimes Mes(G(F))$. Alors cette \'egalit\'e est v\'erifi\'ee pour tout ${\bf f}$ et tout $\boldsymbol{\gamma}\in D_{g\acute{e}om}({\cal O},\omega)\otimes Mes(M(F))^*$.
 
 (ii) Supposons $(G,\tilde{G},{\bf a})$ quasi-d\'eploy\'e et \`a torsion int\'erieure. Soit ${\cal O}$ une classe de conjugaison stable semi-simple dans $\tilde{M}(F)$ form\'ee d'\'el\'ements $\tilde{G}$-\'equisinguliers. Supposons que la distribution
 $${\bf f}\mapsto S_{\tilde{M}}^{\tilde{G}}(\boldsymbol{\delta},{\bf f})$$
 soit stable pour tout $\boldsymbol{\delta}\in D_{g\acute{e}om,\tilde{G}-reg}^{st}(\tilde{M}(F))\otimes Mes(M(F))^*$. Alors elle est stable pour tout $\boldsymbol{\delta}\in D_{g\acute{e}om}^{st}({\cal O})\otimes Mes(M(F))^*$.}
 
 Preuve. Pour (i), on applique le lemme 2.8  en prenant  ${\cal D} =D_{g\acute{e}om,\tilde{G}-reg}(\tilde{M}(F),\omega)\otimes Mes(M(F))^*$. Gr\^ace \`a  2.3(1) et 2.6(8), la condition (ii) de ce lemme est v\'erifi\'ee. Donc aussi  la condition (i) de ce lemme, qui n'est autre que la conclusion de l'\'enonc\'e. Pour le (ii), on applique le lemme 2.9  en prenant  ${\cal D}^{st} =D_{g\acute{e}om,\tilde{G}-reg}^{st}(\tilde{M}(F),\omega)\otimes Mes(M(F))^*$. Gr\^ace au lemme 2.2 et \`a 2.4(2), la condition (ii) de ce lemme est v\'erifi\'ee. Donc aussi  la condition (i) de ce lemme, qui n'est autre que la conclusion de l'\'enonc\'e. $\square$

 \bigskip

\subsection{Une formule d'induction}
Soient $(G,\tilde{G},{\bf a})$ un triplet quelconque, $\tilde{M}$ un espace de Levi de $\tilde{G}$ et  $\tilde{R}$ un espace de Levi de $\tilde{M}$. On rappelle qu'il y a un homomorphisme d'induction
$$\begin{array}{ccc}D_{g\acute{e}om}(\tilde{R}(F),\omega)\otimes Mes(R(F))^*&\to&D_{g\acute{e}om}(\tilde{M}(F),\omega)\otimes Mes(M(F))^*\\ \boldsymbol{\gamma}&\mapsto&\boldsymbol{\gamma}^{\tilde{M}}\\ \end{array}$$
Soit ${\cal O}$ une classe de conjugaison semi-simple par $R(F)$ dans $\tilde{R}(F)$.  On note  ${\cal O}^{\tilde{M}}$ la classe engendr\'ee dans $\tilde{M}(F)$.

\ass{Lemme}{Soit $\boldsymbol{\gamma}\in D_{g\acute{e}om}(\tilde{R}(F),\omega)\otimes Mes(R(F))^*$. On suppose que les \'el\'ements du support de $\boldsymbol{\gamma}^{\tilde{M}}$ sont $\tilde{G}$-\'equisinguliers. Si $\boldsymbol{\gamma}$ est assez voisin de ${\cal O}$, on a l'\'egalit\'e
$$g_{\tilde{M},{\cal O}^{\tilde{M}}}^{\tilde{G}}(\boldsymbol{\gamma}^{\tilde{M}})=\sum_{\tilde{L}\in {\cal L}(\tilde{R})}d_{\tilde{R}}^{\tilde{G}}(\tilde{M},\tilde{L})(g_{\tilde{R},{\cal O}}^{\tilde{L}}(\boldsymbol{\gamma}))^{\tilde{G}}.$$}

{\bf Remarque.} Pour tout $\tilde{L}$ tel que $d_{\tilde{M}}^{\tilde{G}}(\tilde{M},\tilde{L})\not=0$, les \'el\'ements du support de $\boldsymbol{\gamma}$ sont $\tilde{L}$-\'equisinguliers.  Cela r\'esulte de l'assertion suivante. Soit $\eta\in \tilde{R}(F)$. Supposons que $\eta$ soit $\tilde{G}$-\'equisingulier en tant qu'\'el\'ement de $\tilde{M}(F)$ (c'est-\`a-dire $M_{\eta}=G_{\eta}$). Alors $\eta$ est $\tilde{L}$-\'equisingulier. Puisque $d_{\tilde{M}}^{\tilde{G}}(\tilde{M},\tilde{L})\not=0$, les tores $A_{\tilde{M}}$ et $A_{\tilde{L}}$ engendrent $A_{\tilde{R}}$. Un \'el\'ement de $\tilde{M}\cap \tilde{L}$ commute \`a $A_{\tilde{M}}$ et $A_{\tilde{L}}$, donc \`a $A_{\tilde{R}}$, donc appartient \`a $\tilde{R}$. D'o\`u $\tilde{M}\cap \tilde{L}=\tilde{R}$. Un \'el\'ement de $M\cap L$ agissant par multiplication \`a gauche conserve $\tilde{M}\cap \tilde{L}$, donc aussi $\tilde{R}$, donc appartient \`a $R$. D'o\`u $M\cap L=R$. Puisque $L_{\eta}\subset G_{\eta}=M_{\eta}$, on a $L_{\eta}\subset M\cap L=R$, d'o\`u l'\'egalit\'e cherch\'ee $L_{\eta}=R_{\eta}$. 

\bigskip

Preuve. Soit ${\bf f}\in C_{c}^{\infty}(\tilde{G}(F))\otimes Mes(G(F))$. On utilise le lemme 1.7:
 $$I_{\tilde{M}}^{\tilde{G}}(\boldsymbol{\gamma}^{\tilde{M}},{\bf f})=\sum_{\tilde{L}_{1}\in {\cal L}(\tilde{R})}d_{\tilde{R}}^{\tilde{G}}(\tilde{M},\tilde{L}_{1})I_{\tilde{R}}^{\tilde{L}_{1}}(\boldsymbol{\gamma},{\bf f}_{\tilde{L}_{1},\omega}).$$
 On d\'eveloppe en germes les deux membres. A gauche, on obtient
 $$(1) \qquad \sum_{\tilde{L}'\in {\cal L}(\tilde{M})}I_{\tilde{L}'}^{\tilde{G}}(g_{\tilde{M},{\cal O}^{\tilde{M}}}^{\tilde{L}'}(\boldsymbol{\gamma}^{\tilde{M}}),{\bf f}).$$
 A droite, on obtient
 $$\sum_{\tilde{L}_{1}\in {\cal L}(\tilde{R})}d_{\tilde{R}}^{\tilde{G}}(\tilde{M},\tilde{L}_{1})\sum_{\tilde{L}_{2}\in {\cal L}(\tilde{R}), \tilde{L}_{2}\subset \tilde{L}_{1}} I_{\tilde{L}_{2}}^{\tilde{L}_{1}}(g_{\tilde{R},{\cal O}}^{\tilde{L}_{2}}(\boldsymbol{\gamma}),{\bf f}_{\tilde{L}_{1},\omega}).$$
 Consid\'erons l'ensemble $A$ des couples d'espace de Levi $(\tilde{L}_{1},\tilde{L}_{2})$ tels que $\tilde{R}\subset \tilde{L}_{2}\subset \tilde{L}_{1}$ et $d_{\tilde{R}}^{\tilde{G}}(\tilde{M},\tilde{L}_{1})\not=0$. Consid\'erons l'ensemble $B$ des triplets $(\tilde{L}_{1},\tilde{L}_{2},\tilde{L}')$ tels que 
 $\tilde{R}\subset \tilde{L}_{2}\subset \tilde{L}_{1}$, $\tilde{M}\subset \tilde{L}'$, $\tilde{L}_{2}\subset \tilde{L}'$ et $d_{\tilde{R}}^{\tilde{L}'}(\tilde{M},\tilde{L}_{2})d_{\tilde{L}_{2}}^{\tilde{G}}(\tilde{L}',\tilde{L}_{1})\not=0$. On a prouv\'e en 1.7(5) 
que  l'application $(\tilde{L}_{1},\tilde{L}_{2},\tilde{L}')\mapsto (\tilde{L}_{1},\tilde{L}_{2})$ etait une bijection de $B$ sur $A$ et que,  pour $(\tilde{L}_{1},\tilde{L}_{2},\tilde{L}')\in B$, on avait l'\'egalit\'e 
 $$d_{\tilde{R}}^{\tilde{L}'}(\tilde{M},\tilde{L}_{2})d_{\tilde{L}_{2}}^{\tilde{G}}(\tilde{L}',\tilde{L}_{1})=
 d_{\tilde{R}}^{\tilde{G}}(\tilde{M},\tilde{L}_{1}).$$
 En utilisant cela, la somme ci-dessus se r\'ecrit
 $$\sum_{\tilde{L}'\in {\cal L}(\tilde{M})}\sum_{\tilde{L}_{2}\in {\cal L}(\tilde{R}), \tilde{L}_{2}\subset \tilde{L}'}d_{\tilde{R}}^{\tilde{L}'}(\tilde{M},\tilde{L}_{2})\sum_{\tilde{L}_{1}\in {\cal L}(\tilde{L}_{2})}d_{\tilde{L}_{2}}^{\tilde{G}}(\tilde{L}',\tilde{L}_{1}) I_{\tilde{L}_{2}}^{\tilde{L}_{1}}(g_{\tilde{R},{\cal O}}^{\tilde{L}_{2}}(\boldsymbol{\gamma}),{\bf f}_{\tilde{L}_{1},\omega}).$$
 Par le lemme 1.7, la derni\`ere somme en $\tilde{L}_{1}$ devient $I_{\tilde{L}'}^{\tilde{G}}((g_{\tilde{R},{\cal O}}^{\tilde{L}_{2}}(\boldsymbol{\gamma}))^{\tilde{L}'},{\bf f})$. 
 L'expression devient
 $$(2) \qquad \sum_{\tilde{L}'\in {\cal L}(\tilde{M})}I_{\tilde{L}'}^{\tilde{G}}( X_{\tilde{M}}^{\tilde{L}'}(\boldsymbol{\gamma}),{\bf f}),$$
 o\`u
 $$X_{\tilde{M}}^{\tilde{L}'}(\boldsymbol{\gamma})=\sum_{\tilde{L}_{2}\in {\cal L}(\tilde{R}), \tilde{L}_{2}\subset \tilde{L}'}d_{\tilde{R}}^{\tilde{L}'}(\tilde{M},\tilde{L}_{2})(g_{\tilde{R},{\cal O}}^{\tilde{L}_{2}}(\boldsymbol{\gamma}))^{\tilde{L}'}.$$
 Les deux expressions (1) et (2) sont \'egales.
 L'assertion du lemme est que $X_{\tilde{M}}^{\tilde{G}}(\boldsymbol{\gamma})=g_{\tilde{M},{\cal O}^{\tilde{M}}}^{\tilde{G}}(\boldsymbol{\gamma}^{\tilde{M}})$. En raisonnant par r\'ecurrence, on peut supposer que cela est vrai si l'on remplace $\tilde{G}$ par $\tilde{L}'\not=\tilde{G}$.  Par diff\'erence entre (1) et (2), on obtient
 $$I^{\tilde{G}}(g_{\tilde{M},{\cal O}^{\tilde{M}}}^{\tilde{G}}(\boldsymbol{\gamma}^{\tilde{M}})-X_{\tilde{M}}^{\tilde{G}}(\boldsymbol{\gamma}) ,{\bf f})=0.$$
 On en d\'eduit l'\'egalit\'e cherch\'ee $X_{\tilde{M}}^{\tilde{G}}(\boldsymbol{\gamma})=g_{\tilde{M},{\cal O}^{\tilde{M}}}^{\tilde{G}}(\boldsymbol{\gamma}^{\tilde{M}})$. $\square$
 
 \bigskip
 
 \subsection{Une formule d'induction, cas endoscopique}
 On conserve les m\^emes donn\'ees, \`a ceci pr\`es que ${\cal O}$ est maintenant une classe de conjugaison stable semi-simple dans $\tilde{R}(F)$.
 
 \ass{Lemme}{Soit $\boldsymbol{\gamma}\in D_{g\acute{e}om}(\tilde{R}(F),\omega)\otimes Mes(R(F))^*$. On suppose que les \'el\'ements du support de $\boldsymbol{\gamma}^{\tilde{M}}$ sont $\tilde{G}$-\'equisinguliers. Si $\boldsymbol{\gamma}$ est assez voisin de ${\cal O}$, on a l'\'egalit\'e
$$g_{\tilde{M},{\cal O}^{\tilde{M}}}^{\tilde{G},{\cal E}}(\boldsymbol{\gamma}^{\tilde{M}})=\sum_{\tilde{L}\in {\cal L}(\tilde{R})}d_{\tilde{R}}^{\tilde{G}}(\tilde{M},\tilde{L})(g_{\tilde{R},{\cal O}}^{\tilde{L},{\cal E}}(\boldsymbol{\gamma}))^{\tilde{G}}.$$}

   La preuve est la m\^eme que la pr\'ec\'edente, en utilisant la relation 1.15 (1) en lieu et place du lemme  1.7.
   
   \bigskip
   \subsection{Une formule d'induction, cas stable}
   On suppose $(G,\tilde{G},{\bf a})$ qusasi-d\'eploy\'e et \`a torsion int\'erieure. On fixe un syst\`eme de fonctions $B$ comme en 1.9. Soient $\tilde{M}$ un espace de Levi de $\tilde{G}$, $\tilde{R}$ un espace de Levi de $\tilde{M}$ et ${\cal O}$ une classe de conjugaison stable semi-simple dans $\tilde{R}(F)$.
   
 \ass{Lemme}{Soit $\boldsymbol{\delta}\in D_{g\acute{e}om}^{st}(\tilde{R}(F),\omega)\otimes Mes(R(F))^*$. On suppose que les \'el\'ements du support de $\boldsymbol{\delta}^{\tilde{M}}$ sont $\tilde{G}$-\'equisinguliers. Si $\boldsymbol{\delta}$ est assez voisin de ${\cal O}$, on a l'\'egalit\'e
$$Sg_{\tilde{M},{\cal O}^{\tilde{M}}}^{\tilde{G}}(\boldsymbol{\delta}^{\tilde{M}})=\sum_{\tilde{L}\in {\cal L}(\tilde{R})}e_{\tilde{R}}^{\tilde{G}}(\tilde{M},\tilde{L})(Sg_{\tilde{R},{\cal O}}^{\tilde{L}}(\boldsymbol{\delta}))^{\tilde{G}}.$$}

La preuve est la m\^eme qu'en 2.11, en utilisant la proposition 1.14(ii).

   \bigskip
   
   \section{D\'eveloppements des int\'egrales orbitales pond\'er\'ees invariantes}
 
 \subsection{Des espaces associ\'es au couple  $(\tilde{G},\tilde{M})$}
 On consid\`ere un triplet $(G,\tilde{G},{\bf a})$ g\'en\'eral. Soit $\tilde{M}$ un ensemble de Levi de $\tilde{G}$.   Consid\'erons l'ensemble des ensembles  $\{\alpha_{i};i=1,...,n\}$ form\'es d'\'el\'ements lin\'eairement ind\'ependants de $\Sigma(A_{\tilde{M}})$ et de nombre d'\'el\'ements maximal, c'est-\`a-dire tels que $n=a_{\tilde{M}}-a_{\tilde{G}}$ (on consid\`ere ici les racines comme des \'el\'ements de $\mathfrak{a}_{\tilde{M}}^*$). On dit que deux tels ensembles  sont \'equivalents s'ils engendrent le m\^eme ${\mathbb Z}$-module dans $\mathfrak{a}_{\tilde{M}}^*$. On note ${\cal J}_{\tilde{M}}^{\tilde{G}}$ l'ensemble des classes d'\'equivalence. Pour $J\in {\cal J}_{\tilde{M}}^{\tilde{G}}$, on note $R_{J}$ de ${\mathbb Z}$ module engendr\'e par les $\alpha_{i}$ pour n'importe quel \'el\'ement $\{\alpha_{i};i=1,...,n\}\in J$. 
   
  Identifions $\underline{la}$ paire de Borel de $G$ \`a une paire $(B^*,T^*)$ pour laquelle $M$ est standard. Soit $J\in {\cal J}_{\tilde{M}}^{\tilde{G}}$. On consid\`ere l'ensemble des racines $\beta\in \Sigma(T^*)$ qui se restreignent \`a $A_{\tilde{M}}$ en un \'el\'ement de $R_{J}$.  C'est le syst\`eme de racines associ\'e \`a un sous-groupe de $G$, que l'on note $G_{J}$. Il contient $M$. On v\'erifie qu'il  est d\'efini sur $F$ et invariant par $ad_{\gamma}$ pour tout $\gamma\in \tilde{M}(F)$ (parce que $ad_{\gamma}$ induit l'identit\'e sur ${\cal A}_{\tilde{M}}$). Alors l'ensemble $\tilde{G}_{J}=G_{J}\tilde{M}$ est un sous-espace tordu de $\tilde{G}$.   On introduit aussi un sous-espace $U_{J}$ de l'espace des germes au point $1$  de fonctions d\'efinies presque partout sur $A_{\tilde{M}}(F)$. C'est le sous-espace  engendr\'e lin\'eairement par les germes de fonctions
 $$a\mapsto \prod_{i=1,...,n}log(\vert \alpha_{i}(a)-\alpha_{i}(a)^{-1}\vert _{F})$$
 pour les ensembles  $\{\alpha_{i};i=1,...,n\}$ appartenant \`a $J$.  Si $\tilde{M}=\tilde{G}$, ${\cal J}_{\tilde{G}}^{\tilde{G}}$ poss\`ede un unique \'el\'ement $\emptyset$. Alors $U_{\emptyset}$ est la droite form\'ee des germes de fonctions constantes.
 
 {\bf Attention.} Les fonctions ci-dessus ne sont pas lin\'eairement ind\'ependantes en g\'en\'eral. Donnons un contre-exemple. On prend $\tilde{G}=G=SO(5)$ et pour $\tilde{M}=M $ un tore d\'eploy\'e maximal. On peut identifier $\mathfrak{a}_{\tilde{M}}(F)$ \`a $F^2$ de sorte qu'un ensemble de racines positives soit form\'e des quatre applications lin\'eaires
 $$(x,y)\mapsto \left\lbrace\begin{array}{ccc}\alpha(x,y)&=&x-y\\\beta(x,y)&=&y\\(\alpha+\beta)(x,y)&=&x\\(\alpha+2\beta)(x,y)&=&x+y.\\ \end{array}\right.$$
 Il y a six ensembles form\'es de deux racines positives lin\'eairement ind\'ependantes. Cinq d'entre eux sont \'equivalents, le dernier, \`a savoir $\{\alpha,\alpha+2\beta\}$, formant une classe d'\'equivalence \`a lui seul (si on ne tient compte que des racines positives). Prenons pour $J$ la classe des cinq premiers. En identifiant les germes de fonctions sur $A_{\tilde{M}}(F)$  au point $1$ \`a des germes de fonctions sur $\mathfrak{a}_{\tilde{M}}(F)$ au point $0$, l'espace $U_{J}$ contient en particulier les germes des fonctions
 $$(x,y)\mapsto\left\lbrace\begin{array}{c}log(\vert x-y\vert _{F})log(\vert y\vert _{F})\\log(\vert x+y\vert _{F})log(\vert x\vert _{F})\\log(\vert x-y\vert _{F})log(\vert x\vert _{F})\\log(\vert x+y\vert _{F})log(\vert y\vert _{F})\\ \end{array}\right.$$
 (on suppose que la caract\'eristique r\'esiduelle est impaire, sinon il faudrait multiplier toutes les racines par $2$).  Mais la somme des deux premi\`eres fonctions et des oppos\'ees des deux derni\`eres est nulle. En effet, elle s'\'ecrit
 $$(log(\vert x-y\vert _{F})-log(\vert x+y\vert _{F}))(log(\vert y\vert _{F})-log(\vert x\vert _{F})).$$
 Si $x$ et $y$ ont m\^eme valeur absolue, le deuxi\`eme facteur est nul. Si $x$ et $y$ ont des valeurs absolues diff\'erentes, alors $x+y$ et $x-y$ ont des valeurs absolues \'egales au sup de celles de $x$ et de $y$. Alors le premier facteur est nul.
 
 \bigskip

   On a n\'eanmoins la propri\'et\'e ci-dessous. Pour l'\'enoncer, on doit d'abord poser une d\'efinition.
 Appelons domaine admissible dans $A_{\tilde{M}}(F)$ l'intersection d'un voisinage ouvert assez petit de $1$ avec l'ensemble des \'el\'ements $a$ qui v\'erifient la condition $\vert \alpha(a)-1\vert _{F}>cd(a)$ pour tout $\alpha\in \Sigma(A_{\tilde{M}})$, o\`u $c>0$ est un r\'eel fix\'e.   Pour un germe de fonction $u$ d\'efinie presque partout dans un voisinage de $1$ dans $A_{\tilde{M}}(F)$, disons que le germe $u$ est \'equivalent \`a $0$ s'il existe $r>0$ et, pour tout domaine admissible, un r\'eel $C>0$ tel que $\vert u(a)\vert \leq Cd(a)^r$ pour tout $a$ dans le domaine et assez proche de $1$. On dit que deux germes sont \'equivalents si leur diff\'erence est \'equivalente \`a $0$. On note $\simeq$ cette relation d'\'equivalence. Cette d\'efinition d\'epend de l'espace ambiant $\tilde{G}$ (puisqu'elle d\'epend de l'ensemble $\Sigma(A_{\tilde{M}})$), mais on esp\`ere que cela ne cr\'eera pas de difficult\'e. Notons que 
 
 (1) si $u$ est un germe \'equivalent \`a $0$ et si $\{\alpha_{i};i=1,...,n\}$ est un ensemble fini d'\'el\'ements de $\Sigma(A_{\tilde{M}})$, le germe
 $$a\mapsto u(a)\prod_{i=1,...,n}log{\vert \alpha_{i}(a)-\alpha_{i}(a)^{-1}\vert _{F}}$$
 est lui-aussi \'equivalent \`a $0$. 
 
 On a
 
 (2) soit $u\in\sum_{\tilde{L}\in {\cal L}(\tilde{M})} \sum_{J\in {\cal J}_{\tilde{M}}^{\tilde{L}}}U_{J}$; supposons $u$ \'equivalent \`a $0$; alors $u=0$.
 
 Preuve. Comme plus haut, on peut descendre via l'exponentielle le germe $u$ en un germe de fonction sur $\mathfrak{a}_{\tilde{M}}(F)$ au voisinage de $0$.  Pour $\tilde{L}\in {\cal L}(\tilde{M})$, $J\in {\cal J}_{\tilde{M}}^{\tilde{L}}$ et  $\underline{\alpha}=\{\alpha_{i};i=1,...,n\}\in J$, on note $u_{\underline{\alpha}}$ la fonction d\'efinie presque partout sur $\mathfrak{a}_{\tilde{M}}(F)$ par
 $$u_{\underline{\alpha}}(H)= \prod_{i=1,...,n}log(\vert 2\alpha_{i}(H)\vert _{F}).$$
 Alors $u$ est combinaison lin\'eaire de telles fonctions.
  Soit $H\in \mathfrak{a}_{\tilde{M}}(F)$ en position g\'en\'erale. Fixons une uniformisante $\varpi_{F}$ de $F$. Pour $k\in {\mathbb Z}$ et $\alpha\in \Sigma(A_{\tilde{M}})$, on a 
  $$log(\vert 2\alpha(\varpi_{F}^kH)\vert_{F})=-k\,log(q)+log(\vert 2\alpha(H)\vert _{F}),$$
  o\`u $q$ est le nombre d'\'el\'ements du corps r\'esiduel. Il en r\'esulte que $u(\varpi_{F}^kH)$ est un polyn\^ome en $k$. Les \'el\'ements $\varpi_{F}^kH$, pour $k\geq0$, restent dans un domaine admissible. Puisque $u$ est \'equivalent \`a $0$, on a donc $lim_{k\to \infty}u(\varpi_{F}^kH)=0$. Mais un polyn\^ome en $k$ qui tend vers $0$ quand $k$ tend vers l'infini est nul. Donc $u(\varpi_{F}^kH)=0$ pour tout $k$. En particulier, pour $k=0$, $u(H)=0$. L'\'el\'ement $H$ \'etant quelconque dans un ouvert dense, on a $u=0$. $\square$
  
  On a

 (3) supposons $\tilde{M}\not=\tilde{G}$;  soient $J\in {\cal J}_{\tilde{M}}^{\tilde{G}}$   et $u\in U_{J}$; supposons que $u$ soit \'equivalent \`a une constante; alors $u=0$ et cette constante est nulle.  
 
 Preuve. Notons $c$ cette constante. Rappelons que ${\mathbb C}$ s'identifie \`a l'espace $U_{\emptyset}$ associ\'e \`a l'\'el\'ement vide de ${\cal J}_{\tilde{M}}^{\tilde{M}}$. On peut donc appliquer (2) \`a $u-c$, d'o\`u $u-c=0$. On descend les fonctions en des fonctions sur $\mathfrak{a}_{\tilde{M}}(F)$ et on utilise les notations de la preuve de (2). On \'ecrit $u=\sum_{\underline{\alpha}\in J}c_{\underline{\alpha}}u_{\underline{\alpha}}$, avec des coefficients complexes $c_{\underline{\alpha}}$.  
 Fixons   un point $H$ en position g\'en\'erale. Pour $\underline{\alpha}=\{\alpha_{i};i=1,...,n\},\underline{\beta}=\{\beta_{i};i=1,...,n\}\in J$, les $\alpha_{i}(H)$ se d\'eduisent des $\beta_{i}(H)$ par une matrice \`a coefficients entiers. Donc
  $$sup_{i=1,..,n}(\vert \alpha_{i}(H)\vert_{F} )\leq sup_{i=1,...,n}(\vert \beta_{i}(H)\vert _{F}).$$
  On peut \'echanger les r\^oles de $\underline{\alpha}$ et $\underline{\beta}$. Donc les deux sup sont \'egaux. Quitte \`a multiplier $H$ par une puissance convenable de $\varpi_{F}$, on peut supposer que ces sup valent $1$. Mais alors $u_{\underline{\alpha}}(H)=0$ pour tout $\underline{\alpha}\in J$. Donc $u(H)=0$. L'\'egalit\'e $u-c=0$ entra\^{\i}ne $c=0$, puis $u=0$. $\square$

 L'ensemble ${\cal J}_{\tilde{M}}^{\tilde{G}}$ contient un unique \'el\'ement $J$ tel que $R_{J}$ contienne $\Sigma(A_{\tilde{M}})$ tout entier. C'est celui qui contient tout ensemble $\{\alpha_{i};i=1,...,n\}$ formant une base de $\Sigma(A_{\tilde{M}})$ pour un certain ordre. On dit que cet \'el\'ement $J$ est l'\'el\'ement maximal de ${\cal J}_{\tilde{M}}^{\tilde{G}}$.
 
 Soit $J\in {\cal J}_{\tilde{M}}^{\tilde{G}}$.    On a $\Sigma^{\tilde{G}_{J}}(A_{\tilde{M}})\subset \Sigma^{\tilde{G}}(A_{\tilde{M}})$. L'ensemble $ {\cal J}_{\tilde{M}}^{\tilde{G}_{J}}$ s'identifie au sous-ensemble des $J'\in {\cal J}_{\tilde{M}}^{\tilde{G}}$ qui contiennent une famille form\'ee d'\'el\'ements de $\Sigma^{\tilde{G}_{J}}(A_{\tilde{M}})$, ou encore, ce qui revient au m\^eme, dont tous les \'el\'ements sont des familles  form\'ees d'\'el\'ements de $\Sigma^{\tilde{G}_{J}}(A_{\tilde{M}})$.    Pour $J'\in {\cal J}_{\tilde{M}}^{\tilde{G}_{J}}$,   l'espace $U_{J'}$ ne d\'epend pas de l'espace ambiant $\tilde{G}$ ou $\tilde{G}_{J}$. Cas particulier: $J$ s'identifie \`a l'\'el\'ement maximal de $ {\cal J}_{\tilde{M}}^{\tilde{G}_{J}}$.  
 
    Notons $Ann_{\tilde{M}}^{\tilde{G}}$ l'annulateur de l'homomorphisme d'induction
 $$D_{g\acute{e}om}(\tilde{M}(F),\omega)\otimes Mes(M(F))^*\to D_{g\acute{e}om}(\tilde{G}(F),\omega)\otimes Mes(G(F))^*.$$
  
  \ass{Lemme}{ Pour tout $J\in {\cal J}_{\tilde{M}}^{\tilde{G}}$, on a l'inclusion $Ann_{\tilde{M}}^{\tilde{G}_{J}}\subset Ann_{\tilde{M}}^{\tilde{G}}$.} 
  
  Preuve.   On fixe des mesures sur chaque groupe. Il est \'equivalent de prouver que l'image de l'application
  $$\begin{array}{cccc}res_{\tilde{M}}^{\tilde{G}}:&I(\tilde{G}(F),\omega)&\to&I(\tilde{M}(F),\omega)\\&f&\mapsto&f_{\tilde{M},\omega}\\ \end{array}$$
  est contenue dans celle de l'application $res_{\tilde{M}}^{\tilde{G}_{J}}$.  On fixe un espace de Levi minimal $\tilde{M}_{0}\subset \tilde{M}$ et un espace parabolique $\tilde{P}_{0}\in {\cal P}(\tilde{M}_{0})$ de sorte que $\tilde{M}$ soit standard.  D'apr\`es [I] lemme 4.3, l'image de $res_{\tilde{M}}^{\tilde{G}}$ est form\'ee des $\varphi\in I(\tilde{M}(F),\omega)$ qui v\'erifient la condition suivante:
  
  - pour deux espaces de Levi $\tilde{R}, \tilde{R}'\in {\cal L}^{\tilde{M}}(\tilde{M}_{0})$ et pour $w\in W^G(\tilde{M}_{0})$ tel que $w(\tilde{R})=\tilde{R}'$, la fonction $\varphi_{\tilde{R'},\omega}$ est l'image de $\varphi_{\tilde{R},\omega}$ par l'isomorphisme $I(\tilde{R}(F),\omega)\to I(\tilde{R}'(F),\omega)$ d\'eduit de $w$.
  
  En appliquant la m\^eme caract\'erisation pour $\tilde{G}_{J}$, l'assertion r\'esulte simplement de l'inclusion $W^{G_{J}}(\tilde{M}_{0})\subset W^G(\tilde{M}_{0})$. Pour prouver celle-ci,  identifions $\underline{la}$ paire de Borel de $G$ \`a une paire $(B^*,T^*)$ telle que $B^*\subset P_{0}$, $T^*\subset M_{0}$. Notons $W^G$ le groupe de Weyl de $T^*$ et $Norm_{W^G}(M_{0})$ l'ensemble des \'el\'ements de $W^G$ qui conservent $M_{0}$. La paire de Borel d\'etermine un automorphisme $\theta$ de $W^G$ et l'action galoisienne quasi-d\'eploy\'ee sur cette paire d\'etermine une action de $\Gamma_{F}$ sur $W^G$. Alors $W^G(\tilde{M}_{0})$ est isomorphe au sous-groupe des invariants par $\Gamma_{F}$ et $\theta$ dans le quotient $(Norm_{W^G}(M_{0})/W^{M_{0}})$. L'assertion \`a prouver r\'esulte alors des faits que $W^{G_{J}}$ est un sous-groupe de $ W^G$ (parce que le syst\`eme de racines de $G_{J}$ est un sous-syst\`eme de celui de $G$) et que cette inclusion est \'equivariante pour les actions de $\theta$ et de $\Gamma_{F}$. $\square$

  \bigskip
 
 \subsection{Un d\'eveloppement des int\'egrales pond\'er\'ees $\omega$-\'equivariantes}
 
 Soient $(G,\tilde{G},{\bf a}) $ un triplet quelconque,  $\tilde{M}$ un espace de Levi de $\tilde{G}$ et ${\cal O}$ une classe de conjugaison semi-simple dans $\tilde{M}(F)$. On pose $Ann_{{\cal O}}^{\tilde{G}}=Ann_{\tilde{M}}^{\tilde{G}}\cap (D_{g\acute{e}om}({\cal O},\omega)\otimes Mes(M(F))^*$. 
 
  \ass{Proposition}{Pour tout   $J\in {\cal J}_{\tilde{M}}^{\tilde{G}}$, il existe une unique application lin\'eaire $\rho^{\tilde{G}}_{J}:D_{g\acute{e}om}({\cal O},\omega)\otimes Mes(M(F))^*\to U_{J}\otimes (D_{g\acute{e}om}({\cal O},\omega)\otimes Mes(M(F))^*)/Ann_{{\cal O}}^{\tilde{G}}$  de sorte que les propri\'et\'es suivantes soient v\'erifi\'ees.
 
 (i) L'application $\rho_{J}^{\tilde{G}}$  est  la compos\'ee de $\rho_{J}^{\tilde{G}_{J}}$ et de la projection
 $$  U_{J}\otimes (D_{g\acute{e}om}({\cal O},\omega)\otimes Mes(M(F))^*)/Ann_{{\cal O}}^{\tilde{G}_{J}}\to   U_{J}\otimes (D_{g\acute{e}om}({\cal O},\omega)\otimes Mes(M(F))^*)/Ann_{{\cal O}}^{\tilde{G}};$$

 (ii) Pour tout $\boldsymbol{\gamma}\in D_{g\acute{e}om}({\cal O},\omega)\otimes Mes(M(F))^*$ et pour tout ${\bf f}\in I(\tilde{G}(F),\omega)\otimes Mes(G(F))$, le germe en $1$ de la fonction
 $$a\mapsto I_{\tilde{M}}^{\tilde{G}}(a\boldsymbol{\gamma},{\bf f}),$$
 qui est d\'efinie pour tout $a\in A_{\tilde{M}}(F)$ en position g\'en\'erale, est \'equivalent \`a
 $$\sum_{\tilde{L}\in {\cal L}(\tilde{M})}\sum_{J\in {\cal J}_{\tilde{M}}^{\tilde{L}}}I_{\tilde{L}}^{\tilde{G}}(\rho_{J}^{\tilde{L}}(\boldsymbol{\gamma},a)^{\tilde{L}},{\bf f})$$.}
 
 {\bf Remarques.} (1) Le terme $\rho^{\tilde{G}}_{J}(\boldsymbol{\gamma})$  d\'efinit un germe de fonction sur $A_{\tilde{M}}(F)$ \`a valeurs dans $(D_{g\acute{e}om}({\cal O},\omega)\otimes Mes(M(F))^*)/Ann_{{\cal O}}^{\tilde{G}}$. On a not\'e $\rho^{\tilde{G}}_{J}(\boldsymbol{\gamma},a)$ la valeur en $a$ de ce germe.

  (2)  La proposition \'equivaut \`a dire que le germe de la fonction $a\mapsto g_{\tilde{M},{\cal O}}^{\tilde{G}}(a\boldsymbol{\gamma})$ est \'equivalent \`a 
 $$\sum_{J\in {\cal J}_{\tilde{M}}^{\tilde{G}} } \rho_{J}(\boldsymbol{\gamma},a)^{\tilde{G}}.$$
 
 \bigskip
 
 Preuve de l'unicit\'e. On raisonne par r\'ecurrence sur la dimension de $G$. Pour $\tilde{L}\in {\cal L}(\tilde{M})$, $\tilde{L}\not=\tilde{G}$ et $J\in {\cal J}_{\tilde{M}}^{\tilde{L}}$, l'application $\rho_{J}^{\tilde{L}}$ est d\'ej\`a d\'etermin\'ee. Pour $J\in {\cal J}_{\tilde{M}}^{\tilde{G}}$ non maximal, l'application $\rho^{\tilde{G}}_{J}$ est uniquement d\'etermin\'ee par  la condition (i).  Il reste \`a d\'eterminer $\rho^{\tilde{G}}_{J}$ pour l'unique \'el\'ement maximal de ${\cal J}_{\tilde{M}}^{\tilde{G}}$.    Pour ce $J$, fixons une base $(u_{k})_{k=1,...,m}$ de $U_{J}$. On peut \'ecrire $\rho^{\tilde{G}}_{J}(\boldsymbol{\gamma},a)=\sum_{k=1,...,m}u_{k}(a) \boldsymbol{\gamma}_{k}$, avec des \'el\'ements $\boldsymbol{\gamma}_{k}\in D_{g\acute{e}om}({\cal O},\omega)\otimes Mes(M(F))^*$.  L'\'egalit\'e du (ii) d\'etermine \`a \'equivalence pr\`es le germe de la fonction $I^{\tilde{G}}(\rho_{J}^{\tilde{G}}(\boldsymbol{\gamma},a)^{\tilde{G}},{\bf f})$ pour tout ${\bf f}$, c'est-\`a-dire de la fonction
 $$\sum_{k=1,...,m}u_{k}(a)I^{\tilde{G}}(\boldsymbol{\gamma}_{k}^{\tilde{G}},{\bf f}).$$
 D'apr\`es 3.1(3), cela d\'etermine les distributions $ \boldsymbol{\gamma}_{k}^{\tilde{G}}$, ce qui est \'equivalent \`a d\'eterminer les  distributions $\boldsymbol{\gamma}_{k}$ modulo $Ann_{{\cal O}}^{\tilde{G}}$.
 
 Preuve de l'existence. Par lin\'earit\'e, on peut se limiter \`a prouver l'existence des germes $\rho_{J}(\boldsymbol{\gamma},a)$ quand $\boldsymbol{\gamma}$ est l'int\'egrale orbitale associ\'ee \`a un \'el\'ement $\gamma$ dont la partie semi-simple appartient \`a ${\cal O}$. Montrons qu'alors,  pour tout ${\bf f}$, le germe en $1$ de la fonction  $a\mapsto  I_{\tilde{M}}^{\tilde{G}}(a\boldsymbol{\gamma},{\bf f})$ est \'equivalent \`a celui de la fonction qui \`a $a$ associe
 $$(3) \sum_{\tilde{L}\in {\cal L}(\tilde{M})}(-1)^{a_{\tilde{M}}-a_{\tilde{L}}}r_{\tilde{M}}^{\tilde{L}}(\gamma,a)I_{\tilde{L}}^{\tilde{G}}(\boldsymbol{\gamma}^{\tilde{L}},{\bf f}).$$
 Comme on le rappellera ci-dessous, les termes $r_{\tilde{M}}^{\tilde{L}}(\gamma,a)$ sont des combinaisons lin\'eaires de produits de termes $log(\vert\alpha(a)-\alpha(a)^{-1}\vert _{F})$ pour $\alpha\in \Sigma(A_{\tilde{M}})$. En utilisant les relations 1.7(12) et 3.1(1), on voit que la fonction (3) ci-dessus est \'equivalente \`a
 $$\sum_{\tilde{L}\in {\cal L}(\tilde{M})}(-1)^{a_{\tilde{M}}-a_{\tilde{L}}}r_{\tilde{M}}^{\tilde{L}}(\gamma,a)\sum_{\tilde{R}\in {\cal L}(\tilde{L})}r_{\tilde{L}}^{\tilde{R}}(\gamma,a)I_{\tilde{R}}^{\tilde{G}}(a\boldsymbol{\gamma},{\bf f}),$$
 ou encore \`a
 $$\sum_{\tilde{R}\in {\cal L}(\tilde{M})}t_{\tilde{M}}^{\tilde{R}}(\gamma,a)I_{\tilde{R}}^{\tilde{G}}(a\boldsymbol{\gamma},{\bf f}),$$
 o\`u
 $$t_{\tilde{M}}^{\tilde{R}}(\gamma,a)=\sum_{\tilde{L}\in {\cal L}^{\tilde{R}}(\tilde{M})}(-1)^{a_{\tilde{M}}-a_{\tilde{L}}}r_{\tilde{M}}^{\tilde{L}}(\gamma,a)r_{\tilde{L}}^{\tilde{R}}(\gamma,a).$$
 Pour en d\'eduire l'assertion ci-dessus, il suffit de prouver que 
 $$(4) \qquad t_{\tilde{M}}^{\tilde{R}}(\gamma,a)=\left\lbrace\begin{array}{cc}0,&\text{  si }\tilde{R}\not=\tilde{M},\\ 1,& \text{ si }\tilde{R}=\tilde{M}.\\ \end{array}\right.$$
 D\'efinissons la $(\tilde{G},\tilde{M})$-famille $(r'_{\tilde{P}}(\gamma,a;\lambda))_{\tilde{P}\in {\cal P}(\tilde{M})}$ par 
 $$r'_{\tilde{P}}(\gamma,a;\lambda)=r_{\tilde{P}}(\gamma,a;-\lambda)=r_{\tilde{P}}(\gamma,a;\lambda)^{-1}.$$
 La premi\`ere \'egalit\'e entra\^{\i}ne que $r_{\tilde{M}}^{_{'}\tilde{L}}(\gamma,a)=(-1)^{a_{\tilde{M}}-a_{\tilde{L}}}r_{\tilde{M}}^{\tilde{L}}(\gamma,a)$. Donc, par une formule usuelle, $t_{\tilde{M}}^{\tilde{R}}(\gamma,a)$ est la fonction associ\'ee \`a la famille produit $(r'_{\tilde{P}}(\gamma,a;\lambda)r_{\tilde{P}}(\gamma,a;\lambda))_{\tilde{P}\in {\cal P}(\tilde{M})}$. Or celle-ci est une famille de fonctions constantes de valeur $1$, d'o\`u l'\'egalit\'e (4).
 
 La famille $(r_{\tilde{P}}(\gamma,a;\lambda))_{\tilde{P}\in {\cal P}(\tilde{M})}$ est d'une forme particuli\`ere qui permet le calcul des fonctions $r_{\tilde{M}}^{\tilde{L}}(\gamma,a)$. Pr\'ecis\'ement, pour toute base $\underline{\alpha}=\{\alpha_{i};i=1,...,n\}$  de ${\cal A}_{\tilde{M}}^{\tilde{L}}$ form\'ee d'\'el\'ements de $\Sigma^{\tilde{L}}(A_{\tilde{M}})$, notons $m(\underline{\alpha};\gamma)$ le volume du quotient de ${\cal A}_{\tilde{M}}^{\tilde{L}}$ par le ${\mathbb Z}$-module engendr\'e par les $\rho(\alpha_{i},\gamma)$ (cf. 1.5), avec la convention $m(\underline{\alpha},\gamma)=0$ si l'un de ces \'el\'ements est nul. Notons aussi $sgn(\underline{\alpha},\gamma)$ le produit des signes des nombres r\'eels $<\alpha_{i},\rho(\alpha_{i},\gamma)>$ (avec la m\^eme convention). Et notons $u_{\underline{\alpha}}$ la fonction sur $A_{\tilde{M}}(F)$ d\'efinie par 
 $$u_{\underline{\alpha}}(a)=\prod_{i=1,...,n}log(\vert \alpha_{i}(a)-\alpha_{i}(a)^{-1}\vert _{F}).$$
 Dans le cas particulier o\`u $\tilde{L}=\tilde{M}$, auquel cas $\underline{\alpha}$ est l'ensemble vide, on admet par convention que les trois termes que l'on vient de d\'efinir valent $1$.
 Alors le lemme  7.1 de [A5]  entra\^{\i}ne l'\'egalit\'e
 $$r_{\tilde{M}}^{\tilde{L}}(\gamma,a)=\sum_{\underline{\alpha}}m(\underline{\alpha},\gamma)sgn(\underline{\alpha},\gamma)u_{\underline{\alpha}}(a),$$
 o\`u on somme sur les ensembles $\underline{\alpha}$ d\'ecrits ci-dessus. L'ensemble de ces ensembles n'est autre que la r\'eunion des $J\in {\cal J}_{\tilde{M}}^{\tilde{L}}$. Pour tout $J\in {\cal J}_{\tilde{M}}^{\tilde{G}}$, posons
 $$(5) \qquad \rho_{J}^{\tilde{G}}(\boldsymbol{\gamma},a)=\sum_{\underline{\alpha}\in J} m(\underline{\alpha},\gamma)sgn(\underline{\alpha},\gamma)u_{\underline{\alpha}}(a)\boldsymbol{\gamma}.$$
 C'est un \'el\'ement de $U_{J}\otimes D_{g\acute{e}om}({\cal O},\omega)\otimes Mes(M(F))^*$.
 Alors la somme (3) devient
 $$\sum_{\tilde{L}\in {\cal L}(\tilde{M})}\sum_{J\in {\cal J}_{\tilde{M}}^{\tilde{L}}}I_{\tilde{L}}^{\tilde{G}}(\rho_{J}^{\tilde{L}}(\boldsymbol{\gamma},a)^{\tilde{L}},{\bf f}).$$
 Puisque la somme (3) est, \`a \'equivalence pr\`es, le germe de $I_{\tilde{M}}^{\tilde{G}}(a\boldsymbol{\gamma},{\bf f})$, on a obtenu l'assertion (ii) de l'\'enonc\'e. L'assertion (i)  est imm\'ediate d'apr\`es la d\'efinition de $\rho^{\tilde{G}}_{J}(\boldsymbol{\gamma},a)$ et 1.5(4). $\square$
 
 On a
 
 (6) si $\tilde{M}=\tilde{G}$, $\rho_{\emptyset}^{\tilde{G}}$ est l'identit\'e, modulo l'isomorphisme $U_{\emptyset}\simeq {\mathbb C}$.
 
 C'est imm\'ediat sur la d\'efinition (5). 
 
 \bigskip
 
 \subsection{D\'eveloppement des int\'egrales orbitales pond\'er\'ees invariantes et  fonction $B$}
 On suppose $ G=\tilde{G}$ et ${\bf a}=1$. On fixe une fonction $B$ comme en 1.8. Soit $M$ un Levi de $G$. On dispose alors de l'ensemble $\Sigma(A_{M},B)$. On peut reprendre les d\'efinitions de 3.1 en rempla\c{c}ant partout l'ensemble $\Sigma(A_{\tilde{M}})$ par cet ensemble $\Sigma(A_{M},B)$. En particulier, on note ${\cal J}_{M}^{G}(B)$ l'analogue de ${\cal J}_{\tilde{M}}^{\tilde{G}}$. Toutes les propri\'et\'es \'enonc\'ees en 3.1 restent vraies avec ces d\'efinitions modifi\'ees. La seule diff\'erence est que, pour $J\in {\cal J}_{M}^{G}(B)$, le groupe $G_{J}$ n'est plus en g\'en\'eral un sous-groupe de $G$ (les groupes $G_{\alpha}$ de 1.8 en sont des cas particuliers). Cela ne cr\'ee pas de perturbations. De m\^eme, l'analogue de la proposition 3.2 reste vraie, avec la m\^eme preuve. Enon\c{c}ons-la, avec des notations \'evidentes.
  
  \ass{Proposition}{Pour tout   $J\in {\cal J}_{M}^{G}(B)$, il existe une unique application lin\'eaire $\rho^{G}_{J}:D_{unip}(M(F))\otimes Mes(M(F))^*\to U_{J}\otimes (D_{unip}(M(F))\otimes Mes(M(F))^*)/Ann_{unip}^{G}$ de sorte que les propri\'et\'es suivantes soient v\'erifi\'ees.
 
 (i) L'application $\rho^G_{J}$  est la  compos\'ee de $\rho_{J}^{G_{J}}$ et de la projection
 $$U_{J}\otimes (D_{unip}(M(F))\otimes Mes(M(F))^*)/Ann_{unip}^{G_{J}}\to U_{J}\otimes (D_{unip}(M(F))\otimes Mes(M(F))^*)/Ann_{unip}^{G}.$$

 (ii) Pour tout $\boldsymbol{\gamma}\in D_{unip}(M(F))\otimes Mes(M(F))^*$ et pour tout ${\bf f}\in I(G(F))\otimes Mes(G(F))$, le germe en $1$ de la fonction
 $$a\mapsto I_{M}^{G}(a\boldsymbol{\gamma},{\bf f}),$$
 qui est d\'efinie pour tout $a\in A_{M}(F)$ en position g\'en\'erale, est \'equivalent \`a
 $$\sum_{L\in {\cal L}(M)}\sum_{J\in {\cal J}_{M}^{G}(B)}I_{L}^{G}(\rho_{J}^{L}(\boldsymbol{\gamma},a)^{L},B,{\bf f})$$.}
  \bigskip
  
  \subsection{D\'eveloppement des int\'egrales orbitales pond\'er\'ees invariantes et syst\`eme de fonctions $B$}
  On suppose $(G,\tilde{G},{\bf a})$ quasi-d\'eploy\'e et \`a torsion int\'erieure. On fixe un syst\`eme de fonctions $B$ comme en 1.9. Soient $\tilde{M}$ un espace de Levi de $\tilde{G}$ et ${\cal O}$ une classe de conjugaison semi-simple dans $\tilde{M}(F)$. Pour $\eta\in {\cal O}$, on a d\'efini l'ensemble $\Sigma(A_{M},B_{\eta})$ en 1.9. Il ne d\'epend pas du choix de $\eta$, on le note plut\^ot $\Sigma(A_{M},B_{{\cal O}})$. De nouveau, on peut reprendre les d\'efinitions de 3.1 en rempla\c{c}ant partout l'ensemble $\Sigma(A_{\tilde{M}})$ par cet ensemble $\Sigma(A_{M},B_{{\cal O}})$. En particulier, on note ${\cal J}_{\tilde{M}}^{\tilde{G}}(B_{{\cal O}})$ l'analogue de ${\cal J}^{\tilde{G}}_{\tilde{M}}$. Remarquons que cet ensemble peut \^etre vide. En effet, puisque $\Sigma(A_{M},B_{{\cal O}})$ est par d\'efinition l'ensemble des restrictions d'\'el\'ements de $\Sigma^{G_{\eta}}(A_{M_{\eta}},B_{\eta})$, il n'existe pas de sous-ensemble lin\'eairement ind\'ependant et de rang $a_{M}-a_{G}$ si  l'ensemble $\Sigma^{G_{\eta}}(A_{M_{\eta}})$ est trop petit.  
    
   Il y a une diff\'erence cruciale avec la situation de 3.1: pour $J\in {\cal J}_{\tilde{M}}^{\tilde{G}}(B_{{\cal O}})$, on ne peut plus d\'efinir le groupe $G_{J}$ car, pour $\eta\in {\cal O}$, la fonction $B_{\eta}$ n'est d\'efinie que sur un sous-ensemble de l'ensemble de racines de $G$. L'assertion (i) de la proposition 3.2 n'a pas d'analogue dans notre situation. On peut toutefois d\'efinir une application
 $$\rho^{\tilde{G}}_{J}:D_{g\acute{e}om}({\cal O})\otimes Mes(M(F))^*\to U_{J}\otimes (D_{g\acute{e}om}({\cal O})\otimes Mes(M(F))^*)/Ann_{{\cal O}}^{G}$$ 
 de la fa\c{c}on suivante.  Par lin\'earit\'e, il suffit de la d\'efinir sur une int\'egrale orbitale $\boldsymbol{\gamma}$ associ\'ee \`a un \'el\'ement $\gamma$ dont la partie semi-simple appartient \`a ${\cal O}$. On reprend la d\'efinition 3.2(5), en y rempla\c{c}ant l'ensemble $\Sigma(A_{\tilde{M}})$ par $\Sigma(A_{M},B_{{\cal O}})$. C'est-\`a-dire qu'avec des notations adapt\'ees, on pose
 $$(1) \qquad \rho^{\tilde{G}}_{J}(\boldsymbol{\gamma},a)=\sum_{\underline{\alpha}\in J}m(\underline{\alpha},\gamma,B_{{\cal O}})sgn(\underline{\alpha},\gamma,B_{{\cal O}})u_{\underline{\alpha}}(a)\boldsymbol{\gamma}.$$
 Plus exactement, $\rho^{\tilde{G}}_{J}(\boldsymbol{\gamma},a)$ est l'image de ce terme modulo $Ann_{{\cal O}}^{G}$. Remarquons que, si $\tilde{M}=\tilde{G}$,  on a $\rho^{\tilde{G}}_{\emptyset}(\boldsymbol{\gamma},a)=\boldsymbol{\gamma}$. La partie "existence" de la proposition 3.2 reste valable, ce qui conduit \`a l'\'enonc\'e suivant.
  
   \ass{Proposition}{ Pour tout $\boldsymbol{\gamma}\in D_{g\acute{e}om}({\cal O})\otimes Mes(M(F))^*$ et pour tout ${\bf f}\in I(\tilde{G}(F))\otimes Mes(G(F))$, le germe en $1$ de la fonction
 $$a\mapsto I_{\tilde{M}}^{\tilde{G}}(a\boldsymbol{\gamma},{\bf f}),$$
 qui est d\'efinie pour tout $a\in A_{\tilde{M}}(F)$ en position g\'en\'erale, est \'equivalent \`a
 $$\sum_{\tilde{L}\in {\cal L}(\tilde{M})}\sum_{J\in {\cal J}_{\tilde{M}}^{\tilde{L}}(B_{{\cal O}})}I_{\tilde{L}}^{\tilde{G}}(\rho_{J}^{\tilde{L}}(\boldsymbol{\gamma},a)^{\tilde{L}},B,{\bf f})$$.}
 
 {\bf Remarques.}  (2) On peut remplacer ${\cal O}$ par une r\'eunion finie $\cup_{i=1,...,n}{\cal O}_{i}$ de classes de conjugaison semi-simples, pour peu que   l'ensemble $\Sigma(A_{M},B_{{\cal O}_{i}})$ ne d\'epende pas de $i$. C'est le cas si ${\cal O}$ est une classe de conjugaison stable.
 
 (3) Supposons que le syst\`eme de fonctions $B$ soit le syst\`eme "trivial", c'est-\`a-dire que la fonction $B_{\eta}$ soit constante de valeur $1$ pour tout $\eta\in \tilde{G}_{ss}(F)$. Les ensembles $\Sigma(A_{M})$ et $\Sigma(A_{M},B_{{\cal O}})$ ne co\"{\i}ncident pas pour autant car, par d\'efinition, ce dernier est l'ensemble des \'el\'ements du premier qui sont restrictions d'\'el\'ements de $\Sigma^{G_{\eta}}(A_{M_{\eta}})$ pour $\eta\in {\cal O}$. On a toutefois une injection ${\cal J}_{\tilde{M}}^{\tilde{G}}(B_{{\cal O}})\subset {\cal J}_{\tilde{M}}^{\tilde{G}}$. Il r\'esulte des d\'efinitions que, pour $J\in {\cal J}_{\tilde{M}}^{\tilde{G}}(B_{{\cal O}})$, les deux d\'efinitions possibles de $\rho_{J}^{\tilde{G}}$ co\"{\i}ncident, tandis que, pour $J\in {\cal J}_{\tilde{M}}^{\tilde{G}}-{\cal J}_{\tilde{M}}^{\tilde{G}}(B_{{\cal O}})$, le terme $\rho_{J}^{\tilde{G}}$ d\'efini en 3.2 est nul. 
 
  \bigskip
 
 \subsection{Termes d'un d\'eveloppement stable}
 On conserve la m\^eme situation que dans le paragraphe pr\'ec\'edent. L'ensemble ${\cal O}$ est maintenant une classe de conjugaison stable semi-simple dans $\tilde{M}(F)$. On note $Ann_{\tilde{M}}^{\tilde{G},st}$, resp. $Ann_{{\cal O}}^{\tilde{G},st}$, l'intersection de $Ann_{\tilde{M}}^{\tilde{G}}$ et de $D_{g\acute{e}om}^{st}(\tilde{M}(F))\otimes Mes(M(F))^*$, resp. de $Ann_{{\cal O}}^{\tilde{G}}$ et de $D_{g\acute{e}om}^{st}({\cal O})\otimes Mes(M(F))^*$. 
 
 Soit $s\in Z(\hat{M})^{\Gamma_{F}}/Z(\hat{G})^{\Gamma_{F}}$. Comme on l'a dit en 1.10, il se d\'eduit du syst\`eme de fonctions $B$ un tel syst\`eme sur $\tilde{G}'(s;F)$ que l'on note encore $B$. Il r\'esulte des d\'efinitions que $\Sigma^{G'(s)}(A_{M},B_{{\cal O}})\subset \Sigma^G(A_{M},B_{{\cal O}})$.  Supposons ${\bf G}'(s)$ elliptique. Alors $a_{G'(s)}=a_{G}$ et de cette inclusion se d\'eduit une inclusion ${\cal J}_{\tilde{M}}^{\tilde{G}'(s)}(B_{{\cal O}})\subset {\cal J}_{\tilde{M}}^{\tilde{G}}(B_{{\cal O}})$.  Si $J\in {\cal J}_{\tilde{M}}^{\tilde{G}'(s)}(B_{{\cal O}})$, la preuve du lemme 3.1 s'applique: on a l'inclusion $Ann_{{\cal O}}^{\tilde{G}'(s)}\subset Ann_{{\cal O}}^{\tilde{G}}$. On voit aussi que l'espace $U_{J}$ est unique, sa d\'efinition ne d\'ependant pas de l'espace ambiant $\tilde{G}$ ou $\tilde{G}'(s)$.
 
 Soit $J\in {\cal J}_{\tilde{M}}^{\tilde{G}}(B_{{\cal O}})$. Dans le paragraphe pr\'ec\'edent, on a d\'efini  une application $\rho^{\tilde{G}}_{J}$ sur $D_{g\acute{e}om}({\cal O})\otimes Mes(M(F))^*$. Notons $\rho^{\tilde{G}}_{J,st}$ sa restriction \`a $D_{g\acute{e}om}^{st}({\cal O})\otimes Mes(M(F))^*$. On d\'efinit une application 
 $$\sigma_{J}:D_{g\acute{e}om}^{st}({\cal O})\otimes Mes(M(F))^*\to U_{J}\otimes (D_{g\acute{e}om}({\cal O})\otimes Mes(M(F))^*)/Ann_{{\cal O}}^{\tilde{G}}$$
 ou plus pr\'ecis\'ement $\sigma_{J}^{\tilde{G}}$
 par la formule de r\'ecurrence
 $$(1) \qquad \sigma^{\tilde{G}}_{J}=\rho_{J,st}-\sum_{s\in Z(\hat{M})^{\Gamma_{F}}/Z(\hat{G})^{\Gamma_{F}}; s\not=1, J\in {\cal J}_{\tilde{M}}^{\tilde{G}'(s)}(B_{{\cal O}})}i_{\tilde{M}}(\tilde{G},\tilde{G}'(s))\sigma_{J}^{\tilde{G}'(s)}.$$
 Plus exactement, les $\sigma_{J}^{\tilde{G}'(s)}$ prennent leurs valeurs dans 
 $$U_{J}\otimes (D_{g\acute{e}om}({\cal O})\otimes Mes(M(F))^*)/Ann_{{\cal O}}^{\tilde{G}'(s)}$$
 mais, gr\^ace \`a ce que l'on a dit ci-dessus, on les pousse en des applications \`a valeurs dans
 $$U_{J}\otimes (D_{g\acute{e}om}({\cal O})\otimes Mes(M(F))^*)/Ann_{{\cal O}}^{\tilde{G}}.$$

 \ass{Proposition (\`a prouver)}{Pour tout $J\in {\cal J}_{\tilde{M}}^{\tilde{G}}(B_{{\cal O}})$, $\sigma^{\tilde{G}}_{J}$ prend ses valeurs dans 
 $$U_{J}\otimes (D_{g\acute{e}om}^{st}({\cal O})\otimes Mes(M(F))^*)/Ann_{{\cal O}}^{\tilde{G},st}.$$}

 \bigskip
 
 \subsection{Quelques formalit\'es}
  On consid\`ere un triplet $(G,\tilde{G},{\bf a})$ quasi-d\'eploy\'e et \`a torsion int\'erieure, un syst\`eme de fonctions $B$ comme en 1.9, un espace de Levi $\tilde{M}$ de $\tilde{G}$ et une classe de conjugaison stable semi-simple ${\cal O}$ dans $\tilde{M}(F)$. Consid\'erons des extensions compatibles
  $$1\to C_{\natural}\to G_{\natural}\stackrel{q}{\to} G\to 1 \text{ et } \tilde{G}_{\natural}\to \tilde{G}$$
  o\`u $C_{\natural}$ est un tore central induit et $\tilde{G}_{\natural}$ est encore \`a torsion int\'erieure. On fixe un caract\`ere $\lambda_{\natural}$ de $C_{\natural}(F)$. On note $\tilde{M}_{\natural}$ l'image r\'eciproque de $\tilde{M}$ dans $\tilde{G}_{\natural}$. On fixe une classe de conjugaison stable semi-simple ${\cal O}_{\natural}$ dans $M_{\natural}(F)$ se projetant sur ${\cal O}$. Le syst\`eme de fonctions $B$ se rel\`eve \`a $G_{\natural}(F)$. L'application $\alpha\mapsto \alpha\circ q$ est une bijection de $\Sigma(A_{M},B_{{\cal O}})$ sur $\Sigma(A_{M_{\natural}},B_{{\cal O}_{\natural}})$. Via cette bijection, les ensembles ${\cal J}_{\tilde{M}}^{\tilde{G}}(B_{{\cal O}})$ et ${\cal J}_{\tilde{M}_{\natural}}^{\tilde{G}_{\natural}}(B_{{\cal O}_{\natural}})$ s'identifient. Pour un \'el\'ement $J$ de cet ensemble, on a un espace $U_{J}$ de germes de fonctions sur $A_{M}(F)$ et un autre, notons-le $U_{\natural,J}$, de germes de fonctions sur $A_{M_{\natural}}(F)$. Il est clair que $U_{\natural,J}$ est form\'e des compos\'es $u\circ q$ pour $u\in U_{J}$. 
   On peut ainsi identifier ces deux espaces. Rappelons que l'on dispose d'un homomorphisme
 $$D_{g\acute{e}om}(\tilde{M}_{\natural}(F))\to D_{g\acute{e}om,\lambda_{\natural}}(\tilde{M}_{\natural}(F)),$$
 cf. 1.10 (3). Fixons une mesure de Haar sur $C_{\natural}(F)$, qui permet d'identifier $Mes(M_{\natural}(F))$ \`a $Mes(M(F))$. 
 On v\'erifie  sur sa d\'efinition que l'application $\rho_{J}^{\tilde{G}_{\natural}}$ se quotiente en un homomorphisme $\rho_{J,\lambda_{\natural}}^{\tilde{G}_{\natural}}$ de sorte que le diagramme suivant soit commutatif
 $$\begin{array}{ccc}D_{g\acute{e}om}({\cal O}_{\natural})\otimes Mes(M_{\natural}(F))^*&\stackrel{\rho_{J}^{\tilde{G}_{\natural}}}{\to}&U_{\natural,J}\otimes (D_{g\acute{e}om}({\cal O}_{\natural})\otimes Mes(M_{\natural}(F))^*)/Ann_{{\cal O}_{\natural}}^{\tilde{G}_{\natural}}\\ \downarrow&&\downarrow\\ D_{g\acute{e}om,\lambda_{\natural}}({\cal O})\otimes Mes(M(F))^*&\stackrel{\rho_{J,\lambda_{\natural}}^{\tilde{G}_{\natural}}}{\to}&U_{J}\otimes (D_{g\acute{e}om,\lambda_{\natural}}({\cal O})\otimes Mes(M(F))^*)/Ann_{{\cal O},\lambda_{\natural}}^{\tilde{G}_{\natural}}\\ \end{array}$$
 avec une d\'efinition naturelle du dernier annulateur. On voit que l'application $\rho_{J,\lambda_{\natural}}^{\tilde{G}_{\natural}}$ ne d\'epend pas des choix de ${\cal O}_{\natural}$ et de la mesure sur $C_{\natural}(F)$. Par des calculs analogues \`a ceux de la preuve de 1.10, on montre que l'application $\sigma_{J}^{\tilde{G}_{\natural}}$ se quotiente de m\^eme en une application
 $$\sigma_{J,\lambda_{\natural}}^{\tilde{G}_{\natural}}:D_{g\acute{e}om,\lambda_{\natural}}^{st}({\cal O})\otimes Mes(M(F))^* \to U_{J}\otimes (D_{g\acute{e}om,\lambda_{\natural}}({\cal O})\otimes Mes(M(F))^*)/Ann_{{\cal O},\lambda_{\natural}}^{\tilde{G}_{\natural}}.$$
 
 Consid\'erons d'autres extensions
$$1\to C_{\flat}\to G_{\flat}\to G\to 1 \text{ et } \tilde{G}_{\flat}\to \tilde{G}$$
un caract\`ere $\lambda_{\flat}$ de $C_{\flat}(F)$ et une classe de conjugaison stable ${\cal O}_{\flat}$ v\'erifiant des conditions similaires. On renvoie \`a [II] 1.10 pour les notations utilis\'ees ci-dessous.  Supposons donn\'ee un caract\`ere $\lambda_{\natural,\flat}$ du produit fibr\'e $G_{\natural,\flat}(F)$ dont la restriction \`a $C_{\natural}(F)\times C_{\flat}(F)$ soit $\lambda_{\natural}\times \lambda_{\flat}^{-1}$. Supposons donn\'e une fonction non nulle $\tilde{\lambda}_{\natural,\flat}$ sur le produit fibr\'e $\tilde{G}_{\natural,\flat}(F)$ qui se transforme selon le caract\`ere $\lambda_{\natural,\flat}$.  A l'aide de cette fonction, on peut identifier comme en [II] 1.10 tous les espaces intervenant ci-dessus  relatifs  aux donn\'ees index\'ees par $\natural$ avec les espaces analogues relatifs aux donn\'ees index\'ees par $\flat$. On v\'erifie comme en [II] 1.10 que les applications $\rho_{J,\lambda_{\natural}}^{\tilde{G}_{\natural}}$ et $\sigma_{J,\lambda_{\natural}}^{\tilde{G}_{\natural}}$ s'identifient \`a $\rho_{J,\lambda_{\flat}}^{\tilde{G}_{\flat}}$ et $\sigma_{J,\lambda_{\flat}}^{\tilde{G}_{\flat}}$. 

On peut aussi remplacer ${\cal O}$ par une r\'eunion finie de classes de conjugaison stable semi-simples, pourvu qu'elles v\'erifient la condition de la remarque 3.4(2)

 Comme application, consid\'erons un triplet $(G,\tilde{G},{\bf a})$, un espace de Levi $\tilde{M}$ de $\tilde{G}$, une classe de conjugaison stable semi-simple ${\cal O}$ dans $\tilde{M}(F)$, une donn\'ee endoscopique ${\bf M}'=(M',{\cal M}',\tilde{\zeta})$ de $(M,\tilde{M},{\bf a})$ qui est elliptique et relevante et un \'el\'ement $\tilde{s}\in \tilde{\zeta}Z(\hat{M})^{\Gamma_{F},\hat{\theta}}/Z(\hat{G})^{\Gamma_{F},\hat{\theta}}$. On pose simplement ${\bf G}'={\bf G}'(\tilde{s})$, $G'=G'(\tilde{s})$ etc... On suppose ${\bf G}'$ elliptique. On note ${\cal O}'$ la r\'eunion finie des classes de conjugaison stable semi-simples dans $\tilde{M}'(F)$ qui correspondent \`a ${\cal O}$.

 On consid\`ere plus pr\'ecis\'ement les trois cas suivants:
 
(1) $(G,\tilde{G},{\bf a})$ est quelconque, on munit $\tilde{G}'(F)$ du syst\`eme de fonctions $B^{\tilde{G}}$ que l'on note simplement $B$;
 
 (2) $(G,\tilde{G},{\bf a})$ est quasi-d\'eploy\'e et \`a torsion int\'erieure; on suppose donn\'e un syst\`eme de fonctions $B$ sur $\tilde{G}(F)$, dont on d\'eduit un tel syst\`eme sur $\tilde{G}'(F)$ que l'on note encore $B$;
 
 (3) $G=\tilde{G}$ et ${\bf a}=1$; on suppose ${\cal O}=\{1\}$; on suppose donn\'ee une fonction $B$ sur $G(F)$ comme en 1.8, dont on d\'eduit une telle fonction sur $G'(F)$ que l'on note encore $B$.
 
 Fixons des donn\'ees auxiliaires $G'_{1},...,\Delta_{1}$ pour ${\bf G}'$. Pour $J\in {\cal J}_{\tilde{M}'}^{\tilde{G}'}(B_{{\cal O}'})$, on d\'efinit comme ci-dessus les applications
 $$\rho_{J,\lambda_{1}}^{\tilde{G}'_{1}}:D_{g\acute{e}om,\lambda_{1}}({\cal O}')\otimes Mes(M'(F))^* \to U_{J}\otimes (D_{g\acute{e}om,\lambda_{1}}({\cal O}')\otimes Mes(M'(F))^*)/Ann_{{\cal O}',\lambda_{1}}^{\tilde{G}'_{1}}$$
 et
  $$\sigma_{J,\lambda_{1}}^{\tilde{G}'_{1}}:D_{g\acute{e}om,\lambda_{1}}^{st}({\cal O}')\otimes Mes(M'(F))^* \to U_{J}\otimes (D_{g\acute{e}om,\lambda_{1}}({\cal O}')\otimes Mes(M'(F))^*)/Ann_{{\cal O}',\lambda_{1}}^{\tilde{G}'_{1}}.$$
  Quand on fait varier les donn\'ees auxiliaires, ces application se recollent en ces applications
  $$\rho_{J}^{{\bf G}'}:D_{g\acute{e}om}({\bf M}',{\cal O}')\otimes Mes(M'(F))^*\to U_{J}\otimes (D_{g\acute{e}om}({\bf M}',{\cal O}')\otimes Mes(M'(F))^*)/Ann_{{\cal O}'}^{{\bf G}'}$$
 et
 $$ \sigma_{J}^{{\bf G}'}:D_{g\acute{e}om}^{st}({\bf M}',{\cal O}')\otimes Mes(M'(F))^*\to U_{J}\otimes (D_{g\acute{e}om}({\bf M}',{\cal O}')\otimes Mes(M'(F))^*)/Ann_{{\cal O}'}^{{\bf G}'}.$$
 On a adapt\'e les notations de fa\c{c}on \'evidente.

 Pla\c{c}ons-nous sous les hypoth\`eses de (2) et supposons de plus que ${\bf M}'$ est la donn\'ee "maximale" ${\bf M}$. Dans ce cas, $D_{g\acute{e}om}({\bf M}',{\cal O}')$ s'identifie \`a $D_{g\acute{e}om}({\cal O}')$. On v\'erifie en reprenant les d\'efinitions que $\rho_{J}^{{\bf G}'}$ s'identifie \`a $\rho_{J}^{G'}$. Cette propri\'et\'e se propage formellement: $ \sigma_{J}^{{\bf G}'}$ s'identifie \`a $\sigma_{J}^{G'}$.
La formule (1) du paragraphe 3.5 se r\'ecrit
 $$(4) \qquad \sigma^{\tilde{G}}_{J}=\rho_{J,st}-\sum_{s\in Z(\hat{M})^{\Gamma_{F}}/Z(\hat{G})^{\Gamma_{F}}; s\not=1, J\in {\cal J}_{\tilde{M}}^{\tilde{G}'(s)}(B_{{\cal O}})}i_{\tilde{M}}(\tilde{G},\tilde{G}'(s))transfert(\sigma_{J}^{{\bf G}'(s)}).$$
 Le transfert est ici l'isomorphisme  naturel de $D^{st}_{g\acute{e}om}({\bf M})\otimes Mes(M(F))^*$ sur $D^{st}_{g\acute{e}om}(\tilde{M}(F))\otimes Mes(M(F))^*$.

 Dans le cas (3), on a mieux. On a ${\cal O}'=\{1\}$ et on peut choisir pour rel\`evement l'orbite ${\cal O}'_{1}=\{1\}$. Parce $C_{1}$ est induit, l'application $M'_{1}(F)\to M'(F)$ est surjective. On en d\'eduit ais\'ement que les homomorphismes naturels
 $$D_{unip,\lambda_{1}}(M'_{1}(F))\leftarrow D_{unip}(M'_{1}(F))\to D_{unip}(M'(F))$$
 sont des isomorphismes. On en d\'eduit un isomorphisme $D_{unip}({\bf M}')\simeq D_{unip}(M'(F))$.
 Ainsi, les applications $\rho_{J}^{{\bf G}'}$ et $\sigma_{J}^{{\bf G}'}$ s'identifient \`a des homomorphismes
 $$\rho_{J}^{{\bf G}'}:D_{unip}(M'(F))\otimes Mes(M'(F))^*\to U_{J}\otimes (D_{unip}(M'(F))\otimes Mes(M'(F))^*)/Ann_{unip}^{G'}$$
 et
  $$\sigma_{J}^{{\bf G}'}:D_{unip}^{st}(M'(F))\otimes Mes(M'(F))^*\to U_{J}\otimes (D_{unip}(M'(F))\otimes Mes(M'(F))^*)/Ann_{unip}^{G'}.$$
  En reprenant les d\'efinitions, on voit que $\rho_{J}^{{\bf G}'}$ s'identifie \`a $\rho_{J}^{G'}$. Cette propri\'et\'e se propage: $\sigma_{J}^{{\bf G}'}$ s'identifie \`a $\sigma_{J}^{G'}$.
 
 \bigskip
 
 \subsection{D\'eveloppement des int\'egrales orbitales pond\'er\'ees stables}
 On suppose $(G,\tilde{G},{\bf a})$ quasi-d\'eploy\'e et \`a torsion int\'erieure. On fixe un syst\`eme de fonctions $B$ comme en 1.9. Soient $\tilde{M}$ un espace de Levi de $\tilde{G}$ et ${\cal O}$ une classe de conjugaison stable semi-simple dans $\tilde{M}(F)$. 
 
 \ass{Proposition}{(i) Pour tout $\boldsymbol{\delta}\in D^{st}_{g\acute{e}om}({\cal O})\otimes Mes(M(F))^*$ et pour tout ${\bf f}\in I(\tilde{G}(F))\otimes Mes(G(F))$, le germe en $1$ de la fonction
 $$a\mapsto S_{\tilde{M}}^{\tilde{G}}(a\boldsymbol{\delta},{\bf f}),$$
 qui est d\'efinie pour tout $a\in A_{M}(F)$ en position g\'en\'erale, est \'equivalent \`a
 $$\sum_{J\in {\cal J}_{\tilde{M}}^{\tilde{G}}(B_{{\cal O}}) }I^{\tilde{G}}(\sigma_{J}^{\tilde{G}}(\boldsymbol{\delta},a)^{\tilde{G}},{\bf f})$$
 $$+\sum_{\tilde{L}\in {\cal L}(\tilde{M}), \tilde{L}\not=\tilde{G}}\sum_{J\in {\cal J}_{\tilde{M}}^{\tilde{L}}(B_{{\cal O}}) }S_{\tilde{L}}^{\tilde{G}}(\sigma^{\tilde{L}}_{J}(\boldsymbol{\delta},a)^{\tilde{L}},B,{\bf f}).$$
 
 (ii) Supposons v\'erifi\'ee la proposition 3.5. Alors le d\'eveloppement pr\'ec\'edent prend la forme
 $$\sum_{\tilde{L}\in {\cal L}(\tilde{M})}\sum_{J\in {\cal J}_{\tilde{M}}^{\tilde{L}}(B_{{\cal O}})}S_{\tilde{L}}^{\tilde{G}}(\sigma^{\tilde{L}}_{J}(\boldsymbol{\delta},a)^{\tilde{L}},B,{\bf f}).$$}
 
 Preuve. Notons qu'en vertu de 3.5(2) et de nos hypoth\`eses de r\'ecurrence, les termes $\sigma_{J}^{\tilde{L}}(\boldsymbol{\delta},a)^{\tilde{L}}$ sont stables si $\tilde{L}\not=\tilde{G}$. Les termes de la formule du (i) ont donc un sens. Evidemment, si les  termes $\sigma^{\tilde{G}}_{J}(\boldsymbol{\delta},a)^{\tilde{G}}$ sont stables eux-aussi, on peut remplacer les int\'egrales orbitales figurant dans cette formule par des int\'egrales orbitales stables. Donc (ii) r\'esulte imm\'ediatement de (i).
 
On part de la d\'efinition
$$(1)   \qquad S_{\tilde{M}}^{\tilde{G}}(a\boldsymbol{\delta},{\bf f})=I_{\tilde{M}}^{\tilde{G}}(a\boldsymbol{\delta},f)-\sum_{s\in Z(\hat{M})^{\Gamma_{F}}/Z(\hat{G})^{\Gamma_{F}};s\not=1}i_{\tilde{M}}(\tilde{G},\tilde{G}'(s))S_{{\bf M}}^{{\bf G}'(s)}(a\boldsymbol{\delta},{\bf f}^{{\bf G}'(s)}).$$
 La proposition 3.4 nous fournit le d\'eveloppement du premier terme: $I_{\tilde{M}}^{\tilde{G}}(a\boldsymbol{\delta},f)$ est \'equivalent \`a
 $$(2) \qquad\sum_{\tilde{L}\in {\cal L}(\tilde{M})} \sum_{J\in {\cal J}_{\tilde{M}}^{\tilde{L}}(B_{{\cal O}})}I_{\tilde{L}}^{\tilde{G}}(\rho^{\tilde{L}}_{J}(\boldsymbol{\delta},a)^{\tilde{L}},B,{\bf f}).$$
 Pour tout $s\in Z(\hat{M})^{\Gamma_{F}}/Z(\hat{G})^{\Gamma_{F}}$,  avec  $s\not=1$, on peut d\'evelopper le terme $S_{{\bf M}}^{{\bf G}'(s)}(a\boldsymbol{\delta},{\bf f}^{{\bf G}'(s)})$ par la proposition que l'on cherche \`a prouver, appliqu\'ee \`a ${\bf G}'(s)$. On passe sur les formalit\'es permettant d'appliquer cette proposition  \`a une telle donn\'ee plut\^ot qu'\`a un espace $\tilde{G}'(s)$. On obtient que la somme du membre de droite de (1) est \'equivalente \`a
 
$$\sum_{s\in Z(\hat{M})^{\Gamma_{F}}/Z(\hat{G})^{\Gamma_{F}},s\not=1}i_{\tilde{M}}(\tilde{G},\tilde{G}'(s))\sum_{\tilde{L}'_{s}\in {\cal L}^{\tilde{G}'(s)}(\tilde{M})}$$
$$\sum_{J\in {\cal J}_{\tilde{M}}^{\tilde{L}'_{s}}(B_{{\cal O}}) }S_{{\bf L}'(s)}^{{\bf G}'(s)}(\sigma^{{\bf L}'(s)}_{J}(\boldsymbol{\delta},a)^{{\bf L}'(s)},B,{\bf f}^{{\bf G}'(s)}).$$ 
La notation ${\bf L}'(s)$ est la m\^eme que dans la preuve de 2.5. 
   Comme dans la preuve de la proposition 2.5, on regroupe les couples $(s,\tilde{L}'_{s})$ intervenant selon l'espace de Levi $\tilde{L}$ qu'ils d\'eterminent par l'\'egalit\'e ${\cal A}_{\tilde{L}}={\cal A}_{\tilde{L}'_{s}}$. L'expression pr\'ec\'edente devient
$$\sum_{\tilde{L}\in {\cal L}(\tilde{M})}\sum_{s\in Z(\hat{M})^{\Gamma_{F}}/Z(\hat{L})^{\Gamma_{F}}, {\bf L}'(s)\text{ elliptique}}\sum_{J\in {\cal J}_{\tilde{M}}^{\tilde{L}'(s)}(B_{{\cal O}}) }$$
$$\sum_{t\in sZ(\hat{L})^{\Gamma_{F}}/Z(\hat{G})^{\Gamma_{F}},t\not=1} i_{\tilde{M}}(\tilde{G},\tilde{G}'(t))S_{{\bf L}'(s)}^{{\bf G}'(t)}(\sigma_{J}^{{\bf L}'(s)}(\boldsymbol{\delta},a)^{{\bf L}'(s)},B,{\bf f}^{{\bf G}'(t)}).$$
 Comme en 2.5, on a l'\'egalit\'e 
$$ i_{\tilde{M}}(\tilde{G},\tilde{G}'(t))= i_{\tilde{M}}(\tilde{L},\tilde{L}'(s)) i_{\tilde{L}'(s)}(\tilde{G},\tilde{G}'(t))$$
et l'expression devient
$$(3) \qquad \sum_{\tilde{L}\in {\cal L}(\tilde{M})}\sum_{s\in Z(\hat{M})^{\Gamma_{F}}/Z(\hat{L})^{\Gamma_{F}}} i_{\tilde{M}}(\tilde{L},\tilde{L}'(s))\sum_{J\in {\cal J}_{\tilde{M}}^{\tilde{L}'(s)}(B_{{\cal O}}) }$$
$$\sum_{t\in sZ(\hat{L})^{\Gamma_{F}}/Z(\hat{G})^{\Gamma_{F}}, t\not=1} i_{\tilde{L}'(s)}(\tilde{G},\tilde{G}'(t))S_{{\bf L}'(s)}^{{\bf G}'(t)}(\sigma^{{\bf L}'(s)}_{J}(\boldsymbol{\delta},a)^{{\bf L}'(s)},B,{\bf f}^{{\bf G}'(t)}).$$
Fixons $\tilde{L}$, $s$ et $J$ et \'etudions la somme int\'erieure en $t$. Supposons d'abord $\tilde{L}\not=\tilde{G}$ et $\tilde{L}\not=\tilde{M}$.      On sait  par r\'ecurrence que les termes $\sigma^{{\bf L}'(s)}_{J}(\boldsymbol{\delta},a)^{{\bf L}'(s)}$ sont stables. Si $s\not=1$, la somme n'est autre que $I_{\tilde{L}}^{\tilde{G},{\cal E}}({\bf L}'(s),\sigma^{{\bf L}'(s)}_{J}(\boldsymbol{\delta},a)^{{\bf L}'(s)},B,{\bf f})$ ou encore \`a
$$I_{\tilde{L}}^{\tilde{G},{\cal E}}(transfert(\sigma^{{\bf L}'(s)}_{J}(\boldsymbol{\delta},a)^{{\bf L}'(s)}),B,{\bf f}).$$
   Le transfert commute \`a l'induction. Donc 
 $$transfert(\sigma^{{\bf L}'(s)}_{J}(\boldsymbol{\delta},a)^{{\bf L}'(s)})=(transfert(\sigma_{J}^{{\bf L}'(s)}(\boldsymbol{\delta},a)))^{\tilde{L}}.$$
 Puisque $\tilde{L}\not=\tilde{M}$, on peut appliquer le th\'eor\`eme 1.16 et l'expression ci-dessus devient $I_{\tilde{L}}^{\tilde{G}}((transfert(\sigma^{{\bf L}'(s)}_{J}(\boldsymbol{\delta},a)))^{\tilde{L}},{\bf f})$. 
  Si $s=1$, la somme en $t$ n'est pas tout-\`a-fait 
  $$I_{\tilde{L}}^{\tilde{G},{\cal E}}({\bf L}'(s),\sigma^{{\bf L}'(s)}_{J}(\boldsymbol{\delta},a)^{{\bf L}'(s)},B,{\bf f})$$
   car il manque le terme $t=1$.  On a $L'(s)=L$ puisque $s=1$. Le terme manquant est par d\'efinition $S_{\tilde{L}}^{\tilde{G}}(\sigma_{J}^{\tilde{L}}(\boldsymbol{\delta},a)^{\tilde{L}},B,{\bf f})$ .  On obtient que la contribution de $\tilde{L}$ \`a l'expression (3) est
  $$\sum_{s\in Z(\hat{M})^{\Gamma_{F}}/Z(\hat{L})^{\Gamma_{F}}} i_{\tilde{M}}(\tilde{L},\tilde{L}'(s))\sum_{J\in {\cal J}_{\tilde{M}}^{\tilde{L}'(s)}(B_{{\cal O}}) }I_{\tilde{L}}^{\tilde{G}}(transfert(\sigma_{J}^{{\bf L}'(s)}(\boldsymbol{\delta},a))^{\tilde{L}},B,{\bf f})$$
  $$-\sum_{J\in {\cal J}_{\tilde{M}}^{\tilde{L}}(B_{{\cal O}}) }S_{\tilde{L}}^{\tilde{G}}(\sigma^{\tilde{L}}_{J}(\boldsymbol{\delta},a)^{\tilde{L}},B,{\bf f}).$$
  On peut r\'ecrire la premi\`ere somme sous la forme
  $$\sum_{J\in {\cal J}_{\tilde{M}}^{\tilde{L}}(B_{{\cal O}})} \sum_{s\in Z(\hat{M})^{\Gamma_{F}}/Z(\hat{L})^{\Gamma_{F}}; J\in {\cal J}_{\tilde{M}}^{\tilde{L}'(s)}(B_{{\cal O}})} i_{\tilde{M}}(\tilde{L},\tilde{L}'(s))I_{\tilde{L}}^{\tilde{G}}((transfert(\sigma^{{\bf L}'(s)}_{J}(\boldsymbol{\delta},a)))^{\tilde{L}},B,{\bf f}).$$
  Or 
 $$\sum_{s\in Z(\hat{M})^{\Gamma_{F}}/Z(\hat{L})^{\Gamma_{F}}; J\in {\cal J}_{\tilde{M}}^{\tilde{L}'(s)}(B_{\cal O})} i_{\tilde{M}}(\tilde{L},\tilde{L}'(s))transfert( \sigma_{J}^{{\bf L}'(s)}(\boldsymbol{\delta},a))=\rho_{J}^{\tilde{L}}(\boldsymbol{\delta},a)$$
 d'apr\`es la d\'efinition  3.5(1).  Donc la contribution de $\tilde{L}$ \`a l'expression (3) se r\'eduit \`a
 $$(4) \qquad \sum_{J\in {\cal J}_{\tilde{M}}^{\tilde{L}}(B_{{\cal O}}) }\left(I_{\tilde{L}}^{\tilde{G}}(\rho_{J}^{\tilde{L}}(\boldsymbol{\delta},a)^{\tilde{L}},B,{\bf f})-S_{\tilde{L}}^{\tilde{G}}(\sigma^{\tilde{L}}_{J}(\boldsymbol{\delta},a)^{\tilde{L}},B,{\bf f})\right).$$
 Supposons maintenant $\tilde{L}=\tilde{G}$. La somme en $t$ est vide si $s=1$ et est r\'eduite au terme $t=s$ si $s\not=1$. On a dans ce cas $i_{\tilde{G}'(s)}(\tilde{G},\tilde{G}'(s))=1$. La contribution de $\tilde{G}$ s'\'ecrit plus simplement
 $$\sum_{s\in Z(\hat{M})^{\Gamma_{F}}/Z(\hat{G})^{\Gamma_{F}},s\not=1}i_{\tilde{M}}(\tilde{G},\tilde{G}'(s))\sum_{J\in {\cal J}_{\tilde{M}}^{\tilde{G}'(s)}(B_{{\cal O}}) }S^{{\bf G}'(s)}(\sigma^{{\bf G}'(s)}_{J}(\boldsymbol{\delta},a)^{{\bf G}'(s)},B,{\bf f}^{{\bf G}'(s)}).$$
 Les int\'egrales stables n'\'etant plus pond\'er\'ees, on n'a plus besoin de faire appel au th\'eor\`eme 1.16 pour obtenir les \'egalit\'es
$$ S^{{\bf G}'(s)}(\sigma^{{\bf G}'(s)}_{J}(\boldsymbol{\delta},a)^{{\bf G}'(s)},B,{\bf f}^{{\bf G}'(s)})=I^{\tilde{G}}(transfert(\sigma^{{\bf G}'(s)}_{J}(\boldsymbol{\delta},a)^{{\bf G}'(s)}),B,{\bf f})$$
$$=I^{\tilde{G}}((transfert(\sigma^{{\bf G}'(s)}_{J}(\boldsymbol{\delta},a)))^{\tilde{G}},B,{\bf f}).$$
La contribution de $\tilde{G}$ devient
$$\sum_{s\in Z(\hat{M})^{\Gamma_{F}}/Z(\hat{G})^{\Gamma_{F}},s\not=1}i_{\tilde{M}}(\tilde{G},\tilde{G}'(s))\sum_{J\in {\cal J}_{\tilde{M}}^{\tilde{G}'(s)}(B_{{\cal O}}) }I^{\tilde{G}}((transfert(\sigma^{{\bf G}'(s)}_{J}(\boldsymbol{\delta},a)))^{\tilde{G}},B,{\bf f}).$$
On la r\'ecrit comme ci-dessus
$$\sum_{J\in {\cal J}_{\tilde{M}}^{\tilde{G}}(B_{{\cal O}}) }\sum_{s\in Z(\hat{M})^{\Gamma_{F}}/Z(\hat{G})^{\Gamma_{F}},s\not=1, J\in {\cal J}_{\tilde{M}}^{\tilde{G}'(s)}(B_{{\cal O}})}i_{\tilde{M}}(\tilde{G},\tilde{G}'(s))I^{\tilde{G}}((transfert(\sigma^{{\bf G}'(s)}_{J}(\boldsymbol{\delta},a)))^{\tilde{G}},B,{\bf f}).$$
On a 
 $$\sum_{s\in Z(\hat{M})^{\Gamma_{F}}/Z(\hat{G})^{\Gamma_{F}}; J\in {\cal J}_{\tilde{M}}^{\tilde{G}'(s)}(B_{\cal O})} i_{\tilde{M}}(\tilde{G},\tilde{G}'(s))transfert( \sigma^{{\bf G}'(s)}_{J}(\boldsymbol{\delta},a))=\rho^{\tilde{G}}_{J}(\boldsymbol{\delta},a)-\sigma_{J}^{\tilde{G}}(\boldsymbol{\delta},a)$$
 d'apr\`es la d\'efinition  3.5(1).  La contribution de $\tilde{G}$ \`a l'expression (3) est donc
$$(5) \qquad \sum_{J\in {\cal J}_{\tilde{M}}^{\tilde{G}}(B_{{\cal O}}) }\left(I^{\tilde{G}}(\rho_{J}^{\tilde{G}}(\boldsymbol{\delta},a)^{\tilde{G}},B,{\bf f})-I^{\tilde{G}}(\sigma^{\tilde{G}}_{J}(\boldsymbol{\delta},a)^{\tilde{G}},B,{\bf f})\right).$$
Consid\'erons enfin l'espace $\tilde{L}=\tilde{M}$. La somme en $s$ dispara\^{\i}t de l'expression (3). La somme en $J$ se r\'eduit au terme $J=\emptyset$.  Les termes $\sigma_{\emptyset}^{{\bf M}}(\boldsymbol{\delta},a)$ sont tous \'egaux \`a $\boldsymbol{\delta}$.
La contribution de $\tilde{M}$ \`a l'expression (3) se r\'eduit \`a la somme en $t$, qui est alors
$$(6)\qquad I_{\tilde{M}}^{\tilde{G}}(\sigma_{\emptyset}^{\tilde{M}}(\boldsymbol{\delta},a),B,{\bf f})-S_{\tilde{M}}^{\tilde{G}}(\sigma_{\emptyset}^{\tilde{M}}(\boldsymbol{\delta},a),B,{\bf f})$$
par d\'efinition de ce dernier terme. 

Le membre de droite de (1) est \'equivalent \`a la diff\'erence entre  (2) et la somme de (4), (5) et (6). On voit que c'est l'expression du (i) de l'\'enonc\'e. Cela ach\`eve la preuve. $\square$

\bigskip

\subsection{Termes d'un d\'eveloppement endoscopique}
Soient $(G,\tilde{G},{\bf a})$ un triplet quelconque, $\tilde{M}$ un espace de Levi de $\tilde{G}$, ${\cal O}$ une classe de conjugaison stable semi-simple dans $\tilde{M}(F)$ et ${\bf M}'=(M',{\cal M}',\tilde{\zeta})$ une donn\'ee endoscopique   de $(M,\tilde{M},{\bf a})$, elliptique et relevante. On note ${\cal O}'$ la r\'eunion finie des classes de conjugaison stable semi-simples dans $\tilde{M}'(F)$ qui correspondent \`a ${\cal O}$. Rappelons que l'on a un homomorphisme $\xi:A_{\tilde{M}}\to A_{M'}$ dont se d\'eduit un isomorphisme $\xi:\mathfrak{a}_{\tilde{M}}\to \mathfrak{a}_{M'}$.

Soit $\tilde{s}\in \tilde{\zeta}Z(\hat{M})^{\Gamma_{F},\hat{\theta}}/Z(\hat{G})^{\Gamma_{F},\hat{\theta}}$. Posons simplement ${\bf G}'={\bf G}'(\tilde{s})$. On dispose de l'ensemble   $\Sigma^{G}(A_{\tilde{M}})$ de racines, que l'on peut voir comme des formes lin\'eaires sur $\mathfrak{a}_{\tilde{M}}$,  et on a d\'efini l'espace $\Sigma^{G'}(A_{M'},B^{\tilde{G}}_{{\cal O}'})$ de formes lin\'eaires sur $\mathfrak{a}_{M'}$.   Montrons que

(1) on a $\beta\circ \xi\in \Sigma^G(A_{\tilde{M}})$ pour tout $\beta\in \Sigma^{G'}(A_{M'},B^{\tilde{G}}_{{\cal O}})$.

  On fixe $\epsilon\in {\cal O}'$ et on applique la construction de la fonction $B_{\epsilon}^{\tilde{G}}$ faite en 1.11. Pour simplifier, on pose $B=B_{\epsilon}^{\tilde{G}}$. On suppose que les paires de Borel $(B,T)$ et $(B',T')$ sont telles que $\tilde{M}$ et $\tilde{M}'$ soient standard. L'isomorphisme $\xi$ ci-dessus se d\'eduit d'un isomorphisme $\xi:\mathfrak{t}^{\theta}\simeq \mathfrak{t}'$. Un \'el\'ement $\beta\in \Sigma^{G'}(A_{M'},B)$ est la restriction \`a $\mathfrak{a}_{M'}$ d'un \'el\'ement $\beta'\in \Sigma^{G'_{\epsilon}}(T',B)$. Alors $\beta\circ \xi$ est la restriction \`a $\mathfrak{a}_{\tilde{M}}$ de $\beta'\circ \xi$. L'\'el\'ement $\beta'$ est de la forme $\alpha'/B(\alpha')$, o\`u $\alpha'\in \Sigma^{G'_{\epsilon}}(T')$. On a rappel\'e en 1.11 la description de cet ensemble, que l'on a d\'ecompos\'e en cas (a) \`a (d). Dans le cas (a), on a $\alpha'=N\alpha$, $\alpha'\circ\xi=(N\alpha)\circ\xi=n_{\alpha}\alpha_{res}$, o\`u $\alpha_{res}$ est la restriction de $\alpha$ \`a $\mathfrak{t}^{\theta}$. Puisque $B(\alpha')=n_{\alpha}$, on obtient $\beta'\circ \xi=\alpha_{res}$. Dans le cas (b), on a $\alpha'=2N\alpha$, $\alpha'\circ\xi=2n_{\alpha}\alpha_{res}$. Puisque $B(\alpha')=2n_{\alpha}$, on a encore $\beta'\circ\xi=\alpha_{res}$. Dans le cas (c), on a encore $\alpha'\circ\xi=2n_{\alpha}\alpha_{res}$. Cette fois, $B(\alpha')=n_{\alpha}$, d'o\`u $\beta'\circ\xi=2\alpha_{res}$. Mais $\alpha$ est de type 2 donc $n_{\alpha}$ est pair et l'\'el\'ement $\bar{\alpha}=\alpha+\theta^{n_{\alpha}/2}(\alpha)$ est de type 3. On a $\bar{\alpha}_{res}=2\alpha_{res}$, donc $\beta'\circ\xi=\bar{\alpha}_{res}$. Dans le cas (d), on a $\alpha'=N\alpha$, $\alpha'\circ=n_{\alpha}\alpha_{res}$ et $B(\alpha')=2n_{\alpha}$, d'o\`u $\beta'\circ\xi=\alpha_{res}/2$. Mais $\alpha$ est de type 3 et il existe une racine $\underline{\alpha}$ de type 2 telle que $\alpha=\underline{\alpha}+\theta^{n_{\underline{\alpha}}/2}(\underline{\alpha})$. On a $\alpha_{res}=2\underline{\alpha}_{res}$, d'o\`u $\beta'\circ\xi=\underline{\alpha}_{res}$. Ainsi $\beta'\circ\xi$ est toujours la restriction \`a $\mathfrak{t}^{\theta}$ d'un \'el\'ement de $\Sigma^G(T)$. Il s'ensuit que $\beta\circ\xi$ est la restriction \`a $\mathfrak{a}_{\tilde{M}}$ d'un tel \'el\'ement. Donc $\beta\circ\xi\in \Sigma^G(A_{\tilde{M}})$. $\square$
 
Supposons ${\bf G}'$ elliptique.  L'application $\beta\mapsto \beta\circ\xi$ d\'efinit une injection $\Sigma^{G'}(A_{M'},B^{\tilde{G}}_{{\cal O}})\to \Sigma^G(A_{\tilde{M}})$. Il s'en d\'eduit une injection ${\cal J}_{\tilde{M}'}^{\tilde{G}'}(B^{\tilde{G}}_{{\cal O}})\to {\cal J}_{\tilde{M}}^{\tilde{G}}$, que nous noterons simplement $J'\mapsto J$.  Consid\'erons deux tels \'el\'ements $J'$ et $J$ tels que $J'\mapsto J$. Pour  $u'\in U_{J'}$, la fonction $u'\circ\xi$ appartient \`a $U_{J}$, d'o\`u une injection $U_{J'}\to U_{J}$.   L'application de transfert
 $$D^{st}_{g\acute{e}om}({\bf M}',{\cal O}')\otimes Mes(M'(F))^*\to D_{g\acute{e}om}({\cal O})\otimes Mes(M(F))^*$$
 commute \`a l'induction. On en d\'eduit qu'elle se factorise en une application
 $$(D^{st}_{g\acute{e}om}({\bf M}',{\cal O}')\otimes Mes(M'(F))^*)/Ann_{{\cal O}'}^{\tilde{G}',st}\to (D_{g\acute{e}om}({\cal O})\otimes Mes(M(F))^*)/Ann_{{\cal O}}^{\tilde{G}}.$$
 On l'appelle  encore transfert.
  
 Soit $J\in {\cal J}_{\tilde{M}}^{\tilde{G}}$. On va d\'efinir une application
 $$\rho_{J}^{\tilde{G},{\cal E}}({\bf M}'):D^{st}_{g\acute{e}om}({\bf M}',{\cal O}')\otimes Mes(M'(F))^*\to U_{J}\otimes (D_{g\acute{e}om}({\cal O})\otimes Mes(M(F))^*)/Ann_{{\cal O}}^{\tilde{G}}.$$
 Soit $\boldsymbol{\delta}\in D^{st}_{g\acute{e}om}({\bf M}',{\cal O}')\otimes Mes(M'(F))^*$. On peut consid\'erer la valeur $\rho_{J}^{\tilde{G},{\cal E}}({\bf M}',\boldsymbol{\delta})$ comme un germe d'application de $A_{\tilde{M}}(F)$ dans $(D_{g\acute{e}om}({\cal O})\otimes Mes(M(F))^*)/Ann_{{\cal O}}^{\tilde{G}}$ dont on note $\rho_{J}^{\tilde{G},{\cal E}}({\bf M}',\boldsymbol{\delta},a)$ la valeur en un point $a\in A_{\tilde{M}}(F)$ en position g\'en\'erale et proche de $1$.   On pose 
 $$\rho_{J}^{\tilde{G},{\cal E}}({\bf M}',\boldsymbol{\delta},a)=\sum_{\tilde{s}\in \tilde{\zeta}Z(\hat{M})^{\Gamma_{F},\hat{\theta}}/Z(\hat{G})^{\Gamma_{F},\hat{\theta}}}i_{\tilde{M}'}(\tilde{G},\tilde{G}'(\tilde{s}))$$
 $$\sum_{J'\in {\cal J}_{\tilde{M}'}^{\tilde{G}'(\tilde{s})}(B^{\tilde{G}}_{{\cal O}'}); J'\mapsto J}transfert(\sigma_{J'}^{{\bf G}'(\tilde{s})}(\boldsymbol{\delta},\xi(a))).$$
 Notons que la somme en $J'$ est vide ou r\'eduite \`a un seul \'el\'ement.  Les consid\'erations qui pr\'ec\'edent montrent que  $\rho_{J}^{\tilde{G},{\cal E}}({\bf M}')$ prend ses valeurs dans l'espace indiqu\'e. D'apr\`es nos hypoth\`eses de r\'ecurrence, les termes $\sigma_{J'}^{{\bf G}'(\tilde{s})}(\boldsymbol{\delta},\xi(a))$ sont stables et on peut bien les transf\'erer, sauf dans le cas o\`u $(G,\tilde{G},{\bf a})$ est quasi-d\'eploy\'e et \`a torsion int\'erieure et o\`u ${\bf M}'={\bf M}$. Dans ce cas,  le terme correspondant \`a $\tilde{s}=1$ pose probl\`eme. On le remplace simplement par $\sigma_{J}^{\tilde{G}}(\boldsymbol{\delta},a)$. Par d\'efinition de ce dernier terme, on a dans ce cas $\rho_{J}^{\tilde{G},{\cal E}} ({\bf M},\boldsymbol{\delta},a)=\rho_{J}^{\tilde{G}}(\boldsymbol{\delta},a)$. 
 
 \ass{Proposition (\`a prouver)}{Pour tout $J\in {\cal J}_{\tilde{M}}^{\tilde{G}}$, tout $\boldsymbol{\delta}\in D^{st}_{g\acute{e}om}({\bf M}',{\cal O}')\otimes Mes(M'(F))^*$ et tout $a\in A_{\tilde{M}}(F)$ en position g\'en\'erale et proche de $1$, on a l'\'egalit\'e
 $$\rho_{J}^{\tilde{G},{\cal E}}({\bf M}',\boldsymbol{\delta},a)=\rho^{\tilde{G}}_{J}(transfert(\boldsymbol{\delta}),a).$$}

 {\bf Variante.} Supposons $(G,\tilde{G},{\bf a})$ quasi-d\'eploy\'e et \`a torsion int\'erieure. Fixons un syst\`eme de fonction $B$. Il s'en d\'eduit un tel syst\`eme sur les espaces endoscopiques $\tilde{G}'(s)$ intervenant. En utilisant ces syst\`emes \`a la fois sur $\tilde{G}$ et sur ses espaces endoscopiques, on d\'efinit la variante $\rho_{J}^{\tilde{G},{\cal E}}({\bf M}')$ pour $J\in {\cal J}_{\tilde{M}}^{\tilde{G}}(B_{{\cal O}})$.
 
 {\bf Variante.} Supposons $G=\tilde{G}$, ${\bf a}=1$ et ${\cal O}=\{1\}$. Fixons une fonction $B$ comme en 1.8. On d\'efinit de m\^eme la variante $\rho_{J}^{G,{\cal E}}({\bf M}')$ pour $J\in {\cal J}_{M}^G(B)$. 
 
 \bigskip
 
 \subsection{D\'eveloppement des int\'egrales orbitales pond\'er\'ees endoscopiques}
 Les donn\'ees sont les m\^emes que dans le paragraphe pr\'ec\'edent.
 
 \ass{Proposition}{Pour tout $\boldsymbol{\delta}\in D^{st}_{g\acute{e}om}({\bf M}',{\cal O}')\otimes Mes(M'(F))^*$ et tout ${\bf f}\in I(\tilde{G}(F),\omega)\otimes Mes(G(F))$, le germe en $1$ de la fonction
 $$a\mapsto I_{\tilde{M}}^{\tilde{G},{\cal E}}({\bf M}',\xi(a)\boldsymbol{\delta},{\bf f}),$$ 
 qui est d\'efinie pour tout $a\in A_{\tilde{M}}(F)$ en position g\'en\'erale, est \'equivalent \`a
 $$\sum_{\tilde{L}\in {\cal L}(\tilde{M})}\sum_{J\in {\cal J}_{\tilde{M}}^{\tilde{L}}}I_{\tilde{L}}^{\tilde{G},{\cal E}}(\rho_{J}^{\tilde{L},{\cal E}}({\bf M}',\boldsymbol{\delta},a)^{\tilde{L}},{\bf f}).$$}
 
 Preuve. Supposons d'abord $(G,\tilde{G},{\bf a})$ quasi-d\'eploy\'e et \`a torsion int\'erieure et ${\bf M}'={\bf M}$. On v\'erifie qu'aucun terme ne change si on supprime les exposants ${\cal E}$. L'assertion est alors le (ii) de la proposition 3.2. On exclut ce cas.

Posons pour simplifier  $a'=\xi(a)$. On a l'\'egalit\'e 
  $$I_{\tilde{M}}^{\tilde{G},{\cal E}}({\bf M}',a'\boldsymbol{\delta},{\bf f})= \sum_{\tilde{s}\in \tilde{\zeta}Z(\hat{M})^{\Gamma_{F},\hat{\theta}}/Z(\hat{G})^{\Gamma_{F},\hat{\theta}}}i_{\tilde{M}'}(\tilde{G},\tilde{G}'(\tilde{s}))S_{{\bf M}'}^{{\bf G}'(\tilde{s})}(a'\boldsymbol{\delta},B^{\tilde{G}},{\bf f}^{{\bf G}'(\tilde{s})}).$$
 Calculons le germe en $1$ de cette expression, \`a \'equivalence pr\`es. D'apr\`es nos hypoth\`eses de r\'ecurrence, la proposition 3.7 est d\'emontr\'ee pour tous les termes intervenant ici.  En l'utilisant, on obtient
  $$\sum_{\tilde{s}\in \tilde{\zeta}Z(\hat{M})^{\Gamma_{F},\hat{\theta}}/Z(\hat{G})^{\Gamma_{F},\hat{\theta}}}i_{\tilde{M}'}(\tilde{G},\tilde{G}'(\tilde{s}))\sum_{\tilde{L}'_{\tilde{s}}\in {\cal L}^{\tilde{G}'(\tilde{s})}(\tilde{M}')}$$
  $$\sum_{J'\in {\cal J}_{\tilde{M}'}^{\tilde{G}'(\tilde{s})}(B^{\tilde{G}}_{{\cal O}'}) }S_{{\bf L}'_{\tilde{s}}}^{{\bf G}'(\tilde{s})}(\sigma_{J'}^{{\bf L}'(\tilde{s})}(\boldsymbol{\delta},a')^{{\bf L}'_{\tilde{s}}},B^{\tilde{G}},{\bf f}^{{\bf G}'(\tilde{s})}).$$
 On regroupe les $(\tilde{s},\tilde{L}'_{\tilde{s}})$ selon l'espace de Levi $\tilde{L}$ d\'etermin\'e par ${\cal A}_{\tilde{L}}={\cal A}_{\tilde{L}'_{\tilde{s}}}$. Comme dans les d\'emonstrations pr\'ec\'edentes, $\tilde{L}'_{\tilde{s}}$ devient $\tilde{L}'(\tilde{s})$. On obtient
  $$\sum_{\tilde{L}\in {\cal L}(\tilde{M})}\sum_{\tilde{s}\in \tilde{\zeta}Z(\hat{M})^{\Gamma_{F},\hat{\theta}}/Z(\hat{L})^{\Gamma_{F},\hat{\theta}}, {\bf L}'(\tilde{s})\text{ elliptique}}\sum_{J'\in {\cal J}_{\tilde{M}'}^{\tilde{L}'(\tilde{s})}(B^{\tilde{G}}_{{\cal O}'}) }$$
  $$\sum_{\tilde{t}\in \tilde{s}Z(\hat{L})^{\Gamma_{F},\hat{\theta}}/Z(\hat{G})^{\Gamma_{F},\hat{\theta}}}i_{\tilde{M}'}(\tilde{G},\tilde{G}'(\tilde{t}))S_{{\bf L}'(\tilde{s})}^{{\bf G}'(\tilde{t})}(\sigma^{{\bf L}'(\tilde{s})}_{J'}(\boldsymbol{\delta},a')^{{\bf L}'(\tilde{s})},B^{\tilde{G}},{\bf f}^{{\bf G}'(\tilde{t})}).$$
  On peut d\'ecomposer la somme en $J'$ en une somme en les  $J\in {\cal J}_{\tilde{M}}^{\tilde{L}}$  et une somme en les $J'$ tels que $J'\mapsto J$. D'autre part, on a encore
  $$i_{\tilde{M}'}(\tilde{G},\tilde{G}'(\tilde{t}))=i_{\tilde{M}'}(\tilde{L},\tilde{L}'(\tilde{s}))i_{\tilde{L}'(\tilde{s})}(\tilde{G},\tilde{G}'(\tilde{t})),$$
  et l'expression devient
   $$\sum_{\tilde{L}\in {\cal L}(\tilde{M})}\sum_{J\in {\cal J}_{\tilde{M}}^{\tilde{L}}}\sum_{\tilde{s}\in \tilde{\zeta}Z(\hat{M})^{\Gamma_{F},\hat{\theta}}/Z(\hat{L})^{\Gamma_{F},\hat{\theta}}}i_{\tilde{M}'}(\tilde{L},\tilde{L}'(\tilde{s}))\sum_{J'\in {\cal J}_{\tilde{M}'}^{\tilde{L}'(\tilde{s})}(B^{\tilde{G}}_{{\cal O}'}),  J'\mapsto J} $$
   $$\sum_{\tilde{t}\in \tilde{s}Z(\hat{L})^{\Gamma_{F},\hat{\theta}}/Z(\hat{G})^{\Gamma_{F},\hat{\theta}}}i_{\tilde{L}'(\tilde{s})}(\tilde{G},\tilde{G}'(\tilde{t}))S_{{\bf L}'(\tilde{s})}^{{\bf G}'(\tilde{t})}(\sigma_{J'}^{{\bf L}'(\tilde{s})}(\boldsymbol{\delta},a')^{{\bf L}'(\tilde{s})},B^{\tilde{G}},{\bf f}^{{\bf G}'(\tilde{t})}).$$
   La somme en $\tilde{t}$ n'est autre que $I_{\tilde{L}}^{\tilde{G},{\cal E}}({\bf L}'(\tilde{s}),\sigma_{J'}^{{\bf L}'(\tilde{s})}(\boldsymbol{\delta},a')^{{\bf L}'(\tilde{s})},{\bf f})$, ou encore 
   $$I_{\tilde{L}}^{\tilde{G},{\cal E}}(transfert(\sigma_{J'}^{{\bf L}'(\tilde{s})}(\boldsymbol{\delta},a')^{{\bf L}'(\tilde{s})}),{\bf f}).$$
    On a l'\'egalit\'e
 $$transfert(\sigma_{J'}^{{\bf L}'(\tilde{s})}(\boldsymbol{\delta},a')^{{\bf L}'(\tilde{s})})=\left(transfert(\sigma_{J'}^{{\bf L}'(\tilde{s})}(\boldsymbol{\delta},a'))\right)^{\tilde{L}}.$$
 L'expression devient
 $$  \sum_{\tilde{L}\in {\cal L}(\tilde{M})}\sum_{J\in {\cal J}_{\tilde{M}}^{\tilde{L}} }\sum_{\tilde{s}\in \tilde{\zeta}Z(\hat{M})^{\Gamma_{F},\hat{\theta}}/Z(\hat{L})^{\Gamma_{F},\hat{\theta}}}i_{\tilde{M}'}(\tilde{L},\tilde{L}'(\tilde{s}))$$
 $$\sum_{J'\in {\cal J}_{\tilde{M}'}^{\tilde{L}'(\tilde{s})}(B^{\tilde{G}}_{{\cal O}'}),  J'\mapsto J}I_{\tilde{L}}^{\tilde{G},{\cal E}}( \left(transfert(\sigma_{J'}^{{\bf L}'(\tilde{s})}(\boldsymbol{\delta},a'))\right)^{\tilde{L}},{\bf f}).$$
 Il suffit d'appliquer la d\'efinition de 3.8 pour obtenir  transformer cette expression  en celle de l'\'enonc\'e. $\square$

  {\bf Variante.} Supposons $(G,\tilde{G},{\bf a})$ quasi-d\'eploy\'e et \`a torsion int\'erieure. Fixons un syst\`eme de fonction $B$.  On a une proposition analogue. La formule de l'\'enonc\'e prend la forme
  $$\sum_{\tilde{L}\in {\cal L}(\tilde{M})}\sum_{J\in {\cal J}_{\tilde{M}}^{\tilde{L}}(B_{{\cal O}})}I_{\tilde{L}}^{\tilde{G},{\cal E}}(\rho_{J}^{\tilde{L},{\cal E}}({\bf M}',\boldsymbol{\delta},a)^{\tilde{L}},B,{\bf f}).$$
 
 {\bf Variante.} Supposons $G=\tilde{G}$, ${\bf a}=1$ et ${\cal O}=\{1\}$. Fixons une fonction $B$ comme en 1.8.  On a une proposition analogue. La formule de l'\'enonc\'e prend la forme
  $$\sum_{L\in {\cal L}(M)}\sum_{J\in {\cal J}_{M}^{L}(B)}I_{L}^{G,{\cal E}}(\rho_{J}^{L,{\cal E}}({\bf M}',\boldsymbol{\delta},a)^{L},B,{\bf f}).$$

 \bigskip

 \subsection{Termes $\rho_{J}$ et induction}
 Soient $(G,\tilde{G},{\bf a})$ un triplet quelconque, $\tilde{M}$ un espace de Levi de $\tilde{G}$, $\tilde{R}$ un espace de Levi de $\tilde{M}$ et ${\cal O}$ une classe de conjugaison semi-simple dans $\tilde{R}(F)$. On note ${\cal O}^{\tilde{M}}$ la classe de conjugaison dans $\tilde{M}(F)$ qui contient ${\cal O}$. On rappelle l'homomorphisme d'induction
 $$\begin{array}{ccc}D_{g\acute{e}om}({\cal O})\otimes Mes(R(F))^*&\to&D_{g\acute{e}om}({\cal O}^{\tilde{M}})\otimes Mes(M(F))^*\\ \boldsymbol{\gamma}&\mapsto& \boldsymbol{\gamma}^{\tilde{M}}\\ \end{array}$$
 
 Rappelons que ${\cal J}_{\tilde{M}}^{\tilde{G}}$  est l'ensemble des classes d'\'equivalence d'ensembles $\underline{\alpha}=\{\alpha_{1},...,\alpha_{n}\}$ o\`u les $\alpha_{i}$ sont des \'el\'ements lin\'eairement ind\'ependants de $\Sigma^G(A_{\tilde{M}})$ et o\`u $n=a_{\tilde{M}}-a_{\tilde{G}}$. Deux ensembles sont \'equivalents s'ils engendrent le m\^eme ${\mathbb Z}$-module. Soit $\tilde{L}\in {\cal L}(\tilde{R})$ tel que
 $$(1) \qquad {\cal A}_{\tilde{R}}^{\tilde{G}}={\cal A}_{\tilde{R}}^{\tilde{M}}\oplus {\cal A}_{\tilde{R}}^{\tilde{L}}.$$
 Alors de l'injection $A_{\tilde{M}}\to A_{\tilde{R}}$ se d\'eduit une application injective $\Sigma^L(A_{\tilde{R}})\to \Sigma^G(A_{\tilde{M}})$. Il s'en d\'eduit une injection ${\cal J}_{\tilde{R}}^{\tilde{L}}\to {\cal J}_{\tilde{M}}^{\tilde{G}}$ par laquelle on identifie le premier ensemble \`a un sous-ensemble du second. A un \'el\'ement $J\in {\cal J}_{\tilde{R}}^{\tilde{L}}$ sont associ\'es deux espaces $U_{J}$, l'un de germes de fonctions sur $A_{\tilde{R}}(F)$, l'autre de germes de fonctions sur $A_{\tilde{M}}(F)$. Ce dernier est l'ensemble des restrictions des \'el\'ements du premier.

 \ass{Lemme}{Soient $J\in {\cal J}_{\tilde{M}}^{\tilde{G}}$, $\boldsymbol{\gamma}\in D_{g\acute{e}om}({\cal O})\otimes Mes(R(F))^*$ et $a\in A_{\tilde{M}}(F)$ en position g\'en\'erale et assez proche de $1$.   On a l'\'egalit\'e
 $$\rho_{J}^{\tilde{G}}(\boldsymbol{\gamma}^{\tilde{M}},a)=\sum_{\tilde{L}\in {\cal L}(\tilde{R}), J\in {\cal J}_{\tilde{R}}^{\tilde{L}}}d_{\tilde{R}}^{\tilde{G}}(\tilde{M},\tilde{L})\rho_{J}^{\tilde{L}}(\boldsymbol{\gamma},a)^{\tilde{M}}.$$}

 Preuve. On fixe des mesures de Haar  sur tous les groupes intervenant. Par lin\'earit\'e, on peut supposer que $\boldsymbol{\gamma}$ est une int\'egrale orbitale associ\'ee \`a un \'el\'ement $u\eta\in \tilde{R}(F)$, o\`u $\eta\in \tilde{R}(F)$ est semi-simple et $u\in R_{\eta}(F)$ est unipotent. Alors $\boldsymbol{\gamma}^{\tilde{M}}$ est une combinaison lin\'eaire
 $$\sum_{j=1,...,k}c_{k}I^{\tilde{M}}(v_{j}\eta,\omega,.)$$
  d'int\'egrales orbitales associ\'ees \`a  des \'el\'ements $v_{j}\eta$,  o\`u $v_{j}\in M_{\eta}(F)$ appartient \`a l'orbite induite de celle de $u$. On pose $\boldsymbol{\gamma}_{j}=I^{\tilde{M}}(v_{j}\eta,\omega,.)$. En appliquant la d\'efinition 3.2(5), on obtient
 $$(2) \qquad \rho_{J}^{\tilde{G}}(\boldsymbol{\gamma}^{\tilde{M}},a)=\sum_{j=1,...,k}c_{j}\sum_{\underline{\alpha}\in J}m(\underline{\alpha},v_{j}\eta)sgn(\underline{\alpha},v_{j}\eta)u_{\underline{\alpha}}(a)\boldsymbol{\gamma}_{j}.$$ 
 Consid\'erons $\underline{\alpha}=\{\alpha_{1},...,\alpha_{n}\}\in J$. Pour tout $i=1,...,n$, fixons une "coracine" $\check{\alpha}_{i}$  que nous normalisons par la condition $<\alpha_{i},\check{\alpha}_{i}>=1$ (sic!).  Notons $m$ le volume du quotient de ${\cal A}_{\tilde{M}}^{\tilde{G}}$ par le ${\mathbb Z}$-module engendr\'e par ces $\check{\alpha}_{i}$, pour $i=1,...,n$. Le terme $\rho^{\tilde{G}}(\alpha_{i},v_{j}\eta)$ d\'efini en  1.5 est proportionnel \`a $\check{\alpha}_{i}$. Il r\'esulte des d\'efinitions que
$$(3) \qquad m(\underline{\alpha},v_{j}\eta)sgn(\underline{\alpha},v_{j}\eta)=m\prod_{i=1,...,n}<\alpha_{i},\rho(\alpha_{i},v_{j}\eta)>.$$
Pour tout $i$, appliquons la relation 1.7(11):
$$\rho(\alpha_{i},v_{j}\eta)=\sum_{\beta_{i}\in \Sigma^{\tilde{G}}(A_{\tilde{R}}),\beta_{i,\tilde{M}}=\alpha_{i}}\rho(\beta_{i},v\eta)_{\tilde{M}}.$$
Notons $\underline{J}'_{\underline{\alpha}}$ l'ensemble des ensembles  $\underline{\beta}=\{\beta_{1},...,\beta_{n}\}$ d'\'el\'ements de $\Sigma^{\tilde{G}}(A_{\tilde{R}})$ tels que $\beta_{i,\tilde{M}}=\alpha_{i}$ pour tout $i$ (en num\'erotant convenablement les \'el\'ements de cet ensemble).   On obtient
$$m(\underline{\alpha},v_{j}\eta)sgn(\underline{\alpha},v_{j}\eta)=\sum_{\underline{\beta}\in \underline{J}'_{\underline{\alpha}}}m\prod_{i=1,...,n}<\alpha_{i},\rho(\beta_{i},u\eta)>.$$
 Pour chaque ensemble  $\underline{\beta}=\{\beta_{1},...,\beta_{n}\}\in \underline{J}'_{\underline{\alpha}}$, d\'efinissons l'espace de Levi $\tilde{L}_{\underline{\beta}}$ de sorte que ${\cal A}_{\tilde{L}_{\underline{\beta}}}$ soit l'intersection des annulateurs des $\beta_{i}$ dans ${\cal A}_{\tilde{R}}$. Alors $\underline{\beta}$ appartient \`a une unique classe $J_{\underline{\beta}} \in {\cal J}_{\tilde{R}}^{\tilde{L}_{\underline{\beta}}}$.   Il r\'esulte des d\'efinitions que la relation (1) est v\'erifi\'ee pour $\tilde{L}=\tilde{L}_{\underline{\beta}}$ et que la classe $J_{\underline{\beta}}$ s'envoie sur  $J$ par l'injection ${\cal J}_{\tilde{R}}^{\tilde{L}_{\underline{\beta}}}\subset {\cal J}_{\tilde{M}}^{\tilde{G}}$. On introduit des coracines $\check{\beta}_{i}$ comme ci-dessus. On a
$$<\alpha_{i},\rho(\beta_{i},u\eta)>=<\alpha_{i},\check{\beta}_{i}><\beta_{i},\rho(\beta_{i},u\eta)>,$$
d'o\`u
$$m(\underline{\alpha},v_{j}\eta)sgn(\underline{\alpha},v_{j}\eta)=\sum_{\underline{\beta}\in J'_{\underline{\alpha}}}mm'_{\underline{\beta}}\prod_{i=1,...,n}<\beta_{i},\rho(\beta_{i},u\eta)>,$$
o\`u
$$m'_{\underline{\beta}}=\prod_{i=1,...,n}<\alpha_{i},\check{\beta}_{i}>.$$
Le produit $mm'_{\underline{\beta}}$ est le volume du quotient de ${\cal A}_{\tilde{M}}^{\tilde{G}}$ par le ${\mathbb Z}$-module engendr\'e par les $\check{\beta}_{i,\tilde{M}}$ pour $i=1,...,n$. Un calcul simple montre que
$$mm'_{\underline{\beta}}=d_{\tilde{R}}^{\tilde{G}}(\tilde{M}, \tilde{L}_{\underline{\beta}})m_{\underline{\beta}},$$
o\`u
 $m_{\underline{\beta}}$ est  le volume du quotient de ${\cal A}_{\tilde{R}}^{ \tilde{L}_{\underline{\beta}}}$ par le ${\mathbb Z}$-module engendr\'e par les $\check{\beta}_{i}$. D'o\`u
 $$m(\underline{\alpha},v_{j}\eta)sgn(\underline{\alpha},v_{j}\eta)=\sum_{\underline{\beta}\in \underline{J}'_{\underline{\alpha}}}d_{\tilde{R}}^{\tilde{G}}(\tilde{M} \tilde{L}_{\underline{\beta}})m_{\underline{\beta}}\prod_{i=1,...,n}<\beta_{i},\rho(\beta_{i},u\eta)>.$$
  Par une \'egalit\'e similaire \`a (3), cette expression  devient
 $$m(\underline{\alpha},v_{j}\eta)sgn(\underline{\alpha},v_{j}\eta)=\sum_{\underline{\beta}\in \underline{J}'_{\underline{\alpha}}}d_{\tilde{M}}^{\tilde{G}}(\tilde{M}, \tilde{L}_{\underline{\beta}})m(\underline{\beta},u\eta)sgn(\underline{\beta},u\eta).$$
 D'autre part, on a $u_{\underline{\alpha}}(a)=u_{\underline{\beta}}(a)$ pour tout $\underline{\beta}\in \underline{J}'_{\underline{\alpha}}$. 
 L'\'egalit\'e (2) devient
 $$ \rho_{J}^{\tilde{G}}(\boldsymbol{\gamma}^{\tilde{M}},a)=\sum_{j=1,...,k}c_{j}\sum_{\underline{\alpha}\in J}\sum_{\underline{\beta}\in \underline{J}'_{\underline{\alpha}}}d_{\tilde{R}}^{\tilde{G}}(\tilde{M}, \tilde{L}_{\underline{\beta}})m(\underline{\beta},u)sgn(\underline{\beta},u)u_{\underline{\beta}}(a)\boldsymbol{\gamma}_{j}$$
 $$=\sum_{\underline{\alpha}\in J}\sum_{\underline{\beta}\in \underline{J}'_{\underline{\alpha}}}d_{\tilde{R}}^{\tilde{G}}(\tilde{M}, \tilde{L}_{\underline{\beta}})m(\underline{\beta},u)sgn(\underline{\beta},u)u_{\underline{\beta}}(a)\boldsymbol{\gamma}^{\tilde{M}}.$$
 Pour un espace de Levi $\tilde{L}$ v\'erifiant (1), on v\'erifie que la r\'eunion sur les $\underline{\alpha}\in J$ des $\underline{\beta}\in \underline{J}'_{\underline{\alpha}}$ tels que $\tilde{L}_{\underline{\beta}}=\tilde{L}$ est vide si $J\not\in {\cal J}_{\tilde{R}}^{\tilde{L}}$. Sinon, c'est l'ensemble des $\underline{\beta}\in J$, o\`u $J$ est vu comme un \'el\'ement de ${\cal J}_{\tilde{R}}^{\tilde{L}}$. D'o\`u
 $$ \rho_{J}^{\tilde{G}}(\boldsymbol{\gamma}^{\tilde{M}},a)=\sum_{\tilde{L}\in {\cal L}(\tilde{R}), J\in {\cal J}_{\tilde{R}}^{\tilde{L}}}d_{\tilde{R}}^{\tilde{G}}(\tilde{M}, \tilde{L})\sum_{\underline{\beta}\in J}m(\underline{\beta},u)sgn(\underline{\beta},u)u_{\underline{\beta}}(a)\boldsymbol{\gamma}^{\tilde{M}}.$$
 Par une \'egalit\'e similaire \`a (2), on a pour tout $\tilde{L}$ intervenant ci-dessus
 $$\rho_{J}^{\tilde{L}}(\boldsymbol{\gamma},a)=\sum_{\underline{\beta}\in J}m(\underline{\beta},u)sgn(\underline{\beta},u)u_{\underline{\beta}}(a)\boldsymbol{\gamma}.$$
 L'\'egalit\'e pr\'ec\'edente devient
 $$ \rho_{J}^{\tilde{G}}(\boldsymbol{\gamma}^{\tilde{M}},a)=\sum_{\tilde{L}\in {\cal L}(\tilde{R}), J\in {\cal J}_{\tilde{R}}^{\tilde{L}}}d_{\tilde{R}}^{\tilde{G}}(\tilde{M}, \tilde{L})\rho_{J}^{\tilde{L}}(\boldsymbol{\gamma},a)^{\tilde{M}}.$$
 Cela prouve le lemme. $\square$
 
 {\bf Variante.} Supposons $(G,\tilde{G},{\bf a})$ quasi-d\'eploy\'e et \`a torsion int\'erieure et fixons un syst\`eme de fonctions $B$ comme en 1.9. On dispose de l'ensemble ${\cal J}_{\tilde{M}}^{\tilde{G}}(B_{{\cal O}^{\tilde{M}}})$. Les constructions s'adaptent pour cet ensemble et on a un lemme similaire.

 \bigskip
 
 \subsection{Termes $\sigma_{J}$ et induction}
 Soient $(G,\tilde{G},{\bf a})$ un triplet quasi-d\'eploy\'e et \`a torsion int\'erieure, $B$ un syst\`eme de fonctions comme en 1.9, $\tilde{M}$ un espace de Levi de $\tilde{G}$, $\tilde{R}$ un espace de Levi de $\tilde{M}$ et ${\cal O}$ une classe de conjugaison stable semi-simple dans $\tilde{R}(F)$. On note ${\cal O}^{\tilde{M}}$ la classe de conjugaison stable dans $\tilde{M}(F)$ qui contient ${\cal O}$. 
 
  \ass{Lemme}{Soient $J\in {\cal J}_{\tilde{M}}^{\tilde{G}}(B_{{\cal O}^{\tilde{M}}})$, $\boldsymbol{\delta}\in D_{g\acute{e}om}^{st}({\cal O})\otimes Mes(R(F))^*$ et $a\in A_{M}(F)$ en position g\'en\'erale et assez proche de $1$.   On a l'\'egalit\'e
 $$\sigma_{J}^{\tilde{G}}(\boldsymbol{\delta}^{\tilde{M}},a)=\sum_{\tilde{L}\in {\cal L}(\tilde{R}), J\in {\cal J}_{\tilde{R}}^{\tilde{L}}}e_{\tilde{R}}^{\tilde{G}}(\tilde{M},\tilde{L})\sigma_{J}^{\tilde{L}}(\boldsymbol{\delta},a)^{\tilde{M}}.$$}
 
 En utilisant le lemme pr\'ec\'edent, la d\'emonstration est similaire \`a celle du (ii) de la proposition
 1.14. 

 \bigskip
 
 \subsection{Termes $\rho_{J}^{\tilde{G},{\cal E}}({\bf M}',\boldsymbol{\delta},a)$ et induction}
 Soient $(G,\tilde{G},{\bf a})$ un triplet quelconque, $\tilde{M}$ un espace de Levi de $\tilde{G}$, ${\bf M}'=(M',{\cal M}',\tilde{\zeta})$ une donn\'ee endoscopique elliptique et relevante de $(M,\tilde{M},{\bf a})$ et $R'$ un groupe de Levi de $M'$ qui est relevant. On construit comme en [I] 3.4 un espace de Levi $\tilde{R}$ de $\tilde{M}$ qui lui correspond et une donn\'ee endoscopique ${\bf R}'$ de $(R,\tilde{R},{\bf a})$ qui est elliptique et relevante. Soit ${\cal O}$ une classe de conjugaison stable  semi-simple dans $\tilde{R}(F)$. On note  ${\cal O}'$ la r\'eunion des classes de conjugaison stable dans $\tilde{R}'(F)$ qui correspondent \`a ${\cal O}$.
  
 \ass{Lemme}{Soient $J\in {\cal J}_{\tilde{M}}^{\tilde{G}}$, $\boldsymbol{\delta}\in D_{g\acute{e}om}({\bf R}',{\cal O}')\otimes Mes(R'(F))^*$ et $a\in A_{\tilde{M}}(F)$ en position g\'en\'erale et proche de $1$.  On a l'\'egalit\'e
 $$\rho_{J}^{\tilde{G},{\cal E}}({\bf M}',\boldsymbol{\delta}^{{\bf M}'},a)=\sum_{\tilde{L}\in {\cal L}(\tilde{R}), J\in {\cal J}_{\tilde{R}}^{\tilde{L}}}d_{\tilde{M}}^{\tilde{G}}(\tilde{M},\tilde{L})\rho_{J}^{\tilde{L},{\cal E}}({\bf R}',\boldsymbol{\delta},a)^{\tilde{M}}.$$}
 
 La d\'emonstration est similaire \`a celle du (i) de la proposition 1.14.
 
 {\bf Variante.} Supposons $(G,\tilde{G},{\bf a})$ quasi-d\'eploy\'e et \`a torsion int\'erieure. Fixons un syst\`eme de fonctions $B$ comme en 1.9. On a un lemme similaire en rempla\c{c}ant l'ensemble ${\cal J}_{\tilde{M}}^{\tilde{G}}$ par ${\cal J}_{\tilde{M}}^{\tilde{G}}(B_{{\cal O}^{\tilde{M}}})$.

\section{Le cas non ramifi\'e}

\bigskip

\subsection{Int\'egrales orbitales pond\'er\'ees de la fonction caract\'eristique d'un espace hypersp\'ecial}

Dans toute cette section, on suppose  $(G,\tilde{G},{\bf a})$ non ramifi\'e et $p$ grand. Pr\'ecis\'ement, on impose les hypoth\`eses (1) \`a (4) de [I] 6.1 ainsi que l'hypoth\`ese (Hyp) de cette r\'ef\'erence. Le groupe $G(F)$ est muni d'une mesure canonique pour laquelle $mes(K)=1$ pour tout sous-groupe compact hypersp\'ecial $K$ de $G(F)$ (rappelons que deux tels sous-groupes sont conjugu\'es par le groupe $G_{AD}(F))$). On munit $G(F)$ de cette mesure et on se d\'ebarrasse ainsi des espaces de mesures intervenant dans les sections pr\'ec\'edentes. On fera de m\^eme pour les autres groupes non ramifi\'es qui interviendront.

On fixe un sous-espace hypersp\'ecial $\tilde{K}$ de $\tilde{G}(F)$, on note $K$ le sous-groupe hypersp\'ecial de $G(F)$ associ\'e. On note ${\bf 1}_{\tilde{K}}$ la fonction caract\'eristique de $\tilde{K}$. Soit $\tilde{M}$ un espace de Levi de $\tilde{G}$ tel que $M$ soit en bonne position relativement \`a $K$. On d\'efinit une forme lin\'eaire $r_{\tilde{M}}^{\tilde{G}}(.,\tilde{K})$ sur $D_{g\acute{e}om}(\tilde{M}(F),\omega)$ par
$$r_{\tilde{M}}^{\tilde{G}}(\boldsymbol{\gamma},\tilde{K})=J_{\tilde{M}}^{\tilde{G}}(\boldsymbol{\gamma},{\bf 1}_{\tilde{K}})$$
pour tout $\boldsymbol{\gamma}\in D_{g\acute{e}om}(\tilde{M}(F),\omega)$.

Soit $\tilde{Q}=\tilde{L}U_{Q}\in {\cal F}(\tilde{M})$. On a l'\'egalit\'e $({\bf 1}_{\tilde{K}})_{\tilde{Q},\omega}={\bf 1}_{\tilde{K}^{\tilde{L}}}$, o\`u $\tilde{K}^{\tilde{L}}=\tilde{K}\cap \tilde{L}(F)$. En particulier, cette fonction ne d\'epend que de $\tilde{L}$. La formule habituelle de descente des int\'egrales orbitales donne donc la formule suivante. Soient  $\tilde{L}\in {\cal L}(\tilde{M})$ et $\boldsymbol{\gamma}\in D_{g\acute{e}om}(\tilde{M}(F),\omega)$. On a l'\'egalit\'e
$$(1) \qquad r_{\tilde{L}}^{\tilde{G}}(\boldsymbol{\gamma}^{\tilde{L}},\tilde{K})=\sum_{\tilde{L}'\in {\cal L}(\tilde{M} )}d_{\tilde{M}}^{\tilde{G}}(\tilde{L},\tilde{L}')r_{\tilde{M}}^{\tilde{L}'}(\boldsymbol{\gamma},\tilde{K}^{\tilde{L}'}).$$

\bigskip

\subsection{L'avatar stable}
On suppose ici $(G,\tilde{G},{\bf a})$ quasi-d\'eploy\'e et \`a torsion int\'erieure. Pour tout espace de Levi $\tilde{M}$ de $\tilde{G}$, nous allons d\'efinir une forme lin\'eaire $s_{\tilde{M}}^{\tilde{G}}(.,\tilde{K})$ sur $D_{g\acute{e}om}^{st}(\tilde{M}(F))$. Comme les int\'egrales orbitales pond\'er\'ees, elle d\'epend de la mesure sur ${\cal A}_{M}^G$ fix\'ee en 1.2.
 La d\'efinition se faisant par r\'ecurrence, on doit  commencer par quelques formalit\'es.

Notre forme lin\'eaire v\'erifiera la propri\'et\'e 

(1) $s_{\tilde{M}}^{\tilde{G}}(.,\tilde{K})$ ne d\'epend que de la classe de conjugaison de $\tilde{K}$ par $G_{AD}(F)$.

Remarquons que deux sous-groupes hypersp\'eciaux de $G(F)$ sont toujours conjugu\'es par $G_{AD}(F)$ mais ce n'est pas le cas pour deux sous-espaces hypersp\'eciaux de $\tilde{G}(F)$. Pour deux tels sous-espaces $\tilde{K}'$ et $\tilde{K}''$, on a seulement: il existe $g\in G_{AD}(F)$ et $z\in Z(G;F)$ de sorte que $\tilde{K}''=z\,ad_{g}(\tilde{K}')$.

Notre forme lin\'eaire v\'erifiera aussi la propri\'et\'e

(2) $s_{\tilde{M}}^{\tilde{G}}(\boldsymbol{\delta},\tilde{K})=0$ si le support de $\boldsymbol{\delta}$ ne coupe pas $\tilde{K}$.

Consid\'erons des extensions compatibles 
$$1\to C_{1}\to G_{1}\to G\to 1\text{ et }\tilde{G}_{1}\to \tilde{G}$$
o\`u $C_{1}$ est un tore central induit, $G_{1}$ est non ramifi\'e et $\tilde{G}_{1}$ est \`a torsion int\'erieure. Soit $\lambda_{1}$ un caract\`ere non ramifi\'e de $C_{1}(F)$. On fixe un espace hypersp\'ecial $\tilde{K}_{1}$ de $\tilde{G}_{1}(F)$ se projetant sur $\tilde{K}$.   On a ${\cal A}_{M_{1}}^{G_{1}}\simeq {\cal A}_{M}^G$ et on choisit pour mesure sur le premier espace l'image par cet isomorphisme de la mesure fix\'ee sur le second. On suppose d\'efinie la forme lin\'eaire $s_{\tilde{M}_{1}}^{\tilde{G}_{1}}(.,\tilde{K}_{1})$ sur $D_{g\acute{e}om}^{st}(\tilde{M}_{1}(F))$, v\'erifiant la propri\'et\'e (2). On d\'efinit une forme lin\'eaire $s_{\tilde{M}_{1},\lambda_{1}}^{\tilde{G}_{1}}(.,\tilde{K}_{1})$ sur $D_{g\acute{e}om,\lambda_{1}}^{st}(\tilde{M}_{1}(F))$ de la fa\c{c}on suivante. Soit $\boldsymbol{\delta}\in D_{g\acute{e}om,\lambda_{1}}^{st}(\tilde{M}_{1}(F))$. On choisit  un \'el\'ement $\dot{\boldsymbol{\delta}}\in D_{g\acute{e}om}^{st}(\tilde{M}_{1}(F))$ qui s'envoie sur $\boldsymbol{\delta}$ par l'application 1.10(3). On pose
$$s_{\tilde{M}_{1},\lambda_{1}}^{\tilde{G}_{1}}(\boldsymbol{\delta},\tilde{K}_{1})=\int_{C_{1}(F)}s_{\tilde{M}_{1}}^{\tilde{G}_{1}}(\dot{\boldsymbol{\delta}}^c,\tilde{K}_{1})\lambda_{1}(c)^{-1}\,dc.$$
La propri\'et\'e (2) assure que cette int\'egrale est \`a support compact.

On note comme toujours ${\bf M}$ la donn\'ee endoscopique "maximale" de $(M,\tilde{M})$. Soit $s\in Z(\hat{M})^{\Gamma_{F}}/Z(\hat{G})^{\Gamma_{F}}$, avec $s\not=1$. On en d\'eduit une donn\'ee endoscopique ${\bf G}'(s)=(G'(s),{\cal G}'(s),s)$ de $(G,\tilde{G})$, qui est non ramifi\'ee.  Supposons-la elliptique. Alors ${\cal A}_{M}^{G'(s)}\simeq {\cal A}_{M}^G$ et on  choisit pour mesure sur le premier espace l'image par cet isomorphisme de la mesure fix\'ee sur le second. On va d\'efinir une forme lin\'eaire $s_{{\bf M}}^{{\bf G}'(s)}(.,\tilde{K})$ sur $D_{g\acute{e}om}^{st}({\bf M})$.
On associe \`a l'espace $\tilde{K}$ un espace hypersp\'ecial $\tilde{K}'(s)$ de $\tilde{G}'(s)$, dont la classe de conjugaison par $G'(s)_{AD}(F)$ est uniquement d\'etermin\'ee, cf. [I] 6.2. On choisit des donn\'ees auxiliaires $G'_{1}(s)$, $\tilde{G}'_{1}(s)$, $C_{1}(s)$, $\hat{\xi}_{1}(s)$ non ramifi\'ees (cf. [I] 6.3). On choisit un sous-espace hypersp\'ecial $\tilde{K}'_{1}(s)$ de $\tilde{G}'_{1}(s)$ se projetant sur $\tilde{K}'(s)$. On note $\Delta_{1}(s)$ le facteur de transfert associ\'e \`a ce sous-espace, cf. [I] 6.3. Soit $\boldsymbol{\delta}\in D_{g\acute{e}om}^{st}({\bf M})$. Par ce choix de facteur de transfert, cette distribution s'identifie \`a un \'el\'ement  $\boldsymbol{\delta}_{1}(s)\in D_{g\acute{e}om,\lambda_{1}}^{st}(\tilde{M}'_{1}(s))$. Puisqu'on a suppos\'e $s\not=1$, on peut supposer par r\'ecurrence que la forme lin\'eaire $s_{\tilde{M}'_{1}(s)}^{\tilde{G}'_{1}(s)}(.,\tilde{K}'_{1}(s))$ est bien d\'efinie. Un calcul formel utilisant par r\'ecurrence la propri\'et\'e (1) montre que le terme
$$s_{\tilde{M}'_{1}(s)}^{\tilde{G}'_{1}(s)}(\boldsymbol{\delta}_{1}(s),\tilde{K}'_{1}(s))$$
ne d\'epend pas des choix de donn\'ees auxiliaires. On le note $s_{{\bf M}}^{{\bf G}'(s)}(\boldsymbol{\delta},\tilde{K})$, ce qui d\'efinit la forme lin\'eaire $s_{{\bf M}}^{{\bf G}'(s)}(.,\tilde{K})$.

On a d\'efini au paragraphe pr\'ec\'edent la forme lin\'eaire $r_{\tilde{M}}^{\tilde{G}}(.,\tilde{K})$ sur $D_{g\acute{e}om}(\tilde{M}(F))$, sous l'hypoth\`ese que $M$ \'etait en bonne position relativement \`a $K$. Montrons que:

(3) sa restriction \`a $D_{g\acute{e}om}^{st}(\tilde{M}(F))$ ne d\'epend que de la classe de conjugaison de $\tilde{K}$ par $G_{AD}(F)$. 

Soit $g\in G_{AD}(F)$, posons $\tilde{K}'=ad_{g}(\tilde{K})$, $K'=ad_{g}(K)$ et supposons que $M$ est encore en bonne position relativement \`a $K'$. On peut fixer deux sous-tores maximaux   $T$ et $T'$ de $M$, d\'efinis sur $F$ et maximalement d\'eploy\'es,  de sorte que $K$, resp. $K'$, soit le fixateur d'un point hypersp\'ecial dans l'appartement de l'immeuble de $G$ associ\'e \`a $T$, resp. $T'$. Les tores $T$ et $T'$ sont conjugu\'es par  $M(F)$. On peut donc fixer $m\in M(F)$ tel que $ad_{m}(T)=T'$. Posons $K''=ad_{m^{-1}}(K')$. Alors $K$ et $K''$ sont les fixateurs de points hypersp\'eciaux dans l'appartement de l'immeuble de $G$ associ\'e \`a $T$. On sait que deux tels points se d\'eduisent l'un de l'autre par l'action d'un \'el\'ement du normalisateur de $T$ dans $G_{AD}(F)$. Puisque de plus,  le normalisateur de $T$ dans $K$ se projette surjectivement sur le groupe de Weyl de $T$, nos deux points se d\'eduisent en fait l'un de l'autre par l'action d'un \'el\'ement de $T_{ad}(F)$. Soit donc $t\in T_{ad}(F)$ tel que $K''=ad_{t}(K)$. Alors $ad_{mt}(K)=K'$, donc $ad_{g^{-1}mt}$ conserve $K$, donc $g^{-1}mt\in K_{ad}$. Dans notre situation \`a torsion int\'erieure, cela entra\^{\i}ne que $ad_{g^{-1}mt}$ conserve $\tilde{K}$. Donc $\tilde{K}'=ad_{mt}(\tilde{K})$.  Cela montre qu'il existe $x\in M_{ad}(F)$ tel que $\tilde{K}'=ad_{x}(\tilde{K})$. Par simple transport de structure, on a l'\'egalit\'e
$$r_{\tilde{M}}^{\tilde{G}}(\boldsymbol{\gamma},\tilde{K})=r_{\tilde{M}}^{\tilde{G}}(ad_{x}(\boldsymbol{\gamma}),\tilde{K}')$$
pour tout $\boldsymbol{\gamma}\in D_{g\acute{e}om}(\tilde{M}(F))$. Mais l'action par conjugaison de $M_{ad}(F)$ se restreint en l'identit\'e sur les distributions stables. Donc
$$r_{\tilde{M}}^{\tilde{G}}(\boldsymbol{\delta},\tilde{K})=r_{\tilde{M}}^{\tilde{G}}(\boldsymbol{\delta},\tilde{K}')$$
pour tout $\boldsymbol{\delta}\in D_{g\acute{e}om}^{st}(\tilde{M}(F))$. Cela prouve (3). $\square$

Gr\^ace \`a (3), on peut \'etendre la d\'efinition de la restriction de   $r_{\tilde{M}}^{\tilde{G}}(.,\tilde{K})$ \`a $D_{g\acute{e}om}^{st}(\tilde{M}(F))$ au cas o\`u $M$ n'est plus suppos\'e en bonne position relativement \`a $\tilde{K}$: on choisit $g\in G_{AD}(F)$ tel que $M$ soit en bonne position relativement \`a $ad_{g}(K)$; pour $\boldsymbol{\delta}\in D_{g\acute{e}om}^{st}(\tilde{M}(F))$, on pose $r_{\tilde{M}}^{\tilde{G}}(\boldsymbol{\delta},\tilde{K})=r_{\tilde{M}}^{\tilde{G}}(\boldsymbol{\delta},ad_{g}(\tilde{K}))$. L'assertion (3) assure que cela ne d\'epend pas du choix de $g$.

On peut maintenant d\'efinir notre forme lin\'eaire $s_{\tilde{M}}^{\tilde{G}}(.,\tilde{K})$. Pour $\boldsymbol{\delta}\in D_{g\acute{e}om}^{st}(\tilde{M}(F))$, on pose
$$(4) \qquad s_{\tilde{M}}^{\tilde{G}}(\boldsymbol{\delta},\tilde{K})=r_{\tilde{M}}^{\tilde{G}}(\boldsymbol{\delta},\tilde{K})-\sum_{s\in Z(\hat{M})^{\Gamma_{F}}/Z(\hat{G})^{\Gamma_{F}}, s\not=1}i_{\tilde{M}}(\tilde{G},\tilde{G}'(s))s_{{\bf M}}^{{\bf G}'(s)}(\boldsymbol{\delta},\tilde{K}).$$
La v\'erification des propri\'et\'es (1) et (2) est imm\'ediate par r\'ecurrence, et gr\^ace \`a (3).

Notons une autre propri\'et\'e formelle de notre forme lin\'eaire. Le groupe $Z(G;F)$ agit par multiplication sur $\tilde{G}(F)$, donc aussi sur $C_{c}^{\infty}(\tilde{M}(F))$ (pr\'ecis\'ement $f^z(\gamma)=f(z\gamma)$) puis sur $D_{g\acute{e}om}(\tilde{M}(F))$, en conservant l'espace des distributions stables. Pour $z\in Z(G;F)$ et $\boldsymbol{\delta}\in D_{g\acute{e}om}^{st}(\tilde{M}(F))$, on a l'\'egalit\'e
$$s_{\tilde{M}}^{\tilde{G}}(\boldsymbol{\delta}^{z},\tilde{K})=s_{\tilde{M}}^{\tilde{G}}(\boldsymbol{\delta},z\tilde{K}).$$

\bigskip

\subsection{L'avatar endoscopique}

On revient au cas o\`u $(G,\tilde{G},{\bf a})$ est quelconque (mais non ramifi\'e comme dans toute la section). Soient $\tilde{M}$ un Levi de $\tilde{G}$ et ${\bf M}'=(M',{\cal M}',\tilde{\zeta})$ une donn\'ee endoscopique elliptique et non ramifi\'ee de $\tilde{M}$. Pour $\tilde{s}\in \tilde{\zeta}Z(\hat{M})^{\Gamma_{F},\hat{\theta}}/Z(\hat{G})^{\Gamma_{F},\hat{\theta}}$, on dispose de la donn\'ee endoscopique ${\bf G}'(\tilde{s})=(G'(\tilde{s}),{\cal G}'(\tilde{s}),\tilde{s})$ de $(G,\tilde{G},{\bf a})$, qui est non ramifi\'ee. Supposons-la elliptique. Alors ${\cal A}_{M'}^{G'(\tilde{s})}\simeq {\cal A}_{\tilde{M}}^{\tilde{G}}$ et on  choisit pour mesure sur le premier espace l'image par cet isomorphisme de la mesure fix\'ee sur le second. On d\'efinit une forme lin\'eaire $s_{{\bf M}'}^{{\bf G}'(\tilde{s})}(.,\tilde{K})$ sur
$D_{g\acute{e}om}^{st}({\bf M}')$ de la m\^eme fa\c{c}on qu'au paragraphe pr\'ec\'edent. C'est-\`a-dire que l'on choisit des donn\'ees auxiliaires non ramifi\'ees $G'_{1}(\tilde{s})$, $\tilde{G}'_{1}(\tilde{s})$, $C_{1}(s)$, $\hat{\xi}_{1}(s)$. On fixe un sous-espace hypersp\'ecial $\tilde{K}'_{1}(\tilde{s})$ de $\tilde{G}'_{1}(\tilde{s};F)$ se projetant sur un espace $\tilde{K}'(\tilde{s})$ de $\tilde{G}'(\tilde{s};F)$ associ\'e \`a $\tilde{K}$. On utilise le facteur de transfert associ\'e \`a cet espace. Ainsi, un \'el\'ement $\boldsymbol{\delta}'\in D_{g\acute{e}om}^{st}({\bf M}')$ s'identifie \`a un \'el\'ement $\boldsymbol{\delta}'_{1}(\tilde{s})\in D_{g\acute{e}om,\lambda_{1}(\tilde{s})}^{st}(\tilde{M}'_{1}(\tilde{s};F))$. On pose
$$s_{{\bf M}'}^{{\bf G}'(\tilde{s})}(\boldsymbol{\delta}',\tilde{K})=s_{\tilde{M}'_{1}(\tilde{s}),\lambda_{1}}(\boldsymbol{\delta}'_{1}(\tilde{s}),\tilde{K}'_{1}(\tilde{s})).$$
Cela ne d\'epend pas des choix de donn\'ees auxiliaires.

Cela \'etant, on d\'efinit une forme lin\'eaire $r_{\tilde{M}}^{\tilde{G},{\cal E}}({\bf M}',.,\tilde{K})$ sur $D_{g\acute{e}om}^{st}({\bf M}')$ par l'\'egalit\'e
$$(1) \qquad r_{\tilde{M}}^{\tilde{G},{\cal E}}({\bf M}',\boldsymbol{\delta}',\tilde{K})=\sum_{\tilde{s}\in \tilde{\zeta}Z(\hat{M})^{\Gamma_{F},\hat{\theta}}/Z(\hat{G})^{\Gamma_{F},\hat{\theta}}}i_{\tilde{M}'}(\tilde{G},\tilde{G}'(\tilde{s}))s_{{\bf M}'}^{{\bf G}'(\tilde{s})}(\boldsymbol{\delta}',\tilde{K}).$$

\bigskip

\subsection{Le lemme fondamental pond\'er\'e}

\ass{Th\'eor\`eme}{Soit $\tilde{M}$ un espace de Levi de $\tilde{G}$ en bonne position relativement \`a $K$. Soit ${\bf M}'$ une donn\'ee endoscopique elliptique et non ramifi\'ee de $\tilde{M}$.
 Pour $\boldsymbol{\delta}'\in D_{g\acute{e}om}^{st}({\bf M}')$, on a l'\'egalit\'e
 $$r_{\tilde{M}}^{\tilde{G},{\cal E}}({\bf M}',\boldsymbol{\delta}',\tilde{K})=r_{\tilde{M}}^{\tilde{G}}(transfert(\boldsymbol{\delta}'),\tilde{K}).$$}
 
 Supposons le support de $\boldsymbol{\delta}'$ form\'e d'\'el\'ements semi-simples fortement  $\tilde{G}$-r\'eguliers. Dans ce cas, l'assertion est le lemme fondamental pond\'er\'e sous sa forme usuelle. Elle est maintenant prouv\'ee d'apr\`es [W3] th\'eor\`eme 3.8.  Ce th\'eor\`eme \'etait conditionnel, mais les conditions impos\'ees sont lev\'ees par les r\'esultats de Ngo Bao Chau et ceux de  Chaudouard et Laumon, bien que ces derniers ne soient pas encore publi\'es en toute g\'en\'eralit\'e. La suite de la section est consacr\'ee \`a la suppression de l'hypoth\`ese faite ci-dessus sur le support de $\boldsymbol{\delta}'$.  On suppose fix\'es $\tilde{M}$ et ${\bf M}'$ comme dans l'\'enonc\'e.
 
 \bigskip
 
 \subsection{D\'eveloppement en germes}
 
 Soit ${\cal O}\subset \tilde{M}(F)$ une r\'eunion finie de classes de conjugaison semi-simples. La proposition 2.3 d\'efinit des germes $g_{\tilde{M},{\cal O}}^{\tilde{L}}$ pour tout $\tilde{L}\in {\cal L}(\tilde{M})$. D'apr\`es [A1], proposition 9.1, les int\'egrales pond\'er\'ees non $\omega$-\'equivariantes v\'erifient le m\^eme d\'eveloppement que leurs versions $\omega$-\'equivariantes, avec les m\^emes germes. En particulier, on a
 $$(1) \qquad r_{\tilde{M}}^{\tilde{G}}(\boldsymbol{\gamma},\tilde{K})=\sum_{\tilde{L}\in {\cal L}(\tilde{M})}r_{\tilde{L}}^{\tilde{G}}(g_{\tilde{M}}^{\tilde{L}}(\boldsymbol{\gamma}),\tilde{K})$$
 pour tout $\boldsymbol{\gamma}\in D_{g\acute{e}om}(\tilde{M}(F),\omega)$ assez proche de $\tilde{C}$.

 Supposons $(G,\tilde{G},{\bf a})$ quasi-d\'eploy\'e et \`a torsion int\'erieure. Soit ${\cal O}\subset \tilde{M}(F)$ une r\'eunion finie de classes de conjugaison stable. Supposons que ${\cal O}$ soit form\'ee d'\'el\'ements $\tilde{G}$-\'equisinguliers. Montrons que l'on a l'\'egalit\'e
 $$(2) \qquad s_{\tilde{M}}^{\tilde{G}}(\boldsymbol{\delta},\tilde{K})=s_{\tilde{M}}^{\tilde{G}}(g_{\tilde{M},{\cal O}}^{\tilde{M}}(\boldsymbol{\delta}),\tilde{K})$$
 pour tout $\boldsymbol{\delta}\in D_{g\acute{e}om}^{st}(\tilde{M}(F))$ assez proche de ${\cal O}$. Remarquons que, d'apr\`es le lemme 2.2, $g_{\tilde{M},{\cal O}}^{\tilde{M}}(\boldsymbol{\delta})$ est stable, le membre de droite ci-dessus est donc d\'efini.
 
 Preuve. On utilise la d\'efinition 4.2(4). En raisonnant par r\'ecurrence, on peut supposer que, pour $s\in Z(\hat{M})^{\Gamma_{F}}/Z(\hat{G})^{\Gamma_{F}}$, $s\not=1$, on a l'\'egalit\'e
 $$s_{{\bf M}}^{{\bf G}'(s)}(\boldsymbol{\delta},\tilde{K})=s_{{\bf M}}^{{\bf G}'(s)}(g_{\tilde{M},{\cal O}}^{\tilde{M}}(\boldsymbol{\delta}),\tilde{K})$$
 pourvu que $\boldsymbol{\delta}$ soit assez proche de ${\cal O}$. L'hypoth\`ese faite sur ${\cal O}$ et la relation 2.3(1) entra\^{\i}nent que le d\'eveloppement (1) se simplifie en
 $$r_{\tilde{M}}^{\tilde{G}}(\boldsymbol{\delta},\tilde{K})=r_{\tilde{M}}^{\tilde{G}}(g_{\tilde{M},{\cal O}}^{\tilde{M}}(\boldsymbol{\delta}),\tilde{K}).$$
 Mais alors le membre de droite de la relation 4.2(4) pour l'\'el\'ement $\boldsymbol{\delta}$ co\"{\i}ncide avec la m\^eme expression relative \`a l'\'el\'ement $g_{\tilde{M},{\cal O}}^{\tilde{M}}(\boldsymbol{\delta})$, donc avec $s_{\tilde{M}}^{\tilde{G}}(g_{\tilde{M},{\cal O}}^{\tilde{M}}(\boldsymbol{\delta}),\tilde{K})$. $\square$
 
 Revenons au cas g\'en\'eral, soit ${\cal O}'$ une classe de conjugaison stable d'\'el\'ements semi-simples dans $\tilde{M}'(F)$. Il lui correspond une classe de conjugaison par $M(\bar{F})$ dans $\tilde{M}(\bar{F})$. Supposons que cette classe soit form\'ee d'\'el\'ements  $\tilde{G}$-\'equisinguliers.  Alors
 
 (3) l'assertion du th\'eor\`eme 4.4 est v\'erifi\'ee pour tout $\boldsymbol{\delta}'\in D_{g\acute{e}om}^{st}({\bf M}')\cap D_{g\acute{e}om}({\cal O}')$. 
 
 Preuve. Le lemme 2.2 assure que l'on peut trouver $\boldsymbol{\delta}'_{reg}\in D_{g\acute{e}om,\tilde{G}-reg}^{st}({\bf M}')$, aussi proche que l'on veut de ${\cal O}'$, de sorte que $g_{{\bf M}'}^{{\bf M}'}(\boldsymbol{\delta}'_{reg})=\boldsymbol{\delta}'$. En fait, en reprenant les d\'emonstrations, on voit que l'on peut supposer le support de $\boldsymbol{\delta}'_{reg}$ en position g\'en\'erale, en particulier $\tilde{G}$-r\'egulier. Consid\'erons la d\'efinition 4.3(1). Notons que l'hypoth\`ese sur ${\cal O}'$ entra\^{\i}ne que, pour tout $\tilde{s}$ y intervenant, ${\cal O}'$ est form\'e d'\'el\'ements  $\tilde{G}'(\tilde{s})$-\'equisinguliers. La relation (2) ci-dessus entra\^{\i}ne alors l'\'egalit\'e
 $$r_{\tilde{M}}^{\tilde{G},{\cal E}}({\bf M}',\boldsymbol{\delta}',\tilde{K})=r_{\tilde{M}}^{\tilde{G},{\cal E}}({\bf M}',\boldsymbol{\delta}'_{reg},\tilde{K}).$$
 Comme en l'a dit en 4.4, le th\'eor\`eme est d\'ej\`a connu pour $\boldsymbol{\delta}'_{reg}$. Le membre de droite ci-dessus est donc \'egal \`a $r_{\tilde{M}}^{\tilde{G}}(transfert(\boldsymbol{\delta}'_{reg}),\tilde{K})$. En appliquant (1), qui se simplifie gr\^ace \`a l'hypoth\`ese sur ${\cal O}'$, c'est aussi $r_{\tilde{M}}^{\tilde{G}}(g_{\tilde{M}}^{\tilde{M}}(transfert(\boldsymbol{\delta}'_{reg})),\tilde{K})$. On a vu en 2.6 que 
 $$g_{\tilde{M}}^{\tilde{M}}(transfert(\boldsymbol{\delta}'_{reg}))=transfert(g_{{\bf M}'}^{{\bf M}'}(\boldsymbol{\delta}'_{reg})),$$
 d'o\`u
 $$g_{\tilde{M}}^{\tilde{M}}(transfert(\boldsymbol{\delta}'_{reg}))=transfert(\boldsymbol{\delta}').$$
 On obtient $r_{\tilde{M}}^{\tilde{G},{\cal E}}({\bf M}',\boldsymbol{\delta}',\tilde{K})=r_{\tilde{M}}^{\tilde{G}}(transfert(\boldsymbol{\delta}'),\tilde{K})$, comme on le voulait. $\square$

\bigskip

\subsection{Un espace de germes sous hypoth\`eses sur $p$}
On a d\'efini en 3.1 un espace de germes $U_{J}$ pour tout $J\in {\cal J}_{\tilde{M}}^{\tilde{G}}$. On pose
$$U_{\tilde{M}}^{\tilde{G}}=\sum_{\tilde{L}\in {\cal L}(\tilde{M})}\sum_{J\in {\cal J}_{\tilde{M}}^{\tilde{L}}}U_{J} \text{ et }U_{\tilde{M}}^{\tilde{G},+}=\sum_{\tilde{L}\in {\cal L}(\tilde{M}),\tilde{L}\not=\tilde{M}}\sum_{J\in {\cal J}_{\tilde{M}}^{\tilde{L}}}U_{J}.$$
Rappelons que, pour $J=\emptyset$, on a $U_{\emptyset}={\mathbb C}$. Donc $U_{\tilde{M}}^{\tilde{G}}={\mathbb C}+U_{\tilde{M}}^{\tilde{G},+}$.

\ass{Lemme}{Supposons $p\not=2,3,5$. Alors  les espaces ${\mathbb C}$ et $U_{\tilde{M}}^{\tilde{G},+}$ sont en somme directe.}

{\bf Remarque.} L'hypoth\`ese faite sur $p$ en 5.1 entra\^{\i}ne $p\not=2,3,5$. Mais cette derni\`ere condition est suffisante ici.

\bigskip

Preuve. On descend les fonctions \`a l'alg\`ebre de Lie $\mathfrak{a}_{\tilde{M}}(F)$ comme en 3.1  et on utilise les notations de ce paragraphe.  Soit
$$u=\sum_{\underline{\alpha}}c_{\underline{\alpha}}u_{\underline{\alpha}}$$
une combinaison lin\'eaire de fonctions $u_{\underline{\alpha}}$, o\`u $\underline{\alpha}=\{\alpha_{i};i=1,...,n\}$ d\'ecrit les ensembles form\'es d'\'el\'ements lin\'eairement ind\'ependants de $\Sigma(A_{\tilde{M}})$. Supposons $u=0$. On doit alors prouver que $c_{\emptyset}=0$ (rappelons que $u_{\emptyset}$ est la fonction constante de valeur $1$). On raisonne par r\'ecurrence sur $a_{\tilde{M}}-a_{\tilde{G}}$. L'assertion est triviale si ce nombre est nul, autrement dit si $\tilde{M}=\tilde{G}$, puisqu'alors $u$ se r\'eduit \`a $c_{\emptyset }u_{\emptyset}$. Supposons $a_{\tilde{M}}-a_{\tilde{G}}>0$, fixons un \'el\'ement $\tilde{P}\in {\cal P}(\tilde{M})$ et notons $\Delta_{\tilde{P}}$ l'ensemble de racines simples associ\'e. Fixons $\beta\in \Delta_{\tilde{P}}$. Il lui est associ\'e un espace de Levi $\tilde{L}$, maximal parmi les espaces de Levi propres de $\tilde{G}$, de sorte que $\Delta_{\tilde{P}}\cap \tilde{L}=\Delta_{\tilde{P}}-\{\beta\}$.  Notons $u^{\tilde{L}}$ la sous-somme de $u$, o\`u on ne conserve que les $\underline{\alpha}=\{\alpha_{i};i=1,...,n\}$ tels que $\alpha_{i}\in \Sigma^{\tilde{L}}(A_{\tilde{M}})$ pour tout $i$. Elle contient le terme constant $ c_{\emptyset }u_{\emptyset}$ et appartient \`a $U_{\tilde{M}}^{\tilde{L}}$.  On va prouver que $u^{\tilde{L}}=0$.  L'hypoth\`ese de r\'ecurrence permettra alors de conclure que $c_{\emptyset}=0$. Puisque $u^{\tilde{L}}$ est invariante par translations par $\mathfrak{a}_{\tilde{L}}(F)$, il suffit de prouver qu'elle est nulle sur $\mathfrak{a}_{\tilde{M}}^{\tilde{L}}(F)$. Fixons un \'el\'ement $H^{\tilde{L}}$ de cet espace, ainsi qu'un \'el\'ement $H_{\beta}$ de $\mathfrak{a}_{\tilde{L}}(F)$ tel que $\beta(H_{\beta})=1$. Soit $k\in {\mathbb Z}$, posons $H=\varpi_{F}^kH_{\beta}+H^{\tilde{L}}$. Soit $\alpha\in \Sigma(A_{\tilde{M}})$. On peut \'ecrire $\alpha=\sum_{\beta'\in \Delta_{\tilde{P}}}n_{\beta'}\beta'$, avec des coefficients entiers $n_{\beta'}$. Ainsi qu'on l'a d\'ej\`a dit, on voit en consid\'erant tous les syst\`emes de racines possibles que ces entiers appartiennent \`a l'ensemble $\{0,\pm 1,...,\pm 6\}$. L'hypoth\`ese sur $p$ entra\^{\i}ne que $\vert n_{\beta'}\vert_{F}=1$ pour tout $\beta'\in \Delta_{\tilde{P}}$. Si $\alpha\in \Sigma^{\tilde{L}}(A_{\tilde{M}})$, on a $n_{\beta}=0$ et $\alpha(H)=\alpha(H^{\tilde{L}})$. Si $\alpha\not\in \Sigma^{\tilde{L}}(A_{\tilde{M}})$, on a $n_{\beta}\not=0$. Consid\'erons $H^{\tilde{L}}$ comme fix\'e et $k$ comme variable, on a alors $\vert \alpha(H)\vert _{F}=\vert n_{\beta}\beta(\varpi_{F}^kH_{\beta})\vert _{F}$ si $k$ est assez n\'egatif, donc $\vert \alpha(H)\vert _{F}=-klog(q)$ si $k$ est assez n\'egatif. On en d\'eduit que, pour tout ensemble  $\underline{\alpha}=\{\alpha_{i};i=1,...,n\}$, le terme $u_{\underline{\alpha}}(H)$ est un mon\^ome en $k$ pour $k$ assez n\'egatif, dont le degr\'e est nul si et seulement si tous les $\alpha_{i}$ appartiennent \`a $\Sigma^{\tilde{L}}(A_{\tilde{M}})$. Dans ce dernier cas, on a simplement $u_{\underline{\alpha}}(H)=u_{\underline{\alpha}}(H^{\tilde{L}})$. Ainsi  $u(H)$ est un polyn\^ome en $k$ pour $k$ assez n\'egatif, dont le terme constant est $u^{\tilde{L}}(H^{\tilde{L}})$. Puisque $u=0$, le polyn\^ome est identiquement nul. Cela entra\^{\i}ne $u^{\tilde{L}}(H^{\tilde{L}})=0$. On a d\'ej\`a dit que cela permettait de conclure. $\square$

Le lemme permet de d\'efinir le terme constant d'un \'el\'ement de $U_{\tilde{M}}^{\tilde{G}}$: c'est sa projection sur le sous-espace ${\mathbb C}$.

\bigskip

\subsection{D\'eveloppement des fonctions $r_{\tilde{M}}^{\tilde{G}}(.,\tilde{K})$ et $s_{\tilde{M}}^{\tilde{G}}(.,\tilde{K})$}

La proposition 3.2 d\'efinit des applications lin\'eaires $\rho^{\tilde{G}}_{J}$ pour tout $J\in {\cal J}_{\tilde{M}}^{\tilde{G}}$. On a

(1) pour tout $\boldsymbol{\gamma}\in D_{g\acute{e}om}(\tilde{M}(F),\omega)$, le germe en $1$ de la fonction
$$a\mapsto r_{\tilde{M}}^{\tilde{G}}(a\boldsymbol{\gamma},\tilde{K})$$
est \'equivalent \`a
$$\sum_{\tilde{L}\in {\cal L}(\tilde{M})}\sum_{J\in {\cal J}_{\tilde{M}}^{\tilde{L}}}r_{\tilde{L}}^{\tilde{G}}(\rho^{\tilde{L}}_{J}(\boldsymbol{\gamma},a)^{\tilde{L}},\tilde{K}).$$

Il suffit de reprendre la preuve de 3.2. La seule propri\'et\'e que l'on utilisait des int\'egrales orbitales $\omega$-\'equivariantes \'etait la relation 1.7(12). Une relation analogue vaut pour notre fonction $r_{\tilde{M}}^{\tilde{G}}(.,\tilde{K})$ gr\^ace \`a 4.1(1). 

Vu comme fonction de $a$, chaque terme $r_{\tilde{L}}^{\tilde{G}}(\rho_{J}^{\tilde{L}}(\boldsymbol{\gamma},a)^{\tilde{L}},\tilde{K})$ appartient \`a $U_{J}$. On en d\'eduit

(2)  pour tout $\boldsymbol{\gamma}\in D_{g\acute{e}om}(\tilde{M}(F),\omega)$, le germe en $1$ de la fonction
$$a\mapsto r_{\tilde{M}}^{\tilde{G}}(a\boldsymbol{\gamma},\tilde{K})$$
est \'equivalent \`a un \'el\'ement de $U_{\tilde{M}}^{\tilde{G}}$ dont le terme constant est \'egal \`a $r_{\tilde{M}}^{\tilde{G}}(\boldsymbol{\gamma},\tilde{K})$. 

Rappelons que, d'apr\`es 3.1(2), cet \'el\'ement de $U_{\tilde{M}}^{\tilde{G}}$ est uniquement d\'etermin\'e.

Supposons $(G,\tilde{G},{\bf a})$ quasi-d\'eploy\'e et \`a torsion int\'erieure.  

(3) pour tout $\boldsymbol{\delta}\in D_{g\acute{e}om}^{st}(\tilde{M}(F))$, le germe en $1$ de la fonction
$$a\mapsto s_{\tilde{M}}^{\tilde{G}}(a\boldsymbol{\delta},\tilde{K})$$
est \'equivalent \`a un \'el\'ement de $U_{\tilde{M}}^{\tilde{G}}$ dont le terme constant est \'egal \`a $s_{\tilde{M}}^{\tilde{G}}(\boldsymbol{\delta},\tilde{K})$.

Cela r\'esulte par r\'ecurrence de la d\'efinition 3.2(4),  de (2) ci-dessus et de la propri\'et\'e \'evidente suivante: pour $s\in Z(\hat{M})^{\Gamma_{F}}/Z(\hat{G})^{\Gamma_{F}}$, on a l'inclusion $U_{\tilde{M}}^{\tilde{G}(s)}\subset U_{\tilde{M}}^{\tilde{G}}$.

On note $\xi:A_{\tilde{M}}(F)\to A_{\tilde{M}'}(F)$ l'homomorphisme naturel. On a

(4) pour tout $\boldsymbol{\delta}'\in D_{g\acute{e}om}^{st}({\bf M}')$, le germe en $1$ de la fonction
$$a\mapsto r_{\tilde{M}}^{\tilde{G},{\cal E}}({\bf M}',\xi(a)\boldsymbol{\delta}',\tilde{K})$$
est \'equivalent \`a un \'el\'ement de $U_{\tilde{M}}^{\tilde{G}}$ dont le terme constant est \'egal \`a $r_{\tilde{M}}^{\tilde{G},{\cal E}}({\bf M}',\boldsymbol{\delta}',\tilde{K})$.

Preuve.  Soit $\tilde{s}\in \tilde{\zeta}Z(\hat{M})^{\Gamma_{F},\hat{\theta}}/Z(\hat{G})^{\Gamma_{F},\hat{\theta}}$. Soit $u'\in U_{\tilde{M}'}^{\tilde{G}'(\tilde{s})}$ et $v'$ un germe de fonctions d\'efini presque partout au voisinage de $1$ dans $A_{\tilde{M}'}(F)$. Notons $u$, resp. $v$, le germe de fonctions $a\mapsto u'(\xi(a))$, resp. $a\mapsto v'(\xi(a))$, d\'efini presque partout sur $A_{\tilde{M}}(F)$ au voisinage de $1$. On a

(5) $u$ appartient \`a $U_{\tilde{M}}^{\tilde{G}}$ et a m\^eme terme constant que $u'$;

(6) si $v'$ est \'equivalent \`a $u'$, alors $v$ est \'equivalent \`a $u$. 

La preuve de (6) est imm\'ediate. Prouvons (5).  Par lin\'earit\'e, il suffit de prouver cette assertion quand $u'$ appartient \`a un ensemble de g\'en\'erateurs de $ U_{\tilde{M}'}^{\tilde{G}'(\tilde{s})}$. On peut donc fixer des \'el\'ements lin\'eairement ind\'ependants $\alpha'_{1},...,\alpha'_{n}$ de $\Sigma^{G'(\tilde{s})}(A_{\tilde{M}'})$ et supposer que $u'$ est la fonction
$$u'(a')=\prod_{i=1,...,n}log(\vert \alpha'_{i}(a')-\alpha'_{i}(a')^{-1}\vert _{F}).$$
Identifions $\underline{les}$ paires de Borel \'epingl\'ees de $G$ et $G'(\tilde{s})$ \`a des paires de Borel \'epingl\'ees de ces groupes pour lesquelles $M$ et $M'$ sont standard. On note $T$ et $T'$ les tores de ces paires et $\theta$ l'automorphisme habituel qui conserve la paire de Borel \'epingl\'ee de $G$. On rappelle que, pour $\beta\in \Sigma(T)$, on note $n_{\beta}$ le plus petit entier $k\geq1$ tel que $\theta^k(\beta)=\beta$.  Notre hypoth\`ese que $p$ est grand implique que tous les entiers $n_{\beta}$ sont premiers \`a $p$. L'homomorphisme $\xi$ ci-dessus est la restriction d'un homomorphisme encore not\'e $\xi:T\to T'\simeq T/(1-\theta)(T)$.   Pour tout $i=1,...,n$, choisissons $\beta'_{i}\in \Sigma(T')$ dont la restriction \`a $A_{\tilde{M}}$ soit $\alpha'_{i}$. D'apr\`es la description d\'ej\`a utilis\'ee plusieurs fois du syst\`eme de racines de $G'(\tilde{s})$, il y a une racine $\beta_{i}\in \Sigma(T)$ de sorte que $\beta'_{i}\circ \xi$ co\"{\i}ncide sur $T^{\theta,0}$ avec $n_{\beta_{i}}\beta_{i}$ ou $2n_{\beta_{i}}\beta_{i}$ (en notation additive). Notons $\alpha_{i}$ la restriction de $\beta_{i}$ \`a $A_{\tilde{M}}$. Alors $\alpha'_{i}( \xi(a))=\alpha_{i}(a)^{m_{i}}$ pour tout $a\in A_{\tilde{M}}(F)$, avec $m_{i}=n_{\beta_{i}}$ ou $m_{i}=2n_{\beta_{i}}$. En tout cas, $m_{i}$ est un entier premier \`a $p$. Mais alors
$\vert \alpha'_{i}(\xi(a))-\alpha'_{i}(\xi(a))^{-1}\vert _{F}$ co\"{\i}ncide avec $\vert \alpha_{i}(a)-\alpha_{i}(a)^{-1}\vert _{F}$ pour $a$ assez voisin de $1$. Le germe de $u$ co\"{\i}ncide donc avec celui de 
$$a\mapsto \prod_{i=1,...,n}log(\vert \alpha_{i}(a)-\alpha_{i}(a)^{-1}\vert _{F}).$$
Les racines $(\alpha_{i})_{i=1,...,n}$ sont encore lin\'eairement ind\'ependantes. D'apr\`es les d\'efinitions, la fonction ci-dessus appartient \`a $U_{\tilde{M}}^{\tilde{G}}$. Cela prouve la premi\`ere assertion de (5). Si $n>0$, les termes constants de $u'$ et de $u$ sont nuls. Si $n=0$, il est clair que $u=u'$. Cela ach\`eve de prouver (5). 

Revenons \`a la preuve de (4). On utilise la d\'efinition 4.3(1). Pour tout $\tilde{s}\in \tilde{\zeta}Z(\hat{M})^{\Gamma_{F},\hat{\theta}}/Z(\hat{G})^{\Gamma_{F},\hat{\theta}}$, les assertions (3), (5) et (6) impliquent que la fonction 
$$a\mapsto s_{{\bf M}'}^{{\bf G}'(\tilde{s})}(\xi(a)\boldsymbol{\delta}',\tilde{K})$$
est \'equivalent \`a un \'el\'ement de $U_{\tilde{M}}^{\tilde{G}}$ dont le  terme constant est  $s_{{\bf M}'}^{{\bf G}'(\tilde{s})}\boldsymbol{\delta}',\tilde{K})$. On sommant ces r\'esultats sur $\tilde{s}$, on obtient (4). $\square$

\bigskip

\subsection{Preuve du th\'eor\`eme 3.4}

Pour $a\in A_{\tilde{M}}(F)$ en position g\'en\'erale, l'\'el\'ement $\xi(a)\boldsymbol{\delta}'$ est combinaison lin\'eaire d'\'el\'ements  v\'erifiant les hypoth\`eses de 4.5(3) . Donc
$$r_{\tilde{M}}^{\tilde{G},{\cal E}}({\bf M}',\xi(a)\boldsymbol{\delta}',\tilde{K})=r_{\tilde{M}}^{\tilde{G}}(transfert(\xi(a)\boldsymbol{\delta}'),\tilde{K}).$$
On a l'\'egalit\'e $transfert(\xi(a)\boldsymbol{\delta}')=a\,transfert(\boldsymbol{\delta}')$. D'apr\`es les assertions 4.7(2) et (4), les deux membres ci-dessus sont \'equivalents \`a des \'el\'ement de $U_{\tilde{M}}^{\tilde{G}}$  dont les termes constants sont respectivement $r_{\tilde{M}}^{\tilde{G},{\cal E}}({\bf M}',\boldsymbol{\delta}',\tilde{K})$ et $r_{\tilde{M}}^{\tilde{G}}(transfert(\boldsymbol{\delta}'),\tilde{K})$. Ces \'el\'ements de $U_{\tilde{M}}^{\tilde{G}}$ sont forc\'ement les m\^emes d'apr\`es 3.1(2). Donc leurs termes constants sont \'egaux. Cela prouve le th\'eor\`eme.

\bigskip

{\bf Bibliographie}

[A1] J. Arthur: {\it The local behaviour of weighted orbital integrals}, Duke Math. J. 56 (1988), p. 223-293

[A2] -----------: {\it The trace formula in invariant form}, Annals of Math. 114 (1981), p. 1-74

[A3]------------: {\it The invariant trace formula I. Local theory}, J. AMS 1 (1988), p. 323-383

[A4] -----------: {\it On the transfer of distributions: weighted orbital integrals}, Duke Math. J. 99 (1999), p. 209-283

[A5] -----------: {\it On a family of distributions obtained from Eisenstein series II: explicit formulas}, Amer. J. of Math. 104 (1982), p. 1289-1336

[A6]-----------: {\it A stable trace formula I. General expansions}, J. Inst. Math. Jussieu 1 (2002), p. 175-277

[K] R. Kottwitz: {\it Stable trace formula: cuspidal tempered terms}, Duke Math. J. 51 (1984), p. 611-650

[KS] -----------, D. Shelstad: {\it Foundations of twisted endoscopy}, Ast\'erisque 255 (1999)

 [Lan] R. P. Langlands: {\it Stable conjugacy: definitions and lemmas}, Can. J. of Math. 31 (1979), p. 700-725

[LS] -------------------, D. Shelstad: {\it On the definition of transfer factors}, Math. Ann. 278 (1987), p. 219-271

[I] J.-L. Waldspurger: {\it Pr\'eparation \`a la stabilisation de la formule des traces tordue I: endoscopie tordue sur un corps local }, pr\'epublication 2012

[W1] --------------------: {\it La formule des traces locale tordue}, pr\'epublication 2012

[W2]---------------------: {\it L'endoscopie tordue n'est pas si tordue}, Memoirs AMS 908 (2008)

[W3]---------------------:{ \it A propos du lemme fondamental pond\'er\'e tordu}, Math. Ann. 343 (2009), p. 103-174

\bigskip

Institut de Math\'ematiques de Jussieu- CNRS

2 place Jussieu

75005 Paris

e-mail: waldspur@math.jussieu.fr

\end{document}